\documentclass[11pt,letterpaper]{amsart}
\usepackage{blindtext}

\makeatletter
\renewcommand\part{%
   \if@noskipsec \leavevmode \fi
   \par
   \addvspace{4ex}%
   \@afterindentfalse
   \secdef\@part\@spart}

\def\@part[#1]#2{%
    \ifnum \c@secnumdepth >\m@ne
      \refstepcounter{part}%
      \addcontentsline{toc}{part}{\thepart\hspace{1em}#1}%
    \else
      \addcontentsline{toc}{part}{#1}%
    \fi
       \begin{center}\Large \partname\nobreakspace\thepart.\nobreakspace\scshape #2\end{center}%
    \nobreak
    \vskip 1.5ex
    \@afterheading}
\def\@spart#1{%
    {\parindent \z@ \raggedright
     \interlinepenalty \@M
     \normalfont
     \huge \bfseries #1\par}%
     \nobreak
     \vskip 3ex
     \@afterheading}
\makeatother

\usepackage{tikz}
\tikzstyle{every node}=[circle, draw, fill=black!50,
                        inner sep=0pt, minimum width=3pt]

\usetikzlibrary{calc,arrows.meta}

\usetikzlibrary{positioning}

\usepackage{amsfonts,amsmath,latexsym,color,epsfig,hyperref,enumitem, amssymb,bbm}

\usepackage{enumitem}   

\newtheorem{prop}{Proposition}[section]
\newtheorem{definition}{Definition}[section]
\newtheorem{lemma}[prop]{Lemma}

\newtheorem{obs}[prop]{Observation}
\newtheorem{theorem}[prop]{Theorem}

\newtheorem*{prop*}{Proposition}
\newtheorem*{theorem*}{Theorem}
\newtheorem{claim}[prop]{Claim}

\newtheorem{remark}[prop]{Remark}

\DeclareMathOperator{\dir}{dir}

\newcommand{\mc}{\mathcal }
\newcommand{\M}{\operatorname{M} }

\newcommand{\CC}{\mathbb{C}}
\newcommand{\TT}{\mathbb{T}}
\newcommand{\cD}{\mathcal{D}}
\newcommand{\bp}{\mathbf{p}}
\newcommand{\bT}{\mathbf{T}}

\newcommand{\cP}{\mathcal{P}}
\newcommand{\eps}{\epsilon}

\newcommand{\ZZ}{\mathbb{Z}}
\newcommand{\RR}{\mathbb{R}}

\newcommand{\NN}{\mathbb{N}}

\renewcommand{\d}{\delta}

\renewcommand{\le}{\leqslant}
\renewcommand{\ge}{\geqslant}
\renewcommand{\leq}{\leqslant}
\renewcommand{\geq}{\geqslant}

\newcommand{\Mat}{\operatorname{Mat}}

\let\oldtocsection=\tocsection
\let\oldtocsubsection=\tocsubsection
\let\oldtocsubsubsection=\tocsubsubsection

\renewcommand{\tocsection}[2]{\hspace{0em}\oldtocsection{#1}{#2}}
\renewcommand{\tocsubsection}[2]{\hspace{1em}\oldtocsubsection{#1}{#2}}
\renewcommand{\tocsubsubsection}[2]{\hspace{2em}\oldtocsubsubsection{#1}{#2}}

\title{Kakeya and Furstenberg problems for sticky sets of tubes}

\author{Hong Wang, Dmitrii Zakharov}

\address{Courant institute of mathematical sciences, New York University}
\email{hw3639@nyu.edu}

\address{Department of Mathematics, Massachusetts Institute of Technology, Cambridge, USA}
\email{zakhdm@mit.edu}

\date{}

\begin{document}

\begin{abstract}
    We study incidence problems for sticky sets of tubes. In the first part, we prove a sharp Furstenberg set estimate for sticky sets of tubes in $\RR^3$ satisfying convex Wolff axioms. In the second part, we prove the Kakeya set conjecture in $\RR^4$ for sticky sets of tubes satisfying convex Wolff axioms. 
\end{abstract}

\maketitle

\tableofcontents

\section{Introduction}\label{sec:intro}

A Kakeya set is a compact set in $\RR^d$ which contains a unit length line segment in every direction. The Kakeya set conjecture is the assertion that any such set must have full Hausdorff dimension. The two dimensional case was proved by Davies \cite{davies1971some}, see \cite{katz2000recent, zahl2025survey} for an introduction and survey of the Kakeya set problem. 
A Furstenberg set is a generalization of the Kakeya set where instead of containing line segments, the set is required to intersect with many lines in a set of fixed dimension. Furstenberg sets were introduced in Furstenberg's unpublished work and in  Wolff's influential exposition \cite{wolff1999recent}.   
On the discretized level, these are questions about how much can collections of thin tubes overlap with each other and what geometric structure occurs if they do.

A sharp Furtenberg set estimate in the plane was proven in 2023 by  Ren and Wang \cite{ren2023furstenberg} and the Kakeya set conjecture in $\RR^3$ was resolved in 2025 by Wang and Zahl \cite{wang2025volume} (see also surveys by Guth \cite{guth2025introduction}, \cite{guth2025outline}, \cite{guth2026kakeya}, and a streamlined proof by Guth, Wang and Zahl \cite{guth2026streamlined}). In both cases, the proof proceeds by reducing the problem to a more structured subproblem, namely when the set of tubes is `sticky', and then solving this subproblem. The sticky Furstenberg estimate in the plane was proven earlier by Orponen and Shmerkin  using their ABC sum-product theorem \cite{orponen2023projections}. 
The sticky case of the 3-dimensional Kakeya conjecure was resolved by Wang and Zahl \cite{wang2026sticky} by following a strategy proposed by Katz and Tao \cite{tao2014}.

In higher dimensions, some partial progress has been made \cite{guth2018polynomial, katz2018polynomial, katz2020kakeya} but sharp results are still a major challenge. 
In this paper, we study sticky versions of the Fustenberg and Kakeya problems in dimensions three and four, respectively. 
The following two theorems are the main results of this paper. We refer to Section \ref{sec:notation} for standard definitions and notation. 

\begin{theorem}\label{thm:furstenberg-R3}
    Let $s \in (0,1]$, $t \in (0,2]$. For any $\eps>0$ there exists $\eta(s,t, \eps)>0$ such that the following holds for $\d$ sufficiently small. 
    Let $\TT$ be a $(2t, \delta^{-\eta})$-AD-regular set of $\delta$-tubes in $\RR^3$ and satisfies the Convex Wolff axiom $C_{t-CW}(\TT)\leq  \delta^{-\eta}$. Suppose $Y$ is a $(\delta, s, \delta^{-\eta})$-shading on $\TT$ such that 
    \[
    \sum_{T\in \TT} |Y(T)|_\d \ge \d^{\eta} |\TT| \max_{T\in \TT} |Y(T)|_\d.
    \]
    Then 
    \begin{equation}\label{eq:furstenberg}
    \left|\bigcup_{T\in \TT} Y(T)\right|_{\delta} \geq \delta^{\epsilon} |Y(T)|_\d \delta^{-\min( 2t, t+s,  2  )}.     
    \end{equation}
\end{theorem}

The convex Wolff axiom means that every convex set $U$ contains at most $|U|^t |\TT|$ many tubes from $\TT$. For most values of $s,t$, our proof goes through with a somewhat weaker version of the Convex Wolff axiom while giving the same lower bound \eqref{eq:furstenberg}, see Section \ref{subsec:weaker-wolff} for details. 

This resolves the sticky case of the three dimensional Furstenberg set conjecture, see \cite[Conjecture 0.13]{wang2024restriction}. In \cite{wang2024restriction}, they posed this conjecture as an intermediate step towards the two-ends Furstenberg estimate in $\RR^3$ \cite[Conjecture 0.9]{wang2024restriction}, which has since been disproven by Cohen \cite{cohen2025}, and the Fourier restriction conjecture \cite[Conjecture 0.1]{wang2024restriction}. 

Standard examples show that the exponent of $\d$ in (\ref{eq:furstenberg}) is optimal (see \cite{wang2024restriction}). 
The exponent $\d^{-2t}$ corresponds to the `sparse' case when shadings $Y(T)$ are essentially disjoint and $|\bigcup Y(T)| \approx |Y(T)| |\TT| \approx |Y(T)| \d^{-2t}$. The exponent $\d^{-2}$ corresponds to the `dense' case when the union $E = \bigcup Y(T)$ satisfies $|T\cap E| \approx |T| |E|$ for all $T\in \TT$. 
The exponent $\d^{-s-t}$ corresponds to the most interesting example where $E$ is a $\d$-neighborhood of a 3-dimensional grid and $\TT$ consists of tubes with rational slopes and offsets.
More precisely, pick some $N, r\ge 1$ and let $P = [-1,1]^3 \cap N^{-1} \ZZ^3$ be the $N\times N\times N$ grid. Let $\mc L$ be the set of lines $\ell$ of the form
\[
\ell = p + \RR (a,b,c),\quad p\in [-1/2,1/2]^3 \cap N^{-1} \ZZ^3, \quad a,b,c \in \ZZ \cap [-N/r, N/r]
\]
Every such line satisfies $|\ell \cap P| \gtrsim r$ and moreover the intersection $\ell\cap P$ is an arithmetic progression of diameter $\sim 1$. For each line $\ell \in \mc L$ there are $\sim r$-many ways to select the `central' point $p$, so by double counting we have
\[
r |\mc L| \sim N^3 (N/r)^3,
\]
i.e. $|\mc L| \sim N^6 r^{-4}$. Let $\TT$ be the union of $\d$-tubes of the form $T = N_\d \cap [-1,1]^3$ over all lines $\ell \in \mc L$ and for each such $T\in \TT$ define the shading $Y(T) = N_\d (\ell \cap P)$. To get an example for the Furstenberg set problem, we need to arrange parameters so that $\TT$ is a $(\d, 2t)$-set and $Y(T)$ is a $(\d, s)$-set. This forces the choice
\[
\d^{-2t} \sim N^6 r^{-4}, \quad \d^{-s} \sim r.
\]
For these parameters, we then have
\[
\left|\bigcup Y(T)\right|_\d = |P| = N^3 = (N^6 r^{-4})^{1/2} r^2 \sim \d^{-t-2s} \sim |Y(T)|_\d \d^{-t-s},
\]
matching the middle term in \eqref{eq:furstenberg}. Finally, it is not hard to verify that this example satisfies the required Convex Wolff axiom. 
Note this set of tubes is not sticky, so it's technically not an example for Theorem \ref{thm:furstenberg-R3}.
However one can modify this construction to make it sticky, see \cite{orponen2024projectionsAD}.

\begin{theorem}\label{thm:kakeya-R4}
    For any $\eps>0$ there exists $\eta(\eps)>0$ such that the following holds for $\d$ sufficiently small. 
    Let $\TT$ be a $(3, \d^{-\eta})$-AD-regular set of $\d$-tubes in $\RR^4$ and satisfies the Convex Wolff axiom $C_{CW}(\TT) \le \d^{-\eta}$. Suppose $Y$ is a $\d^{\eta}$-dense shading on $\TT$. Then
    \[
    \left|\bigcup_{T\in \TT}Y(T)\right| \geq \d^{\eps}.
    \]
\end{theorem}

This theorem is a direct generalization of the sticky Kakeya estimate in $\RR^3$ \cite[Theorem 5.2]{wang2025assouad}. An interesting feature of Theorem \ref{thm:kakeya-R4} is that we only impose the Convex Wolff axioms on $\TT$ instead of a more general polynomial Wolff axioms. 
The reason is that for sticky sets of tubes, the Convex Wolff axiom (tubes do not concentrate in convex sets) automatically implies a corresponding Polynomial Wolff axiom (tubes do not concentrate in bounded complexity semi-algebraic sets). 
On the other hand, it is well-known that for general sets of tubes in dimension at least 4, Convex Wolff axioms are not sufficient: there exists a collection $\TT$ of $\approx \d^{-3}$-many $\d$-tubes in $\RR^4$ that satisfies the convex Wolff axiom $C_{CW}(\TT)\lesssim 1$ and at the same time $\TT$ is contained in a union of $\d^{-1/2}$-many $\d$-neighborhoods of quadratic surfaces.
So, in particular, $|\bigcup Y(T)| \lesssim \d^{1/2}$. See Section 9 in \cite{guth2025introduction} for more details.

Our Theorem \ref{thm:kakeya-R4} circumvents this issue of tubes concentrating in neighborhoods of quadrics. It appears that in order to resolve the full Kakeya conjecture in $\RR^4$, one would need to deal with this difficulty is some way. While our result does not address this important issue, we hope that our results and the techniques developed in this paper will be useful for the general case. 

We also hope that our techniques will be useful in higher dimensional sticky problems. In fact, some parts of the proof of Theorem \ref{thm:kakeya-R4} address some special cases of sticky Kakeya in $\RR^d$. On the other hand, new ideas seem to be required to deal with other cases. We leave this direction for future work.

\medskip
    \noindent\emph{Paper organization.}
In Section \ref{sec:proof-overview} we give an overview of the proofs. Sections \ref{sec:notation}--\ref{sec:structure} contain some preliminary results used for both main theorems. 
Theorem \ref{thm:furstenberg-R3} is proven in Part I (Sections \ref{sec:setup-furstenberg}--\ref{sec:very-coplanar}) and Theorem \ref{thm:kakeya-R4} is proven in Part II (Sections \ref{sec:setup-kakeya}--\ref{sec:32planar}). 
The two parts are mostly self-contained and do not require readers familiarity with the contents of each other. 

\medskip
    \noindent\emph{Prior work.}
Let us mention some earlier work studying sticky Kakeya sets in higher dimensions. Kulkarni \cite[Proposition 31]{Kulkarni2024StickyPlanebrush} proved an estimate of the form $|\bigcup_{\TT} Y(T)|_\d \gtrapprox \d^{ -0.603\ldots\cdot d}$, where $\TT$ is a $(d-1)$-dimensional sticky set of tubes in $\RR^d$ (not required to satisfy any Wolff axioms). For large $d$ this improves on the general volume lower bound $|\bigcup_{\TT} T|_\d \gtrapprox \d^{-0.597\ldots \cdot d}$ due to Katz--Tao \cite{katz2000improved} for sets of tubes in $\RR^d$ with $\d$-separated directions. 

In dimension 4, Choudhuri \cite{Choudhuri2026StickyKakeyaR4} recently showed that sticky Kakeya sets have Hausdorff dimension at least 3.25, improving on the lower bound 3.059 for general Kakeya sets due to Katz--Zahl \cite{katz2020kakeya}.
These estimates utilize the planebrush estimate from \cite{katz2020kakeya} (for sets of tubes where tubes through a fixed point span a 2-plane) and a trilinear estimate of Guth--Zahl \cite{guth2018polynomial} (for sets of tubes where tubes through a fixed point span a 3-plane). There are some parallels to these approaches in our work, but ultimately the proofs are quite different. 
We refer to the survey \cite{zahl2025survey} for other recent work on the Kakeya problem.

\subsection{Acknowledgements} We thank Alex Cohen and Larry Guth for helpful conversations. 
This collaboration was started when the second author visited NYU Courant in November 2024; he thanks Mehtaab Sawhney for providing a couch during that visit. The second author was partially supported by a Simons Dissertation Fellowship. 

\medskip
    \noindent\emph{AI use Statement.}
    All mathematical content of the paper is solely due to the authors.
LLMs were used at a later stage of manuscript preparation to improve the organization of the Appendix, $L$/$W$-tuple sections and for checking proofs.  

\section{Proof overview.}
\label{sec:proof-overview}

\subsection{Sticky Kakeya in 3D} \label{subsec:sticky-kakeya-r3-outline} Proofs of Theorems \ref{thm:furstenberg-R3} and \ref{thm:kakeya-R4} follow a similar general route and share several important ideas. To illustrate the overall approach, it is convenient to start with a review of the proof of the sticky Kakeya conjecture in $\RR^3$. We will then discuss modifications needed in the Furstenberg setting and in the higher dimensions.


The sticky Kakeya theorem in $\RR^3$ is the following statement proved by Wang--Zahl in \cite{wang2026sticky}, \cite{wang2025assouad}.

\begin{theorem}\label{thm:sticky-R3}
    For any $\eps>0$ the following holds for all sufficiently small $\eta,\d$. 
    Let $\TT$ be a $(2, \d^{-\eta})$-AD-regular set of $\d$-tubes satisfying the Convex Wolff axiom with error $\d^{-\eta}$ and let $Y$ be a $\d^\eta$-dense shading on $\TT$. Then $|\bigcup_\TT Y(T)|_\d \ge \d^{-3+\eps}$. 
\end{theorem}

Let us briefly introduce the notation that we use in this section. We are intentionally slightly informal here and postpone rigorous definitions to Section \ref{sec:notation}. 

\begin{itemize}
    \item For a set $U$ let $\TT[U] = \{T\in \TT: T\subset U\}$. 
\medskip 
    \item For $T\in \TT$ and $\varrho \in [\d,1]$ let $T_\varrho = N_\varrho T$ be a $\varrho$-tube containing $T$ (there is at most $O(1)$ such distinct $\varrho$-tubes). 
\medskip 
    \item $(2, \d^{-\eta})$-AD-regular means that $|\TT[T_\varrho] | \in [ \d^{\eta} (\varrho/\d)^2, \d^{-\eta} (\varrho/\d)^2 ]$ for all $\varrho \in [\d,1]$ and $T \in \TT$. In \cite{wang2026sticky} this property is referred to as {\em sticky}. In the following discussion, we refer to this property as {\em almost AD-regular}.
\medskip 
    \item Convex Wolff axiom with error $\d^{-\eta}$ means that $|\TT[U]| \le \d^{-\eta} |U| |\TT|$ for any convex set $U \subset \RR^3$. 
\medskip 
    \item $\d^\eta$-dense shading $Y$ on $\TT$ is a collection of sets $Y(T) \subset T$ with $\sum_{T\in \TT} |Y(T)|_\d  \ge \d^\eta |\TT|\d^{-1}$. We denote $U(\TT, Y) = \bigcup_{T\in \TT} Y(T)$. 
\end{itemize}

Write $\mu(\TT)=\frac{\sum_{\TT}|Y(T)|}{|U(\TT,Y)|}$ for the {\em multiplicity} of $(\TT, Y)$, i.e. the average number of tubes passing through a point $p \in U(\TT,Y)$. By assumption and double counting, we have $\mu(\TT) \approx \d^{-\kappa}$. Suppose that $\kappa$ is the largest possible constant for which there are sticky Kakeya configurations with multiplicity $\d^{-\kappa}$ for arbitrarily small $\d$. For the sake of contradiction, assume that $\kappa>0$.

\subsubsection{Initial reduction.}
In the first step of the proof of Theorem \ref{thm:sticky-R3} we use rescaling and thickening operations to reduce the problem to a more structured case. Suppose that the statement is false, i.e. for some constant $\kappa>0$ (think of $\kappa=0.1$) and any $\eta>0$ and arbitrarily small $\d$ there is a configuration of tubes $(\TT, Y)$ as in Theorem \ref{thm:sticky-R3} with error $\d^{-\eta}$ and $|U(\TT, Y)|_\d \lessapprox \d^{\kappa-3}$. The error $\d^{-\eta}$ is not essential for the following discussion, so we will subsume it in the `$\approx $' notation. 

The properties of the sticky Kakeya configuration $(\TT, Y)$ are preserved under the following two operations.

\medskip\noindent{\em Thickening.} For $\varrho \in (\d, 1)$ consider the set of essentially distinct $\varrho$-tubes $\TT_\varrho$ covering $\TT$ and define a shading on $\TT_\varrho$ by $Y(T_\varrho) = N_\varrho(\bigcup_{T\in \TT[T_\varrho]} Y(T))$. This gives a new sticky Kakeya configuration $(\TT_{\varrho}, Y_{\varrho})$ on scale $\varrho$, i.e. $\TT_\varrho$ is a collection of $\varrho$-tubes that is almost AD-regular and satisfies the Convex Wolff axiom.

\medskip\noindent{\em Rescaling.}
For some $\varrho \in (\d, 1)$ fix a $\varrho$-tube $T_\varrho \in \TT_\varrho$ and consider an affine linear map $\psi_{T_\varrho}: \RR^3\to \RR^3$ such that $\psi_{T_\varrho}(T_\varrho) \approx B_1(0)$. Consider the rescaled set of tubes 
\[
\TT^{T_\varrho} = \{\psi_{T_\varrho}(T), ~~ T \in \TT[T_\varrho]\}
\]
with the corresponding shading $Y^{T_\varrho}(\psi_{T_\varrho}(T)) = \psi_{T_\varrho}(Y(T))$. This gives a sticky Kakeya configuration $(\TT^{T_\varrho}, Y^{T_\varrho})$ of tubes on scale $\d/\varrho$, i.e. $\TT^{T_\varrho}$ is a collection of $\d/\varrho$-tubes that is almost AD-regular and satisfies the Convex Wolff axiom.

Here it is crucial that $\TT$ is assumed to be almost AD-regular: if we start with an arbitrary set of tubes satisfying Convex Wolff axioms, then $\TT^{T_\varrho}$ might fail to have the same property. Indeed, $\TT$ could be a {\em well-spaced} collection of tubes in the sense that every $\varrho$-tube contains at most one tube from $\TT$ (see \cite{guth2019incidence}). So the rescaled collection $\TT^{T_\varrho}$ consists of a single tube which violates the convex Wolff axiom.\footnote{In the notation of \cite{wang2025volume}, we work with the {\em Convex Frostman Wolff axiom} here. If we instead consider the Convex Katz--Tao Wolff axiom, then in the well-spaced example $\TT^{T_\varrho}$ is also Katz--Tao, but the thickening $\TT_\varrho$ might fail to be Katz--Tao. } The stickiness assumption makes the problem invariant under both operations and it is used crucially in the proof.

The basic property of rescaling and thickening is the following multiplicity bound:
\begin{equation}\label{eq:multiplicity-1}
    \mu(\TT) \lessapprox \mu(\TT_{\varrho}) \mu(\TT^{T_\varrho}).
\end{equation}
Indeed, fix $p \in U(\TT,Y)$. Assuming appropriate regularity properties, there are $\approx \mu(\TT_{\varrho})$ many $\varrho$-tubes through $p$. For each such tube $T_\varrho$, we can estimate the number of $\d$-tubes through $p$ and which are contained in $T_\varrho$ by $\mu(\TT^{T_\varrho})$. Summing over $T_\varrho$ gives \eqref{eq:multiplicity-1}. 

Now by the maximality of $\kappa$, we can estimate $\mu(\TT_{\varrho})\lessapprox \varrho^{-\kappa}$ and $\mu(\TT^{T_\varrho}) \lessapprox (\d/\varrho)^{-\kappa}$ and so applying \eqref{eq:multiplicity-1} matches the assumed lower bound $\mu(\TT) \gtrapprox \d^{-\kappa}$. So if $\TT$ is a sticky Kakeya set with maximum multiplicity, then so are $\TT_{\varrho}$ and $\TT^{T_\varrho}$ for every intermediate scale $\d \ll\varrho \ll 1$. This means that we can replace $\TT$ with the thickened or rescaled set of tubes whenever that makes the problem more structured.

The first key observation (dating back to the work of Katz--Laba--Tao \cite{katz2000improved}, see also Tao's blogpost \cite{tao2014}) is that several thickening and rescaling moves let us reduce to the case when the set $E = U(\TT, Y)$ has an extra `global grain' structure. Namely, we claim that after a sequence of thickening and rescaling moves one can write:
\begin{equation}\label{eq:E-structure-kakeya-r3}
E = \{ (x, y + f(z) x, z), ~ z \in Z, ~ y \in Y_z,~ x \in X_{y, z} \},    
\end{equation}
where:
\begin{itemize}
    \item[(i)] $Z \subset [-1,1]$ is $\approx 1$-dense set, $Y_z \subset [-1,1]$ is $(\d, 1-\kappa, \approx 1)$-AD-regular for each $z\in Z$ and $X_{y, z}$ is $\approx 1$-dense for each $z\in Z$, $y\in Y_z$,
    \item[(ii)] $f: Z \to [-1,1]$ is an $\approx 1$-Lipschitz function.
\end{itemize} 
Furthermore, the set of tubes $\TT$ is compatible with the structure of $E$ in the following sense. For every $p = (x, y+f(z)x, z) \in E$, let $\TT_p = \{T\in \TT:~ p \in Y(T)\}$ be the set of tubes through $p$. Then: 
\begin{itemize}
    \item[(iii)]  $\TT_p$ is a $(\d, \kappa, \approx 1)$-AD-regular set of tubes. Moreover, for every $\varrho \in (\d, 1)$, $|\bigcup_{p' \in B_\varrho(p)\cap E} \TT_{p'}|_\varrho \approx \varrho^{-\kappa}$. This means that up to resolution $\varrho$, sets $\TT_{p'}$ are roughly the same for $p'$ within the same $\varrho$-ball. 
    \item[(iv)] $\TT_p$ is contained in a $O(\d)$-neighborhood of the plane $\Pi_p=p + \langle (1, f(z), 0), ~ (0, \xi_p, 1) \rangle$ for some $\xi_p \in [-1,1]$. The function $p \mapsto \xi_p$ is a $\approx 1$-Lipschitz function on $E$.
\end{itemize}

Let us briefly outline how one can obtain such description of $E$. Start with a sticky Kakeya configuration $(\TT, Y)$ of maximum multiplicity $\d^{-\kappa}$. Let $E = U(\TT,Y)$.

Let $\varrho = \d^{1/2}$ and consider a $\varrho$-ball $B_\varrho$ centered at one of the points $p\in E$. The number of $\varrho$-tubes $T_\varrho$ intersecting with $B_\varrho$ is upper bounded by $\mu(\TT_\varrho) \lesssim \varrho^{-\kappa}$. Let $\TT_{\varrho,B_\varrho}$ denote this collection of $\varrho$-tubes. 
On the other hand, for typical $p' \in B_\varrho \cap E$ we claim that
\begin{equation}\label{eq:perfect-overlap}
    \#\{T_\varrho \in \TT_{\varrho,B_\varrho}:~ \exists T \in \TT_{p'}, T \subset T_\varrho\} \gtrapprox \varrho^{-\kappa}.    
\end{equation}
The left hand side counts the number of thick tubes through the ball $B_\varrho$ which contain a thin tube passing through a fixed point $p'$ inside the ball. To see \eqref{eq:perfect-overlap} we use double counting: there are $\d^{-\kappa}$ many tubes through $p'$ and for each $T_\varrho$ there are at most $\mu(\TT^{T_\varrho}) \lessapprox (\d/\varrho)^{-\kappa}$ many tubes through $p'$ that are contained in $T_\varrho$.
Since the right hand side of \eqref{eq:perfect-overlap} is roughly $|\TT_{\varrho,B_\varrho}|$, this gives a very strong local structural property: for typical tube $T_\varrho \in \TT_{\varrho, B_\varrho}$ we have 
\begin{equation} \label{eq:parallel-tubelets}
    E \cap B_{\varrho} \approx \text{union of parallel }\d\times\d\times \varrho \text{ tubes parallel to }T_\varrho.
\end{equation}
Indeed, for a fixed $T_{\varrho}$, a large fraction of points $p'\in E\cap B_{\varrho}$ has a tube $T$ such that $p' \in T \subset T_\varrho$. But then $T\cap B_\varrho$ is a $\d\times\d\times \varrho$ tube and the direction of this tube in within $O(\varrho)$ of the direction of $T_\varrho$. Since we chose $\varrho$ so that $\d/\varrho = \rho/1$, the direction vectors of these tubes are both defined up to precision $\sim \varrho$, giving \eqref{eq:parallel-tubelets}. 

Pick a pair of tubes $T_{\varrho}, T_{\varrho}' \in \TT_{\varrho, B_{\varrho}}$ that make angle $\gtrsim 1$ with each other (if all tubes make angle $\le \theta\ll 1$ then we can get a gain in \eqref{eq:multiplicity-1} by decomposing at scale $\theta$). Then applying \eqref{eq:parallel-tubelets} we see that $E\cap B_{\varrho}$ is roughly a union of parallel $\d\times \varrho\times \varrho$ slabs in the direction spanned by the directions of $T_{\varrho}, T_{\varrho}'$. Any other tube  $T_\varrho'' \in \TT_{\varrho, B_\varrho}$ must be coplanar with $T_\varrho, T_\varrho'$ up to resolution $\varrho$, otherwise \eqref{eq:parallel-tubelets} implies $|N_{\tilde\delta} E|\approx |E|$ for some $\tilde \delta \gg \delta$,  which would contradict the maximality of $\kappa$. 

In particular, for every $p'$, the set of tubes $\TT_{p'}$ passing through $p'$ must make angle at most $\varrho = \d/\varrho$ with the $\d\times \varrho \times \varrho$ slab containing $p'$. In particular, $\TT_{p'}$ is contained in a common $\varrho \times 1\times 1$ slab centered at $p'$. By thickening $\TT$ to scale $\varrho$ and replacing $(\d,\varrho)$ by $(\varrho, \varrho^{1/2})$, we may assume that $\TT_{p'}$ is contained in a $\d$-neighborhood of a plane $\Pi_{p'}$ passing through $p'$ (for all $p' \in E$).
On this new scale, there are again $\d\times \varrho \times \varrho$ local slabs coming from the argument in the previous paragraph. These slabs must be aligned with the planes $\Pi_{p'}$ that capture $\d$-tubes passing through each point $p'\in E\cap B_{\varrho}$. 

Now fix a $\varrho$-tube $T_\varrho \in \TT_{\varrho}$ and cover $E\cap T_\varrho$ by disjoint $\varrho$-balls $B_{\varrho}^{(1)}, \ldots, B_{\varrho}^{(m)}$ where $m \approx \varrho^{-1}$. The above argument applies inside each of the balls $B_{\varrho}^{(i)}$. 
Consider the rescaled tube configuration $(\tilde \TT, \tilde Y) = (\TT^{T_\varrho}, Y^{T_\varrho})$. After rotating coordinates, we can think of $T_\varrho$ as the vertical $\varrho$-tube around the $Z$-axes, so that the affine rescaling map $\psi_{T_\varrho}: T_\varrho \to B_1$ is given by
\[
\psi_{T_\varrho}:(x,y,z) \mapsto \left(\frac{x}{\varrho}, \frac{y}{\varrho}, z\right).
\]
Under this map, each ball $B_\varrho^{(i)}$ is mapped to a horizontal $\varrho \times 1\times 1$ slab. The set of parallel $\d\times \varrho \times \varrho$ slabs covering $E\cap B_\varrho^{(i)}$ is mapped to a set of parallel $\varrho \times 1 \times \varrho$ tubes contained in the horizontal slab. Note that for each slab, the common direction of these tubes could be different. In other words, we can write the set $\tilde E = U(\tilde \TT, \tilde Y)$ as
\[
\tilde E \approx N_\varrho \{ (x, y+f(z)x,z), ~ z \in \varrho \ZZ \cap [0,1], ~y\in Y_z, ~x \in [0,1]  \}.
\]
Here each fixed value of $z$ corresponds to a horizontal $\varrho$-slab, the vector $(1,f(z),0)$ is the direction vector of parallel tubes in this slab, and $y \in Y_z$ parameterizes the offset of these tubes. So after replacing $\TT$ with $\tilde \TT$ this gives the structure declared in \eqref{eq:E-structure-kakeya-r3}. By tracking through the construction one can verify that it satisfies the properties (i)--(iv) stated above. 

\subsubsection{Analysis of \texorpdfstring{$E$}{E}: first additive relation.} We now study the set $E$ given by \eqref{eq:E-structure-kakeya-r3}. We will use the interaction between the structure of $E$ and the set of tubes $\TT$ to deduce some useful additive relations between the sets $Y_z$, the directions sets of tubes $\TT_p$, the slope function $f$ and the slopes $\xi_p$. After that the problem can be analyzed using sum-product theorems. 
There are actually two different rather reductions to a sum-product problem. In $\RR^3$, either one of these reduction suffices to prove the sticky Kakeya conjecture. We use one of them in the sticky Furstenberg problem. 
For sticky Kakeya in $\RR^4$, we will need to use both.

For a tube $T\in \TT$ let $\theta(T) \in [-1,1]^2$ denote the vector such that the direction vector of $T$ is $(\theta(T),1)$.
By property (iv), namely that the set of tubes $\TT_p$ is contained in a common plane $\Pi_p = p + \langle (1, f(z),0), (0,\xi_p,1)\rangle$, we may write $\theta(\TT_p) = \{ (\varphi, \xi_p + f(z) \varphi),~ \varphi \in \Phi_p \}$ where $\Phi_p \subset [-1,1]$ is a $(\d,\kappa,\approx 1)$-AD-regular set. By (iii), if $|p-p'| \le \varrho$, then we have $N_\varrho \Phi_p \approx N_\varrho \Phi_{p'}$, i.e. the $\varrho$-neighborhoods of these sets are essentially the same.
Fix some intermediate scale $\varrho = \d^\mu \in (\d,1)$ (e.g. $\mu=1/4$ will do), consider a $\varrho$-ball ${\bf q} = B_\varrho(p_0)$ centered at a fixed point $p_0 = (x_0, y_0+f(z_0)x_0, z_0) \in E$. 
We study the set $E({\bf q}) = E \cap {\bf q}$. Let $\Phi_{\bf q} \subset \varrho \ZZ \cap [-1,1]$  be the set so that $N_\varrho \Phi_{\bf q} \approx N_\varrho \Phi_{p}$ for all $p \in E({\bf q})$. Consider the sets $Z({\bf q}) = Z\cap N_\varrho z_0$ and $Y({\bf q}) = N_{\varrho^3} Y_{z_0} \cap N_{\varrho^2} y_0 $. (Later on we use $Y({\bf q})$ to denote a slightly different set that is basically $Y_{z_0} \cap N_{\varrho}y_0$; the $Y({\bf q})$ defined here is a piece of that set that is important in the end, so we skip this distinction in this discussion.)
The set $Z({\bf q})$ is an $\approx 1$ dense subset in a $\varrho$-interval around $z_0$ and $Y({\bf q})$ is a rescaled $(\varrho, 1-\kappa, \approx 1)$-AD-regular set contained in the $\varrho^2$-neighborhood of $y_0$.

The first additive relation is as follows:
\begin{equation}\label{eq:Yq-additive-relation}
    Y({\bf q}) \approx Y({\bf q}) + (z - z') (f(z) - f(z')) (\Phi_{\bf q} - \Phi_{\bf q}), \quad z,z' \in Z({\bf q}).
\end{equation}
More precisely, an $\approx 1$ fraction of tuples $(y, z, z', \varphi_1, \varphi_2) \in Y({\bf q}) \times Z({\bf q})^2 \times \Phi_{\bf q}^2$ has the property that $y + (z-z') (f(z) - f(z'))(\varphi_1-\varphi_2) \in Y({\bf q})$.

This relation is new and it is the key observation of this paper. The previous work on sticky Kakeya --- Tao's outline \cite{tao2014} and the subsequent argument of Wang--Zahl \cite{wang2026sticky} --- make use of a related but not equivalent local symmetry. In particular, the flexibility of choosing the angles $\varphi_1,\varphi_2$ in \eqref{eq:Yq-additive-relation} from a $\kappa$-dimensional set $\Phi_{\bf q}$ is a new feature of this and other relations appearing in this paper.

We would like to use \eqref{eq:Yq-additive-relation} in order to prove that $Y({\bf q}) \approx N_{\varrho^2}y_0$ (i.e. $Y({\bf q})$ is 1-dimensional after appropriate rescaling). 
An important property of \eqref{eq:Yq-additive-relation} is that it captures the difference between the fields $\RR$ and $\CC$. Namely, one could run the proof strategy that we outline here in an attempt to prove the sticky Kakeya estimate over complex numbers (replace tubes by neighborhoods of complex lines in $\CC^3$, etc). The reduction so far was insensitive to the base field (except for the fact that we needed an access to sufficiently many scales). However, at some point one needs to see a difference since a direct analogue of Theorem \ref{thm:sticky-R3} fails over $\RR^3$: take $E$ to be a $\d$-neighborhood of a Hermitian quadric in $\CC^3$ and $\TT$ to be a set of complex tubes contained in it. Then the set $\TT$ will satisfy the Convex Wolff axiom with respect to the class of `complex convex sets', i.e. sets of the form $U = \{ z_1 v_1+z_2 v_2 +z_3 v_3\mid z_i\in \CC, |z_i|\le 1 \}$ for $v_1,v_2,v_3\in \CC^3$. 

This example manifests as the following complex example satisfying \eqref{eq:Yq-additive-relation}. 
Take $Y({\bf q}) = \Phi_{\bf q} =[-1,1] \subset \CC$, $Z({\bf q}) = [-1,1]+i[-1,1] \subset \CC$ and set $f(z) = \overline{z}$. Then (\ref{eq:Yq-additive-relation}) holds as:
\[
(z-z') (f(z)-f(z')) (\varphi_1-\varphi_2) = |z-z'|^2 (\varphi_1-\varphi_2) \in \RR.
\]
As ``$\dim_\CC\RR = 1/2$'', this gives a local example over $\CC$ with $\kappa=1/2$. In what follows we need to show that this cannot happen over $\RR$. As is now customary, we need to use a sum-product theorem to do so.

We will give a justification to \eqref{eq:Yq-additive-relation} later, first let us see how it can be used to conclude the proof. 
Let $z_1, z_2,z_1', z_2' \in Z({\bf q})$ be arbitrary and apply (\ref{eq:Yq-additive-relation}) with $(z,z') \in \{(z_1, z_1'),(z_1, z_2'), (z_2, z_1'), (z_2, z_2')\}$. 
Write $Q(z,z') = (z-z')(f(z)-f(z'))$. 
Using the following `polarization identity'
\begin{align*}
Q(z_1, z_1')& - Q(z_1, z_2') - Q(z_2, z_1') + Q(z_2, z_2')\\
&=  - (z_1-z_2)(f(z_1') -f(z_2')) - (z_1'-z_2') (f(z_1) - f(z_2)),
\end{align*}
we obtain a polarized version of (\ref{eq:Yq-additive-relation}):
\begin{equation}\label{eq:Yq-polarized}
    Y({\bf q}) \approx Y({\bf q}) + ((z_1-z_2)(f(z_1') -f(z_2')) + (z_1'-z_2') (f(z_1) - f(z_2))) (\Phi_{\bf q} - \Phi_{\bf q}).
\end{equation}

Let $V(z) = (z,f(z))$ and $V'(z') = (f(z'),z')$.
Similarly to the proof of the ABC sum-product theorem \cite{orponen2023projections}, we view the $z$-expression as a scalar product of vectors $V(z_1)-V(z_2)=( z_1-z_2 , f(z_1)-f(z_2))$ and $V'(z'_1) - V'(z'_2) =(  f(z'_1)-f(z'_2), z'_1-z'_2)$.
For a fixed pair $(z_1', z_2')$ we can then write
the set 
\[
\{(z_1-z_2)(f(z_1') -f(z_2')) + (z_1'-z_2') (f(z_1) - f(z_2)), ~ z_1, z_2\in Z({\bf q}) \}
\] 
as 
\[
(\Gamma-\Gamma) \cdot (V(z_1')-V(z_2')), \quad \Gamma = \{V(z) =(z,f(z)), ~z \in Z({\bf q})\}.
\]
The slope $\lambda$ of the vector $V(z_1')-V(z_2')$ is given by $\lambda = \frac{f(z_1') - f(z_2')}{ z_1'-z_2'}$. Thus the set of slopes $\lambda$ taken over all pairs $(z_1',z_2')$ is the set of slopes of the graph $\Gamma' = \{V'(z')=(z', f(z')), ~ z' \in Z({\bf q})\}$ (which of course is just a reflection of $\Gamma$ but it is convenient to keep them separate). 
To control the set of slopes of $\Gamma'$ we can apply the radial projections theorem of  of Orponen--Shmerkin--Wang \cite{orponen2023hausdorff}: on a discretized level, it states that if $\Gamma'$ is not contained in a tube then the set of slopes has `dimension' $\min(\dim \Gamma', 1) = 1$. 
As long as the function $f(z)$ cannot be approximated by a linear function in $z$ on the set $Z({\bf q})$, the graph $\Gamma'$ is not concentrated in any tube and we can apply radial projections to $\Gamma'$. Thus, the set of slopes $\lambda = \frac{f(z_1')-f(z_2')}{z_1-z_2}$ is 1-dimensional (more precisely, it is a $(\varrho,1,\approx1)$-set). 
Note that this step fails over $\CC$ if $f(z) = \overline{z}$.

By Kaufman's projection theorem, $\pi_\lambda(\Gamma)$ is 1-dimensional for a typical $\lambda$ taken from a 1-dimensional set. Combining this information over all $z_1',z_2'$ and pruning down if needed, we conclude that the set of expressions $(z_1-z_2)(f(z_1') -f(z_2')) + (z_1'-z_2') (f(z_1) - f(z_2))$ is 1-dimensional. Thus, in light of (\ref{eq:Yq-polarized}), the set $Y({\bf q})$ must be 1-dimensional as well. On the other hand, property (i) of $E$ states that $Y({\bf q})$ is $(1-\kappa)$-AD-regular. So $\kappa=0$ and we are done. 

It remains to consider the case when $f(z)$ can be approximated by a linear function on the set $Z({\bf q})$. This corresponds to the $C^2$-smooth case in \cite{wang2026sticky} and in that paper it was dealt with by applying a non-linear restricted projections theorem. Here we take a different route and use (\ref{eq:Yq-additive-relation}) again. Indeed, suppose that $f(z) \approx a z+b$ holds for some $a, b\in \RR$ on the set $Z({\bf q})$. Plugging this in, gives 
\begin{equation}\label{eq:Yq-quadratic-shift}
Y({\bf q}) \approx Y({\bf q}) + a (z-z')^2 (\Phi_{\bf q} -\Phi_{\bf q}), ~ z,z' \in Z({\bf q}).    
\end{equation}
If $a\gtrapprox 1$ then $a(z-z')^2 (\varphi_1-\varphi_2)$ takes a 1-dimensional set of values and so we conclude that $Y({\bf q})$ is 1-dimensional, as desired. So we must have $a \ll 1$. In particular, $|f(z)-f(z')| \lesssim a |z-z'| \ll |z-z'|$.
This means that the slope function $f$ has slower than Lipschitz variation on scale $\varrho$.

We can apply this argument at every scale $\varrho$ and conclude that $|f(z)-f(z')| \lesssim |z-z'|^{1+\varepsilon}$ on $Z$ for some $\varepsilon>0$ (in fact one can get $\varepsilon=1$). If this were the case for the whole interval $[-1,1]$ instead of $Z$, then this would imply that $f$ is a constant function. 
In this situation, one can check that $E$ decomposes into a union of parallel slabs and one quickly reaches a contradiction with the Convex Wolff axiom of $\TT$. 
However, the set $Z$ can and will be quite `porous' and so $f$ can have `H\"older exponent' $1+\varepsilon$ with $\varepsilon>0$ and not be constant. 
Nevertheless, we can use a geometric argument to show that a sticky Kakeya configuration with a slope function like this must always violate the Convex Wolff axiom. Roughly speaking, we fix a tube $T\in \TT$ and a point $p$ on it. Look at all tubes $T'$ that intersect  $T$ within distance $\varrho$ from $p$ and then look at tubes $T''$ that intersect $T'$ within distance $\varrho$ and make angle $O(\varrho)$ with $T$. We show that all such tubes $T''$ are `trapped' in a convex box of dimensions $\varrho^{1+\varepsilon} \times \varrho \times 1$ around $T$. By iterating this construction and applying a version of \eqref{eq:Yq-additive-relation} we show that the number of tubes that accumulate in this box exceeds the amount allowed by the Convex Wolff axiom. This contradiction shows that $\kappa=0$ and completes the proof of Theorem \ref{thm:sticky-R3}.

\subsubsection{Second additive relation.} We now describe a second additive relation that can be used to give another proof of the sticky Kakeya conjecture in $\RR^3$. Recall from (iv) that $\xi_p$ describes the slope of the plane $\Pi_p$ through $p$ which contains all tubes $\TT_p$. 
Fix some $p_0 = (x_0, y_0+f(z_0)x_0, z_0)\in E$, let $X({\bf q}) = X_{y_0, z_0}\cap N_\varrho x_0$, and define a function $\xi: X({\bf q}) \to [-1,1]$ by $\xi(x) = \xi_{(x, y_0+f(z_0)x, z_0)}$. Note that $\xi(x)$ is a $\approx 1$-Lipschitz function of $x$ (since $\xi_p$ is Lipschitz in $p$ by (iv)). Consider the graphs 
\[A = \left\{\begin{pmatrix}
    x \\ 
    \xi(x)
\end{pmatrix}, ~~ x \in X({\bf q})\right\}
\quad \text{and} \quad B= \left\{\begin{pmatrix}
    z \\ 
    -f(z)
\end{pmatrix}, ~~ z \in Z({\bf q})\right\}.
\] 
Then our second key relation states that
\begin{equation}\label{eq:AB-additive-relation}
    A \approx A + (B-B) (\Phi_{\bf q} - \Phi_{\bf q}).
\end{equation}
More explicitly, we have 
\begin{equation}\label{eq:xi-additive-relation}
\xi(x + (z-z') (\varphi_1-\varphi_2)  ) \approx \xi(x) - (f(z)-f(z')) (\varphi_1-\varphi_2)
\end{equation}
for an $\approx 1$-dense set of tuples $(x, z, z', \varphi_1, \varphi_2) \in X({\bf q})\times Z({\bf q})^2 \times \Phi_{\bf q}^2$. The representation (\ref{eq:AB-additive-relation}) is similar to the ABC-sum-product theorem of Orponen--Shmerkin \cite{orponen2023projections} which roughly states that if $A, B, C \subset [-1,1]$ are sets such that $A \approx A+BC$ then $\dim A \ge \min(1, \dim B+\dim C)$. Here we have a two-dimensional version of their theorem: $A, B \subset \RR^2$ and $C = \Phi_{\bf q} \subset [-1,1]$. Since $A,B$ are graphs, we only need the special case when $\dim A = \dim B=1$. 
If $B$ is not concentrated in any tube, then by the radial projections theorem, the set of slopes $\Lambda \subset [-1,1]$ of $B$ is 1-dimensional. For each slope $\lambda = \frac{b_2-b_2'}{b_1-b_1'}$, where $b=(b_1, b_2), b'=(b_1', b_2') \in B$, the identity (\ref{eq:AB-additive-relation}) then implies that $|\pi_\lambda(A)| \lesssim \frac{|A|}{|\Phi_{\bf q}|}$, where $\pi_\lambda(x,y)= x+\lambda y$. Indeed, each fiber of $A$ above the projection $\pi_\lambda$ contains a dilate of the set $\Phi_{\bf q}-\Phi_{\bf q}$. So the set $A$ drops in size along a $\approx1$ dimensional set of orthogonal projections. Kaufman projection theorem then implies that $\dim A \ge 1+\kappa$. So $\kappa=0$, as desired.
The case when $B$ is contained in a tube corresponds to $f$ being approximated by a linear function. This case can be dealt with using (\ref{eq:Yq-additive-relation}) as previously. Thus, \eqref{eq:AB-additive-relation} gives an alternative way to rule out the `non-linear' case. 

\subsubsection{Deducing the additive relations}
We now explain how to arrive at (\ref{eq:Yq-additive-relation}) and (\ref{eq:AB-additive-relation}). 
Let $Y_z({\bf q}) = N_{\varrho^3} Y_z \cap N_{\varrho^2} y_0$. The first step is to show that $Y_z({\bf q}) \approx Y_{z'}({\bf q})$ for any $z, z' \in Z({\bf q})$.

After a change of coordinates, we may assume that $p_0 = 0$, $0 \in \Phi_{\bf q}$ and that $\xi_{p_0}=0$. So by (iv), for typical $p \in E \cap B_{\varrho^2}(p_0)$ there is a `vertical' tube $T \in \TT_p$ through $p$, i.e. $|\theta(T)| \lesssim \varrho^2$. Thus, for typical $p_1 =(x_1,y_1 + f(z_0)x_1, z_0 )$ with $|p_1-p_0|\lesssim \varrho^2$ we have $p + (z-z_0) (\theta(T),1) \in N_\d E$ for many $z \in Z$. By taking $z \in Z({\bf q})$ (i.e. $|z-z_0|\lesssim \varrho$) and using $|\theta(T)|\lesssim \varrho^2$, we get that $y_1 + O(\varrho^3) \in Y_{z}({\bf q})$. 

So the sets $Y_z({\bf q})$ are approximately the same for different $z \in Z({\bf q})$. Thus, in order to prove (\ref{eq:Yq-additive-relation}), it suffices to take $z' = z_0$. For convenience we use $z_0=0$ and $f(z_0)=0$. Consider a triple of points $p_1, p_2,p\in E({\bf q})$ forming a `$V$-shape': $p_i = (x_i, y_i, 0) \in E\cap\{z=0\}$ lie on level $0$,  $p = (x, y+f(z)x, z)$ lies on level $z$ and there is a pair of tubes $T_1, T_2 \in \TT$ connecting them, i.e. $p_1, p \in T_1$ and $p_2, p \in T_2$. So we have
\[
\begin{aligned}
p&=p_1 + z (\theta(T_1), 1) + O(\d)  \\
&= p_2 + z (\theta(T_2), 1) + O(\d).    
\end{aligned}
\]
Write $\theta(T_i) = (\varphi_i, \xi_p + f(z)\varphi_i)+O(\d)$, with $\xi_p$ as in property (iv) and $\varphi_i \in \Phi_p$ for the set $\Phi_p$ parameterizing the slopes of tubes $\TT_p$ inside the plane $\Pi_p$. Plugging in these into the above relation between $p,p_1,p_2$ gives
\[
\begin{aligned}
(x, y+f(z)x) &= (x_1, y_1) + z (\varphi_1, \xi_p +f(z)\varphi_1) +O(\d) \\&= (x_2, y_2) + z (\varphi_2, \xi_p +f(z)\varphi_2) +O(\d).    
\end{aligned}
\]
Rearranging then yields:
\begin{align*}
    x_2-x_1 &= z (\varphi_1-\varphi_2) + O(\d),\\
    y_2-y_1 &= z f(z) (\varphi_1-\varphi_2) + O(\d)
\end{align*}

By (iii) the set of slopes $\Phi_p$ is essentially independent from the choice of $p \in {\bf q} =B_{\varrho}(0)$ up to precision $\varrho$. That is, for some set $\Phi_{\bf q} \subset \varrho \ZZ\cap [-1,1]$ we have $N_\varrho \Phi_p \approx N_\varrho \Phi_{\bf q}$. Thus, we can write $\varphi_i = \overline{\varphi}_i+O(\varrho)$ for some $\overline{\varphi}_i \in \Phi_{\bf q}$. In this computation we can choose $y_1 \in Y({\bf q})$, $z \in Z({\bf q})$ and $\overline{\varphi}_i \in \Phi_{\bf q}$ essentially arbitrarily and so (\ref{eq:Yq-additive-relation}) follows (recall that we chose coordinates so that $z'=0$).

To obtain (\ref{eq:AB-additive-relation}) or (\ref{eq:xi-additive-relation}), we additionally write 
\[
\theta(T_i) = (\varphi_i, \xi_{p_i} + f(z_0) \varphi_i) = (\varphi_i, \xi_p+f(z)\varphi_i)
\] 
This implies that 
\[
\xi_{p_2} - \xi_{p_1} = (f(z)-f(z_0)) (\varphi_2-\varphi_1) + O(\d).
\] 
Since $p_2 = p_1 + z(\varphi_2-\varphi_1, 0, 0) + O(\varrho^2)$, $z_0=0$, and $\xi_p$ is Lipschitz in $p$, this gives the desired relation on $\xi(x_2)$ and $\xi(x_1)$.

Variants of the computations outlined in this part will appear many times throughout our arguments. We now turn to describe the additional ideas required for the sticky Furstenberg and higher dimensional sticky Kakeya.

\subsection{Sticky Furstenberg} In the Furstenberg setting, the set of tubes $\TT$ is $2t$-AD-regular for some $t \in (0,2]$ and the shading $Y(T)$ is assumed to be $(\d, s)$-Frostman. For simplicity, we will assume that $Y(T)$ is in fact $(\d, s)$-AD-regular, so that $Y(T)\cap B_\varrho(p)$ is a rescaled $(\d/\varrho,s)$-Frostman for every $p \in Y(T)$ and $\varrho \in (\d, 1)$. 
This assumption can be removed by carefully choosing the intermediate scale $\varrho$. 

After a sequence of rescaling and thickening moves similarly to the initial reduction above, we can parameterize the set $E$ as follows
\begin{equation}\label{eq:E-structure-furst}
E = \bigcup_{z \in Z} A_z \times \{z\}    
\end{equation}
where $Z \subset [-1,1]$ is a $(\d, s,\approx 1)$-set and $A_z \subset [-1,1]^2$ is a $(\d, 2t-\kappa, \approx 1)$-AD-regular set for every $z \in Z$. Furthermore, for every tube $T\in \TT$ we have $Y(T) \approx E \cap T$, i.e. $Y(T)$ is essentially the same as $Z$ after coordinate projection.
The numerology for the dimension of $A_z$ comes from a double counting identity: $|E| \approx |Z| |A_z| \approx \frac{|\TT| |Y(T)|}{|\TT_p|}$. 
This is similar to the rescaling step in the sticky Kakeya argument, except that slices $A_z$ do not necessarily have a global grains structure. 

Recall that one of the early steps in the sticky Kakeya proof was to reduce to the case when the set of tubes $\TT_p$ through a point $p\in E$ is coplanar. For the Furstenberg problem, this step is more complicated and we need to consider more general sets $\TT_p$. First, by the rescaling symmetry, we may assume that $\TT_p$ is $(\d,\kappa,\approx 1)$-AD-regular, for some $\kappa \ge 0$, and we choose $\kappa$ to be the largest possible exponent. Our goal is to prove an upper bound $\kappa \le \kappa_0(s, t) = \max(0,t-s, 2t-2)$. 

In Section \ref{sec:structure} we prove a useful structural result for almost AD-regular subsets in $[-1,1]^d$. By applying this result to the set of directions $\theta(\TT_p) \subset [-1,1]^2$, we can rescale and thicken the set of tubes and obtain the following structure. For every $p \in E$ the set of directions $\Theta_p=\theta(\TT_p)$ can be written as
\[
\Theta_p = \{ (\varphi, \xi+ \alpha_p \varphi), ~~ \xi \in \Xi_p, ~~ \varphi \in \Phi_{p, \xi} \},
\]
where $\alpha_p \in [-1,1]$ is a slope depending on $p$, and for some parameter $\gamma \in [0, \kappa]$, we have that $\Xi_p \subset [-1,1]$ is a $(\d,\kappa-\gamma,\approx 1)$-AD-regular set and $\Phi_{p, \xi} \subset [-1,1]$ is a $(\d,\gamma, \approx 1)$-AD-regular set. Furthermore, we have that $|\Theta_p \cap T_\varrho| \lesssim \varrho^{\kappa-\gamma} |\Theta_p|$ for any (planar) $\varrho$-tube $T \subset [-1,1]^2$ and any $\varrho \in [\d, 1]$ (and this inequality is almost achieved by tubes at slope $\alpha_p$ that cover $\Theta_p$). In other words, 
the set $\Theta_p$ is contained in a union of parallel tubes of maximum density and all underlying sets are almost AD-regular. The parameter $\gamma \in [0, \kappa]$ quantifies how coplanar the tubes in $\TT_p$ are and we need to consider several cases depending on the value of $\gamma$. 

When $\gamma$ is small, namely, if $\gamma \le s+\kappa-1$, then we obtain the desired estimate using the high-low method. By considering the intersection $E\cap B_{\varrho}(p)$ and rescaling, we obtain an orthogonal projections problem on scale $\varrho$, where the the set of projections is given by the set $N_\varrho\Theta_p$ (as before, this set does not depend on the particular choice of a point in the $\varrho$-ball). The condition that $\Theta$ does not intersect with tubes too much can be used to effectively estimate the overlap on the Fourier side and upper bound the `high' part. On the other hand, if `low' dominates then $E \cap B_\varrho(p)$ is locally a $(2+s)$-dimensional set. Comparing with (\ref{eq:E-structure-furst}) and recalling that $A_z$ is $(2t-\kappa)$-AD-regular, this leads to the desired bound on $\kappa$. We refer to Section \ref{sec:high-low} for a short self-contained argument. The upper bound $\gamma \le s+\kappa-1$ is precisely the right condition to make this estimate work. 

Next, we consider the intermediate case $\gamma \in [s+\kappa-1, 1-s]$. To simplify the notation a bit,  let us assume that for some fixed $p\in E$ we have $\Theta_p = \Phi_p \times \Xi_p$ (i.e. $\alpha_p=0$ and all sets $\Phi_{p, \xi}$ are the same) where $\Phi_p, \Xi_{p} \subset [-1,1]$ are $\gamma$- and $(\kappa-\gamma)$- almost-AD-regular sets, respectively. Let ${\bf q} =  B_\varrho(p)$, we can localize the representation (\ref{eq:E-structure-furst}) to ${\bf q}$ and write $E\cap {\bf q} = \bigcup_{z\in Z({\bf q})} A_z({\bf q}) \times \{z\}$. Since the sets $N_\varrho\Theta_{p'}$ are approximately the same for all $p' \in E\cap {\bf q}$, this gives 
\[
A_{z'}({\bf q}) \approx A_z({\bf q}) + (z'-z) (\varphi, \xi, 1) + O(\varrho^2)
\]
for many $z,z' \in Z({\bf q})$ and $\varphi \in \Phi_p$, $\xi \in \Xi_p$. Indeed, for each point in the $z$-slice, choose a tube through it pointing roughly in the $(\varphi,\xi,1)$ direction and move along this tube to level $z'$. 
By changing coordinates we may assume that $0\in \Phi_p \cap \Xi_p$, meaning that all sets $A_z({\bf q})$ are approximately the same up to resolution $\varrho^2$. Assume for a moment that the set $A_z({\bf q})$ is also a cartesian product $X({\bf q}) \times Y({\bf q})$ for some sets $X({\bf q})$ and $Y({\bf q})$. Then the above additive relation can be split into x- and y-coordinates:
\begin{align*}
    X({\bf q}) &\approx X({\bf q}) + (Z({\bf q}) - Z({\bf q})) \Phi_p \\
    Y({\bf q}) &\approx Y({\bf q}) + (Z({\bf q}) - Z({\bf q})) \Xi_p.
\end{align*}
We can now apply the ABC-sum-product theorem of Orponen and Shmerkin \cite{orponen2023projections} to conclude that
\begin{align*}
    \dim X({\bf q}) &\ge \min(1, \dim Z({\bf q}) + \dim \Phi_p) = \min(1, s + \gamma) = s+\gamma,\\
    \dim Y({\bf q}) & \ge \min(1, \dim Z({\bf q}) + \dim \Xi_p) = \min(1, s + (\kappa-\gamma)) = s+\kappa-\gamma,
\end{align*}
where we used the assumption that $\gamma$ is in the intermediate range $[s+\kappa-1, 1-s]$ to compute the minimums. By combining these estimates we conclude that $\dim A_z({\bf q}) \ge 2s+\kappa$. On the other hand, we know that $A_z$ is $(2t-\kappa)$-AD-regular and so we get $\kappa \le t-s\le \kappa_0(s,t)$, as desired. In this argument we used some cartesian product structure assumptions on the set of directions $\Theta_p$ and the sets $A_z({\bf q})$, these assumptions are not essential and can be removed. 

One caveat in this part of the argument is that we need to assume that both sets $\Phi_p$ and $\Xi_p$ have positive dimension. Indeed, if either of the sets consists of the single element 0, then the additive relation on $X({\bf q})$ or $Y({\bf q})$, respectively, does not give any new information. In case of $\Phi_p$ this is not really a problem: heuristically, by aligning the coordinate system we may arrange $\Phi_p$ to have at least two well-separated points, and so even if $\Phi_p$ is `zero-dimensional', we can still apply (trivial case of) ABC-sum-product. On the other hand, the case of $\Xi_p = \{0\}$ requires extra care. We will comment on this case below. 

It remains to consider the `almost coplanar' case $\gamma \in [1-s, 1]$. We can apply the same ABC-sum-product argument as above to the set $X({\bf q})$ to conclude that $\dim X({\bf q}) \ge \min(1, s +\gamma) =1 $. In other words, the set $A_z({\bf q})$ fully fills out horizontal lines (with the vertical offsets given by $Y({\bf q})$). This means that the set $E\cap {\bf q}$ has global grains, just as in the sticky Kakeya setting. By rescaling and thickening we can then reduce to the case when we can write 
\[
E = \{(x, y+f(z)x, z), z \in Z, y\in Y_z, x\in X_{y, z}\},
\]
where $Z \subset [-1,1]$ is an $s$-dimensional set, $Y_z \subset [-1,1]$ is $(2t-\kappa-1)$-AD-regular and $X_{y,z} \subset [-1,1]$ is $\approx 1$-dense. Furthermore, we have appropriate versions of properties (i), (ii), (iii), (iv) from Section \ref{subsec:sticky-kakeya-r3-outline}.
The same identity \eqref{eq:Yq-additive-relation} then holds in this setting:
\[
Y({\bf q}) \approx Y({\bf q}) + (z - z') (f(z) - f(z')) (\varphi_1-\varphi_2), \quad z,z' \in Z({\bf q}).
\]
Here the angles $\varphi_1,\varphi_2$ range over a $\gamma$-dimensional set $\Phi_{\bf q}$ (if we had cartesian products $\Theta_p = \Phi_p\times \Xi_p$ then we can take $\Phi_{\bf q} \approx N_\varrho \Phi_p$). If the slope function $f$ cannot be approximated by a linear function on $Z({\bf q})$ then the polarization argument and the applications of radial projection / Kaufman projection theorems imply that 
\[
Y({\bf q}) \approx Y({\bf q}) + S(\bf q) (\Phi_{\bf q} -\Phi_{\bf q}),
\]
where $S({\bf q})$ is an at least $s$-dimensional set produced from $Z({\bf q})$ and $f$ (the set of scalar products of difference sets of two graphs $\Gamma$ and $\Gamma'$). So by the ABC-sum-product theorem we obtain
\[
\dim Y({\bf q}) \ge \min(1, \dim S({\bf q}) + \dim \Phi_{\bf q}) \ge \min(1, s+\gamma)\ge 1,
\]
giving $\kappa\le 2t-2\le \kappa_0(s,t)$, as desired. In the case when $f$ is a locally approximated by a linear function with non-zero derivative, a direct application of (\ref{eq:Yq-additive-relation}) again leads to a desired bound. In the case when $f$ has `H\"older exponent' $1+\varepsilon$, we reach a contradiction with the Convex Wolff axiom, very similarly to the sticky Kakeya case in $\mathbb{R}^3$. 

A similar argument applies to the aforementioned edge case $\Xi = \{0\}$. In this case we have $\gamma = \kappa$, the ABC sum-product argument in the $X$-direction gives $\dim X({\bf q}) \ge \min(1, s+\kappa)$ and we need to deal with the case when the minimum is achieved on the second term. Here we can use a slightly simpler additive relation on $Y({\bf q})$:
\begin{equation}\label{eq:third-relation}
Y({\bf q}) \approx Y({\bf q}) + f(z) (x_2-x_1) - z (\xi(x_2)-\xi(x_1))    
\end{equation}
where $z \in Z({\bf q})$, $x_1, x_2 \in X({\bf q})$ and $\xi$ is the slope function defined in the second additive relation. We view the additive term as the scalar product between vectors $(f(z), z)$ and $(x, -\xi(x))$. Using the Furstenberg set estimate we show that the set of such dot products has dimension at least
\[
\frac{\dim X({\bf q}) + \dim Z({\bf q})}{2} \ge \frac{(s+\kappa) + s}{2}.
\]
On the other hand, $\dim Y({\bf q}) \le (2t-\kappa)-\dim X({\bf q})$, which leads to the bound $\kappa \le \frac{4}{5}(t-s) \le t-s$, as desired. The relation \eqref{eq:third-relation} is essentially the one used in the previous sticky Kakeya argument \cite{wang2026sticky}, \cite{tao2014}.

\subsection{Four dimensional Sticky Kakeya}
We state the problem in arbitrary dimension $d$ since the initial reduction and some of the arguments work for all dimensions, later we specialize to $d=4$. 

Let $\TT$ be a $(\d,d-1, \approx 1)$-AD-regular set of tubes in $\RR^d$ and let $Y(T)$ be a $\approx 1$-dense shading on $\TT$. Suppose that $\TT$ satisfies the Convex Wolff axiom. We start by applying rescaling and thickening to reduce to the case when $E = U(\TT,Y)$ has some special form.

We may assume that the set of tubes $\TT_p$ through a point $p \in E = \bigcup_{\TT} Y(T)$ is a $(\d,\kappa,\approx 1)$-AD-regular set, where $\kappa>0$ is the largest possible exponent. By the multilinear Kakeya argument, we may assume that $\TT_p$ is essentially contained in a common hyperplane. In dimension $d\ge 4$ it might happen that $\TT_p$ is contained in a smaller dimensional subspace. 
Thus, we need to split into several cases depending on the dimension of the span of $\TT_p$. So there exists some $\ell \in \{2, \ldots, d-1\}$ such that $\TT_p$ is contained in a $O(\d)$ neighborhood of an $\ell$-plane $\Pi_{p}^{(\ell)}$ and does not concentrate too much in any $\ell-1$-plane. By a multilinear Kakeya argument using \eqref{eq:parallel-tubelets}, $E$ locally fills out 
\[
\underbrace{\d\times \ldots \times \d}_{d-\ell} \times \underbrace{\d^{1/2} \times \ldots \times \d^{1/2}}_{\ell}
\]
parallel convex sets (previously we had this with $\ell=2$).
Rescaling turns these into $(\ell-1)$-dimensional horizontal global grains:
\begin{equation}\label{eq:kakeya-E-shape}
E = \{(x, y + f(z) x, z), ~ z \in Z, y\in Y_z, x\in X_{y, z}\}    
\end{equation}
where now $Z \subset [-1,1]$ is $\approx 1$-dense, $Y_z \subset [-1,1]^{d-\ell}$ is $(d-\ell-\kappa)$-AD-regular and $X_{y,z} \subset [-1,1]^{\ell-1}$ is $\approx 1$-dense. The `slope function' $f(z)$ is now a $(d-\ell)\times(\ell-1)$-matrix. 
There are two special cases that generalize fairly easily from the 3-dimensional sticky Kakeya argument in Section \ref{subsec:sticky-kakeya-r3-outline}. 
The first case is when the grains have codimension 1, i.e. $\ell = d-1$. Then the set $Y({\bf q}) \subset [-1,1]$ is $(1-\kappa)$-dimensional and $f(z)$ is a linear map from $\RR^{\ell-1}$ to $\RR$. A direct analogue of (\ref{eq:Yq-additive-relation}) gives 
\[
Y({\bf q}) \approx Y({\bf q}) + (z'-z)(f(z')-f(z)) (\Phi_{\bf q}-\Phi_{\bf q}).
\]
The important difference is that now $\Phi_{\bf q} \subset [-1,1]^{\ell-1}$ is a $\kappa$-AD-regular set which does not concentrate in any slab. Thus, for $\ell=d-1$, we get dot products between the vector $f(z')-f(z) \in \RR^{\ell-1}$ and $\Phi_{\bf q}$. If the function $f(z)$ is not approximated by a linear function on $Z({\bf q})$, then we can find a pair $\varphi_1, \varphi_2 \in \Phi_{\bf q}$ such that the function $f_{\varphi_1, \varphi_2}:z\mapsto f(z) \cdot (\varphi_1-\varphi_2)$ is also not approximated by a linear function (using that $\Phi_{\bf q}$ spans $\RR^\ell$ in a robust way). We can now use $f_{\varphi_1, \varphi_2}$ to proceed exactly as in the outline for the sticky Kakeya and conclude that $Y({\bf q})$ is locally 1-dimensional. The cases when $f$ is locally linear or constant can be dealt with in a similar manner. 

Let us point out a major difficulty in trying to extend this argument to the case when $\ell \le d-2$. The identity (\ref{eq:Yq-additive-relation}) is still true but it becomes more difficult to use. Indeed, the set $Y({\bf q})$ on the left hand side now $d-\ell-\kappa$ dimensional. On the other hand, $z \in Z$ ranges over a $1$-dimensional set and $\Phi_{\bf q}$ is $\kappa$-dimensional set. Thus, in order to obtain a contradiction, we need to show that a 1-dimensional set and a $\kappa$-dimensional set together produce a $\ge d-\ell-\kappa\ge 2-\kappa$ dimensional set. It is plausible that something like this should be true given a sufficiently strong non-linearity assumption on $f$ and some non-degeneracy condition of $\Phi_{\bf q}$. However, proving this seems difficult. 

Another nice special case is when $\ell=2$. (One can see that when $d=4$, $\ell=2$ and $\ell=d-1$ are exactly the two cases we need to consider.) Then we have $\Phi_p \subset [-1,1]$ and for each $p \in E$ and we can write $\Theta_p = \{(\varphi, \xi_p+f(z) \varphi), ~ \varphi\in \Phi_p\}$ for some $\xi_p \in [-1,1]^{d-2}$. By defining the function $\xi(x) = \xi_{(x, y_0 +f(z_0)x, z_0)}$ for some fixed $y_0,z_0$, we have a direct analogue of (\ref{eq:AB-additive-relation}):
\begin{equation}\label{eq:AB-relation-again}
A \approx A + (B-B) (\Phi_{\bf q} - \Phi_{\bf q}) 
\end{equation}
where 
\[
A = \left\{ \begin{pmatrix}
    x \\
    \xi(x)
\end{pmatrix}, ~x\in X({\bf q})\right\} \quad 
\text{and}\quad
B = \left\{\begin{pmatrix}
    z \\
    -f(z)
\end{pmatrix}, ~z \in Z({\bf q})\right\}.
\] 
Note that as before $\dim A = \dim B = 1$. By taking a generic linear projection $\RR^{d-1}\to \RR^2$ we may reduce to the case when $A, B \subset [-1,1]^2$ are 1-dimensional sets satisfying (\ref{eq:AB-relation-again}). This is exactly the same setup as in the case of the sticky Kakeya in $\RR^3$ and so the same argument applies. 

From this argument it follows that $B$ must be concentrated in a tube, i.e. the function $f$ is essentially a locally linear function (one can verify that this property is preserved under generic projection to $\RR^2$). 
By using a version of the identity (\ref{eq:Yq-quadratic-shift}), we can use this to conclude that 
\begin{equation}\label{eq:Yq-global-grain}
Y({\bf q}) \approx Y({\bf q}) + v \cdot [-\varrho^2, \varrho^2]    
\end{equation}
for some non-zero vector $v \in \RR^{d-2}$ (which depends on ${\bf q}$). By rescaling into a $\varrho^2$-tube, this implies that the set $E$ fills out out 2-dimensional horizontal global grains. This is not a contradiction but it does give us extra information. It turns out that this is a tricky case that requires additional ideas and structures, which we discuss in the next paragraph.

Let us now specialize to $d=4$ and summarize what we have shown so far. First, we proved that if $\TT_p$ is 3-linear (i.e. robustly spans a 3-plane), then $\dim E = 4$ and $\kappa=0$. Then, assuming that $\TT_p$ is contained in a 2-plane, we showed that $E$ must have 2-dimensional global grains (as opposed to the 1-dimensional ones which we get `for free'). Let $\tilde f(z) = (\tilde f_1(z), \tilde f_2(z))$ denote a new slope function parameterizing the 2-grains, then we can write 
\begin{equation}\label{eq:kakeya-E-tilde}
E = \{ (\tilde x, \tilde y + \tilde f(z)\tilde x, z), ~ z \in Z, ~\tilde y \in \tilde Y_z, ~\tilde x \in \tilde X_{\tilde y, z} \},    
\end{equation}
where $\tilde Y_z \subset [-1,1]$ is $(1-\kappa)$-AD-regular and $\tilde X_{\tilde y,z} \subset [-1,1]^2$ is $\approx 1$-dense. This is similar to the representation (\ref{eq:E-structure-kakeya-r3}) with $\ell=3$ except that the grains are no longer spanned by the directions of tubes in $\TT_p$. The 1-grains parameterized by $f$ and the 2-grains parameterized by $\tilde f$ are `compatible' with each other. For example for $\theta =\theta(T)\in \Theta_p$ we can write
\begin{align*}
\theta =& (\varphi, \xi_p + f(z) \varphi) = (\varphi, \xi_1 +f_1(z)\varphi, \xi_2 + f_2(z) \varphi) \\
=& (\varphi, \xi_p+f_1(z)\varphi, \tilde \xi_p + \tilde f(z) \cdot (\varphi, \xi_p+f_1(z)\varphi) ).
\end{align*}
By comparing the last coordinates we get the following relation:
\[
f_2(z) = \tilde f_1(z) + \tilde f_2(z) f_1(z).
\]
We have another relation which follows from the $\ell=2$ argument above. Recall that we showed that the set $B = \{\begin{pmatrix}
    z \\
    -f(z)
\end{pmatrix}\}$ is contained in a tube when $z \in Z({\bf q})$. This roughly means that $f(z)$ is a $C^1$-differentiable function. The vector $v$ in (\ref{eq:Yq-global-grain}) is simply the derivative $f'(z_0)$. By comparing the two different grain representations it follows that
\[
\left\langle 
\begin{aligned}
(1, f_1(z), f_2(z)) \\(0, f'_1(z), f'_2(z))    
\end{aligned}
\right\rangle = \left\langle
\begin{aligned}
    (1, 0, \tilde f_1(z))\\ (0, 1, \tilde f_2(z))
\end{aligned}\right\rangle 
\]
and so $\tilde f_2(z) = f'_2(z) / f'_1(z)$ and $\tilde f_1(z) = f_2(z) - f_1(z) f'_2(z) / f'_1(z)$. 

We use this information to analyze the set $E\cap {\bf q}$. By tracking down how tubes interact with the structure of $E$, we obtain a somewhat more involved version of the key identity (\ref{eq:Yq-additive-relation}) on $Y({\bf q})$ which is as follows:
\begin{equation}\label{eq:Yq-complicated}
    \tilde Y({\bf q}) \approx \tilde Y({\bf q}) + (z-z') (f_2(z) - f_2(z') - \tilde f_2(z') (f_1(z) - f_1(z'))) (\Phi_{\bf q}-\Phi_{\bf q}),
\end{equation}
where $\tilde Y({\bf q})$ is a thickened piece of $\tilde Y_{z'}$ in a $\varrho^3$-interval. 
Let us briefly explain where this formula comes from. By writing $E$ in the form \eqref{eq:kakeya-E-shape} and using \eqref{eq:Yq-additive-relation} again, we have
\begin{equation}\label{eq:Yq-symmetry-again}
    Y({\bf q}) \approx Y({\bf q}) + (z-z') (f(z) - f(z'))  (\varphi_1-\varphi_2),    
\end{equation}
on the other hand, comparing the representations \eqref{eq:kakeya-E-shape} and \eqref{eq:kakeya-E-tilde} we can write 
\[
(x, y_1 + f_1(z') x, y_2 + f_2(z')) = (\tilde x_1, \tilde x_2, \tilde y + \tilde f_1(z') \tilde x_1+\tilde f_2(z') \tilde x_2)
\]
and putting $x=\tilde x_1=0$ gives representation
\[
Y_{z'} = \{  (\tilde x_2, \tilde y + \tilde f_2(z') \tilde x_2 ), ~ \tilde y \in \tilde Y_{z'}  \},
\]
where $\tilde x_2$ ranges over a dense subset in $[-1,1]$. 
Ignoring the slight difference between $Y_{z'}$ and $Y({\bf q})$, we can apply an orthogonal projection along the vector $(1,\tilde f_2(z'))$ (i.e. take dot product with $(-\tilde f_2(z'), 1)$) to both sides of \eqref{eq:Yq-symmetry-again} and arrive at \eqref{eq:Yq-complicated}.

Now our task is to use (\ref{eq:Yq-complicated}) to show that $Y({\bf q})$ is locally 1-dimensional, or to arrive at a contradiction with the Convex Wolff axiom. The analysis follows roughly the same route as before: we use a polarization identity and then use radial and Kaufman projection theorems. Computations become more complicated and we need to use the `$C^1$-smoothness' assumption on $f$.  

Assuming the normalization $f'(z_0) = (1,0)$, the polarization step ultimately shows that 
\[
Y({\bf q}) \approx Y({\bf q}) + (V(z_1) - V(z_2)) \cdot (U(z_1') - U(z_2')) (\Phi_{\bf q} - \Phi_{\bf q})
\]
where $V(z), U(z') \in \RR^3$ are 3-dimensional vectors given by
\[
V(z) = \begin{pmatrix}
    z \\
    z^2 \\
    f_2(z)
\end{pmatrix}, \quad U(z') = \begin{pmatrix}
    f_2(z') - 2 z' \tilde f_2(z') \\
    \tilde f_2(z') \\
    z'
\end{pmatrix}.
\]
Thus, to finish the proof, it suffices to show that the set of dot products $V(z)\cdot U(z')$ spans a 1-dimensional subset in $[-1,1]$. Note that $\{V(z_1)-V(z_2), ~z_1, z_2\in Z({\bf q})\}$ is at least 2-dimensional set (consider the first two coordinates) and unless $f_2$ is close to a quadratic polynomial, this set is not concentrated in any slab. So a 3-dimensional radial projections theorem implies that the set $\{V(z_1)-V(z_2)\}$ defines a 2-dimensional set of slopes. Since $\{U(z)\}$ is a 1-dimensional set (consider the last coordinate), a 3D Kaufman projection theorem implies that the dot product of $\{U(z)\}$ with a typical slope of $\{V(z_1)-V(z_2)\}$ is 1-dimensional. This shows that $Y({\bf q})$ is locally 1-dimensional and so $\kappa=0$, as desired. 
The case when $f_2$ is locally a quadratic function can be dealt in an analogous way to how we dealt with linear cases. 



    

\section{Notation}\label{sec:notation}

In this section, we recall some definitions and lemmas from \cite{orponen2023projections}.

The letter $C$ will usually denote a constant, while the letter $K$ will usually denote a parameter depending on $\delta$ (e.g. $K = C \delta^{-\eps}$).

We say $A \lesssim B$ if $A \le CB$ for some constant $C$, and define $A \gtrsim B$, $A \sim B$ similarly. We write $A \lesssim_n B$ to emphasize the implicit constant may depend on $n$. We also use $A \gtrsim_{\delta^{O(\eta)}} B$ to denote that $A \geq C_1 \delta^{C_2\eta} B$ for some constant $C_1, C_2>0$ and $A\approx_{\delta^{O(\eta)}} B$ to denote $C_1^{-1}\delta^{C_2\eta} B\leq A \leq C_1 \delta^{-C_2\eta} B$.

For a finite set $A$, $|A|$ denotes cardinality of $A$. For a measurable and infinite set $A$ (usually a ball or a rectangle), $|A|$ denotes the Lebesgue measure of $A$.

\begin{definition}\label{def:deltas-set}
    We say that a subset $A$ of $\mathbb{R}^d$ is $(\d, t, C)$-set if the following holds: for every $r \ge \d$ and $r$-ball $B_r$ we have
    \[
    |A \cap B_r| \le C r^t |A|.
    \]
    Additionally, we say that $A$ is {\em rescaled $(\d, s, C)$-set} if there is a homothety $\phi:\RR^d \to \RR^d$ such that $\phi(A) \subset B_1$ and $\phi(A)$ is $(\d, s, C)$-set. 
\end{definition}

\begin{definition}
    We say that a subset $A$ of $\mathbb{R}^d$ 
    is $(\d, t, C)$-Katz-Tao if the following holds: for every $r\ge \d$ and $r$-ball $B_r$ we have
    \[
    |A \cap B_r|_\d \le C (r/\d)^t. 
    \]
\end{definition}

\begin{definition}\label{def: deltarhoADregularset}
    We say that a subset $A$  of $\mathbb{R}^d$ is $([\d, \varrho], t, C)$-AD-regular if the following holds:
    for any $a \in A$ and $r \in [\d, \varrho]$ we have 
    \[
    |A \cap B_r(a)|_\d \in [C^{-1} (r/\d)^t, C (r/\d)^t].
    \]
\end{definition}

Typically, we will have $\varrho=1$ in which case we say that $A$ is $(\d, t, C)$-AD-regular.

For any compact set $A\subset \mathbb{R}^d$ and $r > 0$, let $N_r A:= \{ x \in \mathbb{R}^d : d(x, A) \le r \}$ denote the Euclidean $r$-neighborhood of $A$. 

\begin{definition}[Dyadic cubes]
    If $n\in \mathbb{Z}$ and $\delta=2^{-n}$, denote by $\cD_{\delta}(\mathbb{R}^d) := \{ [a\delta, (a+1)\delta) \times [b\delta, (b+1)\delta) : a, b \in \mathbb{Z} \}$ the family of dyadic cubes  in $\mathbb{R}^d$ of side-length $\delta$. If $\cP\subset \mathbb{R}^d$, denote 
    \[
    \cD_{\delta}(\cP):=\{ Q\in \cD_{\delta} (\RR^d): Q\cap \cP\neq \emptyset\}.
    \]
    We abbreviate $\cD_{\delta}:=\cD_{\delta}([0,1)^d)$ for $n\geq 0$.  If $\delta <\Delta\in 2^{-\mathbb{N}}$ and $p\in \cD_{\delta}$, let $p^{\Delta}$ denote the dyadic cube $\bp \in \cD_{\Delta}$ containing $p$. 
\end{definition}

For a compact set $A \subset \RR^d$ and $r = 2^{-n}$ we denote $A^{(r)} = \bigcup_{a\in A} a^{(r)} = \bigcup \mc D_r(A)$ the union of dyadic $r$-cubes covering $A$. We typically omit the brackets and write $A^r = A^{(r)}$ for the dyadic neighborhood of $A$ (apart from cases in which it might lead to confusion). 

\begin{definition}\label{defn:covering number}
    For a bounded set $P \subset \RR^d$ and $\delta = 2^{-n}$, $n \ge 0$, define the dyadic $\delta$-covering number
    \begin{equation*}
        |P|_\delta := |\cD_\delta (P)|.
    \end{equation*}
\end{definition}
The number $|P|_{\delta}$ is comparable (up to a universal constant factor) to the least number of $\delta$-balls needed to cover $P$. (This latter definition is the more common definition of covering number.) If $\cP\subset \cD_{\delta}$ is a set of dyadic cubes, we say that $\cP$ is a $(\delta, s, C)$-set if the union  of $\delta$-cubes in $\cP$ is a $(\delta, s, C)$-set in the sense of Definition \ref{def:deltas-set} and we use $|\cP|$ to denote the number of dyadic cubes in $\cP$ and $|\cP|_{\Delta}:=|\cD_{\Delta}(P)|$ (hence $|\cP|_{\delta}=|\cP|$). 

For different parts of the proof it will be sometimes convenient to use slightly different definitions of a tube. 
In most places, we define for $(a, b) \in [-1,1]^{d-1}\times \RR^{d-1}$ a $\d$-tube $T$ as
\begin{equation}\label{eq:usual-tube}
T = \{ (a t + b + x, t) , ~ t \in [-1,1], ~ x \in B_\d(0)  \}.    
\end{equation}
Note that such $T$ makes angle $\lesssim 1$ with the vertical axes. 
For the purposes of rescaling and thickening we use the following notion of dyadic $\d$-tubes.

\begin{definition}[Dyadic $\delta$-tubes]\label{def: dyadictube}
    Let $\delta\in 2^{-\mathbb{N}}$. A \emph{dyadic $\delta$-tube} is a set of the form $T=\cup_{x\in p} \mathbf{D}(x)$, where $p\in \cD_{\delta}([-1,1]^{d-1}\times \mathbb{R}^{d-1}),$ and $\mathbf{D}$ is the \emph{point-line duality map}
\[
\mathbf{D}(a,b):=\{ (y,x)\in \mathbb{R}^{d-1}\times \mathbb{R}: y=ax+b\}\subset \mathbb{R}^d 
\]
sending the point $(a,b)\in \mathbb{R}^{2(d-1)}$ to a corresponding line in $\mathbb{R}^d$. Abusing the notation, we abbreviate $\mathbf{D}(p):=\cup_{x\in p} \mathbf{D}(x).$ The collection of all dyadic $\delta$-tubes is denoted 
\[
\TT^{\delta}:=\{ \mathbf{D}(p): p\in \mathcal{D}_{\delta}([-1,1]^{d-1}\times \mathbb{R}^{d-1})\}.
\]
A finite collection of dyadic $\delta$-tubes $\{\mathbf{D}(p)\}_{p\in \cP}$ is called a $(\delta, s, C)$-set if $\cP$ is a $(\delta, s, C)$-set. 

If $\delta<\Delta\in 2^{-\mathbb{N}}$ and $T\in \TT^{\delta}$, let $T^{\Delta}$ denote the unique dyadic tube $\bT\in \TT^{\Delta}$ containing $T$. 
\end{definition}

Note that for $(a, b) \in [-1,1]^{2(d-1)}$ the dyadic $\d$-tube ${\bf D}(p) \cap [-1,1]^d$, where $p = (a,b)^\d \in \mathcal D_\d$, defines an approximately the same set as the $\d$-tube $T$ defined in (\ref{eq:usual-tube}). We can always switch between normal and dyadic tubes and only lose constant factors. 

One drawback of dyadic tubes is that for $T\in \TT^\d$ it might not be the case that the Euclidean neighborhood $N_\d T$ is contained in a `thicker' tube $T^{\varrho}$ for some $\varrho>\d$. This happens if $T$ is very close to the boundary of the bigger dyadic tube. The following proposition is aimed at mitigating this issue by randomly translating the set of dyadic tubes and thickening them slightly.

\begin{prop}\label{prop:regularize-dyadic-tubes}
    Let $\TT \subset \TT^\d$ and suppose that each $T\in \TT$ is of the form ${\bf D}(p)$ for some $p \in \cD_\d([-1,1]^{2(d-1)})$.
    Let $S \ge 2$ and suppose that $\d \in 2^{S\ZZ}$ and $S \ge C\log \log (1/\d)$. Then there exists $\TT' \subset \TT^{2^S\d}$ and a vector $v \in \RR^d$ such that:
    \begin{itemize}
        \item the number of $T \in \TT$ such that $T+v \subset T'$ for some $T'\in \TT'$ is at least $|\TT|/2$,


        \item for every $\varrho \in [2^S\d,2^{-S}] \cap 2^{\ZZ}$ and $T \in \TT'_\varrho$ we have
        \[
        N_\varrho T  \subset T^{2^S \varrho}.
        \]
    \end{itemize}

\end{prop}

\begin{proof}
    For $\tau>\varrho$ and a $\tau$-tube $T_{\tau}$ let $\TT^\varrho(T_{\tau})$ be the set of tubes $T_\varrho$ such that $N_\varrho T_\varrho \subset T_{\tau}$. It is not hard to see that then we have 
    \[
    |\TT^\varrho(T_{2^S\varrho})| \ge (1-C_d2^{-S}) |\TT^\varrho[T_{2^S\varrho}]|,
    \]
    i.e. most of the $\varrho$-tubes contained in $T_{2^S\varrho}$ are $\varrho$-separated from the boundary. 

    Let $v \in [-1,1]^d$ be a uniformly random vector. One can check that for any fixed $T\in \TT$, with probability at least $1-C2^{-S}$, there exists $\tilde T \in \TT^{2^S\d}$ such that $T+v \subset \tilde T$. 
    Define $\tilde \TT(v) \subset \TT^{2^S \d}$ to be the set of tubes $\tilde T$ such that $T+v \subset \tilde T$ for some $T\in \TT$. So, on average, $\tilde \TT(v)$ covers a large fraction of $\TT$. 

    Furthermore, one can verify that for any fixed $T\in \TT$, the $2^S\d$-tube $\tilde T \in \tilde \TT(v)$ covering $T$ has the following property  with probability at least $1- C 2^{-S} \log(1/\d)$: for any $\varrho \in [2^S\d, 2^{-S}]$ we have $N_\varrho \tilde T^\varrho \subset \tilde T^{2^S \varrho}$. Since $S \ge C \log \log (1/\d)$, this probability is at least $1/2$. So using linearity of expectation, we can fix a value of $v$ and define $\TT'$ to be the set of $\tilde T \in \tilde \TT(v)$ such that the property above holds. 
\end{proof}

For a set of dyadic $\delta$-tubes $\TT$ and $\delta< \Delta\in 2^{-\mathbb{N}}$,  define the dyadic $\Delta$-covering number $|\TT|_{\Delta}:=|\cD_{\Delta}(\TT)|$ where $\cD_{\Delta}(\TT):=\{ \bT\in \TT^{\Delta}: \exists T\in \TT, \bT=T^{\Delta}\}$. When $\delta=\Delta$, we clearly have $|\TT|_{\delta} = |\TT|$.

\begin{definition}[Slope set]
    The slope of a line $\ell=\mathbf{D}(a,b)$ is defined to be the vector $a\in \mathbb{R}^{d-1}$, write $\dir(\ell):=a.$ 

    If $T=\mathbf{D}(p)$ is a dyadic $\delta$-tube, we define the slope set $\dir(T) = \cup_{x\in p} \dir ( \mathbf{D}(x))$. Thus $\dir(T)$ is a cube of sidelength $\delta$. If $\TT$ is a collection of dyadic $\delta$-tubes, we write $\dir(\TT):=\{\dir(T): T\in \TT\}$. 

    Moreover, writing $p = (a,b) + [0,\d)^{2(d-1)}$, we denote $\theta(T) = a$. In other words, $\theta(T)$ is the unique element in the intersection $\dir(T) \cap \d \ZZ^{d-1}$. 
    
\end{definition}


\begin{definition}
    Let $\TT$ be a set of $\delta$-tubes. We say that $Y$ is a $(\delta, s, C)$-dense shading on $\TT$ if for any $T\in \TT$,   $Y(T)$ is a union of dyadic $\delta$-cubes intersecting $T$ and for any $\varrho\in [\delta, 1]$ and any $B_\varrho$, 
    \[
    |Y(T)\cap B_\varrho|_{\delta} \leq C\varrho^s |Y(T)|_{\delta}. 
    \]
\end{definition}
If $\TT$ is a set of $\delta$-tubes and $Y$ is a shading on $\TT$, write 
\[
U(\TT, Y)=\bigcup_{T\in \TT} Y(T). 
\]
For an arbitrary set $U \subset \RR^d$ we define
\[
\TT[U ] = \{T \in \TT:~ T \subset U\}.
\]

\begin{definition}[$(t,K)$-AD regular set of tubes.]
Let $\TT$ be a set of dyadic $\d$-tubes. We say that $\TT$ is $(t, K)$-AD-regular if for any $\varrho \in [\delta,1]\cap 2^{-\mathbb{N}}$ and any $T\in \TT$, 
\begin{equation}\label{eq: ADregular}
K^{-1} (\delta/\varrho)^{-t}\leq |\TT[T^\varrho]|_{\delta} \leq K (\delta/\varrho)^{-t}.
\end{equation}
\end{definition}

For a collection of usual (non-dyadic) $\d$-tubes $\TT$ (of the form (\ref{eq:usual-tube})) we can similarly say that $\TT$ is $(t, K)$-AD-regular if (\ref{eq: ADregular}) holds with $T^\varrho = N_\varrho T$. 
In applications we usually take $K=\delta^{-\eta}$ for some very small parameter $\eta>0$.


Let $p$ be a $\varrho$-cube of the form $(a,b)+ [0,\varrho)^{2(d-1)}$. 
For a $\varrho$-tube $T=\mathbf{D}(p)$, we can define the anisotropic rescaling map $\psi_T:\RR^d\to \RR^d$ by
\[
\psi_T( (y,x) ) = \left( \frac{y - b-a x}{\varrho}  ,x \right).
\]
It follows that for any $T' \in \TT^\d$ we have $\psi_T(T') \in \TT^{\d/\varrho}$, i.e. dyadic tubes rescale into dyadic tubes. 

\begin{obs}[Rescaling]\label{obs:rescaling-AD-regular-tubes}
    Suppose that $\TT$ is $(t, K)$-AD-regular set of $\d$-tubes. Let $\varrho \in [\d,1]\cap 2^{-\mathbb{N}}$. Then for any $T\in \TT$ the set of $(\d/\varrho)$-tubes 
    \[
    \TT^{T^{\varrho}}:= \{ \psi_{T^\varrho} (T'), ~ T' \in \TT[T^\varrho]\}
    \]
    is $(t, K)$-AD-regular.
\end{obs}

\begin{obs}[Thickening]\label{obs:thickening-AD-regular-tubes}
    Suppose that $\TT$ is $(t, K)$-AD-regular set of $\d$-tubes. Let $\varrho \in [\d,1]\cap 2^{-\mathbb{N}}$. 
    Then $\TT_\varrho:= \cD_{\varrho}(\TT)$ is a $(t, K^2)$-AD-regular set of $\varrho$-tubes. 
\end{obs}

\begin{proof}
First we note that $\TT[T_\varrho] \cap \TT[T'_\varrho] = \emptyset$ for any distinct $T_{\varrho}, T_{\varrho}'\in \TT_\varrho$, and  $\TT = \bigsqcup_{T_\varrho\in  \TT_\varrho} \TT[T_\varrho]$.


Let $\tau \in (\varrho, 1]\cap 2^{-\mathbb{N}}$ and $T_\varrho\in \TT_{\varrho}$. Let $T_\tau= T_{\varrho}^\tau$, the dyadic $\tau$ tube in $\TT^{\tau}$ (see Definition~\ref{def: dyadictube})  containing $T_\varrho$. 
Then we can compute
\[
|\TT[T_\tau]| = \sum_{T_{\varrho}' \in \TT_{\varrho}} |\TT[T_\tau] \cap \TT[T_\varrho']| = \sum_{T_\varrho' \in\TT_\varrho[T_\tau]} |\TT[T_\tau] \cap \TT[T_\varrho']| = \sum_{T_\varrho' \in\TT_\varrho[T_\tau]} | \TT[T_\varrho']|
\] 
and so using \eqref{eq: ADregular} on both sides we get
\[
|\TT_\varrho[T_\tau]| \in [K^{-2} (\tau/\varrho)^{t}, K^2 (\tau/\varrho)^t].
\]
Therefore $\TT_\varrho$ is $(t, K^2)$-AD-regular.
\end{proof}
Let $\TT\subset \TT^{\delta}$ be a set of  $\delta$-tubes and $Y$ be a shading on $\TT$.   For any $T_\varrho\in \TT_\varrho$ in Observation~\ref{obs:thickening-AD-regular-tubes}, we can define shading on the thick tube $T_\varrho$ as 
\begin{equation}\label{eq: shadingTrho} Y(T_\varrho)= \cD_\varrho(U(\TT[T_\varrho], Y)).
\end{equation} 

By combining the two operations above for any  $\varrho< \tau \in (\d, 1]\cap 2^{-\mathbb{N}}$ we can define a collection $\TT_\varrho^{T_\tau}$ of $\varrho/\tau$-tubes obtained by rescaling $\TT[T_\tau]$ and thickening to $T_\varrho$. It follows that any such $\TT_{\varrho}^{T_\tau}$ is $(t, K^2)$-AD-regular set of $\varrho/\tau$-tubes. 

When we start with a collection of usual $\d$-tubes $\TT$, then we can first replace it by a comparable collection of dyadic $\d$-tubes and then use the construction above. 

\begin{definition}\label{def:convex-wolff-axiom}
    Let $\TT$ be a set of $\d$-tubes in $\RR^{3}$. For numbers $t, t' \in [0,2]$ we define the Convex Wolff constant with exponents $(t, t')$ as
    \[
    C_{t; t'-CW}(\TT) = \max_{U:~ a\times b\times 2\text{-box}}  \,\, \, \frac{|\TT[U]|}{ a^{2t} (b/a)^{t'} |\TT|}
    \]
    where the maximum is over all $a\le b$ and all $a\times b\times 2$ boxes $U \subset \RR^3$. 

    Furthermore, when $t=t'$ we write $C_{t-CW}(\TT) = C_{t;t-CW}(\TT)$
\end{definition}

\begin{obs}\label{obs:wolff-axiom-rescaling}
    Let $\TT$ be a $(2t, K)$-AD-regular set of $\d$-tubes. For $\varrho<\tau \in [\d,1] \cap 2^{-\mathbb{N}}$,  let $\tilde \TT = \psi_{T^{\tau}}(\TT_{\varrho}[T^\tau])$ for any $T\in \TT$ and $\TT_\varrho$ defined in Observation \ref{obs:thickening-AD-regular-tubes}. Then for any $t'\in [0,2]$ we have
    \[
    C_{t; t'-CW}(\tilde \TT) \lesssim K^4 C_{t;t'-CW}(\TT).
    \]
\end{obs}

\begin{proof}
By the definition of $(2t, K)$-AD-regular set, 
\begin{equation}
    |\TT|\in [K^{-1}\delta^{-2t}, K\delta^{-2t}]
\end{equation}
and 
\begin{equation}\label{eq: ADTrho}
|\TT[T_\varrho]|\in [K^{-1}(\varrho/\delta)^{2t}, K(\varrho/\delta)^{2t}]
\end{equation}
for each $T_\varrho\in \TT_\varrho$. 
By Observation~\ref{obs:thickening-AD-regular-tubes}, 
\[
|\TT_\varrho[T^\tau]|\in [K^{-2}(\tau/\varrho)^{2t}, K^2(\tau/\varrho)^{2t}],
\]
and so $|\TT|\leq K^3 (\frac{\varrho}{\tau\delta})^{2t} |\TT_\varrho[T^\tau]|.$

For any $a\times b\times 2$-box $U$, $\phi_{T^\tau}^{-1}(U)$ is contained in a $C_0a\tau \times C_0b\tau \times C_0$-box $U'$ where $C_0$ is an absolute constant and 
\[
|\TT[U']|\lesssim C_{t; t'-CW}(\TT) (a\tau)^{2t} (b/a)^{t'} |\TT|. 
\]
Now
\[
|\TT[U']|\ge\sum_{T_\varrho\in \TT_\varrho[U']} |\TT[T_\varrho]|, 
\]
by \eqref{eq: ADTrho}, 
\begin{align*}
|\tilde\TT[U]| & = |\TT_\varrho[U']|\leq K(\delta/\varrho)^{2t}|\TT[U']|\\
& \lesssim K(\delta/\varrho)^{2t} C_{t; t'-CW}(\TT) (a\tau)^{2t} (b/a)^{t'} |\TT| \\
& \leq K(\delta/\varrho)^{2t} C_{t; t'-CW}(\TT) (a\tau)^{2t} (b/a)^{t'} K^3 (\frac{\varrho}{\tau\delta})^{2t} |\TT_\varrho[T^\tau]|\\
&= K^4 C_{t; t'-CW}(\TT) a^{2t}(b/a)^{t'} |\tilde\TT|
\end{align*}
since $|\tilde\TT|=|\TT_\varrho[T^\tau]|$. 
\end{proof}

We consider the natural choice $t'=t$, so that the bound on $\TT[U]$ is in terms of $|U|^t$. It was conjectured in \cite{wang2024restriction} that if $\TT$ satisfies the Convex Wolff axiom with exponents $(t,t)$ then it satisfies the Furstenberg set estimate. 



\section{Uniform structure}\label{sec:uniform}

In this section we prove some dyadic versions of pigeonholing from \cite{cohen2025lower}.
Let $\Omega_d$ be the set of triples $(x, t, \theta) \in [-1,1]^{d-1}\times [-1,1] \times [-1, 1]^{d-1}$. We think of $\Omega_d$ as the set of point-line pairs $(p\in \ell)$ in $\RR^d$ where $p = (x, t)$ and $\ell = (x,t) + \RR(\theta,1)$. 
Note that not every point-line pair is of this form. The following proposition allows us to reduce to this form by passing to a large subset of point-line pairs. 


\begin{obs}\label{obs:pigeonhole-slope}
    There exists a set of rigid motions $g_1, \ldots, g_{k}$ of $\RR^d$ with $k = O_d(1)$ so that for any $p \in [-1,1]^d$ and any line $\ell$ through $p$ there exists $j$ such that $g_j(p) \in [-1,1]^d$ and $g_j(\ell)$ has the form $p + \RR(\theta, 1)$ for some $\theta \in [-1,1]^{d-1}$.

    As a corollary, for any collection of point-line pairs $(p_i\in \ell_i)_{i=1}^N$ with $p_i\in [-1,1]^d$ there exists $j$ and some $X \subset \Omega_d$ such that $\#\{i:~ \exists (x, t, \theta)\in X:~~ p_i = g_j(x, t), ~ \ell_i = g_j( p_i+\RR(\theta,1) )\} \ge N/k$.
\end{obs}    


For each $(x,t,\theta) \in \Omega_d$ let $T^\varrho(\omega) = {\bf D}(p )$, $p = (\theta,x-t\theta)^\varrho$, be the dyadic $\varrho$-tube containing the line $(x,t) + \RR (\theta,1)$. 
For dyadic scales $v, h, \alpha$ such that $h \ge v \alpha$ and $\omega_0=(x_0, t_0, \theta_0)\in \Omega_d$ define
\[
R_{v,h,\alpha}(\omega_0) = \{ \omega=(x, t, \theta) \in\Omega_d:~ t \in t_0^{(v)}, ~ (x,t) \in T^{h}(\omega_0), ~ \theta \in \theta_0^{(\alpha)} \}
\]
We view $R=R_{v,h,\alpha}(\omega_0)$ as a `rectangle' in the parameter space $\Omega_d$. Geometrically, it describes the following set of point-line pairs: fix a point-line pair $p_0 \in \ell_0$ and consider a $\underbrace{h\times \ldots \times h}_{d-1} \times v$-tube segment $T_{h,v}=T_{h,v}(p_0,\ell_0)$ defined by this pair. Then $R$ consists of point-line pairs $p\in \ell$ such that $p \in T_{h,v}$ and $\angle (\ell,\ell_0) \le \alpha$.
Note that if $h\ge v\alpha$ then this defines a roughly symmetric relation on the space of point-line pairs.
For example, if $h =v\alpha$ then the condition is essentially that the tube segments $T_{h,v}(p_0,\ell_0)$ and $T_{h, v}(p,\ell)$ are essentially the same. 
The variables $v$, $h$, $\alpha$ stand for vertical, horizontal and angular distance, respectively. 
As in \cite{cohen2025lower}, we can consider a net $\Omega^{v,h,\alpha}_d$ of $\omega_0$ and define a finitely overlapping collection of rectangles $\mc R_{v,h,\alpha}$ covering $\Omega_d$. 
For a configuration $X \subset \Omega_d$ we can then define a \emph{concentration number} 
\[
\M_X(v,h,\alpha) = \max_{R\in \mc R_{v,h,\alpha}} |R \cap X|.
\]

For parameters $K>0$ and $\Delta_1>\ldots>\Delta_m$, let us say that a configuration $X \subset \Omega_d$ is $K$-uniform on a sequence of scales $\Delta_1 > \ldots > \Delta_m$ if for all $i, j, k \in \{1, \ldots, m\}$ such that $\Delta_j \ge \Delta_i \Delta_k$ and for any $\omega_0=(p_0, \ell_0)\in X$ we have
\[
|X \cap R_{\Delta_i, \Delta_j, \Delta_k}(\omega_0)| \ge \frac{1}{K} \M_X(\Delta_i, \Delta_j, \Delta_k).
\]

The following lemma was proven (with a slightly different but essentially equivalent definition of uniformity) in \cite{cohen2025lower} in the two-dimensional case. The higher dimensional version is completely analogous, so we skip the proof.

\begin{lemma}\label{lem:uniform-subset}
    Let $d\ge 2$. For any $\d>0$ and $m \ge 1$ there exists some $K \le (C_d \log 1/\d)^{10m^3}$ such that the following holds. 
    For any $X \subset \Omega_d$ and any scales $\Delta_1> \ldots > \Delta_m \ge \d$ there exists $X' \subset X$ of size at least $K^{-1} |X|$ so that $X'$ is $K$-uniform on scales $\Delta_1, \ldots, \Delta_m$.
\end{lemma}

We say that $X$ is $(\d,K)$-uniform if $X$ is $K$-uniform along the sequence of scales $K^{-j}$ with $j=0, \ldots, m$ and $m  = [\log_K (1/\d)]$. It follows from Lemma \ref{lem:uniform-subset} that we can always find a $K$-uniform subset $X'\subset X$ with $K = \d^{-o(1)}$ (where $o(1)$ is tending to $0$ with $\d$ and does not depend on any other parameters), and this is what we are going to use in applications. 

We can translate between collections of tubes with shadings $(\TT, Y)$ and configurations $X \subset \Omega_d$ in the following way. Given $(\TT,Y)$ define $X$ to be the set of triples $(x,t,\theta)$ where $(x,t) \in U(\TT,Y)^{\d} \cap (\d\ZZ)^d$ and $\theta = \theta(T)$ for some $T\in \TT_Y(p)$, $p = (x,t)$. 
In the other direction, for each $\omega=(x,t,\theta)\in X$ consider the corresponding tube $T=T^{\d}(\omega)$ and add the dyadic $\d$-cube  $(x,t)^\d$ to the shading of $T$.

\section{Structure of AD-regular sets}\label{sec:structure}

The goal of this section is to prove Lemma \ref{lem:finding-parallel-rectangles} which finds a useful structure inside an arbitrary AD-regular set. 

The statement says roughly the following.  Let $d\geq 2$ be an integer. 
Let $A \subset [-1,1]^d$ be a $(\d,t)$-AD-regular set. Then we can find some intermediate scales $\varrho< \d^{\eps}\tau$  
such that (a dense subset of) $A \cap B_\tau$ is contained in a union of parallel $\varrho\times \tau$ tubes (radius $\varrho$, length $\tau$) $T_1, \ldots, T_M$ such that:
\begin{itemize}
\item for every $i$, the set $A\cap T_i$ is $s$-AD-regular for some $s \in [0,t]$,
\item the set of tubes is $(t-s)$-AD-regular, and
\item for any other $\varrho\times \tau$ tube $T$ the set $A\cap T$ is $s$-Katz-Tao between scales $\varrho$ and $\tau$.
\end{itemize}

Let $A \subset [0,1]^d$ be a finite set. For $a \in A$ and scales $ \varrho< \tau \le 1$ define 
\[
M_{\varrho\times \tau}(a, A) = \max_{R} |R \cap A|
\]
where the maximum is over all $\varrho\times \tau$ tubes $R\subset \RR^2$ centered at $a$. 
Let us also define
\[
M_{\varrho\times \tau}(A) = \max_{a \in A} M_{\varrho\times \tau}(a, A).
\]
Up to a constant factor this is the maximum number of points of $A$ contained in a $\varrho\times \tau$ tube. 
We will use the following property later: for any $\varrho \le \varrho'$ and $\tau \le \tau'$ we have
\begin{equation}\label{eq:lipschitz-max}
    M_{\varrho'\times \tau'}(A) \le C (\varrho'/\varrho)^{d-1}(\tau'/\tau) M_{\varrho\times \tau}(A).
\end{equation}

We say that a subset $A$ of some metric space is $(\d,K)$-uniform if for every $r\in [\d, 1]$ and $a\in A$ we have
\[
|A \cap B_r(a)|_\d \ge \frac{1}{K} \max_{a'\in A} |A \cap B_r(a')|_\d.
\]
By a standard dyadic pigeonholing argument, there is some $K = \d^{-o(1)}$ such that any subset $A \subset [0,1]^d$ contains a  $\d$-separated subset $A' \subset A$ of size at least $\frac{1}{K} |A|_\d$ which is $(\d, K)$-uniform.

\begin{prop}\label{prop:pigeonhole-rectangle}
    Let $\d >0$ and let $A \subset [0,1]^d$ be a $\d$-separated set. Then for some $K = \d^{-o(1)}$, there exists a $(\delta, K)$-uniform subset $A' \subset A$ such that $|A'| \ge \frac1K |A|$ and
    \[
    \#\left\{ a\in A'~|~ \forall \varrho\le\tau \in [\d, 1]:~M_{\varrho\times \tau}(a, A') \ge \frac1K M_{\varrho\times \tau}(A')\right\} \ge \frac1K |A'|.
    \]
\end{prop}

This proposition shows that the maximal rectangular density centered at each point is about average.

\begin{proof}
    Let $\varepsilon>0$ and $m = \lceil\varepsilon^{-1}\rceil$. Let $S = \{(i, j):~ m \ge i>j \ge 0\}$ and for $a\in A$ define a vector ${\bf v}(a, A) \in \ZZ^{S}$ so that ${\bf v}(a, A)_{i,j}$ is equal to the number $k$ for which 
    \[
    M_{\d^{i/m}\times \d^{j/m}}(a, A) \in [\d^{-k/m}, \d^{-(k+1)/m}).
    \]
    If $A \subset [0,1]^d$ is $\d$-separated then ${\bf v}(a, A)$ can take at most $(2m+1)^{m^d} = O_\varepsilon(1)$ distinct values. By the pigeonhole principle we can then find some ${\bf v}_1 \in \NN^{S}$ so that ${\bf v}(a, A) = {\bf v}_1$ for at least $c_\varepsilon |A|$ elements $a\in A$. Let $\tilde A_1 \subset A$ be the set of such elements.  By pigeonholing, we can find a $(\delta, O_{\varepsilon}(1))$-uniform subset  $A_1\subset \tilde A_1$ such that $|A_1|\geq c_\varepsilon |\tilde A_1| \geq c_{\varepsilon}^2 |A|$.  Now let ${\bf v}_2$ be so that ${\bf v}(a, A_1) = {\bf v}_2$ for at least $c_\varepsilon |A_1|$ elements $a \in A_1$ and define the corresponding set $\tilde A_2 \subset A_1$. Find a $(\delta, O_{\varepsilon}(1))$-uniform subset $A_2\subset \tilde A_2$ with $|A_2|\geq c_{\varepsilon} |\tilde A_2|\geq c_{\varepsilon}^2 |A_1|$.  Continuing in this manner we define a sequence of vectors ${\bf v_1}, {\bf v_2}, \ldots$ and sets $A_1 \supset A_2 \supset \ldots$ so that $|A_{t+1}| \ge c_\varepsilon^2 |A_t|$ and ${\bf v}(a, A_t) ={\bf v}_{t+1}$ for any $a \in A_{t+1}$. 

    Now note that by $A_{t+1} \subset A_t$ we have ${\bf v}_{t+2} \le {\bf v}_{t+1}$ coordinate-wise. So since the vectors ${\bf v}_t$ take $O_\varepsilon(1)$ distinct values, there exists some $t \le O_\varepsilon(1)$ for which ${\bf v}_{t} = {\bf v}_{t+1}$ holds. This means that ${\bf v}(a, A_t) = {\bf v}_{t+1} = {\bf v}_t$ for any $a \in A_{t+1}$. On the other hand, for each $a\in A_t$ we have ${\bf v}(a, A_t) \le {\bf v}(a, A_{t-1}) = {\bf v}_t$. Unpacking definitions, we conclude that
    \[
    \#\left\{ a\in A_t~|~ \forall (i,j)\in S:  M_{\d^{i/m}\times \d^{j/m}}(a, A_t) \ge \d^{1/m} M_{\d^{i/m}\times \d^{j/m}}(A_t) \right\} \ge c_\varepsilon |A_t|.
    \]
    This gives the desired property for pairs $\varrho, \tau$ of the form $\d^{i/m}, \d^{j/m}$. Now let $\varrho \in [\d^{i/m}, \d^{(i-1)/m})$ and $\tau \in [\d^{j/m}, \d^{(j-1)/m})$ be arbitrary. 
    
    By (\ref{eq:lipschitz-max}) we have
    \[
    M_{\d^{(i-1)/m}\times \d^{(j-1)/m}}(A_t) \le C \d^{-d/m} M_{\d^{i/m}\times \d^{j/m}}(A_t).
    \]
    So using the monotonicity we get $M_{\varrho\times \tau}(a, A_t) \ge M_{\d^{i/m}\times \d^{j/m}}(a, A_t)$ and so for $c_\varepsilon$-fraction of $a\in A_t$ we get 
    \begin{align*}
    M_{\varrho\times \tau}(a, A_t) \ge M_{\d^{i/m}\times \d^{j/m}}(a, A_t) \ge \d^{1/m}M_{\d^{i/m}\times \d^{j/m}}(A_t) \\\ge c \d^{(d+1)/m} M_{\d^{(i-1)/m}\times \d^{(j-1)/m}}(A_t) \ge c\d^{(d+1)/m} M_{\varrho\times \tau}(A_t).    
    \end{align*}
    This completes the proof with $K = \d^{-(d+1)/m}$ and sufficiently small $\d< \d_0(\varepsilon)$.
\end{proof}

\begin{lemma}\label{lem:katz-tao-rectangles}
    Let $\zeta \in (0,1/2)$ and $\eta < \eta_0(\zeta)$, then the following holds for $\d \le \d_0(\eta)$. Suppose that $A \subset [0,1]^d$ is a $\d$-separated and $(\d, \d^{-\eta})$-uniform set. Then we can find a subset $A' \subset A$ of size at least $\d^{2\eta} |A|$, a pair of scales $\varrho\le \tau \in [\d, 1]$ and some $s \in \RR$ with the following properties:
    \begin{itemize}
        \item[(i)] $\tau/\varrho \ge \d^{-\chi(\zeta)}$, where $\chi(\zeta) =  c\zeta^{3\lceil \zeta^{-1}\rceil}$, and $\max(-s, s-1) \le \frac{C\eta}{\zeta^{3\lceil \zeta^{-1}\rceil}}$ (which means $s$ is roughly in $[0,1]$),
        \item[(ii)] for any $a\in A'$ there exists a $\varrho\times \tau$ tube $R$ containing $a$ such that 
        \[
        |A'\cap R| \ge \d^{O(\eta)} (\tau/\varrho)^{s} |A \cap B_{\varrho}(a)|,
        \]
        \item[(iii)] for any $\varrho\le \varrho' \le \tau'\le \tau$, and any $\varrho'\times \tau'$-tube $R'$ we have 
        \[
        |A'\cap R'|_{\varrho'} \leq C \delta^{-O(\eta)} (\tau/\varrho)^{\zeta}(\tau'/\varrho')^s.
        \]
        In particular, this implies
        \[
        |A' \cap R'|\le C\d^{-O(\eta)}(\tau/\varrho)^{\zeta} (\tau'/\varrho')^s |A \cap B_{\varrho'}(a)|
        \]
        for an arbitrary $a \in A$.
    \end{itemize}
\end{lemma}

\begin{proof}
    Let $A \subset [0,1]^d$ be a $\d$-separated, $(\d,\d^{-\eta})$-uniform set. 
    Let $A' \subset A$ be the $(\delta, K)$-uniform subset given by Proposition \ref{prop:pigeonhole-rectangle} with some $K = \d^{-o(1)}$. By taking $\d < \d_0(\eta)$ we can ensure that $K\le \d^{-\eta}$.
    Define a sequence of pairs $(\d_i, \Delta_i)$ as follows. Let $\d_0=\d, \Delta_0=1$ and suppose that we already defined $\d_i, \Delta_i$ for some $i\ge 0$. Let $s_i\in \RR$ be the number such that
    \begin{equation}\label{eq: si}
    M_{\d_i\times \Delta_i}(A') = (\Delta_i/\d_i)^{s_i} M_{\d_i\times \d_i}(A). 
    \end{equation}
    When $i=0$, this  means for $1/K$-fraction of $a\in A'$, there exists a $\d \times a$-tube $T$ centered at $a$ with $|T\cap A'|\geq cK^{-1} \d^{-s_0}$ (note that since $A$ is $\d$-separated, we have $M_{\d\times \d}(A) \le C$).  
    If there exist $\varrho \le \tau \in [\d_i, \Delta_i]$ so that 
    \[
    M_{\varrho\times \tau}(A') \ge (\Delta_i/\d_i)^{\zeta} (\tau/\varrho)^{s_i} M_{\varrho\times \varrho}(A)
    \]
    then we select $(\d_{i+1}, \Delta_{i+1}) = (\varrho,\tau)$ and continue the process.  Otherwise, we stop. 

    Suppose that we terminate after $t$ steps and consider the final pair $(\d_t, \Delta_t)$. Note that we have $\Delta_{i+1}/\d_{i+1} \le \Delta_i/\d_i$ and so 
    \begin{equation}\label{eq:M-iterate}
    (\Delta_{i+1}/\d_{i+1})^{s_{i+1}} M_{\d_{i+1}\times \d_{i+1}}(A) = M_{\d_{i+1}\times \Delta_{i+1}}(A') \ge (\Delta_i/\d_i)^\zeta (\Delta_{i+1}/\d_{i+1})^{s_{i}}  M_{\d_{i+1}\times \d_{i+1}}(A)
    \end{equation}
    \[
    (\Delta_{i+1}/\d_{i+1})^{s_{i+1}} \ge (\Delta_i/\d_i)^\zeta (\Delta_{i+1}/\d_{i+1})^{s_{i}}
    \]
    implying that $s_{i+1} \ge s_i + \zeta$ for all $i$. On the other hand, by uniformity of $A$, we have
    \begin{align}\label{eq:M-packing-bound}
    M_{\varrho \times \tau}(A') &\ge M_{\varrho\times \varrho}(A') \ge \frac{1}{CK} M_{\varrho\times \varrho}(A) \\ \label{eq:M-packing-bound2}
    M_{\varrho \times \tau}(A') &\le M_{\varrho \times \tau}(A) \le C\d^{-\eta} (\tau/\varrho)  M_{\varrho\times \varrho}(A).
    \end{align}
    Write $S=C\max(K, \d^{-\eta} )$ for some large constant $C$. Then (\ref{eq:M-packing-bound}), (\ref{eq:M-packing-bound2}), and \eqref{eq: si}  give
    \begin{equation}\label{eq:bound-Delta-delta-ratio}
    S^{-1}\le (\Delta_i/\d_i)^{s_i} \le S (\Delta_i/\d_i).    
    \end{equation}
    Note that if $s_{i}\ge 0$ then $(\Delta_{i+1}/\d_{i+1})^{s_{i}} \ge 1$ and if $s_{i}<0$ then $(\Delta_{i+1}/\d_{i+1})^{s_{i}} \ge (\Delta_i/\d_i)^{s_i} \ge S^{-1}$. 
    It follows that 
    \[
    S(\Delta_{i+1}/\d_{i+1}) \ge (\Delta_{i+1}/\d_{i+1})^{s_{i+1}} \ge (\Delta_i/\d_i)^\zeta (\Delta_{i+1}/\d_{i+1})^{s_{i}} \ge S^{-1}  (\Delta_i/\d_i)^\zeta, 
    \]
    where the first and the last inequalities  follow from  \eqref{eq:bound-Delta-delta-ratio} and the middle inequality follows from the choice of $\Delta_{i+1}, \delta_{i+1}$ and the definition of $s_{i+1}$ in \eqref{eq: si}. 
    This
    gives a lower bound
    \begin{equation*}
    \Delta_{i+1}/\d_{i+1} \ge S^{-2} (\Delta_i/\d_i)^\zeta.
    \end{equation*}
    So if $\zeta \le 1/2$ then iterating this $i$ times gives
    \begin{equation}\label{eq:ratio-Delta-lower-bound}
        \Delta_{i}/\d_{i} \ge S^{-4} (1/\d)^{\zeta^{i}}.
    \end{equation}
    Using $S \le C\d^{-\eta}$ and $\d\le \d_0(\eta)$, this implies that for any $i \ge0$ such that $\zeta^{i} \ge 10\eta$, we have $\Delta_{i}/\d_i \ge (1/\d)^{\zeta^i/2} \ge S$ and so by (\ref{eq:bound-Delta-delta-ratio}) we have $s_i \in [-1,2]$ for any such $i$. Now suppose that $\zeta^{2\lceil\zeta^{-1}\rceil} > 10\eta$ holds. Then by using the relation $s_{i+1} \ge s_i+\zeta$ for $i\le 3\lceil\zeta^{-1}\rceil$ we get that $s_{3\lceil\zeta^{-1}\rceil}-s_0 > 3$. This gives a contradiction with the above property that $s_i\in [-1,2]$ which means that the process must have stopped before the threshold $3\lceil\zeta^{-1}\rceil$, i.e. we have an upper bound $t < 3\lceil\zeta^{-1}\rceil$. 

    We claim that $s=s_t$, $\tau = \Delta_t$ and $\varrho=\d_t$ satisfy the conclusions of the lemma. Indeed, we have by the above
    \[
    \tau/\varrho \ge S^{-4} (1/\d)^{\zeta^t} \ge (1/\d)^{\zeta^t/2}\ge \d^{- \zeta^{3\lceil \zeta^{-1}\rceil}/2}
    \]
    provided that $\eta$ is small enough, giving the first part of (i). The bounds on $s=s_t$ then follow from (\ref{eq:bound-Delta-delta-ratio}). 
    Since the process has stopped at step $t$, we know that for any $\varrho\le \varrho'\le \tau'\le \tau$ we have 
    \[
    M_{\varrho'\times \tau'}(A') < (\tau/\varrho)^{\zeta} (\tau'/\varrho')^{s_t} M_{\varrho'\times \varrho'}(A),
    \]
    giving (iii) since both $A'$ and $A$ are $(\delta, \delta^{-\eta})$-uniform and $|A'|\geq \delta^{\eta}|A|$ (recall that the restriction that the tube is centered at a point of $A'$ only changes the bound by a constant factor). Finally, by the conclusion of Proposition \ref{prop:pigeonhole-rectangle} we have
    \[
    \#\left\{a\in A':~ \exists R:~ |A'\cap R| \ge \frac{1}{K} M_{\varrho\times \tau}(A')\right\} \ge \frac{1}{K}|A'| \ge \d^\eta |A'|
    \]
     where $R$ is a $\varrho\times \tau$ tube centered at $a \in A'$. Define $A''$ to be the set of $a\in A'$ for which there exists a $\varrho\times \tau$ tube $R$ containing $a$ and such that $|A' \cap R| \ge \frac{1}{K}M_{\varrho\times \tau}(A')$.  Note that for this $\varrho\times \tau$-tube $R$, all of $A'\cap R$ is included in $A''$. 
    It follows that for any $a \in A''$ there is a $\varrho\times \tau$ tube $R$ with $|R\cap A''| \ge \frac{1}{K}M_{\varrho\times \tau}(A')$. 
    By definition, we have $M_{\varrho\times \tau}(A') = (\tau/\varrho)^{s_t} M_{\varrho\times \varrho}(A)$ so this gives property (ii) for the set $A''$.  Since $A'' \subset A'$, the property (iii) is still true for $A''$ in place of $A'$. This finishes the proof.
\end{proof}


 \begin{definition}\label{def:aligned} 
    We say that a set $A \subset [-1,1]^d$ is $(\d, t, s, C)$-aligned if for some $\alpha \in [-1,1]^{d-1}$ we can write 
    \[
    A  = \{ (x, y + \alpha x), ~y \in Y, x \in X_{y} \},
    \]
    where $Y$ is $(\d,  t-s, C)$-AD-regular and for each $y \in Y$, $X_y\subset[-1,1]$ is $(\d, s, C)$-AD-regular. Furthermore, for any $\varrho\times \tau$ tube $T$ we have $|A \cap T|_\varrho \le C (\tau/\varrho)^{s}$.

    We say that a set $A \subset [-1,1]^2$ is nearly $(\d, t, s, C)$-aligned if there exists a $\d$-separated $(\d, t, s, C)$-aligned set $\tilde A$ such that $A\subset N_\d \tilde A $ and $\tilde A \subset N_\d A$.
\end{definition}

The following lemma requires in addition $A$ to be almost AD-regular. 

\begin{lemma}\label{lem:finding-parallel-rectangles}
    Let $\zeta>0$, $t \in [0,d]$, then the following holds for $\eta < \eta_0(\zeta)$,  $\d< \d_0(\eta)$ and  $\chi(\zeta)=c\zeta^{4\lceil \zeta^{-1}\rceil}$.
    Suppose that $A \subset [-1,1]^d$ is a $\d$-separated, $(\d,t, \d^{-\eta})$-AD-regular set.
    Then there exists $A'\subset A, |A'|\geq  (\varrho/\tau )^{\zeta}|A|$ and $s\in [0,t]$ and scales $\varrho, \tau\in [\delta, 1], \tau/\varrho \geq \delta^{-\chi(\zeta)}$  such that for every $a\in A'$, there exists an affine transformation $\phi_{a, \tau}$ such that $\phi_{a, \tau}(A'\cap B(a,\tau))$ is nearly $(\varrho/\tau, t,s, (\varrho/\tau)^{-\zeta})$-aligned.
\end{lemma}

In fact, we can take $\phi_{a,\tau}$ to be a dilation plus possibly a 90 degree rotation. The rotation is only needed because Definition \ref{def:aligned} only allows for slopes $\alpha$ at most 1.

\begin{proof}
    We choose a sufficiently large constant $C>0$ and sufficiently small  $0<\eta\ll \epsilon$ such that $\varepsilon_0:=  \eps+  C\eta \eps^{-3\lceil \eps^{-1}\rceil}  \lesssim \zeta^2$. 
Let $A' \subset A$ be a subset given by Lemma \ref{lem:katz-tao-rectangles} applied with $\eps$ in place of $\zeta$ (we will take $\eta \le \eta_0(\eps)$) and let  $s, \varrho, \tau$ be the output parameters. By (i) and (ii) of  Lemma \ref{lem:katz-tao-rectangles} the set $A'$ is covered by $\varrho\times \tau$ tubes $T$ such that 
\begin{equation}\label{eq: denseT}
|A' \cap T| \ge (\tau /\varrho)^{s - \varepsilon_0} \varrho^t |A'|.
\end{equation} So we can find a set of essentially distinct $2\varrho\times 2\tau$ tubes $\TT$ which cover $A'$ and $|A' \cap T| \ge (\tau /\varrho)^{s - \varepsilon_0} \varrho^t |A'|$ for every $T\in \TT$. By randomly rotating and passing to a constant density subset we may assume all tubes have slope with $y$-axes at most 1. 

    Consider a point-line configuration $X \subset \Omega_d$ defined as follows: for each $a \in A'$ and a tube $T\in \TT$ covering $a$ we associate the point-line pair $(a\in \ell)$ where $\ell$ is a line through $a$ parallel to the coreline of $T$ and put in $X$ the corresponding pair $(a, \theta)$, where $\ell = a+\RR(\theta,1)$ (abusing the notations, we identify $(a, \theta)$ with its corresponding $(a, \ell)$). 
    By construction and Item (iii) of Lemma~\ref{lem:katz-tao-rectangles}, we have
    \[
    (\tau /\varrho)^{s - \varepsilon_0} \varrho^t |A'| |\TT|\le |X| \le  (\tau /\varrho)^{s + O(\varepsilon_0)} \varrho^t |A'| |\TT|. 
    \]
     Furthermore, for every $(p,\ell) \in X$ there are at least $(\tau/\varrho)^{s-\varepsilon_0} \varrho^t |A'|$ pairs $(p', \ell') \in X$ such that $p'$ lies in the $2\varrho\times 2\tau$ tube associate to the pair $(p, \ell)$ and for which the line $\ell'$ is parallel to the core line of the tube.  This implies that for every $\omega\in X$ we have
    \begin{equation}\label{eq:uniform-rho-tau}
    \#\{ \omega' \in X:~ d_{C\tau, C\varrho, \varrho/\tau}(\omega', \omega) \le 1 \} \gtrsim (\varrho/\tau)^{O(\varepsilon_0)}  \M_{X}(C\tau, C\varrho, \varrho/\tau).    
    \end{equation}
    For some $K = \d^{-o(1)}$ we can find a $K$-uniform point-line configuration $X' \subset X$ with $|X'| \ge \frac1K |X|$. Note that by uniformity of $X'$ and (\ref{eq:uniform-rho-tau}), we have for any $\omega=(x, t, \theta)\in X'$:
    \[
    \M_{X'}(C\tau, C\varrho, \varrho/\tau)\ge K^{-O(1)}  (\varrho/\tau)^{O(\varepsilon_0)} \M_{X}(C\tau, C\varrho, \varrho/\tau)
    \]
    i.e. the concentration numbers for the triple of scales $(C\tau, C\varrho, \varrho/\tau)$ do not change significantly after passing to the refinement. In particular, applying Item (iii) of Lemma~\ref{lem:katz-tao-rectangles} with $\varrho'=\varrho$ and $\tau'\in [\varrho, \tau]$,  every $C\varrho\times C\tau$ tube $T$ associate to  a point-line pair $(p,\ell)\in X'$ contains a $([\varrho, \tau], s, (\tau/\varrho)^{O(\varepsilon_0)})$-AD regular set (see Definition~\ref{def: deltarhoADregularset}) of points from $A'$.  To see this, apply Lemma~\ref{lem:katz-tao-rectangles} (iii) with $\varrho'=\varrho$ and $\tau'\in [\varrho,\tau]$ we have for any $a\in A'\cap T$, 
  \begin{equation}\label{eq: ADregularupper}
  |A'\cap T\cap B_{\tau'}(a)|_\varrho \leq (\tau/\varrho)^{O(\varepsilon_0)} (\tau'/\varrho)^s. 
  \end{equation}

  Apply Lemma~\ref{lem:katz-tao-rectangles} (iii) with $\tau'=\tau, \varrho'\in [\varrho,\tau]$,
  we have 
  \[
  |A'\cap T|_{\varrho'}\leq  (\tau/\varrho)^{O(\varepsilon_0)} (\tau/\varrho')^s. 
  \]
  By \eqref{eq: denseT} and \eqref{eq: ADregularupper} (think $\varrho'=\tau'$),  we also have 
  \[
  |A'\cap T|_{\varrho'} \geq (\varrho/\tau)^{O(\varepsilon_0)} (\tau/\varrho')^s 
  \]
  and 
  \begin{equation}\label{eq: ADregularlower}
  |A'\cap T\cap B_{\tau'}(a)|_\varrho \geq (\varrho/\tau)^{O(\varepsilon_0)} (\tau'/\varrho)^s. 
  \end{equation}
  This implies that $T\cap A'$ contains a $([\varrho, \tau], s, (\tau/\varrho)^{O(\varepsilon_0)})$-AD-regular set of points from $A'$. 

  The directions of tubes in $\TT$ could be different, we are going to pigeonhole to find suitable scales such that the corresponding tubes (usually thicker and shorter) point in roughly the same direction. 
    
    For $m \ge 1$ and $i=0, \ldots, m$ let $\varrho_i = \varrho (\tau/\varrho)^{i/m}$ and $\Delta = (\varrho/\tau)^{1/m}$. Consider the quantities $D_i = \frac{\M_{X'}(\varrho_i, \varrho_i, 1) }{\M_{X'}(\varrho_i, \varrho_i, \Delta)}$ for $i=0,\ldots,m-1$. Up to a factor of $K^{O(1)}$, this is counting the $\Delta$-covering number of the set of lines $\ell$ where $(a, \ell) \in X'$ and $a \in B_{\varrho_i}(a_0)$ for a fixed $(a_0, \ell_0)\in X'$ (by uniformity the result does not depend on the choice of $(a_0, \ell_0)$). Note that this implies that we have $D_{i+1} \ge K^{-O(1)} D_i$ for any $i$ and that we have $D_i \in [K^{-O(1)}, K^{O(1)}\Delta^{-1}]$ for any $i$. By pigeonhole principle this implies that there is $i \in \{0, \ldots, m-1\}$ so that $D_{i+1} \le K^{O(1)} \Delta^{-1/m} D_i$.

    Fix some $(p_0, \ell_0) \in X'$ and consider the $\varrho_{i+1}$ box $Q$ centered at $p_0$. 
    Let $A_Q \subset A'\cap Q$ be the set of $a \in A'$ so that $a\in Q$ and there exists $(a, \ell) \in X'$.
    Let $Q_1, \ldots, Q_N \subset Q$ be a collection of finitely overlapping $\varrho_i$-boxes covering $A_Q$. By uniformity, we can choose these boxes so that $N = K^{O(1)} \delta^{O(\eta)} (\varrho_{i+1}/\varrho_i)^{t} =  K^{O(1)} \delta^{O(\eta)}  \Delta^{-t}$ and $|A_Q \cap Q_j|_\d = K^{O(1)} \delta^{O(\eta)}  (\varrho_i/\d)^t$. This is where we use the assumption that $A$ is $(\d,t,\d^{-\eta})$-AD-regular and both refinements $A'\subset A$ and $X'\subset X$ passed to sufficiently dense subsets. 
    Let $\Theta \subset \mathbb{S}^1$ be a collection of angles whose $\Delta$-neighbourhoods cover all directions $\theta(\ell)$ with $(a, \ell) \in X'$, $a\in A_Q$. By uniformity, we have $|\Theta| \le K^{O(1)} D_{i+1}$. Let $\Theta_j \subset \Theta$ the the subset of angles covering directions of lines coming from the box $Q_j$. By the choice of $i$ and uniformity, we have for every $j$ that $|\Theta_j| \ge K^{-O(1)} D_i \ge K^{-O(1)}\Delta^{1/m} D_{i+1} \ge K^{-O(1)}\Delta^{1/m} |\Theta| $.
    By pigeonhole principle and double counting, there exists $\theta \in \Theta$ so that the number of $j \in \{1, \ldots, N\}$ with $\theta\in \Theta_j$ is at least $ K^{-O(1)}\Delta^{1/m} N$. Let $J \subset \{1, \ldots, N\}$ be the set of such $j$.
    Cover the set of boxes $\{Q_j, ~j\in J \}$ by  a set $\TT$ of parallel finitely overlapping $C\varrho_i\times C\varrho_{i+1}$ tubes in the direction of $\theta$. We have $|\TT| \lesssim \varrho_{i+1}/\varrho_i \sim \Delta^{-d+1}$. Note that since  any $T\in \TT$ contains some $Q_j$, it follows that $C T$ contains a $([\varrho_{i}, \varrho_{i+1}], s, (\tau/\varrho)^{O(\varepsilon_0)})$-AD regular set of points from $A'$ (since $T$ contains a piece of $\varrho\times \tau$ rectangle through $(p, \ell)$ which in turn contained a $([\varrho, \tau], s, (\tau/\varrho)^{O(\varepsilon_0)})$-AD regular subset of $A'$). Using that $A'\cap Q$ is $t$-AD-regular, we have $|\TT| \lesssim  (\tau/\varrho)^{O(\varepsilon_0)}(1/\Delta)^{t-s}$. On the other hand, since all boxes $Q_j$, $j\in J$ are covered, it also follows that $|\TT| \gtrsim \Delta^{1/m} (\tau/\varrho)^{-O(\varepsilon_0)}(1/\Delta)^{t-s}$.

        Let $J' \subset J$ be the subset of indices $j$ for which the box $Q_j$ is covered by a tube $T\in \TT$ with at least $c|J|/|\TT|$-many other boxes $Q_{j'}$, $j'\in J$. Let $\TT'\subset \TT$ denote such subset of tubes.  Then we have that each $T\in \TT'$ also contains at least $c|J|/|\TT|$-many boxes $Q_{j'}$ with $j' \in J'$. Note that we can ensure that $|J'| \gtrsim |J|$. Then $|\TT'|\gtrsim \Delta^{1/m} (\varrho/\tau)^{O(\varepsilon_0)}|\TT|$ since $T$ contains a $([\varrho_i, \varrho_{i+1}], s, (\tau/\varrho)^{O(\epsilon_0)})$-AD regular set of points from $A'$ and its union covers a $\geq \Delta^{1/m}(\tau/\varrho)^{\varepsilon_0}$-fraction of $A'\cap Q$. In fact, from this we know that $c|J|/|\TT| \sim (\varrho_{i+1}/\varrho_i)^{s} (\tau/\varrho)^{O(\varepsilon_0)}$. 

Let $\TT''\subset \TT'$ be a (rescaled ) $(\varrho_{i+1}/\varrho_i, (\tau/\varrho)^{\varepsilon_0})$-uniform subset of tubes such that  the subset of indices $J''\subset J'$ covered by $\TT'$ satisfies $|J''|\gtrapprox |J|$ and so $|\TT''|\gtrapprox  \Delta^{1/m} (\varrho/\tau)^{O(\varepsilon_0)}|\TT|$.
        
    We  show that $\TT''$ is $([\varrho_i, \varrho_{i+1}], t-s, S)$-AD regular with $S \lesssim  (1/\Delta)^{1/m} (\tau/\varrho)^{O(\varepsilon_0)}$. To see this, for each $T\in \TT''$,   $r\in [\varrho_i, \varrho_{i+1}]$ and $\tilde T= N_r T$, since $CT$ contains a $([\varrho_i, \varrho_{i+1}, s, (\tau/\varrho)^{O(\varepsilon_0)})$-AD regular set of points from $A'$, 
    \[
    |N_r T\cap A'|_r \gtrsim (\varrho_{i+1}/r)^s (\varrho/\tau)^{O(\varepsilon_0)}.
    \]
  Since $\TT''$ is uniform and  $|J''|\gtrapprox |J|$, which is at least a $\Delta^{1/m}(\tau/\varrho)^{\varepsilon_0}$-fraction of $A'\cap Q$, we have 
  \[
  |\TT''[N_r T]|\gtrsim  \Delta^{1/m}(\varrho/\tau)^{O(\varepsilon_0)} (\varrho_{i+1}/r)^s (r/\varrho_i)^t (\varrho_{i+1}/\varrho_i)^{-s}\gtrsim \Delta^{1/m}(\varrho/\tau)^{O(\varepsilon_0)}(r/\varrho_i)^{t-s}. 
  \]
By Lemma~\ref{lem:katz-tao-rectangles} (iii), 
\[
 |N_r T\cap A'|_r \lesssim (\tau/\varrho)^{O(\varepsilon_0)} (\varrho_{i+1}/r)^s.
\]
So we have the upper bound
\[
|\TT''[N_rT]|\lesssim (\tau/\varrho)^{O(\varepsilon_0)} (r/\varrho_i)^{t-s}. 
\]

Define $A''_Q \subset A'_Q$ to be the set of points of $A'_Q$ contained in the union of boxes $Q_{j'}$, $j'\in J''$ then any tube $T\in \TT''$ which 
contains an $([\varrho_i, \varrho_{i+1}], s, CS)$-AD-regular subset in $A''_Q$. By taking $m > 10/\zeta$ and $\eps$ small enough  we can get $(\tau/\varrho)^{O(\varepsilon_0)} \le \Delta^{-c\zeta} $, i.e. $\varepsilon_0 \lesssim c\zeta/m$, and so we can get $CS < \Delta^{-\zeta}$. We see that the set $A'' = \bigcup A''_Q$ where the union is over $\varrho_{i+1}$-boxes $Q$ covering the point set $P[X']$ satisfies the conclusions of the lemma. Namely we redefine $A'=A''$, $\varrho_0 = \varrho_{i}$, $\tau_0 = \varrho_{i+1}$ and note that $\tau_0/\varrho_0 \ge (\tau/\varrho)^{1/m} \ge \d^{-\chi(\zeta)}$ with $\chi(\zeta)= c\zeta^{4\lceil \zeta^{-1}\rceil}$. For any $a\in A''$, let $\phi_{a, \tau}$ be an affine transformation that maps $B(a,\tau)$ to the unit ball, then upto composing $\phi_{a, \tau}$ with a possible rotation, which is again an affine transformation, $\phi_{a,\tau}(A'\cap B(a,\tau))$ is $(\varrho_0/\tau_0,t, s, (\varrho_0/\tau_0)^{-\zeta})$-aligned. In particular, for any $\varrho/\tau \leq \varrho' \leq \tau'\leq 1$,  any $\varrho'\times \tau'$-tube $R$ we have  $|\phi_{a,\tau}(A'\cap B(a,\tau))\cap R|_{\varrho'}  \leq  (\tau/\varrho)^{\zeta} (\tau'/\varrho')^s$ following from Lemma~\ref{lem:katz-tao-rectangles} (iii). 
\end{proof}

\part{Sticky Furstenberg}

\section{Setup for sticky Furstenberg}\label{sec:setup-furstenberg}

In the following Sections \ref{sec:setup-furstenberg}--\ref{sec:very-coplanar} we prove Theorem \ref{thm:furstenberg-R3}.

\begin{definition}\label{def: Fstkappa}
    We use $\operatorname{F}_{s,t}(\kappa)$ to denote the following assertion. 
    
    Let $s\in (0, 1], t\in (0, 2]$. For any $\epsilon>0$, there exists $\eta = \eta(s,t, \epsilon)>0$ such that the following holds for $\delta>0$ sufficiently small. Let $\TT$ be a $(2t, \delta^{-\eta})$-AD-regular set of $\delta$-tubes and $C_{t-CW}(\TT)\leq  \delta^{-\eta}$. Suppose $Y$ is a $(\delta, s, \delta^{-\eta})$-dense shading on $\TT$ such that 
    \[
    \sum_{T\in \TT} |Y(T)|_\d \ge \d^{\eta} |\TT| \max_{T\in \TT} |Y(T)|_\d.
    \]
    Then 
    \begin{equation}\label{eq:UTYd}
    |U(\TT, Y)|_{\delta} \geq \delta^{\epsilon} |Y(T)|_\d \delta^{-2t+\kappa}.     
    \end{equation}
\end{definition}

Using this definition, Theorem \ref{thm:furstenberg-R3} is equivalent to $\operatorname{F}_{s,t}(\kappa)$ with $\kappa = \max\{0,t-s, 2t-2\}$. 
To prove Theorem \ref{thm:furstenberg-R3} we will choose minimal $\kappa$ such that $\operatorname{F}_{s,t}(\kappa)$ fails and study configurations of tubes which achieve almost equality in (\ref{eq:UTYd}). In the terminology of \cite{wang2026sticky} these are called {\em extremal families} of tubes. 
Note that by double counting, a typical point $p \in U(\TT, Y)$ is contained in approximately $\frac{\sum_{\TT} |Y(T)|}{|U(\TT, Y)|} \approx \frac{\d^{-2t}|Y(T)|}{|U(\TT, Y)|} \approx \d^{-\kappa}$ many shadings $Y(T)$ (assuming that $Y(T)$ is a union of dyadic $\d$-cubes). Thus, the parameter $\kappa$ quantifies the {\em multiplicity} $\mu(\TT, Y)$ of the pair $(\TT, Y)$.  

The following is a special case of Definition \ref{def:aligned}.

 \begin{definition}\label{def:aligned2}
    We say that a set $\Theta \subset [-1,1]^2$ is $(\d, \kappa, \gamma, C)$-aligned if for some $\alpha \in [-1,1]$ we can write 
    \[
    \Theta  = \{ (\varphi, \xi + \alpha \varphi), ~\xi \in \Xi, \varphi \in \Phi_\xi \},
    \]
    where $\Xi$ is a $(\d,  \kappa-\gamma, C)$-AD-regular and for each $\xi \in \Xi$, $\Phi_\xi$ is $(\d, \gamma, C)$-AD-regular. Furthermore, for any $\varrho\times \tau$ tube $T$ we have $|\Theta \cap T|_\varrho \le C (\tau/\varrho)^{\gamma}$.
\end{definition}

Using induction on scale we reduce the sticky Furstenberg problem to studying the following structured case. 

\begin{definition}\label{def: configuration}
A $(\d, \eta, s, t, \kappa, \gamma)$-Furstenberg configuration is the following collection of data: 
\begin{enumerate}
    \item A set of  $\d$-tubes $\TT$ which is $(2t, \d^{-\eta})$-AD-regular. 
    We also assume each tube $T\in \TT$ has $\theta(T) \in [-1,1]^2$. 
    \item A $\d^{\eta}$-uniform $(\d, s)$-set $Z \subset \d \ZZ \cap [-1,1]$.
    \item Set $E$ of the form $E = \bigcup_{z\in Z} E_z$ where $E_z =A_z\times \{z\}$ and $A_z \subset [-1,1]^2$ is a $(\d, 2t-\kappa, \d^{-\eta})$-AD-regular set.

    \item For every $p \in E$, there is a $\d$-separated $(\d, \kappa, \gamma, \d^{-\eta})$-aligned  set $\Theta_p \subset [-1,1]^2$ such that: for every $\theta\in \Theta_p$ there exists $T\in \TT$ with $p\in T$ and $\theta(T) \in N_\d \theta$. We denote by $\TT_p \subset \TT$ the set of all such tubes over $\theta\in \Theta_p$. 
    \end{enumerate}
\end{definition}

For short, we represent this data by the tuple $\mc C=(\TT, E, \{\TT_p\})$. A $\lambda$-refinement of a Furstenberg configuration is a Furstenberg configuration $\mc C'=(\TT, E', \{\TT'_p\})$ with $E' \subset E$ and $\TT'_p \subset \TT_p$ and
\[
\sum_{p \in E'} |\TT'_p| \ge \lambda \sum_{p\in E} |\TT_p|.
\]
Note that we always keep the set of tubes the same. 

\begin{remark}
For each $T\in \TT$, let $Z(T)$ be the set of $z\in Z$ so that $T\in \TT_p$ for some $p\in E_z$. Then 
\[
\sum_{T\in \TT} |Z(T)|_{\delta} \gtrsim \delta^{4\eta} |\TT|_{\delta} |Z|_{\delta}. 
\]
To see this, by Item (3),   $|E|_{\delta}\gtrsim \d^\eta |Z|_{\delta}\delta^{-(2t-\kappa)}$. By Item (4) and  double counting, \[ \sum_{T\in \TT}|Z(T)|_{\delta} \sim |E|_{\delta} |\Theta_p|_{\delta} \gtrsim |Z|_{\delta} \delta^{-(2t-\kappa)} \cdot \delta^{-\kappa+2\eta} \gtrsim \delta^{4\eta}|\TT|_{\delta}|Z|_{\delta}, \]
where the last inequality is because $\TT$ is a $(2t, \delta^{-\eta})$-AD-regular set. 
\end{remark}

We prove Theorem \ref{thm:furstenberg-R3} using the following four propositions. Note that only Proposition \ref{prop: coplanar} requires the Convex Wolff axiom to obtain a desired bound on $\kappa$. 

The first proposition uses {\em critical}  pairs $(\TT, Y)$, i.e. which nearly achieve equality in (\ref{eq:UTYd}), to construct Furstenberg configurations with the same parameter $\kappa$ and some $\gamma \in [0,\kappa]$. 

\begin{prop} \label{prop: reduction}
Let $\kappa$ be the infimum of the set of  $\kappa'$ such that $\operatorname{F}_{s,t}(\kappa')$ holds. 
    Then for any $\d_0, \eta$ there exists $\d< \d_0$ and $\gamma \in [0,\kappa]$ so that there exists a $(\d, \eta, s, t, \kappa, \gamma)$-Furstenberg configuration whose set $\TT$ of $\delta$-tubes satisfies $C_{t-CW}(\TT)\leq \delta^{-\eta}$. 
\end{prop}

The following propositions give upper bounds on $\kappa$ for different ranges of the parameter $\gamma$. 

\begin{prop}\label{prop: high-low}
    For any $\epsilon_0>0$, there exists $\eta>0$ such that the following holds for $\d>0$ sufficiently small. 
    
    Suppose that there is a $(\d, \eta, s, t, \kappa, \gamma)$-Furstenberg configuration. 
    If $\gamma < s+\kappa-1-\eps_0$ then $\kappa \le \max\{2t-2,0\}+2\epsilon_0$.
\end{prop}

\begin{prop}\label{prop: trilinearABC} For any $\epsilon_0>0$, $\epsilon_1>0$, there exists $\eta>0$ such that the following holds for $\d>0$ sufficiently small. 

    Suppose there is a $(\d, \eta, s, t, \kappa, \gamma)$-Furstenberg configuration. Suppose that $\gamma \in [s+\kappa-1-\eps_0, 1-s+\eps_0]$  and $\gamma \in [\eps_1, \kappa-\eps_1]$ then $\kappa \le \max\{t-s,0\}+4\epsilon_0$.
\end{prop}

\begin{prop}\label{prop: coplanar} For any $\epsilon_0>0$, there exists $\eta>0$ such that the following holds for $\d>0$ sufficiently small. 

    Suppose there is a $(\d, \eta, s, t, \kappa, \gamma)$-Furstenberg configuration whose set of $\d$-tubes $\TT$ satisfies $C_{t-CW}(\TT) \le \d^{-\eta}$. If $\gamma \ge \max\{s+\kappa-1-\eps_0, 1-s+\eps_0\}$ then $\kappa \le \max\{2t-2, t-s,0\}+4\eps_0$. 
\end{prop}

The three propositions above roughly consider different cases of how $\gamma$ relates to the numbers $s+\kappa-1$ and $1-s$, except for the fact that Proposition \ref{prop: trilinearABC} in addition requires that $\gamma$ is separated from the edges of the interval $[0,\kappa]$. The following two propositions are aimed at covering those two edge cases. 

\begin{prop}\label{prop: very-non-coplanar} For any $\epsilon_0>0$, there exists $\eta>0$ such that the following holds for $\d>0$ sufficiently small. 

    Suppose there is a $(\d, \eta, s, t, \kappa, \gamma)$-Furstenberg configuration whose set of $\d$-tubes $\TT$ satisfies $C_{t-CW}(\TT) \le \d^{-\eta}$. If $\gamma \le \eps_0$ and $s+\kappa\le 1+\eps_0$ then $\kappa \le \max\{2t-2, t-s,0\}+2\eps_0$. 
\end{prop}

\begin{prop}\label{prop: very-coplanar} For any $\epsilon_0>0$ and sufficiently small $\eps_1>0$, there exists $\eta>0$ such that the following holds for $\d>0$ sufficiently small. 

    Suppose there is a $(\d, \eta, s, t, \kappa, \gamma)$-Furstenberg configuration whose set of $\d$-tubes $\TT$ satisfies $C_{t-CW}(\TT) \le \d^{-\eta}$. If $\gamma \ge \kappa-\eps_1$ and $s+\kappa\le 1+\eps_0$ then $\kappa \le \max\{2t-2, t-s,0\}+2\eps_0$. 
\end{prop}

Proposition \ref{prop: reduction} is proved in Section \ref{sec:configuration-furstenberg}. Proposition \ref{prop: high-low} is proved in Section \ref{sec:high-low}. Proposition \ref{prop: trilinearABC} is proved in Section \ref{sec:intermediate}. Proposition \ref{prop: coplanar} is proved in Sections \ref{sec:global-grains}--\ref{sec:proof-of-coplanar-case}. Proposition \ref{prop: very-non-coplanar} is proved in Section \ref{sec:very-non-coplanar} and Proposition \ref{prop: very-coplanar} is proved in Section \ref{sec:very-coplanar}.

\begin{proof}[Proof of Theorem \ref{thm:furstenberg-R3}]
Note that the statement of Theorem \ref{thm:furstenberg-R3} is equivalent to the assertion $\operatorname{F}_{s,t}(\kappa_0)$ where $\kappa_0 = \max(2t-2, t-s,0)$. Let $\kappa$ be the infimum of all $\kappa'$ such that $\operatorname{F}_{s,t}(\kappa')$ holds. It suffices to show that $\kappa \le \kappa_0$. Let $\eps_0 >0$ be an arbitrarily small constant, let $\eps_1 = \eps_1(2\eps_0)$ be constant from Proposition \ref{prop: very-coplanar} and let $\eta$ be sufficiently small depending on $\eps_0,\eps_1$. 

By Proposition \ref{prop: reduction}, for any $\d_0$ there is $\d <\d_0$ and $\gamma\in [0, \kappa]$ so that there exists a $(\d, \eta, s, t, \kappa, \gamma)$-Furstenberg configuration whose set of $\d$-tubes $\TT$ satisfies $C_{t-CW}(\TT)\le \d^{-\eta}$. Then 
\begin{itemize}
    \item if $\gamma \le s+\kappa-1-\eps_0$, then by Proposition \ref{prop: high-low} we have $\kappa \le \max(2t-2,0)+2\eps_0 \le \kappa_0+2\eps_0$,
    \item if $\gamma \in [s+\kappa-1-\eps_0, 1-s+\eps_0] \cap [\eps_1,\kappa-\eps_1]$ then by Proposition \ref{prop: trilinearABC} we have $\kappa \le \max(t-s,0) +4\eps_0 \le \kappa_0 +4\eps_0$,
    \item if $\gamma \ge \max\{s+\kappa-1-\eps_0, 1-s+\eps_0\}$ then by Proposition \ref{prop: coplanar} we have $\kappa \le \max\{2t-2, t-s,0\}+4\eps_0 \le \kappa_0 +4\eps_0$,
    \item if $\kappa \le 1-s+2\eps_0$ and $\gamma \le \eps_1$ then by Proposition \ref{prop: very-non-coplanar} we have $\kappa \le \max\{2t-2, t-s,0\}+4\eps_0$.
    \item if $\kappa \le 1-s+2\eps_0$ and $\gamma \ge\kappa- \eps_1$ then by Proposition \ref{prop: very-coplanar} we have $\kappa \le \max\{2t-2, t-s,0\}+4\eps_0$. 
\end{itemize}

This exhausts all possibilities for $\gamma$: the first three bullet points cover everything except $\gamma \in [s+\kappa-1-\eps_0, 1-s+\eps_0] \cap [0,\eps_0]$ and $\gamma \in [s+\kappa-1-\eps_0, 1-s+\eps_0] \cap [\kappa-\eps_0, \kappa]$. If the first interval intersection is non-empty, then $s+\kappa -1-\eps_0 \le \eps_0$ and the 4-th bullet applies; if the second interval intersection is non-empty, then $1-s+\eps_0 \ge \kappa-\eps_0$ and the 5-th bullet applies. 

We conclude that in all cases we have $\kappa \le \kappa_0+4\eps_0$. Since $\eps_0>0$ is arbitrary, we conclude that $\kappa \le \kappa_0$ holds, concluding the proof.
\end{proof}

\subsection{Using weaker Convex Wolff axioms.}\label{subsec:weaker-wolff}

By remarks at the end of Sections \ref{sec:proof-of-coplanar-case} and \ref{sec:very-coplanar}, the sharp sticky Fustenberg estimate holds with a slightly relaxed convex Wolff axiom $C_{t;t'-CW}(\TT) \le \d^{-\eta}$ for some $t' \le t'(s,t)$ that can be defined as follows:
\[
t'=
\begin{cases}
    \min(2t,s), \quad t\in [0,s],\\
    2t-s, \quad t \in [s,1],\\
    2-s, \quad t \in [1, 2-s],\\
    \min(2t-2+s, 2), \quad t\in [2-s,2].
\end{cases}
\]
Recall that in the definition of $C_{t;t'-CW}$ we require that $t' \in [0,2]$ and $t' \le 2t$ (so that both exponents in $a^{2t-t'}b^{t'}$ are non-negative), hence the second terms in the two minimums above.
It turns out that all other restrictions on $t'$ are in fact necessary, as we show on the following two examples. 

For some $t_1 \le 2t$ let $(\TT_1, Y_1)$ be a sharp example for Furstenberg in the plane with $\TT_1$ a $(\d, t_1)$-AD-regular and $Y_1(T)$ a $(\d, s)$-set. 
Take a $A \subset [-1,1]$ a $(\d, 2t-t_1)$-AD-regular set and let $\TT$ be the union of copies of $\TT_1$ inside the planes $[-1,1]^2 \times \{a\}$, $a\in A$. Then $\TT$ is almost $2t$-AD-regular and
\[
\mu(\TT)\gtrsim \mu(\TT_1) \gtrsim \d^{-\max(0, (t_1-s)/2, t_1-1)}.
\]
It is straightforward to check that $C_{t;t_1-CW}(\TT) \lessapprox 1$: indeed, the worst case is $U = [-1,1]^2 \times [-\d,\d]$ for which $|\TT[U]| \approx |\TT_1| \approx \d^{-t_1}$ and this should be upper bounded by $\d^{2t} (1/\d)^{t_1} |\TT| \approx \d^{-t_1}$. 

Consider different cases:
\begin{itemize}
    \item If $t \in [0,s]$, put $t_1 = s+\zeta$, so $\mu(\TT) \gg 1 = \d^{-\kappa_0(s,t)}$. Thus, we need $t' \le s$. 

    \item If $t \in [s,1]$, put $t_1 = 2t-s+\zeta$, so $\mu(\TT) \gg \d^{-(t-s)} =  \d^{-\kappa_0(s,t)}$. Thus, we need $t' \le 2t-s$.

    \item If $t \in [1,2-s]$, put $t_1 = t+1-s+\zeta$, so $\mu(\TT) \gg \d^{-(t-s)} =  \d^{-\kappa_0(s,t)}$. Thus, we need $t' \le t+1-s$. 
\end{itemize}

For $t \in [2-s,2]$ we will give a better example. Again let $(\TT_1, Y_1)$ be a planar Furstenberg pair so that $\TT_1$ is almost $(2t-t_1)$-AD-regular for some $t_1$ and $Y$ is $s$-set. Suppose that $t_1 \in [2-s, 2]$. For each $T_1 \in \TT_1$ let $\TT_{T_1}$ be a $(\d, t_1)$-AD-regular set of tubes in $T \times [-1,1]$. Define the shading by $Y(T) = T\cap (Y_1(T_1)\times [-1,1])$ and put $\TT =\bigcup_{T_1\in \TT_1} \TT_{T_1}$. 
Note that for $t_1\in [2-s,2]$ we have $U(\TT_{T_1}, Y) \approx Y(T_1) \times [-1,1]$ and $\mu(\TT_{T_1}) \approx \d^{-(t_1-1)}$.

Then $\TT$ is almost $2t$-AD-regular and 
\[
\mu(\TT) \gtrsim \mu(\TT_1)  \mu(\TT_{T_1}) \gtrsim \d^{-\max(0, (2t-t_1-s)/2, 2t-t_1-1)} \d^{-(t_1-1)}.
\]
As previously, we can verify that $C_{t;t_1}(\TT) \lessapprox 1$. So

\begin{itemize}
    \item If $t\in [1, 2-s]$, put $t_1 = 2-s+\zeta$, so $\mu(\TT) \gg \d^{-(t-s)} =  \d^{-\kappa_0(s,t)}$. Thus, we need $t' \le 2-s$.
    \item If $t \in [2-s,2]$, put $t_1 = 2t-2+s+\zeta$, so $\mu(\TT) \gg \d^{-(2t-2)} =  \d^{-\kappa_0(s,t)}$. Thus, we need $t' \le 2t-2+s$.
\end{itemize}

\section{Preliminaries}\label{sec:prelim-furstenberg}

In the following sections we study the structure of  $(\d, \eta, s, t, \kappa, \gamma)$-Furstenberg configurations inside a $\varrho$-cube. The following two lemmas allow us to refine the set $E$ to obtain useful structure inside $\varrho$-cubes.

\begin{lemma}\label{lem:same-size-in-rho-ball}
     Let $(\TT, E, \{\TT_p\})$ be an $(\d, \eta, s, t, \kappa, \gamma)$-Furstenberg configuration. Let $\varrho \in [\d, 1]\cap 2^{-\mathbb{N}}$. Then there exists a subconfiguration $(\TT, E', \{\TT_p\})$ such that $|E'|_\delta \ge |E|_\delta/2$ and $|E' \cap \bp|_\delta \sim_{\d^{O(\eta)}} |Z \cap \mathbf{z}|_\delta  (\varrho/\d)^{2t-\kappa}$ for any $p \in E'_z$ and the $\varrho$-dyadic cube $\bp$ containing $p$ and $\varrho$-dyadic interval $\mathbf{z}$ containing $z$. 

     Furthermore, for $Z' =\{z: E'_z\neq\emptyset\}$ and $z \in Z'$ we have $|\bigcup_{z'\in Z' \cap {\bf z}} E'_{z'}|_{\varrho} \sim_{\d^{O(\eta)}} \varrho^{-(2t-\kappa)}$.
\end{lemma}

\begin{proof}
    Note that 
    \begin{equation}\label{eq: EBprhoupper}
    |E\cap \bp|_\delta \le \d^{-O(\eta)}|Z \cap \mathbf{z}|_\delta (\varrho/\d)^{2t-\kappa}
    \end{equation} by $\d^\eta$-uniformity of $Z$ and AD-regularity of sets $E_z$. 
    We claim that $|E|_{\varrho} \le \d^{-O(\eta)} |Z|_\varrho \varrho^{-(2t-\kappa)}$. Indeed, by double counting we have 
    \[
    |\TT|_\varrho \gtrsim \frac{1}{|Z|_\varrho} \sum_{\bp \in \cD_\varrho(E)} |\TT_{\bp}|_\varrho \gtrsim \frac{\d^{\eta} \varrho^{-\kappa}}{|Z|_\varrho} |E|_\varrho.
    \]
    
    Let $\mc D \subset \mc D_{\varrho}(E)$ be the subset of dyadic boxes $\bp$ so that $|\bp\cap E|_\delta \ge\d^{C_1\eta} |Z \cap \mathbf{z}|_\delta (\varrho/\d)^{2t-\kappa}$ for a sufficiently large constant $C_1$. By \eqref{eq: EBprhoupper} we have
    \[
    \sum_{\bp \in \mc D_{\varrho}(E) \setminus \mc D} |\bp\cap E|_\delta \le \d^{C_1\eta}|Z \cap \mathbf{z}|_\delta (\varrho/\d)^{2t-\kappa} |E|_\varrho \le \d^{C_1\eta - O(\eta)} |E|_\delta.
    \]
    So for large enough $C_1$, this sum is less than $|E|/4$. So the set $E' = \bigsqcup_{\bp \in \mc D} E\cap \bp$ satisfies the conclusion of the lemma. We obtain the second conclusion by a similar pruning argument. 
\end{proof}

\begin{lemma}\label{lem:pigeonhole-rho-ball}
    Let $(\TT, E, \{\TT_p\})$ be a $(\d, \eta, s, t, \kappa, \gamma)$-Furstenberg configuration. Let $\varrho \in [\d, 1]$. Then there exists a subset $E' \subset E$ such that $|E'| \ge \d^{O(\eta)} |E|$ and 
    for every $p_0 \in E'$ and $\mathbf{q}=p_0^\varrho$  (the dyadic $\varrho$-cube containing $p_0$) there exists $p_0' \in {\bf q}$ such that the following holds:
    \begin{itemize}
        \item[(a)] For every $p \in E'\cap \mathbf{q}$ we have $|\Theta_{p} \cap \Theta^\varrho_\mathbf{q}|_\d \ge \d^{O(\eta)} |\Theta_{p}|_\d$, where $\Theta_\mathbf{q} = \Theta_{p'_0}$. 
        \item[(b)] Let $\TT(\mathbf{q})$ be the set of $T\in \TT$ such that $\theta(T) \in \Theta_\mathbf{q}^{\varrho}$ and $T \in \TT_p$ for some $p \in E'\cap \mathbf{q}$. 
        Then $|\TT(\mathbf{q})| \sim_{\d^{O(\eta)}} \varrho^{-\kappa} (\varrho/\d)^{2t}$.
        Let $Z(\mathbf{q}) = \{z\in Z: E'_z \cap \mathbf{q}\neq\emptyset\}$.
        For $T\in \TT$ let $Z(\mathbf{q}, T)$ be the set of $p \in Z({\bf q})$ such that $T\in \TT_p$ for some $p \in E'_z\cap {\bf q}$. Then for any $\tilde Z\subset Z({\bf q})$ such that $|E'_{\tilde Z} \cap {\bf q}| \ge \d^{C\eta} |\tilde Z| (\varrho/\d)^{2t-\kappa}$ we have
        \[
        \sum_{T\in \TT(\mathbf{q})} |\tilde Z\cap Z(\mathbf{q}, T)|_\delta \sim_{\d^{O(\eta)}} |\tilde Z|_\delta |\TT(\mathbf{q})|.
        \]
    \end{itemize}
\end{lemma}

\begin{proof}
By Lemma \ref{lem:same-size-in-rho-ball}, we can pass to a refinement and assume that $|E \cap {\bf p}|_\d \sim_{\d^{O(\eta)}} |Z \cap {\bf z}|_\d (\varrho/\d)^{2t-\kappa}$ holds for every ${\bf p} = p^\varrho$, $p \in E$. 

For $p \in E$ let $\TT(p, \varrho)$ denote the set of $T\in \TT$ such that $T \in \TT_{p'}$ for some $p'\in E \cap p^\varrho$. We have by double counting and (3) of Definition~\ref{def: configuration}, 
\[
\sum_{ p\in E} |\TT(p, \varrho)|_\delta \le \sum_{T\in \TT} |E \cap N_{C\varrho }T|_\delta \lesssim_{\d^{O(\eta)}} |\TT|_\delta |Z|_\delta (\varrho/\d)^{2t-\kappa} \lesssim_{\d^{O(\eta)}} |\TT|_\delta |E|_\delta  \varrho^{2t-\kappa}
\]
and so by $(2t, \d^{-\eta})$-AD-regularity of $\TT$,
\[
\sum_{ p\in E} |\TT(p, \varrho)|_\varrho \lesssim_{\d^{O(\eta)}} |E|_\delta \varrho^{-\kappa}.
\]
So we can find a subset $E' \subset E$ such that $|E'| \ge \d^{O(\eta)} |E|$ and $|\TT(p, \varrho)|_\varrho  \sim_{\d^{O(\eta)}} \varrho^{-\kappa}$ for every $p \in E'$, where the lower bound comes from (4) of Definition~\ref{def: configuration}. It follows that after replacing $E$ with $E'$ we still have $|\TT(p, \varrho)|_\varrho  \sim_{\d^{O(\eta)}} \varrho^{-\kappa}$. 
Since $|\TT_p|_\varrho \gtrsim \d^{\eta} \varrho^{-\kappa}$ for all $p \in E$, it follows by Cauchy--Schwarz that for every $p_0 \in E$ we have
\[
\sum_{p_1, p_2 \in E \cap p_0^\varrho} |\Theta_{p_1} \cap \Theta_{p_2}^{\varrho}|_\d \gtrsim_{\d^{O(\eta)}} |E \cap p_0^\varrho|^2 \d^{-\kappa}.
\]
So we can find some $p_0' \in E\cap p^\varrho_0$ so that the set
\[
E(p^{\varrho}_0) =\{  p \in E \cap p_0^\varrho: |\Theta_{p} \cap \Theta_{p_0'}^{\varrho}|_\d \ge \d^{C_0\eta} |\Theta_p|_\d \}
\]
has size at least $\d^{O(\eta)} |E \cap p_0^\varrho|$ (for an appropriately large constant $C_0$). Taking $E'$ to be the union of sets $E(p^{\varrho}_0)$ over all $p_0$ representing distinct dyadic $\varrho$ boxes gives (a). 
Part (b) follows from (a) by simple double counting.
\end{proof}

The following lemma will be applied to the shading $Z$ to find scales $\beta$ so that $Z \cap z^{(\beta)}$ is a rescaled $(\beta, s, \beta^{-\varepsilon})$-set. Essentially the same statement appears in \cite[Proposition 4.32]{orponen2023projections}.

\begin{lemma}\label{lem:good-scales} 
    For any $\epsilon>0$ and $s \in (0,1]$ there exists $\mu, \eta_0, \delta_0>0$ such that the following holds for all $\eta< \eta_0$ and $\delta < \delta_0$. 
    Let $A \subset [0,1]$ be a $(\delta, s, \delta^{-\eta})$-set and suppose that $A$ is $(\delta, \delta^{-\eta})$-uniform. Then there exists $\beta \in [\delta^{1/2}, \delta^{\mu}]$ such that $A \cap B(a, \beta)$ is a rescaled $(\beta, s, \beta^{-\epsilon})$-set for every $a\in A$. That is, for every $a\in A$ and $\varrho \in [\beta^2, \beta]$ we have
    \[
    |A\cap B(a, \varrho)|_{\beta^2} \le \beta^{-\epsilon} (\varrho/\beta)^s |A\cap B(a, \beta)|_{\beta^2}.
    \]
\end{lemma}

\begin{proof}
    Let $f: [0,1] \rightarrow [0,1]$ be the branching function of $A$, i.e. 
    \[
    |A|_{\delta^x} = \delta^{-f(x)}.
    \]
    Note that since $A$ is $(\d, \d^{-\eta})$-uniform, $f(x)$ is a meaningful function up to precision $\eta$. 
    Since $A$ is $(\delta, s, \delta^{-\eta})$-set we have
    \[
    f(x) \ge s x - \eta -o(1)
    \]
    where $o(1)$ tends to 0 with $\delta \rightarrow 0$ (in fact, it decays as $\frac{c}{\log(1/\delta)}$). This term comes from a constant factor loss which occurs when we compare the definition of covering number using dyadic intervals and the Frostman set condition using neighborhoods. 
    
    Since for any $\varrho \le \tau$ we have $|A|_{\varrho} \le |A|_{\tau} \le C (\tau/\varrho) |A|_\varrho$, we also know that $f(x) \le f(x+h) \le f(x)+h+o(1)$ for any $x+ h\le 1$. 
    Let $g(x) = f(x) - sx$, so that $g(x) \ge -\eta+o(1)$ for all $x\in [0,1]$ and $g$ is a 1-Lipschitz function (up to precision $o(1)$). We say that a scale $b \in [0, 1/2]$ is {\em $\epsilon$-good} if
    \[
    f(b+x) \ge f(b) + xs - \epsilon b
    \]
    for all $x \in [0,b]$. Equivalently, $b$ is $\eps$-good if
    \[
    g(x) \ge g(b) - \epsilon b/2
    \]
    for every $x \in [b, 2b]$. 

    Let $\mu = 2^{-1-M}$, where $M = \lceil 4\epsilon^{-1} \rceil$ and suppose that $\eta \le \mu/2$.
    Suppose there is no $\epsilon$-good scale $b \in [\mu, 1/2]$. Then we can let $b_0 = \mu$ and construct a sequence $b_i$, $i=0, 1, \ldots$, such that for $i\ge 0$:
    \begin{itemize}
        \item $b_{i+1} \le 2 b_i$,
        \item  $g(b_{i+1}) \le g(b_i) - \epsilon b_i/2$.
    \end{itemize}

    Thus, after $M$ steps we obtain some
    $b_M \le 2^M \mu \le 1/2$ and
    \[
    g(b_{M}) \le g(b_{M-1}) - (\epsilon/2) b_{M-1} \le g(b_{M-2}) - (\epsilon/2)(b_{M-1}+b_{M-2}) \le \ldots \le g(b_0) - (\epsilon/2) \sum_{j=0}^{M-1} b_j
    \]
    Now we use $g(b_0) \le b_0+o(1)$ and $b_j \ge b_0$ to conclude that
    \[
    -\eta+o(1) \le g(b_M) \le b_0 - (\epsilon/2) M b_0 + o(1) \le - b_0 +o(1) \le -\mu+o(1) 
    \]
    and so for sufficiently small $\delta$ we obtain a contradiction with the condition $\eta \le \mu/2$.

    Thus, there is a good scale $b \in [\mu,1/2]$. Unraveling the definition, it follows that for every $x \in [0,b]$:
    \[
    |A|_{\d^{b+x}} \ge \d^{\eps b/2}\d^{-xs}|A|_{\d^{b}}.
    \]
    Let $\beta = \d^{b}$. Since $A$ is $(\d, \d^{-\eta})$-uniform, this implies that $A\cap B(a,\beta)$ is a rescaled $(\beta, s, C\beta^{-\eps/2} \d^{-\eta})$-set. For small enough $\eta$ and $\d\le \d_0$, we have $C\d^{-\eta} < \d^{-\eps \mu/2} \le \beta^{-\eps/2}$ giving the result. 
\end{proof}

\section{Construction of Furstenberg configurations: Proof of Proposition \ref{prop: reduction}}\label{sec:configuration-furstenberg}

Let $\kappa$ be the infimum of all $\kappa'$ such that $\operatorname{F}_{s,t}(\kappa')$ holds. Then for any $\nu>0$, $\operatorname{F}_{s,t}(\kappa-\nu)$ does not hold. 
Therefore, for every $\nu>0$ there exists $\eps>0$ such that for all sufficiently small $\eta>0$ and $\delta>0$, there exists  a $(2t, \delta^{-\eta})$-AD-regular set of $\delta$-tubes with $C_{t-CW}(\TT)\leq \delta^{-\eta}$ and a $(\delta, s, \delta^{-\eta})$-shading $Y$ on $\TT$ with $|Y(T)|_{\delta}\sim M$ for every $T\in \TT$ and 
\begin{equation*}
|U(\TT, Y)|_{\delta}\leq \delta^{\epsilon-\nu}M \delta^{-2t+\kappa}\leq M\delta^{-2t+\kappa-\nu}.
\end{equation*}
Since the properties of Furstenberg configurations are monotone in $\eta$, we may replace the parameter $\eta$ by $\max(\eta, \nu)$, i.e. without loss of generality,
\begin{equation}\label{eq: extremal}
|U(\TT, Y)|_{\delta}\leq \d^{-\eta} M\delta^{-2t+\kappa}.
\end{equation}
Equivalently, 
\begin{equation}\label{eq: extremalmu}
    \mu(\TT, Y):=\frac{\sum_{T\in \TT}|Y(T)|_{\delta}}{|U(\TT, Y)|_{\delta}}\gtrsim \delta^{-\kappa+\eta}.
\end{equation}
We are going to show that for any $\delta_0$ and $\eta_0$, there exists $\delta_1<\delta_0$ and $\gamma\in  [0, \kappa]$ so that there exists a $(\delta_1, \eta_0, s, t, \kappa, \gamma)$-Furstenberg configuration. We are going to take $\eta, \delta$ sufficiently small compared to $\eta_0, \delta_0$, respectively.

\medskip 
\noindent \textit{Step 1. Define point line configuration and obtain uniform refinement.}

By using Proposition \ref{prop:regularize-dyadic-tubes} we can replace $\TT$ by a large subset and increase $\d$ slightly so that the resulting collection of dyadic tubes has the property that $N_\varrho T^\varrho \subset T^{2^S \varrho}$ for all $T^\varrho \in \TT_\varrho$ and $\varrho \in [\d,2^{-S}]$, where $S \sim \log \log (1/\d)$. After some further pigeonholing, we can ensure that $\TT$ is still $(2t, \d^{-2\eta})$-AD-regular set (assuming that $\d$ is sufficiently small in terms of $\eta$).

Define $X\subset \Omega_d$ as the  set of $\delta$-separated  triples $(x,t, \theta)\in [-1,1]^{d-1}\times [-1, 1]\times [-1, 1]^{d-1}$ such that $(x,t)\in U(\TT, Y)^\d \cap (\delta\mathbb{Z})^{d}$ and $\theta=\theta(T)$ for some $T\in \TT_Y(p)$, $p=(x,t)$. By Observation~\ref{obs:pigeonhole-slope}, we can assume all pairs $(p, T), T\in \TT_Y(p)$ can be represented in $X$. 

Apply Lemma~\ref{lem:uniform-subset} to $X$ with  $\Delta_j =\delta^{j\eta}$, $j=1, \dots, \eta^{-1}$ to find a subset $X'\subset X$, $|X'|_{\delta}\gtrapprox |X|_{\delta}$ and $X'$ is $\gtrapprox 1$-uniform on scales $\Delta_j$, $j=1\dots, \eta^{-1}$, which means that for all $j \leq i+k$ and any $\omega_0=(p_0, \ell_0)\in X'$, we have 
\[
\M_{X'}(\Delta_i, \Delta_j, \Delta_k)\geq |R_{\Delta_i, \Delta_j, \Delta_k}(\omega_0) \cap X'|
\gtrapprox \M_{X'}(\Delta_i, \Delta_j, \Delta_k)
\]
where $R_{v,h,\alpha}(\omega_0)$ is the phase-space rectangle around $\omega_0$ defined in Section \ref{sec:uniform}.
Here and in the following formulas, the notation $\approx $, $\gtrapprox$ hides logarithmic factors of the form $(\log(1/\d))^{O_\eta(1)}$, coming from the bounds in Lemma \ref{lem:uniform-subset}.

\medskip
\noindent \textit{Step 2. Obtain a refined shading $Y'$ and show it is uniform.}

For each $T\in \TT$, define shading $Y'$ on $T$ as the union of $\delta$-balls centered at points $p$ satisfying $(p, \theta(T))\in X'$. 
Let $\TT'$ be the set of $T\in \TT$ such that there exists $p$, $(p, T)\in X'$. The uniformity of $X'$ yields a lot of uniformity information on $(\TT', Y')$. 

For any $\delta \leq a \leq 1$ and $T\in \TT', x\in Y'(T),$ $T_a=T^a$, we claim that
\begin{equation}\label{eq:Mdda}
|\TT'[T_a]_{Y'}(x)|_{\delta}\approx  \M_{X'}(\delta, \delta,a).    
\end{equation}
Indeed, the rectangle $R_{\d,\d,a}(\omega_0)$ for $\omega_0 = (x_0,t_0,\theta_0)$, $p_0 = (x_0,t_0)$ captures point-line pairs $(p,\ell)$ such that $d(p,p_0) \lesssim \d$ and $\angle(\ell, \ell_0) \lesssim a$. Thus, $|R_{\d,\d,a}\cap X'|$ is approximately the number of tubes $T \in \TT'_{Y'}(p_0)$ making angle $\lesssim a$ with $T$. 
Use the conclusion of Proposition \ref{prop:regularize-dyadic-tubes} which makes the dyadic structure and the Euclidean distance compatible with each other, we get \eqref{eq:Mdda}. In what follows, we use concentration numbers $\M_{X'}$ to compute other quantities related to $(\TT',Y')$ in a similar manner. 

To ease the notation, let $N_{X'}(v, h, \alpha):=\frac{\M_{X'}(v, h, \alpha)}{\M_{X'}(\delta, \delta, \alpha)}$. Observe that for any $p\in Y'(T)$ and $\varrho\in (\delta,1)$, 
\begin{equation*}\label{eq: shadingsegment}
|Y'(T)\cap B(p,\varrho)|_{\delta} \approx \frac{ \M_{X'}(\varrho, \delta, \delta/\varrho)}{\M_{X'}(\delta, \delta, \delta/\varrho)}=:N_{X'}(\varrho, \delta, \delta/\varrho). 
\end{equation*}
Indeed, by uniformity, every $\d$-cube $p$ contained in  $Y'(T)\cap B(p,\varrho)$ contributes $|\TT'[T_{\d/\varrho}]_{Y'}(p)|\approx \M_{X'}(\delta, \delta, \delta/\varrho)$ many cube-tube incidences to the count $\M_{X'}(\varrho, \delta, \delta/\varrho)$. 
Therefore $Y'(T)$ is $\approx \delta^{\eta}$-uniform. 

By the definition of $\TT'$ and $Y'$, for any $T\in \TT$, if $Y'(T)\neq \emptyset$, then $T\in \TT'$. For any $\varrho \in [\delta, 1]\cap 2^{-\mathbb{N}}$, define $\TT_\varrho'= \cD_\varrho (\TT')$.  Then for any $T_\varrho\in \TT_\varrho'$, 
\[
\sum_{T\in \TT'[T_\varrho]} |Y'(T)|_{\delta} \approx \M_{X'}(1, \varrho, \varrho).
\]
Notice that by definition $\frac{|X'|}{\M_{X'}(1,\varrho,\varrho)} \approx |\TT'|_\varrho$.
Since $Y'$ is a $\gtrapprox 1$-refinement of $Y$, and $\TT$ is $( 2t, \delta^{-2\eta})$-AD-regular,   we have 
\begin{equation*}\label{eq: TT'AD}
 \delta^{2\eta} M (\delta/\varrho)^{-2t} \lessapprox   \sum_{T\in \TT'[T_\varrho]} |Y'(T)|_{\delta} \lessapprox \delta^{-2\eta} M (\delta/\varrho)^{-2t}. 
\end{equation*}
In particular, $\TT'$ is $(2t, \delta^{-4\eta})$-AD-regular. 



Recall that we define $Y'(T_\varrho)= \cD_\varrho(U(\TT'[T_\varrho], Y'))$. Then for any $b\in [\varrho, 1]$ of the form $\Delta_j$,  
\begin{equation*}
     |Y'(T_\varrho)|_b \approx \frac{\sum_{T\in \TT'[T_\varrho]}|Y'(T)|_{\delta}}{\M_{X'}(b, \varrho, \varrho)} \approx \frac{ \M_{X'}(1,\varrho, \varrho)}{\M_{X'}(b, \varrho, \varrho)}
\end{equation*}
and for $x \in Y'(T_\varrho)$
\begin{equation}\label{eq:Y-rho-shading}
|Y'(T_\varrho) \cap B(x,b)|_\varrho \approx \frac{\M_{X'}(b, \varrho, \varrho)}{\M_{X'}(\varrho, \varrho, \varrho)}.    
\end{equation}
This shows that $Y'(T_\varrho)$ is $\approx \delta^{\eta}$-uniform for any $\varrho \in [\delta, 1]$ and $T_\varrho \in \TT'^\varrho$.

\medskip

\noindent \textit{Step 3. $Y'(T_\varrho)$ is a $(\varrho, s, \delta^{-4\eta})$-set.}
For any $\varrho<r<1$ of the form $\Delta_j$ and $x\in Y'(T_\varrho)$, we want to show  \begin{equation}\label{eq: YTrhospacing} |Y'(T_\varrho)\cap B(x,r)|_{\varrho} \lessapprox r^s \delta^{-\eta} |Y'(T_\varrho)|_\varrho.
\end{equation} 
By (\ref{eq:Y-rho-shading}), it suffices to show 
\begin{equation}
    \M_{X'}(r, \varrho, \varrho) \lessapprox r^s \delta^{-\eta}  M_{X'}(1, \varrho, \varrho). 
\end{equation}
Indeed, since $Y'(T)$ is uniform,  $|Y'(T)|_{\delta}\gtrapprox |Y(T)|_{\delta}$,  and $Y(T)$ is a $(\delta, s, \delta^{-\eta})$-set, we get that $|Y'(T)|_r\gtrapprox \delta^{\eta} r^{-s}$. Therefore, 
\begin{align*}
    \M_{X'}(r, \varrho, \varrho) & \approx \sum_{T\in \TT'[T_\varrho]} \frac{|Y'(T)|_\delta}{|Y'(T)|_{r}}\\
    &\lessapprox r^s \delta^{-\eta}  \sum_{T\in \TT'[T_\varrho]} |Y'(T)|_{\delta} \approx r^s \delta^{-\eta} \M_{X'}(1, \varrho, \varrho). 
\end{align*}
For any $r\in [\varrho, 1]$, we approximate it with the closest $\Delta_j$ to show \eqref{eq: YTrhospacing} with $4\eta$ in place of $\eta$.

\medskip
\noindent \textit{Step 4. Apply $\operatorname{F}_{s,t}(\kappa)$ on $(\TT[T_\varrho], Y')$ and $(\TT_\varrho, Y')$ to show extremal multiplicities.}

Let $\eps_1\in (0, \eta_0)$ be a small parameter to be determined. Since $\operatorname{F}_{s,t}(\kappa+\eps_1/2)$ holds, the following inequalities hold for $\eta>0$ sufficiently small and $\varrho>0$ sufficiently small, 
\begin{equation}
    |U(\TT_\varrho, Y')|_{\varrho} \geq  \varrho^{\eps_1} |Y'(T_\varrho)|_\varrho \varrho^{-2t+\kappa}, 
\end{equation}
or equivalently 
\begin{equation}\label{eq: muthicktubes0}
\mu(\TT_\varrho, Y') \leq \varrho^{-\eps_1-\eta-\kappa}.
\end{equation}
Similarly, apply $\operatorname{F}_{s,t}(\kappa+\eps_1/2)$ for $\delta/\varrho>0$ sufficiently small, 
\begin{equation}\label{eq: muthintubes0}
    \mu(\TT[T_\varrho], Y')\leq (\delta/\varrho)^{-\eps_1-\eta-\kappa}.
\end{equation}
Since $Y'$ is a $\gtrapprox 1$-refinement of $Y$, 
\[
\mu(\TT, Y')\gtrapprox \mu(\TT, Y)\gtrsim \delta^{-\kappa+\eta}.
\]
Since $\mu(\TT, Y')\lesssim \mu(\TT[T_\varrho], Y') \mu(\TT_\varrho, Y')$, combining with \eqref{eq: muthicktubes0} and  \eqref{eq: muthintubes0}, we have 
\begin{equation}\label{eq: muthicktubes}
     \delta^{\epsilon_1+2\eta} \varrho^{-\kappa} \leq \mu(\TT_\varrho, Y') \leq \varrho^{-\eps_1-\eta-\kappa}
\end{equation}
and 
\begin{equation}\label{eq: muthintubes}
        \delta^{\epsilon_1+2\eta} (\delta/\varrho)^{-\kappa} \leq \mu(\TT[T_\varrho], Y')\leq (\delta/\varrho)^{-\eps_1-\eta-\kappa}.
\end{equation}
Recall that by the definition of $Y'$ and $\TT'$, if $T\in \TT$ has nonempty shading $Y'(T)$, then $T\in \TT'$. 
For any $\delta\leq r<\varrho  \leq 1$, and $T_\varrho\in \TT_\varrho'$, 
\begin{equation}\label{eq: mutubes}
   \delta^{2\eps_1+4\eta}  (r/\varrho)^{-\kappa} \leq \mu(\TT_r'[T_\varrho], Y')=\mu(\TT_r[T_\varrho], Y') \leq    \delta^{-2\eps_1-4\eta}  (r/\varrho)^{-\kappa}. 
\end{equation}

\medskip
\noindent \textit{Step 5. Pigeonhole to find $Z$ at suitable scales.}

Let $\eta_1\in (\eta,1)\cap 2^{-\mathbb{N}}$ be a small parameter and let $\delta_j = \delta^{j\eta_1}, j=1, \dots, \eta_1^{-1}$. By pigeonholing, there exists $j$ such that 
\begin{equation}\label{eq: Z}
    |Y'(T_{\delta_{j+1}})|_{\delta_1} \gtrapprox \delta_1^{\eta_1} |Y'(T_{\delta_j})|_{\delta_1}.
\end{equation}

Let $\tilde\TT = \TT_{\delta_{j+1}}^{T_{\delta_j}}$ be the anisotropic rescaling of $\TT_{\delta_{j+1}}[T_{\delta_j}]$ and $\tilde Y$ be the rescaling of $Y'$. Let $Z$ be the set of $z$ such that $\mathbb{R}^2 \times \{z\}\cap U(\tilde\TT, \tilde Y)\neq \emptyset$. Since $Y'(T_{\delta_j})$ is $\approx \delta^{\eta}$-uniform, $Z$ is $\approx \delta^{\eta}$-uniform. 

\medskip
\noindent \textit{Step 6. AD-regularity of horizontal slices.}

The set $E=U(\tilde\TT, \tilde Y)$ is of the form $E=\cup_{z\in Z} E_z$, $E_z=A_z\times \{z\}$ for some $A_z\subset [-1, 1]^2$. We would like to show $A_z$ is a $(2t-\kappa, \delta_1^{-\eta_0})$-AD regular set.  To see this, for any $z\in Z, x\in A_z$ and $r\in [\delta_1, 1]$,
\begin{equation} \label{eq: Azregular}
|A_z\cap B(x, r)|_{\delta_1} 
\approx \frac{ \M_{X'}(\delta_1, r\delta_j, \delta_j) }{\M_{X'}(\delta_1, \delta_{j+1}, \delta_j)}.
\end{equation}

We are not going to use the RHS of \eqref{eq: Azregular} to calculate $|A_z\cap B(x,r)|_{\delta_1}$, what is important for us is the fact that \eqref{eq: Azregular} holds for any $z\in Z, x\in A_z$ and $r\in [\delta_1, 1]$. 

 Recall $\TT'$ is $(2t, \delta^{-2\eta})$-AD-regular and $|Y'(T_\varrho)|_\varrho$ is approximately the same for each $T_\varrho \in \TT_\varrho'$, $\varrho$ of the form $\Delta_j$. By \eqref{eq: mutubes},  and $\mu(\tilde \TT, \tilde Y)=\mu(\TT_{\delta_{j+1}}'[T_{\delta_j}], Y')$,  
 \begin{equation}\label{eq: Azupper}
 |A_z|_{\delta_1} \approx \frac{|E|_{\delta_1}}{|Z|_{\delta_1}} \lesssim   \frac{|E|_{\delta_1}}{|\tilde Y(\tilde T)|_{\delta_1}}\approx  \frac{ \sum_{\tilde T\in \tilde\TT} |\tilde Y(\tilde T)|_{\delta_1} }{ \mu(\tilde \TT, \tilde Y)  |\tilde Y (\tilde T)|_{\delta_1}} \lessapprox  \delta^{-2\epsilon_1-4\eta} \delta_1^{-2t+\kappa}. 
 \end{equation}

By \eqref{eq: Z}, we have $|Z|_{\delta_1}\lessapprox \delta_1^{-\eta_1} |\tilde Y (\tilde T)|_{\delta_1}$. Combining with  \eqref{eq: mutubes}, 
 \begin{equation*}
      |A_z|_{\delta_1} \approx \frac{|E|_{\delta_1}}{|Z|_{\delta_1}} \gtrapprox \delta_1^{\eta_1}  \frac{|E|_{\delta_1}}{|\tilde Y(\tilde T)|_{\delta_1}}\approx \delta_1^{\eta_1} \frac{ \sum_{\tilde T\in \tilde\TT} |\tilde Y(\tilde T)|_{\delta_1} }{ \mu(\tilde \TT, \tilde Y)  |\tilde Y (\tilde T)|_{\delta_1}} \gtrapprox \delta_1^{\eta_1} \delta^{2\epsilon_1+4\eta} \delta_1^{-2t+\kappa}. 
 \end{equation*}

Similarly, for any $r\in [\delta_1, 1]$, 
\begin{equation}\label{eq: lowerbdAzr}
|A_z|_{r} \gtrapprox \delta_1^{\eta_1}\frac{|E|_{r}}{|\tilde Y (\tilde T)|_r} \gtrapprox \delta_1^{\eta_1} \delta^{2\epsilon_1+4\eta} r^{-2t+\kappa}. 
\end{equation}

For any $T_r\in \tilde \TT_r$, since $Z$  and $\tilde Y(\tilde T)$ are $\approx \delta^{\eta}$-uniform, $|Z|_{\delta_1/r} \lessapprox \delta_1^{-\eta_1} \delta^{-2\eta} |\tilde Y(\tilde T)|_{\delta_1/r}$, 
\begin{equation*}
    |A_z\cap B(x,r)|_{\delta_1} \gtrsim  \delta_1^{\eta_1} \delta^{2\eta} \frac{\sum_{\tilde T\in \tilde \TT[T_r]} |\tilde Y (\tilde T)|_{\delta_1/r}}{\mu(\tilde \TT[T_r], \tilde Y) |\tilde Y (\tilde T)|_{\delta_1/r}} \gtrapprox \delta_1^{\eta_1} \delta^{2\epsilon_1+6\eta} (\delta_1/r)^{-2t+\kappa}. 
\end{equation*}
Combining with \eqref{eq: lowerbdAzr} and \eqref{eq: Azupper}, we have 
\begin{equation}\label{eq: Azxr}
    \delta_1^{\eta_1} \delta^{2\epsilon_1+6\eta} (\delta_1/r)^{-2t+\kappa}\lessapprox     |A_z\cap B(x,r)|_{\delta_1} \lessapprox \delta_1^{-2\eta_1} \delta^{-4\epsilon_1-12\eta} (\delta_1/r)^{-2t+\kappa}. 
\end{equation}
This shows that $A_z$ is $(2t-\kappa, \delta_1^{-\eta_0})$-AD regular if we choose $\eta_1>0, \eps_1, \eta>0$  sufficiently small depending on $\eta_0$ and $\delta>0$ sufficiently small depending on all previous parameters.

\medskip

\noindent \textit{Step 7. Apply structure of AD-regular set.}

So far we have shown (1), (2), (3) from Definition \ref{def: configuration}. Now we obtain (4). 
Let $\tilde \d = \d^{\eta_1} = \d_{j+1}/\d_j$. 
For $p \in E$ let $\Theta_p = \{ \theta(T):~ T \in \tilde\TT_{\tilde Y}(p) \}$. By (\ref{eq: muthicktubes}) we have that $\Theta_p$ is $(\tilde \d, \kappa, \d^{-2\eps_1-4\eta})$-AD-regular set of directions. Apply Lemma \ref{lem:finding-parallel-rectangles} to each set $\Theta_p$ with some $\zeta>0$ to be determined. We obtain some scales $\varrho, \tau \in [\tilde \d, 1]$ and $\Theta'_p \subset \Theta_p$, so that $\tau/\varrho \ge \tilde\d^{-\chi(\zeta)}$, $|\Phi_p'|_{\tilde \d} \ge (\varrho/\tau)^\zeta |\Phi_p|$ and some $\gamma \in [0, \kappa]$ such that $\phi_{a,\tau}(\Theta'_p \cap B(\theta,\tau))$ is nearly $(\varrho/\tau, \kappa, \gamma, (\varrho/\tau)^{-\zeta})$-aligned for every $\theta \in \Theta'_p$. By pigeonhole and further refinement, we may assume that parameters $\varrho, \tau, \gamma$ do not depend on the choice of $p \in E$. We may also assume that $\phi_{a, \tau}$ is a dilation. By incurring a $\d^{-C\eta}$ loss, we may assume that $\varrho, \tau \in \d^{\eta \ZZ}$. 

Define a refinement $\tilde Y' \subset \tilde Y$ by setting $\tilde Y'(T)$ be the set of $p \in E$ such that $\theta(T) \in \Theta'_p$. Then $\tilde Y'$ is $\gtrapprox (\varrho/\tau)^{\zeta}$-refinement of $\tilde Y$. 
Let $\hat \TT = \tilde \TT_{\varrho}^{\tilde T_\tau}$ for some $\tilde T\in \tilde \TT$ and let $\hat Y$ be the rescaling of $\tilde Y'$. Then for every $p \in U(\hat \TT, \hat Y)$ we have a corresponding set $\hat\Theta_p$ which is nearly $(\varrho/\tau, \kappa, \gamma, (\varrho/\tau)^{-\zeta} \d^{-C\eta})$-aligned and such that for every $\theta\in \hat\Theta_p$ there is $T \in \hat \TT_{\hat Y}(p)$ with $|\theta(T) - \theta| \lesssim \varrho/\tau$. By slightly modifying the set $\hat\Theta_p$ we may assume that it is $(\varrho/\tau, \kappa, \gamma, (\varrho/\tau)^{-\zeta} \d^{-C\eta})$-aligned (not just nearly). This gives (4). 

Let $\hat\d = \varrho/\tau$. 
Since $\hat\d \ge \tilde \d^{\chi(\zeta)} = \d^{\eta_1 \chi(\zeta)}$, we can choose $\d$ small enough so that $\hat\d \le \d_0$. 
By taking $\zeta>0$ small enough, we also get $(\varrho/\tau)^{-\zeta} \d^{-C\eta} < \hat\d^{-\eta_0}$. By similar considerations, we verify that $(\hat \TT, U(\hat \TT, \hat Y), \{\hat \TT_p\})$ has the properties (1), (2), (3), (4) in Definition \ref{def: configuration}. 

The conclusion about $C_{t-CW}(\hat \TT)\le \hat \d^{-\eta_0}$ follows from Observations~\ref{obs:rescaling-AD-regular-tubes}, \ref{obs:thickening-AD-regular-tubes} and Observation~\ref{obs:wolff-axiom-rescaling} provided that $\epsilon_1,\eta \ll \eta_0$ are small enough.

\section{High-low estimate: Proof of Proposition \ref{prop: high-low}}\label{sec:high-low}

In the following proposition we use high-low method to prove a restricted projection estimate in the case when the set of directions $\Theta$ has strong non-concentration in rectangles. 
This result is similar to the incidence estimate \cite[Proposition 2.1]{guth2019incidence}. The main difference is that the restriction on the direction set $\Theta$ allows us to more effectively estimate the high frequency part.

\begin{prop}\label{prop:small-gamma-range}
    Let $\alpha \in (0,1)$. Then for any $\varepsilon>0$ the following holds for $\eta \le \eta_0(\varepsilon, \alpha)$ and $\d \le \d_0(\eta,\alpha)$.

    Let $\Theta \subset [-1,1]^2$ be a $\d$-separated set such that for any $\d\times 1$ tube $T$ we have $|\Theta \cap T| \le \nu |\Theta|$ for some $\nu \in (0,1)$. 
    For each $\theta\in \Theta$ let $A_\theta \subset [-1,1]^2$ be a $\d$-separated set of size $\sim M$ for some $M\ge 1$ and such that $|A_\theta \cap B_{S\d}| \le \d^{-\eta} S^{2-\alpha}$ for any $S\ge 1$.
    For each $a\in A_\theta$ let $E(a, \theta) \subset ((a,0) + \RR (\theta,1)) \cap [-1,1]^3$ be a $\d$-separated set of size $\sim \lambda \d^{-1}$ for some $\lambda \in (0,1)$. Let $E$ be the union of $E(a, \theta)$ over all $a\in A_\theta$ and all $\theta\in \Theta$ and suppose that $|E|_\d \le \d^{-\eta} M \lambda \d^{-1}$.

    Then we have $\lambda \le \d^{-\varepsilon}\nu$.
\end{prop}

\begin{proof}
    Let $S = \d^{-\eta'}$ for some $\eta'\gg \eta$ to be chosen later. For $a\in A_\theta$ and $\theta \in \Theta$ let $f_{\theta, a}$ be a smooth bump function essentially supported on the $\d \times 2$ tube $T_{\theta,a}$ around the line $(a,0)+\RR(\theta,1)$ and suppose that $f_{\theta,a}(x) \sim 1$ for $x \in T_{\theta,a}$. If we let $E_\d$ denote the $\d$-neighborhood of the set $E$, then we have 
    \[
    \int_{E_\d} f_{\theta,a} \ge \int_{E(\theta,a)_\d} f_{\theta,a} \gtrsim \operatorname{Vol}(E(\theta,a)_\d) \gtrsim \lambda \d^2
    \]
    and so we get 
    \[
    \int_{E_\d} \sum_{a,\theta} f_{\theta,a} \ge \d^{O(\eta)} (\lambda \d^2)|\Theta| M \ge \d^{O(\eta)} |\Theta| \operatorname{Vol}(E_\d).
    \]
    On the other hand, we can decompose $f_{\theta, a} = f^{high}_{\theta, a}+ f^{low}_{\theta, a}$ where $\widehat{f}^{low}_{\theta,a}$ is essentially supported on $B_{(S\d)^{-1}}(0)$ and $\widehat{f}^{high}_{\theta,a}$ is essentially supported on $B_{\d^{-1}}\setminus B_{(S\d)^{-1}}$.

    First, we estimate the low contribution to the integral as follows. Note that $f^{low}_{\theta,a}$ is essentially supported on the $\d^{-o(1)} S\d$-tube around the line $(a,0)+\RR(\theta,1)$ (where $o(1)$ is much smaller than $\eta$ and comes from the choice of the smooth bump function) and has value $\lesssim \d^{o(1)} S^{-2}$ at every point.
    Fix $\theta \in \Theta$ and cover $A_\theta$ by finitely overlapping $ \d^{-o(1)} S\d$-balls $B_j$. By assumption, $|A_\theta \cap B_j|\lesssim \d^{-\eta} S^{2-\alpha}$ for each ball $B_j$. Let $T_j$ be the $\d^{-o(1)} S\d\times 4$ tube above $B_j$ and direction $(\theta,1)$. 
    Then we can estimate
    \begin{align*}
    \int_{E_\d} \sum_{a\in A_\theta} |f^{low}_{\theta, a} | &\lesssim \sum_j \int_{E_\d \cap T_j} \sum_{a\in A_\theta \cap B_j} |f^{low}_{\theta, a} | \\&\lesssim \d^{-O(\eta)} \sum_j \operatorname{Vol}(E_\d \cap T_j)  \d^{o(1)} S^{-\alpha} \lesssim \d^{-O(\eta)}  \operatorname{Vol}(E_\d) S^{-\alpha}.    
    \end{align*}
    By choosing $\eta' \ge \frac{C\eta}{\alpha}$ we can ensure that this is less than $\frac12 \int_{E_\d} \sum_{a,\theta} f_{\theta,a}$. It follows that 
    \begin{equation}\label{eq:lower-bound-high}
    \left|\int_{E_\d} \sum_{a,\theta} f^{high}_{\theta,a}\right| \ge \frac{1}{2}\int_{E_\d} \sum_{a,\theta} f_{\theta,a} \gtrsim \d^{O(\eta)} |\Theta| \operatorname{Vol}(E_\d).    
    \end{equation}
    Denote $f_\theta = \sum_{a\in A_\theta} f_{\theta, a}$. 
    By Cauchy--Schwarz we have
    \[
    \left|\int_{E_\d} \sum_{\theta} f^{high}_{\theta}\right|^2 \le \operatorname{Vol}(E) \int_{E_\d} \left|\sum_{\theta} f_{\theta}^{high}\right|^2 \le \operatorname{Vol}(E) \int_{\RR^3} \left|\sum_{\theta} f_{\theta}^{high}\right|^2 = \operatorname{Vol}(E) \int_{\RR^3} \left|\sum_{\theta} \widehat{f}_{\theta}^{high}\right|^2.
    \]
    The function $\widehat{f}_\theta^{high}$ is essentially supported on the region $\Pi_\theta =\{x \in \RR^3:~ |x| \in (\frac{1}{S\d}, \frac{1}{\d}), ~ |x \cdot (\theta,1)| \lesssim 1 \}$. So $\Pi_\theta$ is essentially a $1\times \d^{-1}\times \d^{-1}$ slab with the middle $1\times (S\d)^{-1}\times (S\d)^{-1}$ subslab removed. Consider the overlap pattern of the slabs $\Pi_\theta$ for $\theta \in \Theta$: note that the set of directions $\theta$ such that $x \in \Pi_\theta$ for a fixed $x \in \RR^3$ is contained in a strip of width $\sim S\d$. Using the restriction of $\Theta$, we get that the overlap between the slabs $\Pi_\theta$ is at most $\lesssim S \nu |\Theta|$, giving
    \[
    \int_{\RR^3}\left|\sum_{\theta} \widehat{f}_{\theta}^{high}\right|^2  \lesssim  S \nu |\Theta| \sum_\theta \int_{\RR^3} \left|\widehat{f}_{\theta}^{high}\right|^2 + \text{Rap.Dec.}
    \]
    So since $f^{high}_{\theta} = \sum_a f^{high}_{\theta, a}$ and functions $f^{high}_{\theta, a}$, $a \in A_\theta$, are physically separated, we conclude that
    \[
    \left|\int_{E_\d} \sum_{\theta, a} f^{high}_{\theta, a}\right|^2 \lesssim S \nu |\Theta| \sum_{\theta, a} \int_{\RR^3} \left|f^{high}_{\theta, a}\right|^2 \lesssim \d^{-O(\eta)} S\nu |\Theta|^2 M \operatorname{Vol}(T_{\theta, a}) \lesssim \d^{-O(\eta)} S\nu |\Theta|^2\lambda^{-1} \operatorname{Vol}(E_\d)
    \]
    So combining with (\ref{eq:lower-bound-high}) gives
    \[
    \d^{O(\eta)} |\Theta|^2 \operatorname{Vol}(E_\d)^2\lesssim \d^{-O(\eta)} S\nu|\Theta|^2\lambda^{-1} \operatorname{Vol}(E_\d)^2
    \]
    for sufficiently small $\eta$, we can take the error term to be less than $\d^{-\varepsilon}$, giving the desired upper bound on $\lambda$.
\end{proof}

\begin{proof}[Proof of Proposition \ref{prop: high-low}]
    Let $\eta_0>0$ be a small parameter to be determined and $\eta>0$ sufficiently small. 
Apply Lemma~\ref{lem:good-scales} to $Z$ and find a good scale $\varrho\in [\delta^{1/2}, \delta^{\mu}]$ such that $Z\cap {\bf z}$ is a rescaled $(\varrho, s, \varrho^{-\eta_0})$-set for every $z\in Z$ and ${\bf z} = z^\varrho$.     
    Pick $p_0 \in E$, let ${\bf q} = p_0^\varrho$ and apply Lemma \ref{lem:pigeonhole-rho-ball} to find $E({\bf q}) \subset E\cap {\bf q}$ and the direction set $\Theta_{\bf q} = \Theta_{p_0'}$.

    Using property (b), by dyadic pigeonholing we can find some $M$ and a subset $\TT' \subset \TT$ such that $|Z({\bf q}, T)| \sim M$ for every $T\in \TT'$ and we have
    \[
    \sum_{T\in \TT'} |Z({\bf q}, T)| \sim M |\TT'| \gtrsim \d^{O(\eta)} |Z({\bf q})| \varrho^{-\kappa} (\varrho/\d)^{2t}.
    \]
    We know by (b) also that $|\TT'| \le \d^{-O(\eta)} \varrho^{-\kappa} (\varrho/\d)^{2t}$ and so it follows that $M \ge \d^{O(\eta)} |Z({\bf q})|$. Since $Z({\bf q}, T) \subset Z$ and $Z$ is $\d^{\eta}$-uniform, this implies that $|Z({\bf q}, T)|_{\varrho^2} \gtrsim \d^{O(\eta)} \frac{M}{|Z \cap N_{\varrho^2} z|} \gtrsim \d^{O(\eta)} |Z({\bf q})|_{\varrho^2}$.

    For $\theta \in \Theta_{\bf q}$ consider a maximal set $\mathbb U_\theta$ of essentially distinct $\varrho^2 \times \varrho^2\times \varrho$ tubes $U$ with $\theta(U) = \theta$ and such that $|U \cap E({\bf q})| \ge \d^{C_0\eta} M (\varrho^2/\d)^{2t-\kappa} $. Note that for small fixed constant $c$, each $T\in \TT$ contributes at least one such segment $U$ to every $\mathbb U_\theta$ with $|\theta - \theta(T)| \le c \varrho$. On the other hand, each tube segment $U$ is being contributed to by at most $\d^{-O(\eta)}\varrho^{-\kappa} (\varrho^2/\d)^{2t}$-many tubes $T\in \TT'$ (there are $\approx \varrho^{-\kappa}$ many $\varrho^2$-tubes $T_{\varrho^2} \in \TT_{\varrho^2}$ going through $U$ and each such tube contains $\approx (\varrho^2/\d)^{2t}$ tubes $T\in \TT'$). 
    Let $\Theta \subset\Theta_{\bf q}$ be a maximal $\frac{c}{2}\varrho$-separated set. From the above we then get
    \[
    \sum_{\theta\in \Theta} |\mathbb U_\theta| \gtrsim \d^{O(\eta)} \frac{|\TT'| }{\varrho^{-\kappa} (\varrho^2/\d)^{2t}} \ge \d^{O(\eta)} \varrho^{-2t} = \d^{O(\eta)} |\Theta| \varrho^{\kappa-2t}.
    \]
    For some $z' \in Z({\bf q})$ and every $\theta\in \Theta$ define $A_\theta\subset  [-1,1]^2 \cap c\varrho^2\ZZ^2$ be the set of points $a$ such that $a^{(\varrho^2)} \cap A_{z'} \neq\emptyset$ and $(a, z') \in U$ for some $U\in \mathbb U_\theta$. By averaging over $Z({\bf q})$, we can choose $z'$ so that
    \[
    \sum_{\theta\in \Theta} |A_\theta| \ge \d^{O(\eta)} |\Theta| \varrho^{\kappa-2t}.
    \]
    Let $\Theta' \subset \Theta$ be the subset of directions $\theta$ so that $|A_\theta| \ge \d^{C_1 \eta} \varrho^{\kappa-2t}$. By  choosing $C_1$ large enough (larger than the implied constants above) we get that $|\Theta'| \ge \d^{C_1\eta} \varrho^{-\kappa}$. In particular, $A_\theta$ is $(\varrho^2, 2t-\kappa, \d^{-O(\eta)})$-Katz--Tao set for every $\theta\in \Theta'$. 
    We may assume $2t-\kappa \le 2-2\epsilon_0$, as otherwise we get the desired bound on $\kappa$. 
    So for $S\ge 1$ we have 
    \begin{equation}\label{eq:Atheta-KT}
    |A_\theta \cap B_{S\varrho^2}| \le     \d^{-O(\eta)} S^{2t-\kappa} \le  \d^{-O(\eta)} S^{2-2\eps_0}.
    \end{equation}
    Define for $a\in A_\theta$ the set $E(a, \theta) \subset ((a,z') +\RR(\theta,1))$ to be a $c\varrho^2$-separated collection of points whose $\varrho^2$-neighborhood covers   $E({\bf q}) \cap U$ where $U\in \mathbb U_\theta$ goes through $(a,z')$. By construction,
    \begin{equation}\label{eq:E-small}
    \left|\bigcup_{\theta\in \Theta', a\in A_\theta} E(a, \theta)\right|_{\varrho^2} \lesssim |E({\bf q})|_{\varrho^2} \le \d^{-O(\eta)} |A_\theta|_{\varrho^2} |E(a, \theta)|_{\varrho^2}.    
    \end{equation}
    
    Since $\Theta_{\bf q} = \Theta_{p_0'}$ for some $p_0' \in E\cap {\bf q}$ by the definition of a $(\d,\kappa,\gamma,\d^{-\eta})$-aligned set (Definition \ref{def:aligned2}) we have $|\Theta_{\bf q} \cap R|_{\varrho} \le \d^{-O(\eta)} \varrho^{-\gamma}$ for every $\varrho \times 1$ rectangle $R$. 
    Thus, for any $\varrho\times 1$ rectangle $R$ we have
    \begin{equation}\label{eq:Theta-rectangle}
    |\Theta' \cap R|_{\varrho} \le \d^{-O(\eta)} \varrho^{-\gamma} \le \d^{-O(\eta)}\varrho^{\kappa-\gamma} |\Theta'|_{\varrho}.    
    \end{equation}
    By rescaling ${\bf q}$ to a unit cube and applying Proposition \ref{prop:small-gamma-range} using \eqref{eq:Atheta-KT}, \eqref{eq:E-small} and \eqref{eq:Theta-rectangle} we obtain
    \[
    \lambda \le \varrho^{-\eta_0} \nu,
    \]
    where
    \[
    \lambda = |E(a, \theta)|_{\varrho^2}\varrho, \quad   \nu = \d^{-O(\eta)} \varrho^{\kappa-\gamma}
    \]
    and we take $\eta$ small enough so that the $\varepsilon$ in Proposition \ref{prop:small-gamma-range} is smaller than $\eta_0$. 
    On the other hand we have 
    \[
    |E(a, \theta)|_{\varrho^2} \ge \d^{O(\eta)} |Z({\bf q})|_{\varrho^2} \ge \d^{O(\eta)} \varrho^{\eta_0} \varrho^{-s}
    \]
    using that $Z({\bf q})$ is a rescaled $(\varrho, s, \varrho^{-\eta_0})$-set. So since $\varrho \le \d^{\mu}$ for some $\mu = \mu(\eta_0)>0$ we get 
    \[
    s-1 \le \gamma-\kappa + O(\eta/\mu+\eta_0).
    \]
    In the limit $\eta, \eta_0\to 0$, this gives $\gamma \ge s+\kappa-1$. This contradicts the assumption $\gamma \le s+\kappa-1 -\eps_0$ and so we are done. 
\end{proof}

\section{Intermediate tube density: Proof of Proposition \ref{prop: trilinearABC}}\label{sec:intermediate}

In this section, we prove  Proposition~\ref{prop: trilinearABC}. 
The main tool is Orponen-Shmerkin discretized ABC-sum production theorem.
\begin{theorem}[Orponen-Shmerkin Theorem 1.7 \cite{orponen2023projections}]\label{thm: ABC}
    Let $0< \beta\leq \alpha <1.$ Then, for every $\gamma \in (\alpha-\beta, 1]$, 
    there exist $\chi, \delta_0\in (0, 1/2]$ such that the following holds. Let $\d\in 2^{-\mathbb{N}}$ with $\d\in (0, \delta_0]$, and let $A, B\subset \delta \cdot \mathbb{Z}\cap [0,1]$ be sets satisfying the following hypotheses: 
    \begin{enumerate}
        \item $|A|\leq \delta^{-\alpha}$, 
        \item $B\neq \emptyset, $ and $|B\cap B(x,r)|\leq \delta^{-\chi} r^{\beta}|B|$ for all $x\in \mathbb{R}$ and $r\in [\delta, 1]$. 
    \end{enumerate}
    Further, let $C\subset \delta\cdot \mathbb{Z}\cap [-1,1]$ be a non-empty set satisfying $|C\cap B(x,r)|\leq \delta^{-\chi}r^{\gamma}|C|$ for $x\in \mathbb{R}$ and $r\in [\delta, 1]$. 
    Then there exists $c\in C$ such that 
    \begin{equation}\label{eq:ABC-conclusion}
    |\{a+cb: (a,b)\in G\}|_{\delta}\geq \delta^{-\chi}|A|, \quad G\subset A\times B, |G|\geq \delta^{\chi}|A||B|.     
    \end{equation}
\end{theorem}

\begin{remark}
Note that our statement has $C \subset [-1,1]$ instead of $C\subset [1, 2]$ in \cite{orponen2023projections}. This version follows from exactly the same proof and a version with slightly different error terms in place of various $\d^{-\chi}$ follows directly from their result. 
\end{remark}

\begin{proof}[Proof of Proposition~\ref{prop: trilinearABC}]

Let $\eta_0>0$ be a small parameter to be determined and $\eta>0$ sufficiently small. 
Apply Lemma~\ref{lem:good-scales} to $Z$ and find a scale $\varrho\in [\delta^{1/2}, \delta^{\mu}]$ such that $Z\cap {\bf z}$ is a rescaled $(\varrho, s, \varrho^{-\eta_0})$-set for every $z\in Z$ and ${\bf z} = z^\varrho$.  By passing through a subset of $Z\cap {\bf z}$, we may assume that $|Z\cap {\bf z}|_{\varrho^2} \sim \varrho^{-s}$. 

Without loss of generality we may assume that $E$ is $\d$-separated. 
Let $p_0\in E$, consider ${\bf q} =p_0^\varrho$.  Apply Lemma \ref{lem:pigeonhole-rho-ball} and let $E({\bf q}) \subset E\cap {\bf q}$, $Z({\bf q}) \subset Z\cap {\bf z}$ and $\Theta_{\bf q} = \{ (\varphi, \xi+\alpha_{\bf q} \varphi), ~ \xi\in \Xi_{\bf q}, ~ \varphi \in \Phi_{{\bf q}, \xi}\}$ be the resulting sets. By a change of coordinates we may assume that $\alpha_{\bf q} = 0$. 

For $z \in Z({\bf q})$ let $A_z({\bf q}) \subset A_z$ be the set of $a \in A_z$ so that $(a,z)\in E({\bf q})$. Note that $\sum_{z \in Z({\bf q})} |A_z({\bf q})|_\d \sim |E({\bf q})|_\d \gtrsim \d^{O(\eta)}|Z({\bf q})|_\d (\varrho/\d)^{2t-\kappa}$. For $z\in Z({\bf q})$ let $Y_z$, $X_{y,z}$ be $c\varrho^2$-separated sets such that 
\[
A_z({\bf q}) \subset N_{\varrho^2} \{ (x, y),~ x\in X_{y, z}, y\in Y_z  \}
\]
and $\sum_{y\in Y_z} |X_{y,z}| \sim |A_z({\bf q})|_{\varrho^2} \le \d^{-O(\eta)} \varrho^{\kappa-2t}$. After pigeonholing we may further assume that all $X_{y, z}$ have the same size, up to a constant factor. 
After passing to subsets in $Z({\bf q})$ and $E({\bf q})$ of density $\ge \d^{O(\eta)}$ we may also assume that each set $A_z({\bf q})$ has the property that $|A_z({\bf q}) \cap B(a, \varrho^2/100)| \gtrsim \d^{O(\eta)} (\varrho^2/\d)^{2t-\kappa}$ for every $a \in A_z({\bf q})$. Note that after this refinement we still have the lower bound from (b) in Lemma \ref{lem:pigeonhole-rho-ball}. 
Let $\Xi \subset \Xi_{\bf q}$ and $\Phi_\xi \subset \Phi_{{\bf q}, \xi}$ be maximal $\varrho/100$-separated subsets and define $\Theta = \{ (\varphi, \xi), \xi\in \Xi, ~\varphi\in \Phi_\xi \} = \bigsqcup_{\xi\in \Xi} \Phi_\xi \times \{\xi\}$. Note that $\Xi$ is $(\varrho, \kappa-\gamma, \d^{-O(\eta)})$-AD-regular and $\Phi_{\xi}$ is $(\varrho, \gamma, \d^{-O(\eta)})$-AD-regular for every $\xi \in \Xi$. Note that by the assumption of the Proposition, we have $\gamma,\kappa-\gamma\ge \eps_1$ for some fixed $\eps_1 >0$ and other parameters can be taken much smaller depending on $\eps_1$, so both of these sets are positive dimensional. This will allow us to apply Theorem \ref{thm: ABC}.

Consider the following set:
\[
H = \{ (p, p', \theta) \in E({\bf q}) \times E({\bf q}) \times \Theta: |\pi_{\theta}(p) - \pi_{\theta}(p')| \le \varrho^2 \}.
\]
Let us lower bound the size of $H$. We start with (b) from Lemma \ref{lem:pigeonhole-rho-ball} to get:
\[
\sum_{T\in \TT({\bf q})} |Z({\bf q}, T)|^2 \ge \d^{O(\eta)} |Z({\bf q})|^2 \varrho^{-\kappa} (\varrho/\d)^{2t}.
\]
Note that each tube $T\in \TT({\bf q})$ contributes $\gtrsim |Z({\bf q}, T)|^2 \left(\d^{O(\eta)} (\varrho^2/\d)^{2t-\kappa}\right)^2$ triples to $H$: first, pick any $\tilde p, \tilde p' \in Z({\bf q}, T)$ and then choose $p \in E_z({\bf q}) \cap N_{\varrho^2/100}\tilde p$ and $p' \in E_z({\bf q}) \cap N_{\varrho^2/100}\tilde p'$ and let $\theta \in \Theta$ be an angle such that $\theta(T) \in N_{\varrho/10} \theta$. With these choices we get $|\pi_{\theta}(p) - \pi_{\theta}(p')| \le \varrho^2$. By double counting this implies a lower bound
\[
|H| \ge \d^{O(\eta)} |E({\bf q})| \cdot |\Theta| \cdot |Z({\bf q})|.
\]
Now consider 
\[
H_{z, z'} = \{ (x, y, x', y', \varphi, \xi): ~ \exists (p, p', \theta)\in H:~ p \in N_{\varrho^2} (x,y), ~ p' \in N_{\varrho^2} (x',y'), ~ \theta = (\varphi, \xi)  \}
\]
where we consider pairs $(x,y)$ with $y \in Y_z$ and $x \in X_{y, z}$ and similarly for $(x',y')$. By double counting we conclude that 
\begin{equation}\label{eq:graph-large-average}
\frac{1}{|Z({\bf q})|^2}\sum_{z,z'\in Z({\bf q})} |H_{z,z'}| \ge \d^{O(\eta)} \varrho^{\kappa-2t} |\Theta|.
\end{equation}
Crucially, each tuple $(x, y, x', y', \varphi, \xi) \in H_{z,z'}$ satisfies the property $|\pi_{\theta}(x, y, z)-\pi_\theta(x',y',z')| \le C\varrho^2$, i.e.
\begin{equation}\label{eq:linear-relations}
|x' - x - (z'-z) \varphi| \le C\varrho^2,\quad |y'-y - (z'-z)\xi| \le C\varrho^2.    
\end{equation}
Define $G_{z, z'} \subset Y_z \times \Xi$ to be the set of pairs $(y, \xi)$ such that $(x, y, x', y', \varphi, \xi) \in H_{z, z'}$ for some choice of $x, x', y', \varphi$. Note that $y' \in Y_{z'}$ is essentially determined by $y, \xi$ (since $Y_{z'}$ is $c\varrho^2$-separated and (\ref{eq:linear-relations})) and $x'$ is essentially determined by $x$ and $\varphi$. This shows that there are at most $C|X_{y,z}| |\Phi_{\xi}| $ ways to extend $y, \xi$ to a tuple $(x, y, x', y', \varphi, \xi)$. So using $|X_{y, z}| |Y_z| \sim_{\d^{O(\eta)}} \varrho^{\kappa-2t}$ and $|\Phi_{\xi}| |\Xi| \sim_{\d^{O(\eta)}} |\Theta|$ it follows from (\ref{eq:graph-large-average}) that
\[
\frac{1}{|Z({\bf q})|^2}\sum_{z,z'\in Z({\bf q})} |G_{z,z'}| \ge \d^{O(\eta)} |Y_{z}| |\Xi|.
\]
For $z \in Z({\bf q})$ consider $C_z$ to be the set of $z' \in Z({\bf q})$ such that $|z'-z| \ge \varrho^{2\eta_0/s}\d^{C_1 \eta} \varrho$ and $|G_{z,z'}| > \d^{C_1 \eta} |Y_z| |\Xi|$. Since $Z({\bf q})$ is rescaled $(\varrho, s, \varrho^{-\eta_0} \d^{-O(\eta)})$-set, for a large enough constant $C_1$, by (\ref{eq:graph-large-average}) there is some $z$ such that we have $|C_z| \ge \d^{C_1 \eta} |Z({\bf q})|$. 

Suppose that $|Y_z|_{\varrho^2} \sim \varrho^{-\alpha}$ for every $z \in Z({\bf q})$ and some $\alpha \in \RR$.  
We apply Theorem \ref{thm: ABC} to the rescaled versions of sets $A = Y_z$, $B = \Xi$ and $C = C_z-z$: we have $|Y_z| \lesssim \varrho^{-\alpha}$, $B$ is rescaled $(\varrho, \kappa-\gamma, \d^{-O(\eta)})$-AD-regular set and $C$ is $(\varrho, s, \varrho^{-\eta_0}\d^{-O(\eta)})$-set. 
So assuming that $\alpha \le \min(s+\kappa-\gamma, 1)-\epsilon_0$ there is $\chi, \d_0$ depending on $\epsilon_0$, $s$, $\kappa-\gamma$ such that if $\eta_0/s, \eta/\mu \ll \chi$ then there exists some $c=z'-z \in C = C_z$ such that
\[
| \{ y + (z'-z) \xi : (y,\xi) \in G_{z, z'} \} |_{\varrho^2} \ge \varrho^{-\chi} |Y_z| \gtrsim \varrho^{-\chi - \alpha}.
\]
On the other hand, we know from (\ref{eq:linear-relations}) that $\{ y + (z'-z) \xi : (y,\xi) \in G_{z, z'} \}  \subset N_{C\varrho^2} Y_{z'}$ and so by our assumption it has covering number $\lesssim \varrho^{-\alpha}$. This is a contradiction and we get 
\[ 
|Y_z|_{\varrho^2} \ge \varrho^{\epsilon_0} \varrho^{-\min(s+\kappa-\gamma, 1)}
\]
for every $z \in Z({\bf q})$.

By a completely analogous argument, applying Theorem \ref{thm: ABC} with rescaled versions of sets $A = X_{y,z}$, $B =\Phi_\xi$ and $C \approx Z({\bf q})-z$ we get
\[
|X_{y, z}|_{\varrho^2} \ge \varrho^{\epsilon_0} \varrho^{-\min(s+\gamma, 1)}
\]
More precisely, we define a new graph $G_{z,z', y, \xi} \subset X_{y, z} \times \Phi_{\xi}$ consisting of pairs $(x, \varphi)$ which extend to a tuple $(x, y, x', y', \varphi, \xi) \in H_{z,z'}$. By averaging we again get that $|G_{z, z', y, \xi}| \gtrsim_{\d^{O(\eta)}} |X_{y,z}| |\Phi_\xi|$ for many choices of $z, z', y, \xi$. We let $C'_{z, y, \xi} \subset Z({\bf q})$ be the set of $z'$ so that $|G_{z, z', y, \xi}| \gtrsim_{\d^{O(\eta)}} |X_{y,z}| |\Phi_\xi|$ and $|z'-z| \ge \varrho^{2\eta_0/s}\d^{C_2 \eta} \varrho$ and choose $z$ so that $|C'_{z, y, \xi}| \ge \d^{C_2 \eta} |Z({\bf q})|$. This allows us to apply Theorem \ref{thm: ABC} with $C = C'_{z, y, \xi} - z$ and get a contradiction unless $|X_{y, z}|_{\varrho^2} \ge \varrho^{\epsilon_0} \varrho^{-\min(s+\gamma, 1)}$ holds. 

So by combining the bounds on $|Y_z|$ and $|X_{y,z}|$ (recall that these sets are $\sim \varrho^2$-separated) we see that for $z \in Z({\bf q})$ we have
\[
|A_z({\bf q})|_{\varrho^2} \gtrsim \sum_{y\in Y_z} |X_{y, z}| \gtrsim \varrho^{2\epsilon_0}  \varrho^{-\min(s+\kappa-\gamma, 1)} \varrho^{-\min(s+\gamma, 1)}
\]
Using the assumption that $\gamma \in [s+\kappa-1-\eps_0, 1-s+\eps_0]$ this gives $|A_z({\bf q})|_{\varrho^2} \gtrsim \varrho^{4\eps_0}  \varrho^{-2s-\kappa} $. Comparing this to the upper bound $|A_z({\bf q})|_{\varrho^2} \le \d^{-O(\eta)} \varrho^{\kappa-2t}$ we get the desired bound on $\kappa$.
\end{proof}

\section{Global grains reduction} \label{sec:global-grains}

In this section and the following sections we prove Proposition \ref{prop: coplanar}. Recall that in this proposition we assume that $\gamma \ge \max(s+\kappa-1-\eps_0, 1-s+\eps_0)$ for a small fixed $\eps_0>0$.
In this section we apply rescaling and thickening operations to reduce to the case when the set $E$ has `global grains' structure (see Section \ref{sec:proof-overview} to compare with the 3D sticky Kakeya case). 

\begin{definition}
    Let $Z \subset [-1,1]$ and let $f: Z\to [-1,1]^m$ be a function. We say that $(Z, f)$ is $(\d, K)$-uniform if for any $\varrho, \sigma \in [\d, 1]$ there is some number $M_{\varrho, \sigma}$ such that for every $z\in Z$:
    \[
    \#\{z' \in Z:~ |z-z'| \le \varrho, ~ |f(z)-f(z')|\le \sigma \} \in [K^{-1} M_{\varrho, \sigma}, M_{\varrho, \sigma}].
    \]
\end{definition}

By standard pigeonholing, for any $\d$ there is some $K=\d^{o(1)}$ such that for any $Z$ and $f$ we can find a $\d$-separated subset $Z' \subset Z$ such that $|Z'| \ge K^{-1} |Z|_\d$ and $(Z', f|_{Z'})$ is $(\d, K)$-uniform. 

The next definition is a modification of Definition \ref{def: configuration}, where the set $E$ has additional structure called `global grains'.

\begin{definition}\label{def:global-grains}
    Let $t \in [0,2]$, $s \in [0,1]$, $\kappa, \gamma \ge 0$ so that $2t\ge 1+\kappa$. 
    
    A global grains $(\d,\eta, s, t, \kappa, \gamma)$-Furstenberg configuration is the following collection of data:
    \begin{enumerate}
        \item A $(2t, \d^{-\eta})$-AD-regular set of tubes $\TT$.
        \item A $\d$-separated $(\d, s)$-set $Z \subset [-1,1]$ and a function $f: Z\to [-1,1]$ such that $(Z, f)$ is $(\d, \d^{-\eta})$-uniform.
        \item A $\d$-separated set $E$ of the form
        \[
        E = \bigsqcup_{z\in Z} E_z=\{ (x, y + f(z) x, z),~ z \in Z,~ y\in Y_z,~ x \in X_{y, z} \} \subset [-1,1]^3,
        \]
        where $Y_z \subset [-1,1]$ is a $(\d, 2t-1-\kappa, \d^{-\eta})$-AD-regular set and $X_{y,z} \subset [-1,1]$ is such that $|X_{y, z}|_\d \ge \d^{\eta -1}$.
        \item For each $p \in E_z$, some
        $\d$-separated $(\d, \kappa-\gamma, \d^{-\eta})$-AD-regular set $\Xi_{p} \subset [-1,1]$ and $\d$-separated $(\d, \gamma, \d^{-\eta})$-AD-regular sets $\Phi_{p, \xi} \subset [-1,1]$ such that for any $\xi \in \Xi_p$ and $\varphi \in \Phi_{p, \xi}$ there is a tube $T\in \TT$ through $p$ with $\theta(T) \in N_\d  (\varphi, \xi + f(z) \varphi)$.
    \end{enumerate}
\end{definition}

The following lemma is a variant of Lemma \ref{lem:pigeonhole-rho-ball} adapted to the global grains configurations. 

\begin{lemma}\label{lem:pigeonhole-rho-ball-2}
    Let $(\TT, E, \{\TT_p\})$ be a global grains $(\d, \eta, s, t, \kappa, \gamma)$-Furstenberg configuration. 
    Let $\varrho \in [\d,1]$. Then there exists $E' \subset E$ such that $|E'| \ge \d^{O(\eta)}|E|$ and the following holds.
    Let $z_0 \in Z$, $p_0 \in E_{z_0}$ and write ${\bf z} = z_0^\varrho$, ${\bf q} = p_0^\varrho$. Then for some $Z({\bf q}) \subset Z \cap {\bf z}$, $E({\bf q}) = E' \cap {\bf q}$, and some $p_0' \in E\cap {\bf q}$ we have the following. 
    \begin{itemize}
        \item[(a)]  Write $\Theta_{\bf q} = \Theta_{p_0'}$ then for each $p \in E({\bf q})$ we have $|\Theta_p \cap  \Theta_{\bf q}^{\varrho}| \ge \d^{O(\eta)} |\Theta_p|$.

        \item[(b)] Let $\TT({\bf q})$ be the set of $T\in\TT$ such that $\theta(T) \in \Theta_{\bf q}$ and $T\in \TT_p$ for some $p \in E({\bf q})$. Then $|\TT({\bf q})| \sim_{\d^{O(\eta)}} \varrho^{-\kappa} (\varrho/\d)^{2t}$. For $T\in \TT({\bf q})$ let $Z({\bf q}, T)$ be the set of $z \in Z({\bf q})$ such that $T \in \TT_p$ for some $p \in E_z({\bf q})$. Then for any $\tilde Z\subset Z({\bf q})$ we have
        \[
        \sum_{T\in \TT({\bf q})} |\tilde Z\cap Z({\bf q}, T)| \sim_{\d^{O(\eta)}} |\tilde Z| |\TT({\bf q})|.
        \]

        \item[(c)] We can write 
        \[
        E({\bf q}) = \{  (x, y + f(z) x, z), ~ z \in Z({\bf q}), ~y \in Y_{z}({\bf q}), ~ x \in X_{y, z}({\bf q}) \},
        \]
        where $Y_z({\bf q})$ is a rescaled $(\d/\varrho, 2t-\kappa-1, \d^{-O(\eta)})$-AD-regular set and $X_z({\bf q})$ is a rescaled $(\d/\varrho, 1, \d^{-O(\eta)})$-AD-regular set.
    \end{itemize}
\end{lemma}

\begin{proof}
    Properties (a) and (b) are the content of Lemma \ref{lem:pigeonhole-rho-ball}. To obtain (c) we use property (3) in the defition of a global grains configuration and pigeonhole so that sets appearing in the definition of $E({\bf q})$ are AD-regular. 
\end{proof}

\subsection{Constructing global grains.} Suppose we have a $(\d, \eta, s, t, \kappa, \gamma)$-Furstenberg configuration and that $s+\gamma \ge 1$. We use this structure to construct a global grains configuration $(\tilde \d, \tilde \eta, s, t, \kappa, \gamma)$. We accomplish this in the following two steps. First, we show that after rescaling each set $A_z$ has the global grains structure. Second, after another rescaling we show that the direction sets $\Theta_p$ are aligned with global grains.

\begin{prop}[Global grains, step 1]\label{prop:global-grains-1}
    Suppose that $\gamma \ge 1-s$ and suppose that for any $\eta$ and $\d_0$ there is $\d<\d_0$ so that a $(\d, \eta, s, t, \kappa, \gamma)$-Furstenberg configuration exists.  
    Then for any $\d_0, \eta$ there is $\d<\d_0$ so that there exists a $(\d, \eta, s, t, \kappa, \gamma)$-Furstenberg configuration with the following extra property.
    
    For each $z\in Z$ there exists $f(z) \in [-1,1]$ such that we can write
    \begin{equation}\label{eq:layers-global-grains}
    A_z = \{ (x, y+f(z)x), ~ y \in Y_{z}, x\in X_{y, z} \}    
    \end{equation}
    for some $(\d, 2t-\kappa-1, \d^{-\eta})$-AD-regular set $Y_z \subset [-1,1]$ and $(\d, 1, \d^{-\eta})$-AD-regular set $X_{y, z} \subset [-1,1]$. 
\end{prop}

\begin{proof}
    Let $\eta_0>0$ be a small parameter to be determined and $\eta>0$ sufficiently small. Consider an $(\d, \eta, s, t, \kappa, \gamma)$-Furstenberg configuration.
    Apply Lemma~\ref{lem:good-scales} to $Z$ and find a good scale $\varrho\in [\delta^{1/2}, \delta^{\mu}]$ such that $Z\cap N_\varrho z$ is a rescaled $(\varrho, s, \varrho^{-\eta_0})$-set for every $z\in Z$. Without loss of generality we may take $\varrho$ to be admissible for $\TT$ so that $\TT[T_\varrho]$ is $(2t, \d^{-\eta})$-AD-regular set of tubes for every $T_\varrho \in \TT_\varrho$.

    For any $p_0\in E$ we can consider the dyadic cube ${\bf q} = p_0^\varrho$ and construct $E({\bf q}) \subset E\cap {\bf q}$ and $\alpha_{\bf q}$, $\Theta_{\bf q}$ using Lemma \ref{lem:pigeonhole-rho-ball}. Let $E({\bf q}) =\bigsqcup_{z \in Z({\bf q})} A_z({\bf q}) \times \{z\}$. 
    By repeating the proof in the previous section, after some pigeonholing, we may ensure that sets $A_z({\bf q})$ for $z \in Z({\bf q})$ can be covered by $\varrho^2$-neigborhood of the form
    \[
    A_z({\bf q}) \subset N_\varrho\{ (x, y+ \alpha_{\bf q} x), ~ y\in Y_{z}({\bf q}), x \in X_{y, z}({\bf q}) \}
    \]
    where sets $Y_z({\bf q})$ and $X_{y,z}({\bf q})$ are $c\varrho^2$-separated and satisfy $\sum_{y\in Y_z({\bf q})} |X_{y,z}({\bf q})| \lesssim |A_z({\bf q})|_{\varrho^2} \lesssim \d^{-O(\eta)} \varrho^{\kappa-2t}$. Furthermore we may assume that the set $\{ (x, y+ \alpha_{\bf q} x), ~ y\in Y_{z}({\bf q}), x \in X_{y, z}({\bf q}) \}$ is rescaled $(\varrho, 2t-\kappa, \d^{-O(\eta)})$-AD-regular set. 

    Now by applying Theorem \ref{thm: ABC} in the same way as in the proof of Proposition \ref{prop: trilinearABC} we get
    \[
     |X_{y, z}({\bf q})|_{\varrho^2}\gtrsim \varrho^{\epsilon} \varrho^{-\min(s+\gamma,1)} \ge \varrho^{\epsilon - 1}
    \]
    for every $y \in Y_z({\bf q})$ (after we pigeonhole all sets $X_{y, z}({\bf q})$ to have approximately the same size). 
    
    From the above bound on $X_{y, z}({\bf q})$ it then follows that $Y_{z}({\bf q})$ must be a rescaled $(\varrho, 2t-\kappa-1, \varrho^{-O(\epsilon)}\d^{-O(\eta)})$-AD-regular set. Indeed, let $y \in Y_z({\bf q})$ and let $\tau \in [\varrho^2, \varrho]$. Assume for simplicity that $\alpha_{\bf q}=0$ (we can reduce to this case by a coordinate change). 
    Then by covering $A_z({\bf q}) \cap ([-1,1] \times N_\tau y)$ with $\sim \frac{\varrho}{\tau}$-many $\tau$-balls we have
    \[
    |A_z({\bf q}) \cap ([-1,1] \times N_\tau y)|_{\varrho^2} \lesssim \d^{-O(\eta)} \frac{\varrho}{\tau}  (\tau/\varrho^2)^{2t-\kappa}
    \]
    and on the other hand
    \[
    |A_z({\bf q}) \cap ([-1,1] \times N_\tau y)|_{\varrho^2} \sim \sum_{y' \in Y_z({\bf q}) \cap N_\tau y} |X_{y',z}({\bf q})|_{\varrho^2} \gtrsim |Y_z({\bf q}) \cap N_\tau y|_{\varrho^2} \varrho^{\epsilon-1}
    \]
    giving $|Y_z({\bf q}) \cap N_\tau y|_{\varrho^2} \lesssim \d^{-O(\eta)} \varrho^{-\epsilon} (\tau/\varrho^2)^{2t-\kappa-1} $. On the other hand, if we let $x_i \in X_{y, z}({\bf q})$ for $i = 1, \ldots$ be a maximal $2\tau$-separated collection of points (there are at least $\varrho^{\epsilon} (\varrho/\tau)$-many points by the lower bound on $|X_{y,z}({\bf q})|_{\varrho^2}$), then by the AD-regularity of $A_z({\bf q})$ we get
    \[
    |A_z({\bf q}) \cap ([-1,1] \times N_\tau y)|_{\varrho^2} \gtrsim \sum_{i} |A_z({\bf q}) \cap B( (x_i, y), \tau)|_{\varrho^2} \gtrsim \d^{O(\eta)}\varrho^{\epsilon} \frac{\varrho}{\tau} (\tau/\varrho^2)^{2t-\kappa}.
    \]
    So the estimate $|A_z({\bf q}) \cap ([-1,1] \times N_\tau y)|_{\varrho^2} \lesssim  |Y_z({\bf q}) \cap N_\tau y|_{\varrho^2} \varrho^{-1}$ gives the lower bound $|Y_z({\bf q}) \cap N_\tau y|_{\varrho^2} \gtrsim \d^{O(\eta)} \varrho^{\epsilon} (\tau/\varrho^2)^{2t-\kappa-1} $.

    The above argument applies to any $\varrho$-ball ${\bf q}$. Consider a $\varrho$-tube $T_{\varrho}$ for some $T\in \TT$. By averaging, we can choose $T$ so that $|Z(T)| = \#\{p: T\in \TT_p\} \ge \d^{O(\eta)} |Z|$. Let $z_1, \ldots, z_m \in Z(T)$ be a maximal $2\varrho$-separated set of points and let $p_i \in E_{z_i}$ be points so that $T\in \TT_{p_i}$. Let ${\bf q}_i = p_i^\varrho$.  
    Let $\psi=\psi_{T_\varrho}$ be the anisotropic rescaling map sending $T_\varrho$ to the unit ball and define a rescaled configuration as follows. First let
    \[
    \mc E_i = \{ (x, y+\alpha_{{\bf q}_i} x, z), ~ z \in Z({\bf q}_i), y \in Y_{z}({\bf q}_i), x \in X_{y, z}({\bf q}_i) \} \subset {\bf q}_i \subset T_\varrho
    \]
    and put $\mc E = \bigcup_{i=1}^m \mc E_i$. We set $\tilde E = \psi(\mc E)$ and $\tilde Z = \bigcup Z({\bf q}_i)$. Let $\tilde \TT = \psi(\TT_{\varrho^2}[T_\varrho])$. It follows that $\tilde \TT$ is $(2t, \d^{-2\eta})$-AD-regular set of $\varrho$-tubes and $C_{t;t'-CW}(\tilde\TT) \le \d^{-O(\eta)} C_{t;t'-CW}(\TT)$. 

    For $\tilde p = \psi(p)\in \tilde E$ we let $\tilde \TT_{\tilde p}$ be the set of $\varrho$-tubes $\psi(T'_{\varrho^2})$ where $T'_{\varrho^2}$ is the dyadic thickening of some $T' \in \TT_{p}$. By averaging over the choice of the initial tube $T_\varrho$ we can ensure that we have
    \[
    \frac{1}{|\tilde E|}\sum_{\tilde p \in \tilde E} |\tilde \TT_{\tilde p}| \ge \d^{O(\eta)} \varrho^{-\kappa}. 
    \]
    By construction, we now have $\tilde E = \bigsqcup_{z\in \tilde Z} \tilde A_z \times \{z\}$ and sets $\tilde A_z$ have the global grains structure with $f(z) = \alpha_{{\bf q}_i}$ for $z \in Z({\bf q}_i)$. 
    Now by applying Lemma \ref{lem:uniform-subset} to the point-line configuration $X=\{(\tilde p, \theta(\tilde T))\}$  
    we can find a dense subset $X'\subset X$ and a corresponding subsets $\tilde E'\subset \tilde E$ and $\tilde \TT_{\tilde p}' \subset \tilde \TT_{\tilde p}$, $\tilde p \in \tilde E'$ which are $(\varrho, \kappa, \d^{-O(\eta)})$-AD-regular. After further pigeonholing we ensure that $\tilde A_z$ is $(\varrho, 2t-\kappa, \varrho^{-O(\epsilon)})$-AD-regular and has the form (\ref{eq:layers-global-grains}), $\tilde Z$ is $\varrho^{-\epsilon}$-uniform and $\Theta_{\tilde p}$ is $(\varrho,  \kappa, \gamma, \varrho^{-\epsilon})$-aligned. So we obtain a configuration satisfying Definition \ref{def: configuration} and the global grain property stated in the lemma.
\end{proof}

\begin{prop}[Global grains, step 2]\label{prop:global-grains-2}
    Suppose that $\gamma \ge 1-s$ and that $2t-\kappa \le 2-\epsilon_0$ for some fixed $\eps_0>0$.
    Suppose that for any $\eta$ and $\d_0$ there is $\d<\d_0$ so that a $(\d, \eta, s, t, \kappa, \gamma)$-Furstenberg configuration exists.  
    Then for any $\d_0, \eta$ there is $\d<\d_0$ so that there exists a global grains $(\d, \eta, s, t, \kappa, \gamma)$-Furstenberg configuration.
    
\end{prop}

\begin{proof}
    We start with a $(\d, \eta, s, t, \kappa, \gamma)$-Furstenberg configuration satisfying the global grains property from Proposition \ref{prop:global-grains-1}. Out goal is to produce a global grains Fustenberg configuration with error $\eta_0$. Suppose that $\eta\ll \eta_0$ is sufficiently small. 
    By dyadic pigeonholing we can find some $\tau \in [\d,1]$ and a subset $E' \subset E$ so that $|\alpha_p - f(z)| \sim \tau$  (or $\le \d$ if $\tau=\d$) for every $p\in E' \cap E_z$ and we have $|E'| \gtrsim \frac{1}{\log (1/\d)} |E|$. 

    Fix some small $\varepsilon >0$. First suppose that we have $\tau \le \d^{\varepsilon}$. In this case we thicken the set of tubes $\TT$ and points $E'$ to scale $C\tau$ to obtain a new $(C\tau, O(\eta/\varepsilon), s, t, \kappa, \gamma)$-Furstenberg configuration (which still has global grains from Proposition \ref{prop:global-grains-1}) and redefine $\alpha_p = f(z)$. Since $|\alpha_p-f(z)| \lesssim \tau$, the set of directions $\Theta_p$ can be slightly modified to satisfy property (4) from Definition \ref{def:global-grains}. So in this case we get the desired conclusion provided that $\eta/\varepsilon \ll \eta_0 $.

    Pick some $\varepsilon \ll \eta_0$ and suppose that $\tau \ge \d^{\varepsilon}$ holds. Let $\varrho \in [\d^{1/2}, \d^{\mu}]$ be the scale given by Lemma \ref{lem:good-scales} so that $Z \cap {\bf z}$ is rescaled $(\varrho, s, \varrho^{-\eta_0})$-set for ${\bf z} = z^\varrho$ (here $\mu=\mu(\eta_0)$). Suppose that $\varepsilon$ is small enough in terms of $\eta_0$ so that $\varepsilon \ll \mu$. Choose some $p_0 \in E$ and let ${\bf q} = p_0^\varrho$ and run the argument in Proposition \ref{prop:global-grains-1} inside ${\bf q}$. By repeating the proof, we get that after some pigeonholing the set $A_z({\bf q})$ for $z \in Z({\bf q})$ satisfies
    \[
    A_z({\bf q}) \subset N_{\varrho^2} \{ (x, y+\alpha_{\bf q} x), ~ y\in Y_z({\bf q}), x \in X_{y,z}({\bf q})  \}
    \]
    where we have $|X_{y,z}({\bf q})|_{\varrho^2} \ge \varrho^{\eta_0-1}$ and $Y_{z}({\bf q})$ is rescaled $(\varrho, 2t-\kappa-1, \varrho^{-O(\eta_0)}\d^{-O(\eta)})$-set. On the other hand, we know from global grains that $A_z({\bf q}) \subset A_z$ and 
    \[
    A_z = \{ (x,y+ f(z) x), ~x\in X_{y,z}, y \in Y_z \}
    \]
    where $|X_{y, z}|_\d \ge \d^{\eta-1}$ and $Y_z$ is $(\d, 2t-\kappa-1, \d^{-\eta})$-AD-regular set. By restricting to $E\cap {\bf q}$, we then can construct $c\varrho^2$-separated sets $\tilde X_{y, z}({\bf q})$ and $\tilde Y({\bf q})$ such that
    \[
    A_z({\bf q}) \subset N_{\varrho^2}\{ (x, y+f(z)x),~ x\in \tilde X_{y, z}({\bf q}), y\in \tilde Y_z({\bf q}) \}
    \]
    where we have $|\tilde X_{y,z}({\bf q})| \ge \d^{O(\eta)}\varrho^{-1}$ and $\tilde Y_{z}({\bf q})$ is rescaled $(\varrho, 2t-\kappa-1, \d^{-O(\eta)})$-AD-regular set.

    This gives us upper bounds on orthogonal projections
    \[
    |\pi_{(1,\alpha_{\bf q})}(A_z({\bf q}))|_{\varrho^2}, |\pi_{(1,f(z))}(A_z({\bf q}))|_{\varrho^2}\le \varrho^{-O(\eta_0)} \d^{-O(\eta)} \varrho^{\kappa+1-2t}
    \]
    and so using the lower bound $|\alpha_{\bf q}-f(z)|\ge \tau$ (recall that when we apply Lemma \ref{lem:pigeonhole-rho-ball} we can choose $\alpha_{\bf q}=\alpha_p$ for a typical $p \in E\cap {\bf q}$) this implies that 
    \[
    |A_z({\bf q})|_{\varrho^2} \lesssim \tau^{-1} |\pi_{(1,\alpha_{\bf q})}(A_z({\bf q}))|_{\varrho^2}\cdot  |\pi_{(1,f(z))}(A_z({\bf q}))|_{\varrho^2}
    \]
    \[
    \varrho^{\kappa-2t} \lesssim \tau^{-1}\varrho^{-O(\eta_0)} \d^{-O(\eta)} \varrho^{2(\kappa+1-2t)}
    \]
    Recall that $\tau \ge \varrho^{\varepsilon}$ so this implies that
    \[
    2t-\kappa \le 2 ( 2t-\kappa-1) + O(\varepsilon+\eta/\mu).
    \]
    Since by the assumption we have $2t-\kappa \le 2-\epsilon_0$, this gives a contradiction for small enough $\varepsilon,\eta$. We conclude that the alternative $\tau \le \d^\varepsilon$ must hold and we get the desired global grains Furstenberg configuration.
\end{proof}

In the next lemma we show that we may assume the function $f$ is dyadic Lipschitz on a given scale $\varrho$, provided that $\kappa$ is larger than some threshold. 

\begin{lemma}\label{lemma:f-Lipschitz}
    Consider a global grains $(\d,\eta, s,t,\kappa,\gamma)$-Furstenberg configuration with slope function $f$. Let $\varrho \in [\d^{1/4}, \d^{\mu}]$ be such that $Z \cap {\bf z}$, ${\bf z}=z^\varrho$ is rescaled $(\varrho, s, \varrho^{-\eps})$-set. 
    Suppose that $2t-\kappa-1 \le 1-\zeta$ for some $\zeta>0$. 
    Then there is a subset $\tilde Z \subset Z\cap {\bf z}$ of density $\d^{O(\eta/\zeta)}$ such that $|f(z_1)-f(z_2)| \le \varrho$ for all $z_1,z_2\in \tilde Z$. 
\end{lemma}

\begin{proof}
    Apply Lemma \ref{lem:pigeonhole-rho-ball-2} and fix a cube ${\bf q} = p_0^\varrho$ on the same level as ${\bf z}$. For $z \in Z({\bf q})$ let us write
     \begin{equation}\label{eq:Ezq-shape}
    E_z({\bf q}) = \{ (x, y+ f(z)x, z),~ x \in X_{y,z}({\bf q}), ~ y \in Y_z({\bf q}) \}     
     \end{equation}
    for rescaled $(\d/\varrho, 2t-1-\kappa, \d^{-O(\eta)})$-AD-regular set $Y_z({\bf q})$ and $(\d/\varrho, 1, \d^{-O(\eta)})$-AD-regular sets $X_{y, z}({\bf q})$. 
    We have $|\Theta_p \cap \Theta_{\bf q}^\varrho| \ge \d^{O(\eta)} |\Theta_p|$ for all $p \in E({\bf q})$. By pigeonhole, we can find $\varrho$-dyadic cube $\theta \subset \Theta_{\bf q}^\varrho$ so that $\Theta_p \cap \theta \neq\emptyset$ for at least $\d^{O(\eta)} |E({\bf q})|$-many $p \in E({\bf q})$. Since $\TT_p$ is $(\d, \kappa, \d^{-\eta})$-AD-regular, we get
    \[
    \sum_{T\in \TT({\bf q}): \theta(T) \in \theta} |Z({\bf q}, T)| \ge \d^{O(\eta)}   |Z({\bf q})| (\varrho/\d)^{2t}.
    \]
    It follows from Cauchy--Schwarz that
    \begin{equation}\label{eq:layers-overlap}
    \sum_{z,z' \in Z({\bf q})} |\pi_{(\theta,1)} ( E_z({\bf q}) ) \cap N_{C\varrho^2}\pi_{(\theta,1)}( E_{z'}({\bf q}))| \ge \d^{O(\eta)} |Z({\bf q})|^2 (\varrho/\d)^{2t-\kappa},    
    \end{equation}
    i.e. up to resolution $O(\varrho^2)$, the layers $E_z({\bf q})$ have essentially the same projection in direction $(\theta,1)$ (here we abuse notation and identify the cube $\theta$ with its middle point). 

    On the other hand, by (\ref{eq:Ezq-shape}) we have
    \[
    \pi_{(\theta,1)} ( E_z({\bf q}) ) = \{ (x, y + f(z) x) - \theta z,\quad x \in X_{y,z}({\bf q}), ~y\in Y_z({\bf q})   \}.
    \]
    Denote $\tau = \max(|f(z) - f(z')|, \varrho)$ for some fixed $z,z' \in Z({\bf q})$. 
    Cover the projection $\pi_{(\theta,1)} ( E_z({\bf q}) )$ by $\sim |Y_z({\bf q})|_{\tau\varrho}$-many $\tau\varrho \times \varrho$ tubes $T_j$. Then inside each tube we have a bound
    \[
    |T_j \cap \pi_{(\theta,1)} ( E_z({\bf q}) ) \cap N_{C\varrho^2} \pi_{(\theta,1)} ( E_{z'}({\bf q}) )|_{\varrho^2} \lesssim \d^{-O(\eta)} \tau^{-1} |Y_z({\bf q}) \cap  N_{\tau\varrho} y|_{\varrho^2} |Y_{z'}({\bf q}) \cap N_{\tau\varrho} y'|_{\varrho^2}
    \]
    where $y \in Y_z({\bf q})$ and $y' \in Y_{z'}({\bf q})$ are chosen so that the corresponding lines $(x, y+f(z)x)$ and $(x', y'+f(z')x')$ are contained in the tube $T_j$. 
    Indeed, each $\varrho^2$-ball $B_{\varrho^2}$ intersecting the set on the left hand side is contained in a $\varrho^2 \times \varrho^2/\tau$ rectangle $R$ proportional to $T_j$. This rectangle projects onto $\sim \varrho^2$-length intervals $I_1$, $I_2$ by the maps $(a,b) \to (a, b-f(z)a)$ and $(a,b) \to (a, b- f(z')a)$ and $R$ is essentially uniquely determined by $I_1, I_2$. 
    Now there are $\tau^{-1}$-many ways to select $B_{\varrho^2}\subset R$ and by definition we have $I_1 \subset N_{C\varrho^2} Y_z({\bf q})$ and $I_2 \subset N_{C\varrho^2} Y_{z'}({\bf q})$, giving the bound above.

    So summing over $T_j$ gives
    \begin{align*}
    |\pi_{(\theta,1)} ( E_z({\bf q}) ) \cap N_{C\varrho^2} \pi_{(\theta,1)} ( E_{z'}({\bf q}) )|_{\varrho^2} &\lesssim \d^{-O(\eta)} (\tau^{-1}|Y_z({\bf q})|_{\varrho\tau}) |Y_z({\bf q}) \cap N_{\tau\varrho} y|_{\varrho^2}|Y_{z'}({\bf q}) \cap  N_{\tau\varrho} y'|_{\varrho^2} \\
    &=  \d^{-O(\eta)} \varrho^{-1} |Y_z({\bf q})|_{\varrho^2} (\varrho/\tau) |Y_{z'}({\bf q}) \cap  N_{\tau\varrho} y'|_{\varrho^2} \\
    &= \d^{-O(\eta)} |E_z({\bf q})|_{\varrho^2} (\tau/\varrho)^{2t-2-\kappa}
    \end{align*}
    where we used that $Y_z({\bf q})$ is $(\d/\varrho, 2t-\kappa-1, \d^{O(\eta)})$-AD-regular. It follows from (\ref{eq:layers-overlap}) that
    \[
    \sum_{z, z' \in Z({\bf q})} \min( 1, \varrho / |f(z) - f(z')|  )^{2+\kappa-2t} \ge \d^{O(\eta)} |Z({\bf q})|^2.
    \]
    By assumption, we have $2+\kappa-2t \ge \zeta$. We conclude that some constant $C$ and for any $S\ge 1$ we have
    \begin{align*}
    \#\{ (z, z') \in Z({\bf q})\times Z({\bf q}): ~ |f(z) - f(z')| \le S\varrho \} +  \#\{ (z, z') \in Z({\bf q})\times Z({\bf q}): ~ |f(z) - f(z')| > S\varrho \} S^{-\zeta}\\ > \d^{C\eta} |Z({\bf q})|^2.    
    \end{align*}
    So by taking $S = \d^{ - 2C \eta /\zeta}$ we conclude 
    \[
    \#\{ (z, z') \in Z({\bf q})\times Z({\bf q}): ~ |f(z) - f(z')| \le \d^{-C\eta/\zeta}\varrho \} \gtrsim \d^{C\eta}|Z({\bf q})|^2.
    \]
    Thus, by pigeonholing we can find the desired set $\tilde Z \subset Z({\bf q})$ on which $|f(z)-f(z')|\le \varrho$, completing the proof. 
\end{proof}

\section{Linear and small variation of the slope function.}\label{sec:structured-slope}

In this section we analyze two special types of global grains Furstenberg configuration. The first is the linear slope function case: $f(z) = a z+b+ O(\varrho^{1+\varepsilon})$ for $z$ in a fixed $\varrho$-interval and $a\gtrsim 1$. The second is small variation case: $|f(z) - f(z')| \lesssim \varrho^{1+\varepsilon}$ for $z,z'$ in the same $\varrho$-interval. 

\subsection{\texorpdfstring{$W$}{W}-tuples.}
\label{subsec:W-tuples}

\medskip\noindent\emph{The tube relation.}
Our analysis is based on counting paths in the Furstenberg
configuration $(\TT,E,\{\TT_p\})$.  For tubes $T_1,T_2\in\TT$ and
$p\in E_z$, write
$T_1\underset p\leftrightarrow T_2$ if $T_1,T_2\in\TT_p$ and there
are $\xi\in\Xi_p$ and
$\varphi_1,\varphi_2\in\Phi_{p,\xi}$ such that
\[
\theta(T_i)\in
N_\d(\varphi_i,\xi+f(z)\varphi_i),
\qquad i=1,2.
\]
Thus the two tubes belong to the same $\gamma$-dimensional direction
fiber at $p$.  In particular,
\begin{equation}\label{eq:connected-tubes-bound}
\#\{(p,T_1,T_2):T_1\underset p\leftrightarrow T_2\}
\sim_{\d^{O(\eta)}}|E|\d^{-\kappa-\gamma}.
\end{equation}
Indeed, after choosing $p$, there are
$\d^{-\kappa+O(\eta)}$ choices for $T_1$.  Its direction determines
$\xi$ up to $O(1)$ choices, and the fiber
$\Phi_{p,\xi}$ gives $\d^{-\gamma+O(\eta)}$ choices for $T_2$.

\begin{obs}
    Suppose that $T_1 \underset{p}{\leftrightarrow} T_2$ for some $p\in E_z$ and let $\varphi_i$ be the first coordinate of $\theta(T_i)$. Then 
    \begin{equation}\label{eq:difference-theta}
        \theta(T_1) - \theta(T_2) = (\varphi_1-\varphi_2) (1, f(z)) + O(\d).
    \end{equation}
\end{obs}

\begin{proof}
    By definition, we have some $\xi \in \Xi_p$ and $\tilde \varphi_i \in \Phi_{p, \xi}$ so that $\theta(T_i) \in N_\d (\tilde\varphi_i, \xi +f(z) \tilde \varphi_i)$. It follows that $|\varphi_i - \tilde\varphi_i|\le \d$ and so subtracting gives the result. 
\end{proof}

Fix a scale $\varrho\in[\d,1]$ and a dyadic cube
${\bf q}=p_0^\varrho$.  We use the refinement
$E({\bf q})\subset E\cap{\bf q}$, together with
$Z({\bf q})$ and $\TT({\bf q})$, given by Lemma
\ref{lem:pigeonhole-rho-ball-2}. Using Lemma \ref{lemma:f-Lipschitz}, we can refine $Z({\bf q})$ so that $|f(z)-f(z')| \lesssim \varrho$ for all $z,z '\in Z({\bf q})$ (unless $\kappa \ge 2t-2-\zeta$ in which case we obtain a satisfactory lower bound on $\kappa$).

\medskip\noindent\emph{Basic $V$-tuples.}
For $z,z'\in Z({\bf q})$, let $V_{z,z'}({\bf q})$ consist of the
tuples
\[
(p,p_1,p_2,T_1,T_2)
\in E_z({\bf q})\times E_{z'}({\bf q})^2
\times\TT({\bf q})^2
\]
such that
$T_i\in\TT_p\cap\TT_{p_i}$ for $i=1,2$ and
$T_1\underset p\leftrightarrow T_2$.

\begin{lemma}
We have
\begin{equation}\label{eq:V-tuples-counting-1}
\sum_{z,z'\in Z({\bf q})}|V_{z,z'}({\bf q})|
\sim_{\d^{O(\eta)}}
|E({\bf q})|\d^{-\kappa-\gamma}|Z({\bf q})|.
\end{equation}
\end{lemma}

\begin{proof}
For $p\in E({\bf q})$, $\xi\in\Xi_p$, and
$z'\in Z({\bf q})$, define
\[
\begin{split}
\mathcal S_{z'}(p,\xi)
=\{(p',T):\;&p'\in E_{z'}({\bf q}),\
T\in\TT_p\cap\TT_{p'},\\
&\theta(T)\in
N_\d(\varphi,\xi+f(z)\varphi)
\text{ for some }\varphi\in\Phi_{p,\xi}\}.
\end{split}
\]
The incidence estimate in Lemma
\ref{lem:pigeonhole-rho-ball-2}(b) gives
\begin{equation}\label{eq:S-first-moment-Furstenberg}
\sum_{p,\xi,z'}|\mathcal S_{z'}(p,\xi)|
\sim_{\d^{O(\eta)}}
|E({\bf q})|\d^{-\kappa}|Z({\bf q})|.
\end{equation}
The number of nonempty indices $(p,\xi,z')$ is at most
\[
\d^{-O(\eta)}
|E({\bf q})|\d^{-(\kappa-\gamma)}|Z({\bf q})|.
\]
An ordered pair in the same $\mathcal S_{z'}(p,\xi)$ is precisely a
$V$-tuple, up to bounded multiplicity.  Hence Cauchy--Schwarz and
\eqref{eq:S-first-moment-Furstenberg} give
\[
\sum_{z,z'}|V_{z,z'}({\bf q})|
\gtrsim_{\d^{O(\eta)}}
|E({\bf q})|\d^{-\kappa-\gamma}|Z({\bf q})|.
\]
For the reverse bound, use
$|\mathcal S_{z'}(p,\xi)|
\lesssim_{\d^{O(\eta)}}\d^{-\gamma}$ and
\eqref{eq:S-first-moment-Furstenberg}.
\end{proof}

For a tuple $(p,p_1,p_2,T_1,T_2) \in V_{z,z'}({\bf q})$, write
\[
p=(x,y+f(z)x,z),\qquad
p_i=(x_i,y_i+f(z')x_i,z'),
\]
and choose $\xi\in\Xi_p$ and
$\varphi_i\in\Phi_{p,\xi}$ corresponding to $T_i$.  The incidence
relations
$p_i=p+(z'-z)(\theta(T_i),1)+O(\d)$ give
\[
\begin{aligned}
x_i-x&=(z'-z)\varphi_i+O(\d),\label{eq:V-tuple-x}\\
y_i+f(z')x_i-(y+f(z)x)
&=(z'-z)(\xi+f(z)\varphi_i)+O(\d).
\end{aligned}
\]
Subtracting the equations with $i=1,2$ and eliminating
$x_1-x_2$ gives
\begin{equation}\label{eq:V-tuple-relation-Furstenberg}
y_1-y_2
=(z'-z)(f(z)-f(z'))(\varphi_1-\varphi_2)+O(\d).
\end{equation}

We will use \eqref{eq:V-tuple-relation-Furstenberg} in the linear case (Section \ref{subsec:linear-grains-Furstenberg}).
For the small-variation case (Section \ref{subsec:slow-grains-Furstenberg}) we need to use longer paths, as defined below.

\begin{figure}[ht]
\begin{tikzpicture}[
    scale=1.0,
    line join=round,
    cube/.style={draw=black, thick},
    mainline/.style={draw=black, very thick},
    redline/.style={draw=red!70!black, thick},
    every node/.style={inner sep=0pt, outer sep=2pt, fill=none, draw=none}
]


\coordinate (A) at (0, 0);
\coordinate (B) at (5, 0);
\coordinate (C) at (7, 2);
\coordinate (D) at (2, 2);

\coordinate (E) at (0, 5);
\coordinate (F) at (5, 5);
\coordinate (G) at (7, 7);
\coordinate (H) at (2, 7);


\draw[cube] (A) -- (B) -- (C) -- (D) -- cycle;
\draw[cube] (E) -- (F) -- (G) -- (H) -- cycle;
\draw[cube] (A) -- (E);
\draw[cube] (B) -- (F);
\draw[cube] (C) -- (G);
\draw[cube] (D) -- (H);


\coordinate (Z1-A) at ($(A)!0.33!(E)$);
\coordinate (Z1-B) at ($(B)!0.33!(F)$);
\coordinate (Z1-C) at ($(C)!0.33!(G)$);
\coordinate (Z1-D) at ($(D)!0.33!(H)$);
\draw[cube] (Z1-A) -- (Z1-B) -- (Z1-C) -- (Z1-D) -- cycle;

\coordinate (Z2-A) at ($(A)!0.66!(E)$);
\coordinate (Z2-B) at ($(B)!0.66!(F)$);
\coordinate (Z2-C) at ($(C)!0.66!(G)$);
\coordinate (Z2-D) at ($(D)!0.66!(H)$);
\draw[cube] (Z2-A) -- (Z2-B) -- (Z2-C) -- (Z2-D) -- cycle;


\coordinate (P) at (3.5, 1);
\coordinate (P1) at (3.0, 2.32);
\coordinate (P2) at (4.0, 2.16);
\coordinate (P1p) at (3.4, 4.32);
\coordinate (P2p) at (4.4, 4.16);


\draw[mainline] ($(P)!-0.3!(P1)$) -- ($(P)!4.0!(P1)$);
\draw[mainline] ($(P)!-0.3!(P2)$) -- ($(P)!4.0!(P2)$);
\draw[redline] ($(P1)!-0.6!(P1p)$) -- ($(P1)!2.6!(P1p)$);
\draw[redline] ($(P2)!-0.6!(P2p)$) -- ($(P2)!2.6!(P2p)$);


\filldraw[black] (P) circle (1pt);
\filldraw[black] (P1) circle (1pt);
\filldraw[black] (P2) circle (1pt);
\filldraw[red!70!black] (P1p) circle (1pt);
\filldraw[red!70!black] (P2p) circle (1pt);


\draw (-0.7, 0.3) node {$Z$};
\draw (-0.5, 2.0) node {$Z'$};
\draw (-0.3, 3.8) node {$Z''$};

\draw (3.5, 0.4) node {$P$};
\draw (2.6, 2.32) node {$P_1$};
\draw (4.4, 2.16) node {$P_2$};
\draw (3.1, 4.32) node[red!70!black] {$P_1'$};
\draw (4.1, 4.16) node[red!70!black] {$P_2'$};

\draw (2.6, 6.3) node {$T_1$};
\draw (5.2, 6.0) node {$T_2$};
\draw (4.0, 7.5) node[red!70!black] {$T_1'$};
\draw (5.3, 7.5) node[red!70!black] {$T_2'$};


\draw[{Stealth[length=4mm]}-{Stealth[length=4mm]}, thick] (8.5, 2) -- (8.5, 6);
\draw (8.8, 4) node {$\varrho$};

\end{tikzpicture}

\caption{$W$-tuple $(p, p_1, p_2, p_1', p_2', T_1,T_2, T'_1, T'_2)$.}
\label{fig:W-tuple-furstenberg}
\end{figure}
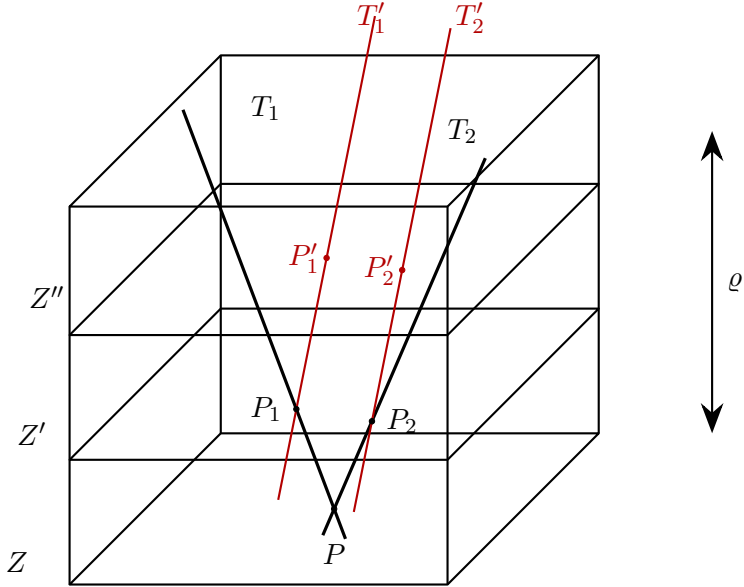

\medskip\noindent\emph{Definition of the $W$-tuples.}
For $z,z',z''\in Z({\bf q})$, let $W_{z,z',z''}({\bf q})$ be the
set of tuples
\[
\begin{aligned}
W_{z,z',z''}({\bf q}) = \{
(p,p_1,p_2,p_1',p_2',T_1,T_2,T_1',T_2')
\in E_z({\bf q})\times E_{z'}({\bf q})^2
\times E_{z''}({\bf q})^2\times\TT({\bf q})^4,   \\
T_1'\underset{p_1}{\leftrightarrow}T_1,\qquad
T_1\underset p\leftrightarrow T_2,\qquad
T_2\underset{p_2}{\leftrightarrow}T_2',\\
T_i'\in\TT_{p_i'}\quad(i=1,2),
\qquad
|\theta(T_1')-\theta(T_2')|\leq\varrho\}
\end{aligned}
\]
See Figure \ref{fig:W-tuple-furstenberg} for an illustration. 

\medskip\noindent\emph{Counting the $W$-tuples.}
The next lemma gives a sharp up to $\d^{O(\eta)}$ bound on the number of $W$-tuples. 

\begin{lemma}\label{lem:counting-W-tuples}
Let $\widetilde Z\subset Z({\bf q})$.  Then
\[
\sum_{z,z',z''\in\widetilde Z}|W_{z,z',z''}({\bf q})|
\gtrsim_{\d^{O(\eta)}}
(\varrho/\d)^{2t-\kappa}
\d^{-\kappa-3\gamma}\varrho^\gamma
|\widetilde Z|^3.
\]
\end{lemma}

\begin{proof} 
    For $(z, z', z'') \in \tilde Z^3$ consider the following set:
    \begin{align*}
    P_{z, z', z''} = \{ (p, p_1, p_1', T_1, T_1') \in E_z({\bf q})\times E_{z'}({\bf q})\times E_{z''}({\bf q}) \times \TT({\bf q})\times \TT({\bf q}):\\
    T_1'\underset{p_1}{\leftrightarrow} T_1, \quad T_1 \in \TT_{p}, T_1' \in \TT_{p_1'} \}. 
    \end{align*}
    Using (\ref{eq:connected-tubes-bound}) we then have
    \[
    \sum_{z, z', z''\in \tilde Z} |P_{z, z', z''}| \sim_{\d^{O(\eta)}} |E_{\tilde Z}({\bf q})| \d^{-\kappa-\gamma} |\tilde Z |^{2}
    \]
    where $E_{\tilde Z}({\bf q}) = \bigcup_{z \in \tilde Z} E_z({\bf q})$. 
    Indeed, first we pick the triple $T_1' \underset{p_1}{\leftrightarrow} T_1$ and then select $z \in \tilde Z({\bf q}, T_1), z'' \in \tilde Z({\bf q}, T'_1)$ and pick the corresponding points $p \in E_z({\bf q})$ and $p'_1 \in E_{z''}({\bf q})$. By (b) in Lemma \ref{lem:pigeonhole-rho-ball-2} we get the desired bound. 
    

    Let $\overline{\Theta}_{\bf q} \subset (\frac\varrho 2 \ZZ)^2$ be a set whose $\varrho/2$-neighborhood covers $\Theta_q$. 
    For $(z, z', z'') \in \tilde Z^3$, $p \in E_z({\bf q})$, $\xi \in \Xi_p$ and $\theta \in \overline{\Theta}_{\bf q}$ consider the set $P_{z, z', z''}(p, \xi, \theta) \subset P_{z, z', z''}$ consisting of tuples $(p_1, p_1', T_1, T_1')$ such that $\theta(T_1) \in N_\d (\varphi, \xi+f(z)\varphi)$ for some $\varphi \in \Phi_{p, \xi}$ and $\theta(T_1') \in N_{\varrho/2} \theta$. Clearly, $P_{z, z', z''}$ is covered by finitely overlapping sets $P_{z, z', z''}(p, \xi, \theta)$. For fixed $p, \xi$, there are at most $\d^{-O(\eta)}\varrho^{-\gamma}$ choices for $\theta \in \overline{\Theta}_{\bf q}$ for which the set $P_{z, z', z''}(p, \xi, \theta)$ is non-empty.
    
     It follows by Cauchy--Schwarz that 
     \begin{align*}
     \sum_{(z,z',z'') \in S} \sum_{p, \xi, \theta} |V_{z, z', z''}(p, \xi, \theta)|^2 &\ge \frac{\d^{O(\eta)}}{|E_{\tilde Z}({\bf q})| \d^{-(\kappa-\gamma)} \varrho^{-\gamma}|\tilde Z|^{2}} \left(  \sum_{(z,z',z'') \in S} \sum_{p, \xi, \theta} |V_{z, z', z''}(p, \xi, \theta)|\right)^2\\
    &\ge \d^{O(\eta)} |E_{\tilde Z}({\bf q})| \d^{-\kappa-3\gamma}\varrho^\gamma |\tilde Z|^{2}.
    \end{align*}
    Any element of the cartesian product $V_{z, z', z''}(p, \xi, \theta) \times V_{z, z', z''}(p, \xi, \theta)$ gives rise to an element of $W_{z,z',z''}({\bf q})$. Indeed, let $(p_1, p_1',T_1, T_1'), (p_2, p_2',T_2, T_2') \in V_{z, z', z''}(p, \xi, \theta)$. Then it follows that $T_1 \underset{p}{\leftrightarrow} T_2$ (since $\theta(T_i) \in N_\d (\varphi_i, \xi+f(z)\varphi_i)$) and $|\theta(T_1') - \theta(T_2')| \le \varrho$ (since $|\theta(T'_i) - \theta| \le \varrho/2$). Each $W$-tuple occurs this way at most constantly many times, so this gives the desired lower bound. 
\end{proof}

\medskip\noindent\emph{Additive relations.}
Fix a tuple $(p,p_1,p_2,p_1',p_2',T_1,T_2,T_1',T_2')$ in $W_{z,z',z''}({\bf q})$, and let
$\varphi_i,\varphi_i'$ be the first direction coordinates of
$T_i,T_i'$.  Tube incidence relations and
\eqref{eq:difference-theta} give
\[
\begin{aligned}
\theta(T_1)-\theta(T_2)
&=(\varphi_1-\varphi_2)(1,f(z))+O(\d),\\
\theta(T_i')-\theta(T_i)
&=(\varphi_i'-\varphi_i)(1,f(z'))+O(\d),
\qquad i=1,2.
\end{aligned}
\]
Consequently,
\begin{equation}\label{eq:diff-theta-prime}
\begin{split}
\theta(T_1')-\theta(T_2')
&=(\varphi_1'-\varphi_2')(1,f(z'))\\
&\quad
+(\varphi_1-\varphi_2)(0,f(z)-f(z'))+O(\d).
\end{split}
\end{equation}

Write
\[
p_i=(x_i,y_i+f(z')x_i,z'),\qquad
p_i'=(x_i',y_i'+f(z'')x_i',z'').
\]
Since $(p,p_1, p_2, T_1,T_2)$ form a $V$-tuple, (\ref{eq:V-tuple-x}) and (\ref{eq:V-tuple-relation-Furstenberg}) give
\begin{equation}\label{eq:y-axes-relation}
\begin{aligned}
x_1-x_2
&=(z'-z)(\varphi_1-\varphi_2)+O(\d),\\
y_1-y_2
&=(z'-z)(f(z)-f(z'))
(\varphi_1-\varphi_2)+O(\d).
\end{aligned}
\end{equation}
Using the relations $p_i' = p_i + (z''-z') \theta(T'_i) + O(\d)$ we similarly obtain
\[
\begin{aligned}
x_1'-x_2'
&=(z'-z)(\varphi_1-\varphi_2)
+(z''-z')(\varphi_1'-\varphi_2')+O(\d),\\
y_1'-y_2'
&=(z'-z)(f(z)-f(z''))
(\varphi_1-\varphi_2)\\
&\quad
+(z''-z')
\bigl((f(z)-f(z'))(\varphi_1-\varphi_2)
+(f(z')-f(z''))(\varphi_1'-\varphi_2')\bigr)
+O(\d).
\end{aligned}
\]
Since a $W$-tuple satisfies
$|\theta(T_1')-\theta(T_2')|\leq\varrho$, we have
$|\varphi_1'-\varphi_2'|\leq\varrho$.  The last display therefore
implies the following two estimates that will be used later:
\begin{align}
x_1'-x_2'
&=(z'-z)(\varphi_1-\varphi_2)
+O(\varrho^2+\d),\label{eq:x-axes-prime-2}\\
y_1'-y_2'
&=O\bigl(
\varrho|f(z)-f(z')|
+\varrho|f(z')-f(z'')|+\d
\bigr).\label{eq:y-axes-prime-2}
\end{align}

\subsection{Linear slope function.} \label{subsec:linear-grains-Furstenberg}

In this section we consider the case when $f$ can be locally approximated by a linear function. 

\begin{lemma}\label{lem:linear-grains}
    For any $\zeta, \mu >0, \eta_0\in (0,1/2)$ the following holds for all $\eta < \eta(\zeta, \mu, \eta_0)$ and $\eps < \eps(\zeta, \eta_0)$. 

    Consider a global grains $(\d, \eta, s, t, \kappa, \gamma)$-Furstenberg configuration with $s+\gamma \ge 1-\zeta/4$ and let $f: Z\to [-1,1]$ be the slope function. Let $\varrho \in [\d^{1/10}, \d^{\mu}]$ and suppose that $\tilde Z\subset Z\cap {\bf z}$ is a rescaled $(\varrho, s, \varrho^{-\epsilon})$-set for every ${\bf z} = z_0^{\varrho}$, $z_0\in \tilde Z$.
    Suppose that for some $|r| \ge \varrho^{5/2}$ and for all $z , z'\in \tilde Z$ we have 
    \begin{equation}\label{eq:linear-approximation-furstenberg}
    f(z') = f(z) + r(z'-z) + O(r \varrho^{1+\eta_0}).    
    \end{equation}
    Then we have $\kappa \le 2t-2+\zeta$. 
\end{lemma}

\begin{proof}

    Replace $E $ by the set $E'$ given by Lemma \ref{lem:pigeonhole-rho-ball-2}. 
Fix a cube ${\bf q} = p_0^{\varrho}$ for some $p_0 \in E$ and let $Z({\bf q})$ and $E({\bf q})$ be the corresponding sets. It is not hard to verify that one can perform this pigeonholing step and guarantee that the set $\tilde Z$ from the statement of the lemma satisfies $|\tilde Z \cap Z({\bf q})| \gtrsim_{\d^{O(\eta)}} |\tilde Z|$. So we can replace $Z({\bf q})$ by $\tilde Z \cap Z({\bf q})$ and by (b) in Lemma \ref{lem:pigeonhole-rho-ball-2} keep the desired double-counting properties. So we may assume that the approximation (\ref{eq:linear-approximation-furstenberg}) holds for every $z,z' \in Z({\bf q})$. Furthermore, it still holds that $Z({\bf q})$ is a rescaled $(\varrho, s, \varrho^{-\eps} \d^{-C\eta})$-set. 

For $z \in Z({\bf q})$ we thus can write 
    \[
    E_z({\bf q}) = \{ (x, y+ f(z)x, z),~ x \in X_{y,z}({\bf q}), ~ y \in Y_z({\bf q}) \}
    \]
    for rescaled $(\d/\varrho, 2t-1-\kappa, \d^{-O(\eta)})$-AD-regular set $Y_z({\bf q})$ and $(\d/\varrho, 1, \d^{-O(\eta)})$-AD-regular set $X_{y, z}({\bf q})$. 
    Let $\overline{\Xi}_{\bf q} \subset \varrho\ZZ$ be the $\varrho/2$-covering of $\Xi_{\bf q}$. For $\xi \in \overline{\Xi}_{\bf q}$ let $\overline{\Phi}_{{\bf q},\xi} \subset \varrho \ZZ$ be the $\varrho/2$-covering of $\bigcup_{\xi' \in N_\varrho\xi}\Phi_{{\bf q}, \xi'}$. Since $\Theta_{\bf q}$ is $(\varrho, \kappa, \d^{-O(\eta)})$-AD-regular, we have $\sum_{\xi \in \overline{\Xi}_{\bf q}} |\overline{\Phi}_{{\bf q}, \xi}| \sim_{\d^{O(\eta)}} \varrho^{-\kappa}$. 
    Let $\overline{\Psi}$ be the set of triples $(\xi, \varphi_1, \varphi_2)$ such that $\varphi_1, \varphi_2 \in \overline{\Phi}_{{\bf q}, \xi}$.
    
    Let $\overline{Z} = Z({\bf q})-z$. We define a graph $G_z \subset Y_z({\bf q}) \times \overline{Z} \times \overline{\Psi}$ as follows. For any tuple
    $(p,p_1,p_2, T_1, T_2) \in V_{z',z}({\bf q})$ we write $p_i = (x_i, y_i + f(z) x_i, z)$ and $\theta(T_i) \in N_\varrho(\varphi_i, \xi+f(z')\varphi_i)$ for some $\xi \in \overline{\Xi}_{\bf q}$ and $\varphi_i \in \overline{\Phi}_{{\bf q}, \xi}$.
    We include the tuple $(y_1, z'-z, (\xi,\varphi_1,\varphi_2))$ in $G_z$.
    By construction and (\ref{eq:V-tuple-relation-Furstenberg}), for any $(y, \overline{z}, (\xi, \varphi_1, \varphi_2)) \in G_z$ there is $y' \in Y_z({\bf q})$ such that 
    \[
    y'= y+ \overline{z} (f(z+\overline{z}) - f(z))(\varphi_1-\varphi_2 + O(\varrho))  + O(\d).
    \]
    So by linearity assumption on $f$, 
    \begin{equation}\label{eq:y-prime-also-in-Y}
    y'= y+ r \overline{z}^2 (\varphi_1-\varphi_2 ) + O(r \varrho^{2+\eta_0}+r \varrho^3+\d).    
    \end{equation}
    Cover $Y_z({\bf q})$ by finitely overlapping intervals $I_{z,j}$ of length $C r \varrho^2$, so that $|Y_z({\bf q}) \cap I_{z,j}| \sim_{\d^{O(\eta)}} (r\varrho^2/\d)^{2t-\kappa-1}$ holds for every index $j$. We want to apply Theorem \ref{thm: ABC} to rescaled sets $Y_z({\bf q}) \cap I_{z,j}$ on scale $\tau := \varrho^{\eta_0}$ with $\alpha = 1-\zeta/2$, $\beta = s$,$\gamma = \gamma$. Let $\chi$ be the constant given by Theorem \ref{thm: ABC}. 
     The set $\{r \overline{z}^2, ~ \overline{z} \in \overline{Z}\}$ is clearly a rescaled $(\varrho, s, C\varrho^{-\epsilon})$-set.
     Let $\overline{\Phi}$ be the multiset of all differences $\varphi_1-\varphi_2$ taken over all triples $(\xi, \varphi_1, \varphi_2) \in \overline{\Psi}$. Then $\overline{\Phi}$ forms a $(\varrho, \gamma, \d^{-O(\eta)})$-set (viewed as a multiset; one can turn it into a proper set by a dyadic pigeonholing argument).
    So if $\eta \ll  \chi \mu \eta_0$ then $\overline{\Phi}$ is $(\tau, \gamma, \tau^{-\chi})$-set and if $\eps  \ll \chi \eta_0$ then  $\{r \overline{z}^2, ~ \overline{z} \in \overline{Z}\}$ is a rescaled $(\tau, s, \tau^{-\chi})$-set (restrict $\overline{Z}$ to be well separated from 0).

    From (\ref{eq:V-tuples-counting-1}) and $|\overline{\Phi}|=|\overline{\Psi}| \sim_{\d^{O(\eta)}} \varrho^{-\kappa-\gamma}$, we get 
    \[
    \sum_{z \in Z({\bf q})} |G_z| \ge \d^{O(\eta)} |\overline{Z}| |\overline{\Psi}| \sum_{z \in Z({\bf q})} |Y_z({\bf q})|.
    \]
    So for some $z \in Z({\bf q})$ we get $|G_z| \ge \d^{O(\eta)} |\overline{Z}| |\overline{\Phi}|  |Y_z({\bf q})|$. Let $G_{z,j}$ be the restriction of $G_z$ onto the set $Y_z({\bf q}) \cap I_{z,j}$. Then for some $j$ we have $|G_{z,j}| \ge \d^{O(\eta)} |\overline{Z}| |\overline{\Phi}|  |Y_z({\bf q}) \cap I_{z,j}|$. For $\eta \ll \chi\mu\eta_0 $ we get that $\d^{O(\eta)} \ge \tau^\chi$. If $|Y_{z, j}({\bf q})|_{r\varrho^2 \tau} < \tau^{\zeta/2-1}$ then by Theorem \ref{thm: ABC} there is $\overline{\varphi} \in \overline{\Phi}$ such that
    \[
    |\{ y+ r \overline{z}^2\overline{\varphi},~~ (y, \overline{z}, \overline{\varphi}) \in G_{z,j} \}|_{r\varrho^2 \tau} > \tau^{-\chi} |Y_{z, j}({\bf q})|_{r\varrho^2 \tau}.
    \]
    However, we know that $|Y_{z, j}({\bf q})|_{r\varrho^2 \tau} \sim_{\d^{O(\eta)}} \tau^{-(2t-\kappa-1)}$ and that by (\ref{eq:y-prime-also-in-Y}) we have
    \[
    \{ y+ r\overline{z}^2\overline{\varphi},~~ (y, \overline{z}, \overline{\varphi}) \in G_{z,j} \} \subset N_{Cr \varrho^2 \tau}Y_z({\bf q}) \cap CI_{z,j}
    \]
    and so by AD-regularity of the set $Y_z({\bf q})$, it has $r\varrho^2\tau$-covering number at most $\d^{-O(\eta)} \tau^{-(2t-\kappa-1)}$. By taking $\eta$ small enough, we obtain a contradiction. We conclude that 
    \[
    \d^{-O(\eta)} \tau^{-(2t-\kappa-1)} \ge |Y_z({\bf q}) \cap I_{z,j}|_{r\varrho^2 \tau} \ge \tau^{\zeta/2-1}
    \]
    giving $2t-\kappa-1 \ge 1 - \zeta/2 + O(\eta / \mu\eta_0)$. Since $\eta$ is sufficiently small, this gives us the desired bound $\kappa \le 2t-2 +\zeta$.  
\end{proof}

\subsection{Small variation of slope function.}\label{subsec:slow-grains-Furstenberg}

In this section we show that if the function $f$ has small variation then this implies that tubes in $\TT$ organize into planks. This will be used to contradict the Convex Wolff axiom. 

\begin{lemma}\label{lem:slowly-rotating}
    For any $\zeta, \mu>0$ the following holds for all sufficiently small $\eps < \eps(\zeta)$ and $\eta < \eta(\zeta, \mu)$.
    Consider a global grains $(\d, \eta, s, t, \kappa, \gamma)$-Furstenberg configuration with $s+\gamma \ge 1$ and let $f: Z\to [-1,1]$ be the slope function. Let $\epsilon>0$ and $\varrho \in [\d^{1/4}, \d^\mu]$ and suppose that $\tilde Z \subset Z\cap {\bf z}$, for some ${\bf z}=z_0\in Z$, is a rescaled $\d^{-\eta}$-uniform $(\varrho, s, \varrho^{-\epsilon})$-set such that  for all $z, z' \in \tilde Z$ we have $|f(z) - f(z')| \le \sigma$ for some $\sigma \in [\max(\d^{1/4}, \varrho^2), \varrho]$.
    
    Then there exists a $\sigma\times \varrho\times 1$ box $U$ such that
    \[
    |\TT[U] |_{\sigma} \ge  \varrho^{\zeta} (\varrho/\sigma)^{1+\gamma}.
    \]
\end{lemma}

\begin{proof}
    Replace $E $ by the set $E'$ given by Lemma \ref{lem:pigeonhole-rho-ball-2}. 
Fix a cube ${\bf q} = p_0^{\varrho}$ for some $p_0 \in E$ and let $Z({\bf q})$ and $E({\bf q})$ be the corresponding sets. It is not hard to verify that one can perform this pigeonholing step and guarantee that the set $\tilde Z$ from the statement of the lemma satisfies $|\tilde Z \cap Z({\bf q})| \gtrsim_{\d^{O(\eta)}} |\tilde Z|$. Thus, by Lemma \ref{lem:pigeonhole-rho-ball-2} (b) we can replace $Z({\bf q})$ by $\tilde Z \cap Z({\bf q})$ and keep the double counting properties. Then it follows that $Z({\bf q})$ is a rescaled $(\varrho, s, \varrho^{-\eps} \d^{-C\eta})$-set and $|f(z)-f(z')| \le \sigma$ for all $z,z' \in Z({\bf q})$. 

For $z \in Z({\bf q})$ we thus can write 
    \[
    E_z({\bf q}) = \{ (x, y+ f(z)x, z),~ x \in X_{y,z}({\bf q}), ~ y \in Y_z({\bf q}) \}
    \]
    for rescaled $(\d/\varrho, 2t-1-\kappa, \d^{-O(\eta)})$-AD-regular set $Y_z({\bf q})$ and $(\d/\varrho, 1, \d^{-O(\eta)})$-AD-regular set $X_{y, z}({\bf q})$. 

    For $T\in \TT$ and $z \in Z$ let
    \[
    U(z, T) = N_{\sigma} T + [-\varrho, \varrho] \cdot (1,f(z),0).
    \]
    Note that $U(z, T)$ is approximately a $\sigma\times \varrho\times 1$ box. 

    Consider the set of tuples $W_{z,z', z''}({\bf q})$ defined for the set $E({\bf q})$ in Section \ref{subsec:W-tuples}. 
    For a pair of tubes $T_1', T_2'$ and $z'' \in Z({\bf q})$ we write $(z'', T_1') \sim (z'', T'_2)$ if they can be extended to a tuple $(p,p_1',p_2', p_1'', p_2'', T_1, T_2, T'_1, T'_2) \in W_{z,z', z''}({\bf q})$ for some $z, z' \in Z({\bf q})$.

    \begin{claim}\label{claim:adjacent-tubes-same-box}
        If $(z'', T_1') \sim (z'', T'_2)$ then $U(z'', T_1') \subset C\cdot U(z'', T_2')$.
    \end{claim}

    \begin{proof}
    Consider a tuple $(p, p_1, p_2, p_1', p_2', T_1,T_2, T'_1, T'_2) \in W_{z, z', z''}({\bf q})$ and write $p = (x, y+ f(z) x, z)$, $p_i = (x_i, y_i + f(z') x_i, z' )$ and $p_i' = (x_i', y'_i+ f(z'') x'_i,z'')$. Let $\varphi_1,\varphi_2, \varphi_1', \varphi_2'$ denote the first coordinates of the slopes $\theta(T_1), \theta(T_2), \theta(T_1'), \theta(T_2')$. Then we have by (\ref{eq:diff-theta-prime}):
    \begin{align*}
    \theta(T_1') - \theta(T_2') = (\varphi_1'-\varphi_2')(1, f(z')) + (\varphi_1-\varphi_2) (0, f(z)-f(z')) + O(\d)     \\
    =  (\varphi_1'-\varphi_2')(1, f(z')) +O(\sigma) = O(\varrho) (1, f(z'')) + O(\sigma).
    \end{align*}
    and by (\ref{eq:x-axes-prime-2}), (\ref{eq:y-axes-prime-2}):
    \begin{align*}
    x_1' - x_2' &= (z'-z) (\varphi_1-\varphi_2) + O(\varrho^2) = O(\varrho),\\
    y_1'-y_2' &=  O(\varrho \sigma+\d) = O(\varrho\sigma).
    \end{align*}
    From this we conclude that $p_2' \in U(z'', T_1')$ and $T_2' \subset C\cdot U(z'', T_1')$ which implies
    $U(z'', T'_1)$ and $U(z'', T'_2)$ define essentially the same box, concluding the proof. 
    \end{proof}

Let $\overline{\Xi}_{\bf q} \subset \varrho \ZZ$ be the set whose $\varrho/2$ neighborhood covers $\Xi_{{\bf q}}$ and let $\overline{\Phi}_{{\bf q},\xi} \subset \varrho \ZZ$ be the set whose $\varrho/2$ neighborhood covers $\Phi_{{\bf q}, \xi}$. 
Let $\mc P$ be the set of all pairs $(z, T)$ with $T \in \TT({\bf q})$ and $z \in Z({\bf q}, T)$. 
For a pair $(z'', T'_1) \in Z({\bf q}) \times \TT({\bf q})$ we define a subset 
\[
G(z'', T'_1) \subset Z({\bf q}) \times Z({\bf q}) \times \overline{\Phi}_{{\bf q}, \xi} \times \overline{\Phi}_{{\bf q}, \xi}
\]
where $\xi \in \overline{\Xi}_{\bf q}$ is the element such that $\theta(T'_1) \in N_{\varrho} (\varphi'_1, \xi+ f(z'') \varphi'_1)$. Namely, for each tuple 
\[
(p,p_1',p_2', p_1'', p_2'', T_1, T_2, T'_1, T'_2) \in W_{z,z', z''}({\bf q})
\]
we consider $\varphi_1, \varphi_2 \in \overline{\Phi}_{{\bf q}, \xi}$ so that we have $\theta(T_i) \in N_{\varrho}(\varphi_i, \xi + f(z) \varphi_i)$. We then include the tuple $(z, z', \varphi_1, \varphi_2)$ into $G(z'', T_1')$. Using the lower bound on $W_{z, z', z''}({\bf q})$ in Lemma \ref{lem:counting-W-tuples} we get
\[
\frac{1}{|E({\bf q})| \d^{-\kappa}}\sum_{z'', T_1'} |G(z'', T_1')| \ge \d^{O(\eta)} |Z({\bf q})|^2 \varrho^{-2\gamma}
\]
where we note that the number of valid pairs $(z'', T_1')$ is approximately $|E({\bf q})| \d^{-\kappa}$. So it follows that on average the graph $G(z'', T_1')$ is $\d^{O(\eta)}$-dense. Let $\mc P' \subset \mc P$ be the subset of all $(z'', T'_1)$ with $|G(z'', T_1')| \ge \d^{C_1 \eta}|Z({\bf q})|^2 \varrho^{-2\gamma}$ for some large constant $C_1$. It follows from the above that most of the $W$-tuples are captured by pairs $(z, T) \sim (z, T')$ with both $(z,T), (z, T') \in \mc P'$. Let $\mc G' \subset \mc P'\times \mc P'$ be the adjacency graph for the relation $(z, T) \sim (z, T')$.

Fix some $(z, T) \in \mc P'$ and for an integer $m \ge 1$ consider the $\ell$-th neighborhood set $\mc N_m(z, T) \subset \mc P'$ in $\mc G$, i.e. this is the set of $(z, T') \in \mc P'$ which are connected to $(z, T)$ by a path with at most $m$ edges (note that the $z$-coordinate doesn't change along the edges of $\mc G'$). By Claim \ref{claim:adjacent-tubes-same-box} all pairs $(z,T') \in \mc N_{\ell}(z, T)$ are all contained in the same convex set $Cm \cdot U(z, T)$. So it suffices to show that for a bounded $m$, the neighborhood is large. 

Now we observe that for $(z'', T'_1) \sim (z'', T_2')$, if we denote by $p_1' = (x_1', y_1'+f(z'') x_1', z'')$ and $p_2'=(x_2', y_2'+f(z'')x_2', z'')$ the points in $T_1' \cap E_{z''}$ and $T_2'\cap E_{z''}$, respectively, then by (\ref{eq:x-axes-prime-2}) we get 
\[
x'_1-x'_2 = (z'-z) (\varphi_1-\varphi_2) + O(\varrho^2)
\]
holds. For $(z, T') \in \mc P'$ let $x(z, T')$ be the $x$-coordinate of the intersection point $T' \cap E_z({\bf q})$.
Now let $(z, T') \in \mc N_m(z, T)$ be arbitrary and let $(z_1,z_2, \varphi_1, \varphi_2) \in G(z, T')$. Then it follows from the above that there exists some $(z, T'') \in \mc N_{m+1}(z, T)$ such that
\[
x(z, T'') - x(z, T') = (z_2-z_1)(\varphi_2-\varphi_1) + O(\varrho^2).
\]
Let $X_{m}$ be a maximal $\varrho^2$-separated subset of points $x(z, T')$ over all $(z, T') \in \mc N_{m}(z, T)$. We obtain a graph 
\[
G_m \subset X_m \times Z({\bf q}) \times Z({\bf q}) \times \overline{\Phi}_{{\bf q}, \xi}\times \overline{\Phi}_{{\bf q}, \xi}
\]
where $\xi \in \overline{\Xi}_{\bf q}$ is such that $\theta(T) \in N_\varrho (\varphi, \xi+f(z) \varphi)$, so that for every $(x, z_1, z_2, \varphi_1, \varphi_2) \in G_m$ there exists $x' \in X_{m+1}$ such that
\[
x' - x = (z_2-z_1) (\varphi_2-\varphi_1)+ O(\varrho^2).
\]
Furthermore, we have the lower bound 
\[
|G_m| \ge \d^{O(\eta)}|X_m|  |Z({\bf q})|^2\varrho^{-2\gamma}.
\]
Recall that $Z({\bf q})$ is rescaled $(\varrho, s, \d^{-O(\eta)}\varrho^{-\epsilon})$-set and $\overline{\Phi}_{{\bf q}, \xi}$ is $(\varrho, \gamma, \d^{-O(\eta)})$-AD-regular. Using the assumption $s+\gamma\ge 1$ and Theorem \ref{thm: ABC}, we either have $|X_{m}| \ge \varrho^{-1+\zeta}$ or we get $|X_{m+1}| \ge \varrho^{-\chi} |X_{m}|$. Here $\zeta>0$ is from the statement of the lemma and $\chi$ is given by Theorem \ref{thm: ABC} and it is a function of $\alpha = 1-\zeta$, $\beta = s$, $\gamma=\gamma$. 

By the pigeonhole principle, we can find some $m = O_{\chi}(1)$ so that $|X_{m+1}| < \varrho^{-\chi} |X_{m}|$ holds. So for this choice of $m$ it must be the case that $|X_{m}| \ge \varrho^{\zeta-1}$. Let $\{x_1, \ldots, x_N\} \subset X_m$ be a maximal $\sigma$-separated subset (recall $\sigma \ge \varrho^2$). For each $x_i$ we have a corresponding pair $(z, T_i) \in \mc N_\ell(z, T)$ and a corresponding point $p_i \in E_z \cap T_i$. We know that $T_i \subset C m \cdot U(z,T)$. Furthermore, if we write $\theta(T_i) \in N_\d (\varphi_i, \xi_i+f(z) \varphi_i)$ for some $\xi_i \in \Xi_{p_i}$ and $\varphi_i \in \Phi_{p_i, \xi_i}$, then using AD-regularity of these sets, we get $\sim_{\d^{O(\eta)}} (\sigma/\varrho)^{-\gamma}$-many $\sigma$-tubes $\tilde T_i \in \TT_{\sigma}$ which intersect the $z$-plane within $\sigma$-ball around $T_i \cap E_z$ and satisfy $\theta(\tilde T_i) = (\tilde\varphi_i, \xi_i + f(z) \tilde \varphi_i+O(\sigma) )$ with some $\tilde \varphi_i = \varphi_i+ O(\varrho)$. Indeed, we choose some $\varphi_i' \in \Phi_{p_i, \xi_i}$ within distance $\varrho$ from $\varphi_i$, then we choose a tube $T'_i \in \TT_{p_i}$ in direction $(\varphi_i', \xi_i+f(z)\varphi_i')$ and let $\tilde T_i = (T'_i)^{\sigma}$. For large enough constant $K = O_{\chi, \eta}(1)$ we get
\[
|\TT_{\sigma}[ K\cdot U(z,T) ]| \gtrsim_{ \d^{O(\eta)}}  N (\varrho/\sigma)^{\gamma} \gtrsim_{\d^{O(\eta)}} \varrho^{\zeta} (\varrho/\sigma) (\varrho/\sigma)^{\gamma}.
\]
This finishes the proof after replacing $\zeta$ with $\zeta/2$ and taking $\eta \ll \mu \zeta$.
\end{proof}

\section{Non-linear slope function}\label{sec:unstructured-slope}

The goal of this section is to deal with the most difficult case when the slope function $f$ is approximately Lipschitz and cannot be approximated by a linear function on any $\varrho$-interval. 

\subsection{\texorpdfstring{$L$}{L}-tuples.}\label{subsec:Lambda-tuple}
This is the second variation on the $V$-tuples from Section
\ref{subsec:W-tuples}. These will be used to analyze the non-linear case in Section \ref{subsec:non-linear-grains-furstenberg} below. 

Fix a global grains configuration $(\TT,E,\{\TT_p\})$ and an
admissible cube ${\bf q}=p_0^\varrho$.  Let
$E({\bf q})$, $Z({\bf q})$, and $\TT({\bf q})$ be given by Lemma
\ref{lem:pigeonhole-rho-ball-2}.

\medskip\noindent\emph{One-arm paths.}
For $z\in Z({\bf q})$ and $p,\widetilde p\in E_z({\bf q})$, write
$p\leftrightarrow\widetilde p$ if
\[
p=(x,y+f(z)x,z),\qquad
\widetilde p=(\widetilde x,y+f(z)\widetilde x,z)
\]
for some $y\in Y_z({\bf q})$ and
$x,\widetilde x\in X_{y,z}({\bf q})$.  Thus the two points lie in
the same global grain.

For $z,z',z''\in Z({\bf q})$, let 
\begin{align*}
    J_{z,z',z''} = \{  (p,p',\widetilde p',p'',T,T')
\in E_z({\bf q})\times E_{z'}({\bf q})^2
\times E_{z''}({\bf q})\times\TT({\bf q})^2 \\
T\in\TT_p\cap\TT_{p'},\quad 
T'\in\TT_{\widetilde p'}\cap\TT_{p''}\\
p'\leftrightarrow\widetilde p' \}
\end{align*}

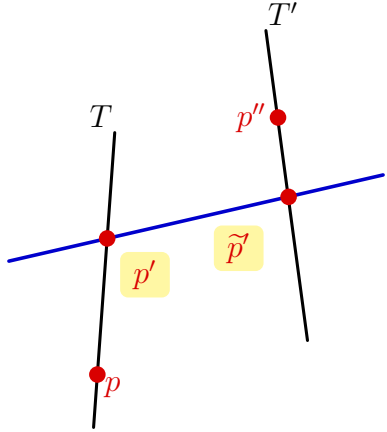
\begin{figure}[ht]
    \centering
\begin{tikzpicture}[
    line cap=round,
    line join=round,
    point/.style={
        circle,
        draw=none,
        fill=red!85!black,
        inner sep=2.2pt,
        outer sep=0pt
    },
    plainlabel/.style={
        rectangle,
        draw=none,
        fill=none,
        inner sep=0pt,
        outer sep=0pt,
        font=\large
    },
    redlabel/.style={
        plainlabel,
        text=red!85!black
    },
    highlight/.style={
        rectangle,
        draw=none,
        fill=yellow!45,
        rounded corners=3pt,
        inner xsep=5pt,
        inner ysep=3pt,
        outer sep=0pt,
        text=red!85!black,
        font=\large
    }
]

\coordinate (pprime)       at (0,0);
\coordinate (ptildeprime)  at (2.40,0.55);
\coordinate (p)            at (-0.13,-1.80);
\coordinate (pdoubleprime) at (2.26,1.60);

\draw[black, line width=1.15pt]
    (-0.18,-2.50) -- (0.10,1.40);

\draw[black, line width=1.15pt]
    (2.65,-1.35) -- (2.10,2.75);

\draw[blue!80!black, line width=1.35pt]
    (-1.30,-0.30) -- (3.65,0.84);

\node[point] at (p) {};
\node[point] at (pprime) {};
\node[point] at (ptildeprime) {};
\node[point] at (pdoubleprime) {};

\node[plainlabel, anchor=south east]
    at ($(0.10,1.40)+(-0.03,0.08)$) {$T$};

\node[plainlabel, anchor=south west]
    at ($(2.10,2.75)+(0.03,0.07)$) {$T'$};

\node[redlabel, anchor=north west]
    at ($(p)+(0.10,-0.04)$) {$p$};

\node[redlabel, anchor=east]
    at ($(pdoubleprime)+(-0.16,0)$) {$p''$};

\node[highlight]
    at ($(pprime)+(0.50,-0.48)$) {$p'$};

\node[highlight]
    at ($(ptildeprime)+(-0.66,-0.68)$) {$\widetilde{p}'$};

\end{tikzpicture}

    \caption{An one-arm path, blue line indicates the horizontal $(k-1)$-grain.}
    \label{fig:one-arm path}
\end{figure}

See Figure \ref{fig:one-arm path}. 
Let $\widetilde Z\subset Z({\bf q})$ and
$E_{\widetilde Z}({\bf q})
=\bigcup_{z\in\widetilde Z}E_z({\bf q})$.  Then
\begin{equation}\label{eq:J-tuples-counting-Furstenberg}
\sum_{z,z',z''\in \widetilde Z}|J_{z,z',z''}|
\sim_{\d^{O(\eta)}}
|E_{\widetilde Z}({\bf q})|
\d^{-2\kappa}|\widetilde Z|^2(\varrho/\d).
\end{equation}
Indeed, we choose elements of a $J$-tuple successively in the order 
$
p, T, (z',p'),\widetilde p',
 T',(z'',p'').
$
The two tube choices contribute $\d^{-2\kappa+O(\eta)}$,
the two level choices contribute
$\d^{O(\eta)}|\widetilde Z|^2$ by Lemma
\ref{lem:pigeonhole-rho-ball-2}(b), and the grain jump contributes
$\varrho/\d$ choices by Lemma
\ref{lem:pigeonhole-rho-ball-2}(c).  

\medskip\noindent\emph{Definition of the $L$-tuples.}
Fix $\tau\in[\d,\varrho]$.  For $p\in E({\bf q})$,
$\xi\in\Xi_p$, and
$\widetilde x\in\tau\ZZ\cap[-1,1]$, let
$J_{z,z',z''}(p,\xi,\widetilde x)$ be the set of paths
$ (p, p', \tilde p', p'', T, T')$ in $J_{z,z',z''}$ where $\theta(T) \in N_\d (\varphi, \xi +f(z) \varphi)$ for some $\varphi \in \Phi_{p, \xi}$ and the $x$-coordinate of $\tilde p'$ lies within $\tau$-neighborhood of $\tilde x$. 

Define
\[
L_{z,z',z''}(p,\xi,\widetilde x)
\subset J_{z,z',z''}(p,\xi,\widetilde x)^2
\]
to be the set of ordered pairs of $J$-tuples $$((p, p'_1, \tilde p'_1, p''_1, T_1, T'_1), (p, p'_2, \tilde p'_2, p''_2, T_2, T'_2))$$ satisfying
$|\theta(T_1')-\theta(T_2')|\leq\tau$. Put
\[
L_{z,z',z''}
=\bigcup_{p,\xi,\widetilde x}
L_{z,z',z''}(p,\xi,\widetilde x).
\]

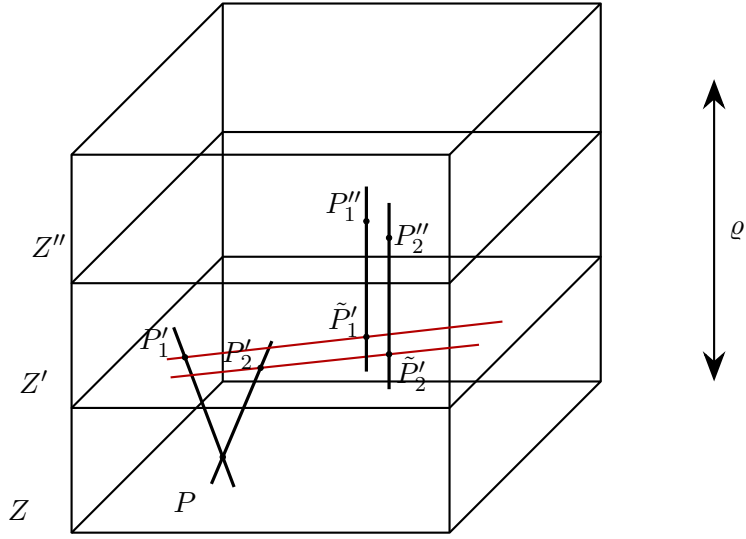
\begin{figure}[ht]
\begin{tikzpicture}[
    scale=1.0,
    line join=round,
    cube/.style={draw=black, thick},
    mainline/.style={draw=black, very thick},
    redline/.style={draw=red!70!black, thick},
    every node/.style={inner sep=0pt, outer sep=2pt, fill=none, draw=none}
]


\coordinate (A) at (0, 0);
\coordinate (B) at (5, 0);
\coordinate (C) at (7, 2);
\coordinate (D) at (2, 2);

\coordinate (E) at (0, 5);
\coordinate (F) at (5, 5);
\coordinate (G) at (7, 7);
\coordinate (H) at (2, 7);


\draw[cube] (A) -- (B) -- (C) -- (D) -- cycle;
\draw[cube] (E) -- (F) -- (G) -- (H) -- cycle;
\draw[cube] (A) -- (E);
\draw[cube] (B) -- (F);
\draw[cube] (C) -- (G);
\draw[cube] (D) -- (H);


\coordinate (Z1-A) at ($(A)!0.33!(E)$);
\coordinate (Z1-B) at ($(B)!0.33!(F)$);
\coordinate (Z1-C) at ($(C)!0.33!(G)$);
\coordinate (Z1-D) at ($(D)!0.33!(H)$);
\draw[cube] (Z1-A) -- (Z1-B) -- (Z1-C) -- (Z1-D) -- cycle;

\coordinate (Z2-A) at ($(A)!0.66!(E)$);
\coordinate (Z2-B) at ($(B)!0.66!(F)$);
\coordinate (Z2-C) at ($(C)!0.66!(G)$);
\coordinate (Z2-D) at ($(D)!0.66!(H)$);
\draw[cube] (Z2-A) -- (Z2-B) -- (Z2-C) -- (Z2-D) -- cycle;


\coordinate (P) at (2.0, 1);
\coordinate (P1p) at (1.5, 2.32);
\coordinate (P2p) at (2.5, 2.18);
\coordinate (P1t) at (3.9, 2.59);
\coordinate (P2t) at (4.2, 2.36);
\coordinate (P1pp) at (3.9, 4.12);
\coordinate (P2pp) at (4.2, 3.9);


\draw[mainline] ($(P)!-0.3!(P1p)$) -- ($(P)!1.3!(P1p)$);
\draw[mainline] ($(P)!-0.3!(P2p)$) -- ($(P)!1.3!(P2p)$);
\draw[redline] ($(P1p)!-0.1!(P1t)$) -- ($(P1p)!1.75!(P1t)$);
\draw[redline] ($(P2p)!-0.7!(P2t)$) -- ($(P2p)!1.7!(P2t)$);
\draw[mainline] ($(P1t)!-0.3!(P1pp)$) -- ($(P1t)!1.3!(P1pp)$);
\draw[mainline] ($(P2t)!-0.3!(P2pp)$) -- ($(P2t)!1.3!(P2pp)$);


\filldraw[black] (P) circle (1pt);
\filldraw[black] (P1p) circle (1pt);
\filldraw[black] (P2p) circle (1pt);
\filldraw[black] (P1t) circle (1pt);
\filldraw[black] (P2t) circle (1pt);
\filldraw[black] (P1pp) circle (1pt);
\filldraw[black] (P2pp) circle (1pt);


\draw (-0.7, 0.3) node {$Z$};
\draw (-0.5, 2.0) node {$Z'$};
\draw (-0.3, 3.8) node {$Z''$};

\draw (1.5, 0.4) node {$P$};
\draw (1.1, 2.52) node {$P'_1$};
\draw (2.2, 2.36) node {$P'_2$};
\draw (3.6, 2.82) node {$\tilde P'_1$};
\draw (4.5, 2.10) node {$\tilde P'_2$};
\draw (3.6, 4.32) node {$P''_1$};
\draw (4.5, 3.9) node {$P''_2$};




\draw[{Stealth[length=4mm]}-{Stealth[length=4mm]}, thick] (8.5, 2) -- (8.5, 6);
\draw (8.8, 4) node {$\varrho$};

\end{tikzpicture}

\caption{$L$-tuple $((p, p'_1, \tilde p'_1, p''_1, T_1, T'_1), (p, p'_2, \tilde p'_2, p''_2, T_2, T'_2))$.}
\label{fig:L-tuple}
\end{figure}

\medskip\noindent\emph{Counting the $L$-tuples.}
The following lemma gives an estimate on the total number of $L$-tuples for an appropriate choice of $\tau$.

\begin{lemma}\label{lem:number-of-Lambda-tuple}
Let $\widetilde Z \subset Z$ be a subset and 
suppose that
$|f(z)-f(z')|\leq c\tau/\varrho$ for all
$z,z'\in\widetilde Z$, where $c>0$ is a sufficiently small absolute
constant.  Then
\[
\sum_{z,z',z''\in \widetilde Z}|L_{z,z',z''}|
\sim_{\d^{O(\eta)}}
|E_{\widetilde Z}({\bf q})|
\d^{-3\kappa-\gamma}\tau^\kappa
(\varrho/\d)(\tau/\d)|\widetilde Z|^2.
\]
\end{lemma}

\begin{proof}
Note that if $((p, p'_1, \tilde p'_1, p''_1, T_1, T'_1), (p, p'_2, \tilde p'_2, p''_2, T_2, T'_2)) \in L_{z,z',z''}$ then we have $T_1 \underset{p}{\leftrightarrow} T_2$. So using (\ref{eq:y-axes-relation}) we get that $p_1'$ and $p_2'$ have their $y$-coordinates differ by $(z'-z) (f(z) - f(z')) (\varphi_1-\varphi_2) +O(\d) = O( \tau)$. Since $\tilde p'_i$ lies in the same global grain as $p_i'$ for $i=1,2$, $\tilde p_1'$ and $\tilde p_2'$ also have their $y$-coordinates differ by $O(\tau)$. So given a choice of $(p, p_1', p_2', T_1, T_2)$, by Cauchy--Schwarz there are typically $\d^{O(\eta)}(\varrho/\d) (\tau/\d)$ ways to select $\tilde p_1'$ and $\tilde p_2'$. Indeed, we first select $\tilde p_1'$ in $|X_{y,z}({\bf q})| \ge \d^{O(\eta)} (\varrho/\d)$ ways and then choose $\tilde p_2'$ on the global grain of $p_2'$ and within distance $\tau$ from $\tilde p_1'$ in $\d^{O(\eta)} (\tau/\d)$ ways (on average). 
    Since both $x$-coordinates of $\tilde p_1'$ and $\tilde p_2'$ are within $\tau$ from some fixed point $\tilde x$, we get $|\tilde p_1' - \tilde p_2'| = O(\tau)$.

    For a given $\tilde p \in E({\bf q})$ let $\TT({\bf q}, \tilde p) \subset \TT({\bf q})$ be the set of tubes $T$ such that $T \in \TT_{\tilde p'}\cap \TT({\bf q})$ for some $\tilde p' \in N_{C\tau}\tilde p$ (where $C$ is a fixed absolute constant). By a double counting argument similar to Lemma \ref{lem:pigeonhole-rho-ball-2}, we have
    \[
    \frac{1}{|E_{\tilde Z}({\bf q})|}\sum_{\tilde p \in E_{\tilde Z}({\bf q})} |\TT({\bf q}, \tilde p)|_\tau \sim_{\d^{O(\eta)}} \tau^{-\kappa}.
    \]
    It follows that when we average over all $(z, z', z'') \in \widetilde Z^3$, all $p, \xi, \tilde x$ and all pairs in  $J_{z,z', z''}(p, \xi, \tilde x)\times J_{z,z', z''}(p, \xi, \tilde x)$,  the intersection $ N_\tau \TT_{\tilde p_1'} \cap \TT_{\tilde p_2'}$ has $\tau$-covering number $\d^{O(\eta)}\tau^{-\kappa}$.
    So, on average, there are $\sim_{\d^{O(\eta)}} \tau^\kappa\d^{-2\kappa}$ ways to choose a pair of tubes $T_1' \in \TT_{\tilde p'_1} \cap \TT({\bf q})$ and $T_2' \in \TT_{\tilde p'_2} \cap \TT({\bf q})$ with $|\theta(T'_1) - \theta(T'_2)| \le \tau$. Once these tubes are selected, we can choose the points $p_1'', p_2''$ on the $z''$-slice in an essentially unique way. 
    Altogether, this gives us the claimed estimate on the sum of sizes of $L_{z,z',z''}$. To make the `on average' type of claims rigorous, we need to use the Cauchy--Schwarz inequality to lower bound the number of pairs of $J$-tuples similarly to how it was done in Lemma \ref{lem:counting-W-tuples}. 
\end{proof}

\medskip\noindent\emph{Additive relations.}
Fix an $L$-tuple and write, for $i=1,2$,
\[
\begin{aligned}
p_i'&=(x_i',y_i'+f(z')x_i',z'),\\
\widetilde p_i'
&=(\widetilde x_i',\widetilde y_i'
  +f(z')\widetilde x_i',z'),\\
p_i''&=(x_i'',y_i''+f(z'')x_i'',z''),\\
\theta(T_i)&=(\varphi_i,\xi+f(z)\varphi_i)+O(\d).
\end{aligned}
\]
The first two edges form the $V$-tuple appearing in
\eqref{eq:y-axes-relation}, so
\begin{equation}\label{eq:L-first-y-Furstenberg}
y_1'-y_2'
=(z'-z)(f(z)-f(z'))(\varphi_1-\varphi_2)+O(\d).
\end{equation}
Since the grain jump preserves the $y$-parameter,
$\widetilde y_i'=y_i'$. We have
\[
p_i''=\widetilde p_i'
+(z''-z')(\theta(T_i'),1)+O(\d).
\]
Subtracting and using
$|\theta(T_1')-\theta(T_2')|\leq\tau$ gives
\[
\begin{aligned}
x_1''-x_2''
&=\widetilde x_1'-\widetilde x_2'+O(\tau\varrho+\d),\\
y_1''-y_2''
&=\widetilde y_1'-\widetilde y_2'
+(f(z')-f(z''))(\widetilde x_1'-\widetilde x_2')
+O(\tau\varrho+\d).
\end{aligned}
\]
Since the two $\widetilde x_i'$ lie in the same $\tau$-ball, this
and \eqref{eq:L-first-y-Furstenberg} yield
\begin{equation}\label{eq:y-axes-another-relation}
\begin{split}
y_1''-y_2''
&=(z'-z)(f(z)-f(z'))(\varphi_1-\varphi_2)\\
&\quad
+O\bigl(\tau|f(z')-f(z'')|+\tau\varrho+\d\bigr).
\end{split}
\end{equation}
In the applications, $|f(z')-f(z'')|\leq\varrho$ and
$\d\leq\tau\varrho$, so the error in
\eqref{eq:y-axes-another-relation} is $O(\tau\varrho)$.

\subsection{Kaufman and radial projections in \texorpdfstring{$\RR^2$}{R2}}

We will need a variant of the radial projection theorem proved by Orponen, Shmerkin and Wang \cite{orponen2024kaufman}. For points $p, p'\in \RR^2$ let us denote $T_{p, p'}(\varrho)$ the $\varrho\times1$-tube around the line $\overline{pp'}$.

\begin{theorem}\label{thm:radial-projections}
    Let $t \in (0,2]$ and $\varepsilon_1, \varepsilon_2, \zeta \in (0,1)$. Suppose that $\varepsilon_1 \le \zeta/8$ .Then there exists $ \delta_0 = \d_0(t, \varepsilon_1, \varepsilon_2, \zeta), \eta_0 = \eta_0(t, \varepsilon_1, \zeta)>0$ such that the following holds for all $\delta < \delta_0$, $\eta < \eta_0$.
    Let $P \subset [0,1]^2$ be a $\delta$-separated  $(\d, \d^{-\eta})$-uniform $(\delta, t, \delta^{-\eta})$-set such that $|P \cap T| \le \d^{\varepsilon_2} |P|$ for any $\d^{\varepsilon_1}\times 1$-tube $T$.

    Then there exists a subset $G \subset P\times P$ such that $|G| \ge (1-\delta^{\min(\eta, \varepsilon_2/2)}) |P|^2$ such that for every $(p,p') \in G$ we have $|T_{p,p'}(\varrho) \cap P| \le \delta^{-\zeta} \varrho^{\min(1, t)} |P|$ for all $\varrho\in [\delta, 1]$.

    In particular, for every $p\in P$ such that $|G|_p| \ge \nu |P|$, where $G|_p:=\{p'\in P: (p, p')\in G\}$,   the set of tubes $\{T_{p, p'}(\delta), ~(p, p')\in G\}$ is a $(\delta, \min(1, t), C\delta^{-\zeta} \nu^{-1})$-set.
\end{theorem}

This version is different from the discretized radial projection theorem stated in \cite[Corollary 2.22]{orponen2024kaufman} and does not appear to follow directly from it. We give a proof of Theorem \ref{thm:radial-projections} in Appendix \ref{appendix:radial}.

\begin{remark}
    Note that in the statement of Theorem \ref{thm:radial-projections} the Frostman constant $\d^{-\eta}$ does not depend on the two-ends parameter $\varepsilon_2$. This is similar to how the quantitative bounds in \cite{orponen2024kaufman} or \cite{wang2024restriction} do not depend on the `2-ends exponents'. This feature is used in an important way in the proof of Lemma \ref{lem:non-quadratic-case} (originally we also used it in the proof of Theorem \ref{thm:furstenberg-R3} but managed to simplify the argument later).
\end{remark}

We will also need Kaufman's projection theorem, we will use a slightly stronger version compared to the one stated in \cite[Theorem 8.2]{wang2026sticky}. For $\lambda \in [-1,1]$ we denote $\pi_\lambda(x,y) = x-\lambda y$ the linear projection with slope $\lambda$.

\begin{theorem}\label{thm:kaufman}
    Let $t\in (0,1]$ and let $\varepsilon>0$. The following holds for $\d\le \d_0$ and $\eta < \eta_0$. 
    
    Let $P \subset [-1,1]^2$ be a $\d$-separated $(\d, t, \d^{-\eta})$-set and let $\Lambda \subset [-1,1]$ be a $\d$-separated $(\d, t, \d^{-\eta})$-set. Then there exists $\Lambda' \subset \Lambda$, $|\Lambda'| \ge (1-\d^{\eta})|\Lambda|$, and $P_\lambda \subset P$ for each $\lambda \in \Lambda'$, $|P_\lambda| \ge (1-\d^{\eta}) |P|$
    such that for any $P_\lambda' \subset P_\lambda$ with $|P_\lambda'|\geq \delta^{\eta}|P_\lambda|$, $\pi_{\lambda}(P_\lambda')$ is a $(\d, t, \d^{-\varepsilon})$-set (if there are overlap, view projection as multiset).
\end{theorem}

\begin{proof}

    Suppose that the statement does not hold. Then we can find a subset $\tilde \Lambda \subset \Lambda$, $|\tilde \Lambda| \ge \d^{\eta} |\Lambda|$ such that for every $\lambda \in \tilde \Lambda$ and any subset $P_\lambda \subset P$ of size at least $(1-\d^\eta)|P|$, there exists $P_\lambda' \subset P_\lambda$, $|P_\lambda'|\geq \delta^{\eta}|P|$ so that $\pi_{\lambda}(P_\lambda')$ is not a $(\d, t, \d^{-\varepsilon})$-set. 

    By a pigeonholing argument, this implies that there exists a dyadic scale $\varrho \in [\d, 1]$ (depending on $\lambda$) and a collection of dyadic $\varrho$-intervals $\mc I_\lambda$ such that 
    \[
    |\pi^{-1}_\lambda (I) \cap P| \gtrsim \d^{-\varepsilon} \varrho^t |P|, \quad I \in \mc I_\lambda
    \]
    and $|\bigcup_{\lambda} \pi^{-1}_\lambda (I) \cap P| \gtrsim \frac{1}{\log 1/\d} \d^{\eta} |P| \ge \d^{2\eta}|P|$ (if $\d\le \d_0$). By passing to a logarithmicly dense subset in $\tilde \Lambda$ we may assume that $\varrho$ is the same for all $\lambda$. After this we still have the lower bound $|\tilde \Lambda| \ge \d^{2\eta} |\Lambda|$.
    
    Suppose $\varepsilon > \eta$ and $\d < \d_0$. For $I \in \mc I_\lambda$ define
    \[
    N_\lambda(I):=\#\{ (p, p') \in P \cap \pi_\lambda^{-1}(I): d(p, p') \ge C \varrho \}.
    \]
    Then we have
    \[
    N_\lambda(I)\gtrsim |P \cap \pi_\lambda^{-1}(I)|^2 \gtrsim \d^{-\varepsilon} \varrho^{t} |P| |P \cap \pi_\lambda^{-1}(I)|
    \]
    and
    \[
    N_\lambda := \sum_{I\in  \mc I_\lambda} N_\lambda(I) \gtrsim \d^{2\eta-\varepsilon} \varrho^{t} |P|^2.
    \]
    So we get $\sum_{\lambda \in \tilde \Lambda} N_\lambda \gtrsim \d^{-\varepsilon+4\eta} \varrho^{t} |P|^2 |\Lambda|$. On the other hand,
    \[
    \sum_{\lambda \in \tilde \Lambda} N_\lambda = \sum_{(p, p') \in P^2: d(p, p') \ge C\varrho} \#\{ \lambda~|~\exists I\in \mc I_\lambda: p, p' \in  \pi_\lambda^{-1}(I)\}.
    \]
    For given $p, p'$ with $d(p, p') \sim r$, the set of such $\lambda$ is contained in an interval of length $\sim \varrho/r$. So since $\Lambda$ is a $(\d, t, \d^{-\eta})$-set, we get
    \[
    \sum_{\lambda \in \tilde \Lambda} N_\lambda \lesssim \sum_{r} \#\{(p, p') \in P^2: d(p, p') \sim r\} \d^{-\eta} (\varrho/r)^t |\Lambda| \lesssim \log (1/\d) \d^{-2\eta} \varrho^{t} |P|^2 |\Lambda|
    \]
    Comparing the upper and lower bounds, we obtain $\varepsilon \le 7\eta$ for small enough $\d$. So by defining $\eta_0 = \varepsilon/10$ we get a contradiction and the theorem is proved.
\end{proof}

\subsection{Non-linearly rotating grains.}\label{subsec:non-linear-grains-furstenberg}

In this section we consider the case when $f$ is not approximated by a linear function and has maximum variation. The main tool in the following ABC sum-product-type statement.

\begin{theorem}\label{thm:ABC-with-F}
    Let $\alpha, \beta \in (0,1], \gamma\geq 0$ such that $\min(\beta+\gamma,1)> \alpha$ and let $\varepsilon_1, \varepsilon_2>0$ such that $\varepsilon_1 \ll_{\alpha, \beta, \gamma}1$. Then there exist $\eta_0 = \eta_0(\alpha, \beta, \gamma, \varepsilon_1)$, $\overline{\eta}_0 = \eta_0(\alpha, \beta, \gamma, \varepsilon_1, \varepsilon_2)$ such that the following holds for all $\eta <\eta_0$, $\overline{\eta} < \overline{\eta}_0$ and small enough $\d$. 
    Let $B \subset [-1,1]$ be a $\d$-separated $(\d, \beta, \d^{-\eta})$-set. Let $C\subset [1,2]$ be a nonempty $\d$-separated $(\d, \gamma, \d^{-\eta})$-set (when $\gamma =0$, $C$ could be a single point). Let $F: B\to [-1,1]$ be a function such that for any $r, u\in \RR$:
    \begin{equation}\label{eq:non-linear-F}
    \#\{b\in B:~ |F(b) - r b- u|\le \d^{\varepsilon_1} \} \le \d^{\varepsilon_2} |B|.    
    \end{equation}
    Let $A \subset [-1,1]$ be a $\d$-separated subset. Let $G \subset A\times B\times B \times C$ be a subset of density at least $\d^{\overline{\eta}}$ and suppose that 
    \begin{equation}\label{eq:closed-under-shifts}
    |\{  a + c( b - b') (F(b) - F(b')), ~ (a, b, b', c) \in G \}|_\d \le \d^{-\overline{\eta}}|A|.
    \end{equation}
    Then we have $|A| \ge \d^{-\alpha}$.

    In particular, if $\beta>\alpha$, then suppose $G\subset A\times B \times B$ has density at least $\delta^{\overline \eta}$ and 
    \[
    |\{a+(b-b')(F(b)-F(b')), (a, b, b')\in G\}|_\delta\leq \delta^{-\overline{\eta}} |A|.
    \]
    Then we have $|A|\geq \delta^{-\alpha}$. 
\end{theorem}

In other words, if a set $A$ is essentially closed under translations by elements $c (b-b') (F(b) - F(b'))$, then $A$ must have `dimension' at least $\min(\beta+\gamma, 1)$.  This is similar to the ABC sum-product theorem except for the expression $(b-b')(F(b)-F(b'))$ appearing in the statement. In the application, we will take $\eta= \overline{\eta}$, so the reader can restrict to this case. In general, having $\eta$ not depend on $\varepsilon_2$ provides some extra flexibility but it will not be needed here.

\begin{remark}
    The analogous statement over $\CC$ is false: take $F(b) = \overline{b}$ (complex conjugation) and take $A = [-1,1] \subset \CC$, $C = [1,2]$ and $B = [-1,1]+ i [-1,1] \subset \CC$. For this choice we have $\dim_\CC A=\dim_\CC C = \gamma=1/2$, $\dim_\CC B = \beta = 1$ and $A$ is essentially invariant under translation by elements $c (b-b') (F(b)-F(b')) = c |b-b'|^2 \in \RR$. 
\end{remark}

\begin{proof}
Consider the graph of the function $F$:
    \[
    \Gamma = \{ (b, F(b)), ~ b \in B \}.
    \]
    Since $B$ is a $(\d, \beta, \d^{-\eta})$-set, so is $\Gamma$.  We claim that for any $\d^{2\varepsilon_1} \times 1$ tube $T$ we have
    \[
    |T\cap \Gamma | \lesssim \d^{\min(\varepsilon_2, \beta \varepsilon_1 - \eta)} |\Gamma| \le \d^{\min(\varepsilon_2, \beta\varepsilon_1/2)} |\Gamma|,
    \]
    where the second inequality holds when we choose  $\eta \le \beta\varepsilon_1/2$. 
    Indeed, if the slope of the tube $T$ is less than $\d^{-\varepsilon_1}$ then we use (\ref{eq:non-linear-F}) and if the slope exceeds $\d^{-\varepsilon_1}$ then we use that $B$ is a $(\d, \beta, \d^{-\eta})$-set. 

    Thus, by Theorem \ref{thm:radial-projections}, for any $\zeta'> 8\varepsilon_1$ there exists $\eta_1 = \eta_1(\zeta', \varepsilon_1, \beta)$ such that for any $\eta \le \eta_1$ we have the following. 
    There exists a subset $H \subset B\times B$ with $|H| \ge (1-\d^{\min(\eta_1, \varepsilon_2/2)}) |B|^2$ such that for every $(b, b')\in H$ and $\varrho \in [\delta, 1]$, any $\varrho$-tube $T_{b, b'}(\rho)$ containing $(b, F(b))$ and $(b', F(b'))$ satisfies 
    \begin{equation}\label{eq: thin-tube}
        |T_{b, b'}(\varrho)\cap \Gamma|\leq \delta^{-\zeta'} \varrho^{ \beta}|\Gamma|.
    \end{equation}
    As a consequence, for each $b$ such that $H|b\neq \emptyset$, the set of $\d\times 1$-tubes containing  $(b, F(b))$ and $(b', F(b'))$ over all $(b, b') \in H$ is a $(\d, \beta, \d^{-\zeta' } \frac{|B|}{|H|_{b}|})$-set.
    Since $B$ is $(\d, \beta, \d^{-\eta})$-set, for any $\varrho \ge \d$ we have $\#\{(b, b') \in B\times B:~ |b-b'| \le \varrho\} \le \d^{-\eta} \varrho^\beta |B|^2$. Write $\overline{\varepsilon}_2 = \min(\eta_1, \varepsilon_2/2)$.
    So by choosing $\varrho =  \d^{2\eta_1/ \beta}$, we can remove from $H$ all pairs at distance less than $\varrho$ and still have the lower bound $|H| \ge (1 -2\d^{\overline{\varepsilon}_2})|B|^2$. 

    Let $B' \subset B$ be the set of $b$ such that \begin{equation}\label{eq: Hb}
    |H|_{b}| \ge (1-\d^{\overline{\varepsilon}_2/2}) |B| 
    \end{equation} (note that $\delta^{\overline{\varepsilon_2}/2}\gg 2\delta^{\overline{\varepsilon_2}}$ when $\delta>0$ is sufficiently small). We have
    \[
    (1-2\d^{\overline{\varepsilon}_2})|B|^2 \le |H| \le |B'| |B| + (1- \d^{\overline{\varepsilon}_2/2})  |B \setminus B'| |B|
    \]
    \[
     (\d^{\overline{\varepsilon}_2/2} - 2\d^{\overline{\varepsilon}_2}) |B|^2 \le \d^{\overline{\varepsilon}_2/2}|B'| |B|
    \]
    and so \begin{equation}\label{eq: B'}
    |B'| \ge (1-2\d^{\overline{\varepsilon}_2/2})|B|.
    \end{equation} Then the above $H$ obtained from Theorem~\ref{thm:radial-projections} says that for every $b \in B'$ the set of slopes 
    \[
    \Lambda_b =\left\{ \frac{F(b') - F(b)}{b'-b} :~ b' \in H|_{b} \right\}
    \]
    is a $(\d, \beta, 2\d^{-\zeta'})$-set. Let $\pi_\lambda(x, y) = x-\lambda y$ denote the linear projection at slope $\lambda$. Let $\Gamma^{t} = \{ (F(b'), b'), ~b' \in B \}$ be the transposed graph of $F$. Let $\zeta>0$ be a parameter depending on $\alpha, \beta, \gamma$ to be determined.  Suppose that $\zeta'$ is sufficiently small in terms of $\zeta$. 
 Decompose $\Lambda_b =\sqcup_{j} \Lambda_{b,j}$ where
 \[ 
 \Lambda_{b,j}:= \left\{\frac{F(b')-F(b)}{b'-b}: ~b'\in H|_b, ~|T_{b,b'}(\delta)\cap \Gamma|\sim 2^{-j}|\Gamma| \right\},
 \]
 then \eqref{eq: thin-tube} implies $2^{-j} \leq \delta^{-\zeta'}\delta^\beta$ and $\Lambda_{b,j}$ is a $(\delta, \beta, \delta^{-10\zeta'})$-set.  
    Since $\Gamma^{t}$ is also $(\d, \beta, \d^{-\eta})$-set, for every $b \in B'$  apply Kaufman's projection theorem (Theorem \ref{thm:kaufman}) with $(\Gamma^t, \Lambda_{b,j})$ in place of $(P, \Lambda)$ and $10\zeta'$ in place of $\eta$, we get a subset $\Lambda'_b \subset \Lambda_b$,  
    \begin{equation}\label{eq: Lambda'b}
    |\Lambda'_{b,j}| \ge (1-\d^{10\zeta'}) |\Lambda_{b,j}|,
    \end{equation} and for every $\lambda = \frac{F(b')-F(b)}{b'-b} \in \Lambda'_{b,j}$ a subset 
    $B_{b, b'} \subset B$, 
    \begin{equation}
    \label{eq: Bbb'}
    |B_{b, b'}| \ge (1-\d^{10\zeta'}) |B|
    \end{equation} such that for any $\tilde\Gamma^t_{b, b'}\subset \Gamma^{t}_{b, b'} := \{(F(b''), b''), ~b'' \in B_{b, b'}\}$  with $|\tilde\Gamma^t_{b, b'}|\geq \delta^{10\zeta'}|\Gamma^t_{b, b'}|$, $\pi_{\lambda}(\tilde \Gamma^{t}_{b, b'})$ is a $(\d, \beta, \d^{-\zeta})$-set. Define $\Lambda'_b= \sqcup_{j} \Lambda'_{b,j}$. 
    For $b \in B'$ let $B'_b$ be the set of $b'$ so that $\frac{F(b')-F(b)}{b'-b} \in \Lambda'_b$. Then \eqref{eq: Lambda'b}, \eqref{eq: Hb}, and $|b-b'|\geq \varrho=\delta^{2\eta_1/\beta}$ implies 
    \begin{equation}\label{eq: Bb}
        |B_b'|\geq (1-2\delta^{\overline{\varepsilon_2}/2})|B|.
    \end{equation}
    Using $|b'-b|\ge \varrho$ for $b' \in H|_b$, it follows that the set 
    \begin{equation}\label{eq: pause}
    \{  (b'-b) F(b'') + (F(b') - F(b))b''~|~  b'' \in B_{b, b'} \}
    \end{equation}
    is a $(\d, \beta, 2\varrho^{-1}\d^{-\zeta})$-set for any $b \in B'$ and $b'\in B'_b$. 

Now let's pause and compare \eqref{eq: pause} with \eqref{eq:closed-under-shifts}: if \eqref{eq: pause} were a set consisting of $(b-b')(F(b)-F(b'))$, then we could apply the Orponen-Shmerkin ABC sum-product theorem. But we are not quite there yet, instead we shall apply the proof ideas of ABC sum-product theorem using Balog-Szemer\'{e}di-Gowers and an arithmetic trick (identity \eqref{eq:identity-Fb} below) to relate  $(b'-b) F(b'') + (F(b') - F(b))b''$ to $(b-b')(F(b)-F(b'))$.

    Define 
    \[
    Q_0 = \{ (b_1, b_2, b'_1, b'_2): b_1 \in B', \quad b_2 \in B'_{b_1}, \quad ~ b_1',b_2' \in B_{b_1,b_2}\} \subset B^4.
    \] 
    It follows from \eqref{eq: B'}, \eqref{eq: Bb}, \eqref{eq: Bbb'} that we have 
    \begin{equation} |Q_0| \ge (1 - 10\d^{\overline{\varepsilon}_2/2} ) |B|^4
    \end{equation}  By construction of sets $B', B'_b$ and $B_{b, b'}$, the set
    \begin{equation}\label{eq:Q-frostman}
    \{  (b_1 - b_2) (F(b_1') - F(b_2')) +(b_1' - b_2')  (F(b_1) - F(b_2))~|~ (b_1, b_2, b_1', b_2') \in Q_0\}    
    \end{equation}
    contains a $(\d, \beta, 2\varrho^{-1}\d^{-\zeta})$-set  and  for any $Q \subset B^4$ with $|Q| \ge 20 \d^{\overline{\varepsilon}_2/2} |B|^4$ if we replace $Q_0$ with $Q\cap Q_0$ in (\ref{eq:Q-frostman}) then  it contains  a $(\d, \beta, 4\varrho^{-1}\d^{-\zeta})$-set.

    We now turn to analyzing (\ref{eq:closed-under-shifts}). We will assume that $\overline{\eta} \ll \overline{\varepsilon}_2$ (since $\overline{\varepsilon}_2=\min(\eta_1, \varepsilon_2/2)$ and $\eta_1=\eta_1(\zeta, \varepsilon_1, \beta)$, this can be done). For shorthand, denote~$S(b, b') = (b-b') (F(b) - F(b'))$. 
    First, by Cauchy--Schwarz, (\ref{eq:closed-under-shifts}) implies an additive energy lower bound
    \begin{equation}\label{eq: energylowerbd}
    \#\{ (a,\tilde a, c, b, b', \tilde b, \tilde b') \in A^2\times C\times B^4: | a+ c S(b, b') -  \tilde a- c S(\tilde b, \tilde b')| \le \d \} \gtrsim \d^{O(\overline{\eta})} |A| |C||B|^4.
    \end{equation}
    Let $\mc S_j \subset \d \ZZ$ be the set of numbers $s = \d j$ such that $\#\{ (b, b') \in B\times B: S(b, b') \in [s, s+\d] \} \sim 2^j$. By dyadic pigeonhole and \eqref{eq: energylowerbd}, we can find an index $j$ such that $2^j |\mc S_j| \gtrsim \d^{O(\overline{\eta})} |B|^2$ and 
    \[
    \#\{ (a,\tilde a, c, s, \tilde s) \in A^2\times C\times \mc S_j^2: | a+ c s - \tilde a- c \tilde s| \le \d \} \\\gtrsim \d^{O(\overline{\eta})} |A| |C| |\mc S_j|^2.
    \]
    So there is a nonempty subset $C' \subset C$ such that $|C'| \ge \d^{O(\overline{\eta})} |C|$  (when $C$ is a single point, we take $C'=C$) and for every $c \in C'$ we have
    \[
    \#\{ (a,\tilde a, s, \tilde s) \in A^2\times \mc S_j^2: | a+ c s - \tilde a- c \tilde s| \le \d \} \\\gtrsim \d^{O(\overline{\eta})} |A| |\mc S_j|^2.
    \]
    
    We will use a corollary of the asymmetric Balog--Szemer\'edi--Gowers theorem (adapted to our discretized setting):
    
    \begin{lemma}[Corollary 2.36 in \cite{tao2006additive}]\label{lem:asymmetric-bsg}
        For any $C \ge 1, \varepsilon>0$ there exists $K$ such that the following holds for all $\d, \nu>0$. 
        Let $X, Y \subset [-1,1]$ be $\d$-separated sets such that $|Y| \le |X|^C$ and suppose that $\#\{(x, x, y, y') \in X^2\times Y^2:~ |x - y-x'+y' | \le \d\} \ge \nu |X|^2 |Y|$. Then there exist subsets $X' \subset X$, $|X'| \gtrsim \d^\varepsilon \nu^K |X|$ and $Y' \subset Y$, $|Y'| \gtrsim \d^{\varepsilon} \nu^K |Y|$ such that
        \[
        |Y' + n \cdot X' - m\cdot X'|_\d \lesssim_{n,m} (\d^{\varepsilon}\nu^{K} )^{-n-m} |Y|_\d.
        \]
        Here $n\cdot X = X+\ldots+X$ denotes the iterated sumset. 
    \end{lemma}

    We will apply this result with $Y = A$ and $X = c \mc S_j$ for each $c \in C'$. 
    First, we need to show that $\log |Y| / \log |X|$ is upper bounded by some constant to satisfy the assumption of Lemma~\ref{lem:asymmetric-bsg}. We have the following crude lower bound.
    \begin{claim}\label{claim:S_j-non-concentrated}
        We have $|\mc S_j| \gtrsim \varrho\d^{\zeta} \d^{-\beta/4}$.
    \end{claim}

    \begin{proof}
        Let $\mc G$ be the set of pairs $(b, b') \in B\times B$ so that $S(b, b') \in N_\d \mc S_j$. By the choice of $j$, we have $|\mc G| \ge \d^{O(\overline{\eta})} |B|^2$. Let $\mc C \subset B^4$ be the set of 4-cycles in $\mc G$, i.e. tuples $(b_1, b_2, b'_1, b'_2) \in B^4$ such that $(b_i, b_{i'}') \in \mc G$ for all $i, i' \in \{1, 2\}$. It follows from a standard application of Cauchy--Schwarz that $|\mc C| \gtrsim \d^{O(\overline{\eta})} |B|^4$. Now observe that we have the following identity for any $b_1, b_2, b_1', b_2' \in B$:
    \begin{equation}\label{eq:identity-Fb}
        \sum_{i, i'=1}^2 (-1)^{i+i'}S(b_i, b_{i'}') = -(b_1 - b_2) (F(b_1') - F(b_2')) - (b_1'-b_2') ( F(b_1)- F(b_2))
    \end{equation}
    We conclude that 
    \[
    \{ -(b_1 - b_2) (F(b_1') - F(b_2')) - (b_1'-b_2') ( F(b_1)- F(b_2) ~|~ (b_1, b_2, b_1', b_2') \in \mc C \} \subset N_{4\d} (\mc S_j + \mc S_j - \mc S_j - \mc S_j).
    \]
    If we take $\overline{\eta} \le c_0 \overline{\varepsilon_2}$ for a small enough constant $c_0$, then the set on the left hand side contains a $(\d, \beta, O(\varrho^{-1}\d^{-\zeta}))$-set by (\ref{eq:Q-frostman}) and the paragraph below it. In particular, it follows that the sum-set on the right hand side has covering number at least $\gtrsim \varrho \d^{\zeta} \d^{-\beta}$. The claim follows.
    \end{proof}
    
    By taking $\zeta < \beta/10$ and $\eta_1$ small enough, we get $|\mc S_j| > \d^{-\beta/8} \ge |A|^{\beta/8}$. Let $\varepsilon' >0$ and apply Lemma \ref{lem:asymmetric-bsg} with parameter $\varepsilon'>0$ and $X=c \mathcal{S}_j$ for each $c\in C'$,  we obtain subsets $A(c) \subset A$ and $\mc S_j(c) \subset \mc S_j$ so that for some $K = K(\varepsilon', \beta)$ we have
    \[
    |A(c)| \ge \d^{O(\overline{\eta} K+\varepsilon')} |A|,\quad |\mc S_j(c)| \ge \d^{O(\overline{\eta} K+\varepsilon')} |\mc S_j|
    \]
    \begin{equation}\label{eq:plunneke-ruzsa}
    |A(c) + c( n \cdot \mc S_j(c) - m \cdot \mc S_j(c))|_\d \lesssim \left(\d^{O(\overline{\eta} K+\varepsilon')}\right)^{-m-n} |A|.    
    \end{equation}
    Let $\mc G(c) \subset B\times B$ be the set of pairs $(b, b')$ so that $S(b, b') \in N_\d \mc S_j(c)$. We have $|\mc G(c)| \ge \d^{O(\overline{\eta} K+\varepsilon')} |B|^2$. Let $\mc C(c)$ be the set of 4-cycles in $\mc G(c)$, i.e. $(b_1, b_2, b_1', b_2') \in B^4$ so that $(b_i, b'_{i'}) \in \mc G(c)$ for all $i, i' \in \{1, 2\}$. By Cauchy--Schwarz we have $|\mc C(c)| \ge  \d^{O(\overline{\eta} K+\varepsilon')} |B|^4$. Taking $n=m=2$ in (\ref{eq:plunneke-ruzsa}) and using (\ref{eq:identity-Fb}) then gives that for 
    \begin{equation}\label{eq:B-of-c}
    \mc B(c) = \{ (b_1 - b_2) (F(b_1') - F(b_2')) + (b_1'-b_2') ( F(b_1)- F(b_2))~|~ (b_1, b_2, b_1', b_2') \in \mc C(c) \} 
    \end{equation}
    we have 
    \begin{equation}\label{eq:reduced-to-ABC}
        |A(c) - c \mc B(c)|_\d \lesssim \d^{O(\overline{\eta} K+\varepsilon')} |A|.   
    \end{equation}
    For any $\overline{\varepsilon}_2>0$ we can choose $\varepsilon'$ and $\overline{\eta}$ sufficiently small so that $\overline{\eta} K + \varepsilon' \ll\overline{\varepsilon_2}$ holds. So by replacing $\mc C(c)$ with $\mc C(c) \cap Q_0$ in (\ref{eq:B-of-c}) and using (\ref{eq:Q-frostman}) we get a subset $\mc B'(c) \subset \mc B(c)$ which is a $(\d, \beta, \d^{-2\zeta})$-set. When $\gamma =0, \beta >\alpha$, this already gives $|A|\geq |\mathcal{B}(c)|_\delta \geq \delta^{-\alpha}$ by taking $\overline\eta, \varepsilon', \zeta>0$ sufficiently small. When $\gamma>0$ and $\min(\beta+\gamma ,1)>\alpha$, we know that $C'$ is a $(\d, \gamma, \d^{-O(\eta)})$-set. Let $\chi = \chi(\alpha, \beta, \gamma)>0$ be the constant from Orponen-Shmerkin ABC sum production theorem (Theorem \ref{thm: ABC}) and choose $\zeta, \eta \ll \chi$. 
    It then follows from Theorem \ref{thm: ABC} and (\ref{eq:reduced-to-ABC}) that $|A| \ge \d^{-\alpha}$, as desired. (Define $\mc B$ to be the union over all $\mc B(c)$ and for each $c$ define the corresponding dense subgraph $\mc G := A(c) \times \mc B(c) \subset A \times \mc B$.)
\end{proof}

\begin{lemma}\label{lem:non-linear-grains}
    Let $\zeta, \mu, \varepsilon_1, \varepsilon_2 >0$ be so that $\varepsilon_1$ is sufficiently small in $\zeta$. Then the following holds for all $\eta < \eta(\zeta, \mu, \varepsilon_1, \varepsilon_2)$ and $\eps < \eps(\zeta, \varepsilon_1)$. 

    Consider a global grains $(\d, \eta, s, t, \kappa, \gamma)$-Furstenberg configuration with $s+\gamma \ge 1$ and let $f: Z\to [-1,1]$ be its slope function. Let $\varrho \in [\d^{1/4}, \d^{\mu}]$ and suppose that $\tilde Z\subset Z\cap {\bf z}$, ${\bf z} = z_0^\varrho$, is a rescaled $(\varrho, s, \varrho^{-\epsilon})$-set for some $z_0\in Z$.
    Suppose that there is $\sigma \in [\varrho^2,\varrho]$ such that for all $z , z'\in \tilde Z$ we have $|f(z') -f(z)| \le \sigma$ and for any $r,u \in \RR$ and $\tau \in [\varrho,1]$ we have
   \begin{equation}\label{eq:f-not-approximated}
   \#\{ z\in \tilde Z:~ |f(z)  -r z -u| \le \varrho^{\varepsilon_1} \sigma \} \le \varrho^{\varepsilon_2} |\tilde Z|.    
   \end{equation}
    Then we have $\kappa \le 2t-2+\zeta$. 
\end{lemma}

\begin{proof}
Replace $E $ by the set $E'$ given by Lemma \ref{lem:pigeonhole-rho-ball-2}. 
Fix a cube ${\bf q} = p_0^{\varrho}$ for some $p_0 \in E$ and let $Z({\bf q})$ and $E({\bf q})$ be the corresponding sets. It is not hard to verify that one can perform this pigeonholing step and guarantee that the set $\tilde Z$ from the statement of the lemma satisfies $|\tilde Z \cap Z({\bf q})| \gtrsim_{\d^{O(\eta)}} |\tilde Z|$. Thus, by Lemma \ref{lem:pigeonhole-rho-ball-2} (b) we can replace $Z({\bf q})$ by $\tilde Z \cap Z({\bf q})$ and keep the double counting properties. Then it follows that $Z({\bf q})$ is a rescaled $(\varrho, s, \varrho^{-\eps} \d^{-C\eta})$-set and (\ref{eq:f-not-approximated}) is satisfied with $Z({\bf q})$ in place of $\tilde Z$ and with $\d^{-C\eta}\varrho^{\varepsilon_2}$ in place of $\varrho^{\varepsilon_2}$.

For $z \in Z({\bf q})$ we thus can write 
    \[
    E_z({\bf q}) = \{ (x, y+ f(z)x, z),~ x \in X_{y,z}({\bf q}), ~ y \in Y_z({\bf q}) \}
    \]
    for rescaled $(\d/\varrho, 2t-1-\kappa, \d^{-O(\eta)})$-AD-regular set $Y_z({\bf q})$ and $(\d/\varrho, 1, \d^{-O(\eta)})$-AD-regular set $X_{y, z}({\bf q})$. 
    Let $\overline{\Xi}_{\bf q} \subset \varrho\ZZ$ be the set who $\varrho/2$-neighborhood covers $\Xi_{\bf q}$. For $\xi \in \overline{\Xi}_{\bf q}$ let $\overline{\Phi}_{{\bf q}, \xi} \subset \varrho \ZZ$ be the set whose $\varrho/2$-neighborhood covers $\bigcup_{\xi' \in N_{\varrho/2}\xi}\Phi_{{\bf q}, \xi'}$. By AD-regularity of $\Theta_{\bf q}$ we have $\sum_{\xi \in \overline{\Xi}_{\bf q}} |\overline{\Phi}_{{\bf q}, \xi}| \sim_{\d^{O(\eta)}} \varrho^{-\kappa}$. 
    Let $\overline{\Psi}$ be the set of triples $(\xi, \varphi_1, \varphi_2)$ such that $\varphi_1, \varphi_2 \in \overline{\Phi}_{{\bf q}, \xi}$.

    Fix an arbitrary $z_0 \in \tilde Z'$ and let
    $\mc Z = \{ \frac{z-z_0}{\varrho}, ~ z \in Z({\bf q})\} \subset [-1,1]$. Define the rescaled function $F: \mc{Z}\to [-1,1]$ by $F(z) = \frac{f(z_0 +\varrho z) - f(z_0)}{\sigma}$. By (\ref{eq:f-not-approximated}) for all $r, u \in \RR$ and $\tau \in [\varrho,1]$ we have
    \begin{equation}\label{eq:F-non-concentrated}
    \#\{z \in \mc{Z}:~ |F(z) - r z - u| \le \varrho^{\varepsilon_1}\} \lesssim \d^{-O(\eta)}\varrho^{\varepsilon_2} |\mc{Z}|.    
    \end{equation}
    We define a graph $G \subset Y_{z_0}({\bf q}) \times Z({\bf q}) \times Z({\bf q}) \times \overline{\Psi}$ as follows. Let $\tau = C \varrho \sigma$ for a large constant $C$ and consider the set of tuples $L_{z, z', z''}$ defined in Section \ref{subsec:Lambda-tuple}. Note that $|f(z)-f(z')| \le \sigma \le C^{-1}\tau/\varrho  $ for any $z,z' \in Z({\bf q})$ and so the assumption of Lemma \ref{lem:number-of-Lambda-tuple} is satisfied. 
    Consider a tuple 
    $((p, p'_1, \tilde p'_1, p''_1, T_1, T'_1), (p, p'_2, \tilde p'_2, p''_2, T_2, T'_2)) \in L_{z,z',z_0}$ and write $p_i'' = (x_i'', y''_i + f(z_0) x_i'', z_0)$ and $\theta(T_i) \in N_\varrho(\varphi_i, \xi+f(z)\varphi_i)$ for some $\xi \in \overline{\Xi}_{\bf q}$ and $\varphi_i \in \overline{\Phi}_{{\bf q}, \xi}$. Then we include the tuple $(y_1'', z,z', (\xi,\varphi_1,\varphi_2))$ in $G$. 
    
    By construction and (\ref{eq:y-axes-another-relation}), for any $(y, z,z', (\xi, \varphi_1, \varphi_2)) \in G$ there is $y' \in Y_{z_0}({\bf q})$ such that 
    \begin{equation}\label{eq:y-y-prime-relation}
    y'= y+ (z'-z) (f(z') - f(z))(\varphi_1-\varphi_2+O(\varrho)) + O(\varrho^2\sigma).
    \end{equation}
    By Lemma \ref{lem:number-of-Lambda-tuple} and averaging over $z_0 \in Z({\bf q})$, we can find some $z_0$ such that 
    \begin{equation}\label{eq:G-lower-bound}
    |G| \gtrsim_{\d^{O(\eta)}} |Y_{z_0}({\bf q})| |Z({\bf q})|^2 |\overline{\Psi}|.    
    \end{equation}
    Indeed, note that we have $|\overline{\Psi}| \sim_{\d^{O(\eta)}} \varrho^{-\kappa-\gamma}$ and $|Y_{z_0}({\bf q})| \sim_{ \d^{O(\eta)}} (\varrho/\d)^{2t-\kappa-1}$. If we fix a tuple in $G$ then we fix the choice of $(\xi, \varphi_1, \varphi_2)$ and $z, z' , z_0\in \tilde Z$ and the choice of $y''_1 \in Y_{z''}({\bf q})$. Note that this determines the directions of tubes $T_1, T_2$ up to $O(\varrho)$. 
    We can recover $x''_1$ in $O(\varrho/\d)$ ways, giving us $p_1''$.
    Then we recover $T_1' \in \TT_{p''_1}$ in $\d^{O(\eta)} \d^{-\kappa}$ ways. This gives us $\tilde p_1'$. Next we can select $p_1'$ in $O(\varrho/\d)$ ways. Then we can choose the triple $T_1 \underset{p}{\leftrightarrow} T_2$ in at most $\d^{O(\eta)} (\varrho/\d)^{\kappa+\gamma}$ ways (since their directions are already given up to precision $\varrho$). This then essentially uniquely recovers $p_2'$. We can then recover $\tilde p_2'$ in $O( \tau/\d )$ ways (since we have the constraint $|\tilde p_1' -\tilde p_2'| \le \tau$). Finally we can recover $T_2' \in \TT_{\tilde p_2'}$ in $\d^{O(\eta)} (\tau/\d)^{\kappa}$ ways. Multiplying these estimates together matches the bound in Lemma \ref{lem:number-of-Lambda-tuple} and so (\ref{eq:G-lower-bound}) follows.

    Now cover $Y_{z_0}({\bf q})$ by finitely overlapping intervals $I_j$ of length $C \varrho\sigma$ so that sets $Y_j =  Y_{z_0}({\bf q}) \cap I_j$ have size $\sim_{\d^{O(\eta)}} (\varrho\sigma/\d)^{2t-\kappa-1}$.
    Let $G_j$ be the graph $G$ restricted to $Y_j$. By pigeonhole, we can find an index $j$ so that $|G_j| \ge \d^{O(\eta)} |Y_j| |Z({\bf q})|^2 |\overline{\Psi}|$. By (\ref{eq:y-y-prime-relation}) we get that
    \[
    \{ y +  (z' - z) (f(z') - f(z)) (\varphi_1-\varphi_2) ~|~ (y, z, z', \xi, \varphi_1,\varphi_2) \in G_j\} \subset N_{C\varrho^2\sigma} Y_{z_0}({\bf q}) \cap C\cdot I_j
    \]
    By pigeonhole, we can find some $\xi, \varphi_1$ such that restricting $G_j$ on $(\xi, \varphi_1)$ produces a subgraph in $Y_j \times \mc Z\times \mc Z \times \overline{\Phi}_\xi$ of density $\d^{O(\eta)}$. 
    We are now in the position to apply Theorem \ref{thm:ABC-with-F} on scale $\varrho$: we let $\mc A$ be the rescaled set $Y_j$, $\mc B = \mc Z$ and $\mc C = \overline{\Phi}_\xi$. Then by (\ref{eq:F-non-concentrated}) we get that (\ref{eq:non-linear-F}) holds with $\varepsilon_2/2$ in place of $\varepsilon_2$ (if $\eta$ is small enough in $\varepsilon_2, \mu$). Let $\mc G \subset \mc A\times \mc B\times \mc B \times \mc C$ be the graph obtained by rescaling.
    So we get
    \[
    |\{ a + c (b'-b) (F(b') - F(b)), ~ (a, b, b', c) \in \mc G\}|_\varrho \le \d^{-O(\eta)} |\mc A|_\varrho.
    \]
     
    We have that $\mc B$ is $(\varrho, s, \varrho^{-\eps}\d^{-O(\eta)})$-set, $\mc C$ is $(\varrho, \gamma, \d^{-O(\eta)})$-set and $|\mc A|_\varrho \le \d^{-O(\eta)} \varrho^{-(2t-\kappa-1)}$. 
    Let $\alpha = 2t-\kappa-1 + \zeta/2$. If $\alpha > 1-\zeta/2$ then we get the desired bound on $\kappa$. Otherwise, for small enough $\eta$, we get $|\mc A|_\varrho \le \varrho^{-\alpha}$. Let $\beta = s$ and let $\eta_0 = \eta_0(\alpha, \beta, \gamma, \varepsilon_1)$, $\overline{\eta_0} = \overline{\eta_0}(\alpha, \beta, \gamma, \varepsilon_1, \varepsilon_2/2)$ be the functions given by Theorem \ref{thm:ABC-with-F}. By taking $ \eps$ to be sufficiently small as a function of $\zeta$ and $\varepsilon_1$, and by taking $\eta$ to be sufficiently small as a function of $\zeta$, $\varepsilon_1, \varepsilon_2$ and $\mu$ (and using that $\varrho \le \d^\mu$), we get that $\mc B$ is $(\varrho, \beta, \varrho^{-\eta_0})$-set, $\mc C$ is $(\varrho, \gamma, \varrho^{-\eta_0})$-set and we have $|\mc G| \ge\varrho^{\overline{\eta_0}} |\mc A||\mc B|^2 |\mc C|$ and 
    \[
    |\{ a + c (b'-b) (F(b') - F(b)), ~ (a, b, b', c) \in \mc G\}|_\varrho \le \varrho^{-\overline{\eta_0}} |\mc A|_\varrho,
    \]
    for all small enough $\d$. With this choice of parameters, we obtain a configuration which contradicts the statement of Theorem \ref{thm:ABC-with-F}. Thus, we must have $\alpha > 1-\zeta/2$ and $\kappa \le 2t-2 +\zeta$, which completes the proof.
\end{proof}

\section{Putting everything together: Proof of Proposition \ref{prop: coplanar}}\label{sec:proof-of-coplanar-case}

Suppose that for any $\d_0, \eta_0$ there are $\d<\d_0$ and $\eta< \eta_0$ for which there exists a $(\d, \eta, s, t, \kappa, \gamma)$-Furstenberg configuration. Here we assume that $\gamma \ge \max\{s+\kappa-1-\eps_0, 1-s+\eps_0\}$. Our goal is to show that $\kappa \le \max(2t-2, t-s)+4\eps_0$ holds, for the sake of contradiction suppose the contrary.

By Propositions \ref{prop:global-grains-1} and \ref{prop:global-grains-2} we can construct a new $(\d, \eta, s, t, \kappa, \gamma)$-Furstenberg configuration (perhaps for different but also arbitrarily small $\d,\eta$) for which the sets $E_z$ and $\Theta_p$ have the global grains structure given by some slope function $f: Z \to [-1,1]$. In other words, we obtain a global grains $(\d, \eta, s, t, \kappa, \gamma)$-Furstenberg configuration $(\TT, E, \{\TT_p\})$ where $\TT$ is a $(2t, \d^{-\eta})$-AD-regular set of tubes satisfying $C_{t-CW}(\TT) \le \d^{-\eta}$.

Now we split into three cases depending on the behavior of the function $f$ and reduce to either Lemma \ref{lem:slowly-rotating}, Lemma \ref{lem:linear-grains} or Lemma \ref{lem:non-linear-grains}. Recall that by Definition \ref{def:global-grains}, the pair $(Z, f)$ is $\d^\eta$-uniform. Fix an arbitrary $z_0 \in Z$. 
Let $\eps>0$ and let $\varrho \in [\d^{1/10}, \d^{\mu}]$ be a scale such that $Z \cap {\bf z}$, ${\bf z} =z_0^\varrho$, is a rescaled $(\varrho, s, \varrho^{-\eps})$-set (such $\varrho$ exists by Lemma \ref{lem:good-scales}).  
Apply Lemma \ref{lem:pigeonhole-rho-ball-2} and replace $E$ by the set $E'$ and $Z$ by $Z' = \{z: E_z'\neq\emptyset\}$. For every $p_0 \in E_{z_0}$ and ${\bf q} = p_0^{\varrho}$ we then have well-defined sets $Z({\bf q}) \subset Z\cap {\bf z}$, ${\bf z} = z_0^{\varrho}$ and $E({\bf q})$.

Let $0<\zeta < \eps_0/2$ and suppose that $\kappa\ge t-s+2\eps_0$. 
Then by Lemma \ref{lemma:f-Lipschitz}, we can refine $Z({\bf q})$ to a subset of density $\d^{O(\eta/\zeta)}$ such that $|f(z) -f(z')|\le \varrho$ for every $z,z'\in Z({\bf q})$. After further pigeonholing we may assume the same uniformity conditions hold for the new $Z({\bf q})$.

Next, we show that $f$ cannot vary much slower than an $\approx 1$-Lischitz function. 

\begin{claim}\label{claim:slowly-rotating}
    For any $0<\varepsilon_1 \le \eps_0/2$ and $\eps >0$ there exists $\varepsilon_2>0$ such that the following holds. Suppose that for every choice of ${\bf q}$ as above and $z_0 \in Z({\bf q})$,
    \[
    \#\{ z \in Z({\bf q}): |f(z)-f(z_0)| \le \varrho^{1+\varepsilon_1} \} \ge \varrho^{\varepsilon_2} |Z({\bf q})|
    \]
    Then $\kappa \le t-s+2\eps_0$. 
\end{claim}

\begin{proof}
    Using $\d^{-\eta}$-uniformity of $(Z({\bf q}),f)$, we can find a subset $\tilde Z \subset Z$ with $|\tilde Z| \ge \d^{O(\eta)} \varrho^\alpha |Z({\bf q})|$ such that $|f(z)-f(z')| \le \varrho^{1+\varepsilon_1}$ whenever $z,z'\in \tilde Z\cap {\bf z}$ for a dyadic $\varrho$-interval ${\bf z}$. Write $\sigma=\varrho^{1+\varepsilon_1}$. After further pigeonholing, we may assume that $(\tilde Z, f)$ is $\d^{-\eta}$-uniform and $\tilde Z \cap {\bf  z}$ is a rescaled $(\varrho, s, \varrho^{-\eps-\varepsilon_2} \d^{-O(\eta)})$-set for $z \in \tilde Z$. 
    By Lemma \ref{lem:slowly-rotating} applied to $\tilde Z$, we conclude that for any $\zeta'>0$ and $\eps, \varepsilon_2$ sufficiently small in $\zeta'$ and $\eta$ sufficiently small in $\zeta', \mu$ we have
    \[
    |\TT[U]|_{\sigma} > \varrho^{\zeta'} (\sigma/\varrho)^{-1-\gamma}
    \]
    for some $\sigma\times \varrho\times 1$ box $U$. On the other hand, by the $t$-Convex--Wolff axiom, we have
    \[
    |\TT[U]| \le \d^{-\eta} \sigma^{2t}(\varrho/\sigma)^{t}  |\TT| 
    \]
    and so, passing $\sigma$-tubes,
    \[
    |\TT[U]|_{\sigma} \lesssim_{\d^{-O(\eta)}} \sigma^{2t}(\varrho/\sigma)^{t} |\TT|_{\sigma} \sim_{\d^{-O(\eta)}} ( \varrho/\sigma)^{ t},
    \]
    since $\TT$ is $2t$-AD-regular. We conclude that 
    \[
    1+\gamma \le t + \zeta'/\varepsilon_1 + O(\eta/\mu\zeta) \le t +\varepsilon_1
    \]
    where we take $\zeta' = \varepsilon_1^2/2$ and $\eta$ small enough in the other parameters.
    Recall that we assume $\gamma \ge \max(s+\kappa-1-\eps_0, 1-s+\eps_0)$ for some $\eps_0>0$. So we deduce the bound $t \ge \max(2-s +\eps_0- \zeta, \kappa+s-\eps_0-\zeta) $.
    The second term in the maximum gives $\kappa \le t-s +\zeta +\eps_0 \le t-s+2\eps_0$, as desired.
\end{proof}

Thus, we may now restrict to the situation when for all $z_0$:
\begin{equation}\label{eq:non-concentration-fz}
    \#\{ z \in Z({\bf q}): |f(z)-f(z_0)| \le \varrho^{1+\varepsilon_1} \} \le \varrho^{\varepsilon_2} |Z({\bf q})|
\end{equation}
for some $\varepsilon_2>0$ depending on $\varepsilon_1$ and $\varepsilon_1$ sufficiently small.  Take $\eta$ sufficiently small in $\varepsilon_1, \mu$ so that we have $\eta/\varepsilon_1 \ll \mu \varepsilon_2$. 

Let $\varepsilon'_1 = 4\varepsilon_1$ and $\varepsilon'_2= \varepsilon_2/4$ be parameters. Fix $z_0 \in \tilde Z$ and ${\bf z} = z_0^\varrho$.  Suppose that for at least half of choices $z_0 \in \tilde Z$, there exist $u, r \in \RR$ such that 
\begin{equation}\label{eq:linear-approximation-fz}
\#\{ z \in Z \cap {\bf z}:~ |f(z) - r z - u| \le \varrho^{\varepsilon'_1} \varrho   \} \ge \varrho^{\varepsilon'_2} | Z \cap {\bf z}|.    
\end{equation} 
We want to apply Lemma \ref{lem:linear-grains} to the set of $z$ as above. First, we need to get an estimate on $|r|$. Consider the case when $|r| \le \varrho^{2\varepsilon_1}$. 
By triangle inequality $|f(z) - rz - u| \le |f(z) - rz_0 -u| + |r|\varrho$ and so by (\ref{eq:non-concentration-fz}) we get that
\begin{align*}
&\#\{ z \in Z \cap {\bf z}:~ |f(z)- r z - u| \le \varrho^{\varepsilon'_1} \varrho  \} \\ 
\le&  \#\{ z \in Z \cap {\bf z}:~ |f(z) - rz_0-u| \le \varrho^{\varepsilon'_1} \varrho + |r|\varrho  \} \le \varrho^{\varepsilon_2/2}| Z \cap {\bf z}|
\end{align*}
where we use $|r|+\varrho^{\varepsilon'_1} < \varrho^{\varepsilon_1}$. This contradicts the choice $\varepsilon'_2 = \varepsilon_2/4$. So we must have $|r| \ge \varrho^{2\varepsilon_1}$.
So (\ref{eq:linear-approximation-fz}) and $\varepsilon'_1\ge 4\varepsilon_1$ implies that
\[
\#\{ z \in Z \cap {\bf z}:~ |f(z) - r z - u| \le \varrho^{\varepsilon'_1/2} |r| \varrho   \} \ge \varrho^{\varepsilon'_2} |\tilde Z \cap {\bf z}|.    
\]
Let $ Z'$ be the set of $z$ satisfying the condition above. Note that $ Z'$ is a rescaled $(\varrho, s, \d^{-O(\eta/\varepsilon_1)}\varrho^{-\eps -\varepsilon'_2})$-set. 
Thus by replacing $Z\cap {\bf z}$ with $Z' \cap {\bf z}$, 
by Lemma \ref{lem:linear-grains} we get $\kappa \le 2t-2 +\zeta$ provided that $\eta \ll_{\zeta, \mu, \varepsilon'_1} 1$ and $\eps, \varepsilon'_2 \ll_{\zeta, \varepsilon'_1} 1$.

Now consider the case when for at least half of the choices $z_0 \in Z$ and ${\bf z} = z_0^\varrho$ we have 
\[
\#\{ z \in Z \cap {\bf z}:~ |f(z) - r z - u| \le \varrho^{\varepsilon'_1} \varrho   \} \le \varrho^{\varepsilon'_2} |Z \cap {\bf z}|
\]
for all $r, u \in \RR$. After a small refinement of $ Z$ we may assume that this holds for all $z_0$. By applying Lemma \ref{lem:non-linear-grains} with $\sigma=\varrho$, $\tilde\eps = \eps+ O(\eta/\varepsilon_1\mu)$ and $\varepsilon'_1, \varepsilon'_2$ as above, we conclude that $\kappa \le 2t-2+\eps_0$ provided that $\eps$ is sufficiently small in $\eps_0, \varepsilon'_1$, and $\eta$ is sufficiently small in $\eps_0, \mu, \varepsilon'_1, \varepsilon'_2$. 

We select parameters as follows. First let $\varepsilon_1$ be small enough in $\eps_0$, then $\varepsilon_2$ is a function depending on $\varepsilon_1$. We then let $\varepsilon'_1 = 4\varepsilon_1$ and $\varepsilon'_2= \varepsilon_2/4$. Then let $\eps$ be sufficiently small and let $\mu$ be a function of $\eps$ given by Lemma \ref{lem:good-scales}. Finally let $\eta$ be sufficiently small in all other parameters. 
This completes the proof of Proposition \ref{prop: coplanar}.

\begin{remark}
    The only place where $t$-Convex Wolff axiom was used is in the proof of Claim \ref{claim:slowly-rotating}. One can verify that the proof of the Claim goes through if we assume that $C_{t;t'-CW}(\TT) \le \d^{-\eta}$ where
    \[
    t' =\begin{cases}
        2-s, \quad t\in (0,2-s],\\
        \min(2t-2+s, 2),\quad t \in [2-s,2].
    \end{cases}
    \]
    Thus, for Proposition \ref{prop: coplanar}, i.e. the case when $\gamma \ge \max(s+\kappa-1, 1-s)$ it suffices to assume this weaker $(t;t')$-Convex Wolff axiom. 
\end{remark}

\section{Fully non-coplanar case: proof of Proposition \ref{prop: very-non-coplanar}}\label{sec:very-non-coplanar}

In this section we consider the case when $\gamma \le \eps_0$ and $s+\kappa \le 1+\eps_0$. The proof in Section \ref{sec:intermediate} does not apply because the direction sets $\Phi_\xi$ might be too small to apply the ABC sum-product theorem. In this section we present a slightly different argument that fixes this issue. We will use Ren--Wang's  Furstenberg set estimate \cite{ren2023furstenberg} in the following form.

\begin{theorem}\label{thm:ren-wang}
    For any $s \in (0,1]$, $t \in (0,2]$, $\zeta>0$ the following holds for $\eta<\eta_0$ and $\d<\d_0$. 
    
    Let $P \subset [-1,1]^2$ be a $(\d, t, \d^{-\eta})$-set  and for every $p \in P$ let $\TT_p$ be a $(\d, s, \d^{-\eta})$-set of tubes through $p$. Suppose $|\TT_p|\sim M$ for all $p$, then 
    \[
    \left|\bigcup_{T\in \TT}\TT_p\right| \gtrsim \d^\zeta\d^{-\min(t, \frac{s+t}{2}, 1)}M.
    \]
\end{theorem}

Let $\eta_0>0$ be a small parameter to be determined and $\eta>0$ sufficiently small. 
Apply Lemma~\ref{lem:good-scales} to $Z$ and find a scale $\varrho\in [\delta^{1/2}, \delta^{\mu}]$ such that $Z\cap {\bf z}$ is a rescaled $(\varrho, s, \varrho^{-\eta_0})$-set for every $z\in Z$ and ${\bf z} = z^\varrho$.  By passing through a subset of $Z\cap {\bf z}$, we may assume that $|Z\cap {\bf z}|_{\varrho^2} \sim \varrho^{-s}$. 

Without loss of generality we may assume that $E$ is $\d$-separated. 
Let $p_0\in E$, consider ${\bf q} =p_0^\varrho$.  Apply Lemma \ref{lem:pigeonhole-rho-ball} and let $E({\bf q}) \subset E\cap {\bf q}$, $Z({\bf q}) \subset Z\cap {\bf z}$ and $\Theta_{\bf q}$ be the resulting sets. By a change of coordinates we may assume that $\alpha_{\bf q} = 0$. 
For $z \in Z({\bf q})$ let $A_z({\bf q}) \subset A_z$ be the set of $a \in A_z$ so that $(a,z)\in E({\bf q})$. Note that $\sum_{z \in Z({\bf q})} |A_z({\bf q})|_\d \sim |E({\bf q})|_\d \gtrsim \d^{O(\eta)}|Z({\bf q})|_\d (\varrho/\d)^{2t-\kappa}$. So after refinement, $A_z({\bf q})$ is $([\d, \varrho], 2t-\kappa, \d^{-O(\eta)})$-AD-regular set for every $z \in Z({\bf q})$. Let $\overline{A}_z \subset {(\varrho^2\ZZ)^2}$ be the set such that $\overline{A}_z^{\varrho^2} = A_z({\bf q})^{\varrho^2}$. 

Let $\overline{\Theta} \subset \varrho \ZZ^2 \cap [-1,1]^2$ be the set such that $\overline{\Theta}^{\varrho} = \Theta_{\bf q}^{\varrho}$. By the definition of a $(\d, \eta,s,t,\kappa, \gamma)$-Furstenberg configuration, the set $\overline{\Theta}$ is $(\varrho, \gamma, \d^{-\eta})$-tube Katz-Tao, i.e. $|\Theta_{\bf q} \cap T|_{\sigma} \le \d^{-\eta} (\tau/\sigma)^\gamma$ for any $\sigma \times \tau$ tube $T$ and any $\varrho\le \sigma\le \tau$. 

Consider the following set
\[
G_{z,z'} = \{ (a,\theta, \theta') \in \overline{A}_z \times \overline{\Theta}\times \overline{\Theta}: a + (z'-z) (\theta-\theta') \in N_{C\varrho^2} \overline{A}_z \}.
\]
Using Lemma \ref{lem:pigeonhole-rho-ball} and averaging over pairs of tubes $T_1,T_2 \in \TT({\bf q})$ and levels $z,z'$ we obtain
\[
\frac{1}{|Z({\bf q})|^2}\sum_{z,z' \in Z({\bf q})} |G_{z,z'}| \gtrsim_{\d^{O(\eta)}} \varrho^{-(2t-\kappa)} |\overline{\Theta}|^2.
\]
Fix $z_0 \in Z({\bf q})$ and $\theta_0'$ so that $\frac{1}{|Z({\bf q})|}\sum_{z} |G_{z_0,z}|_{\theta_0'}| \gtrsim_{\d^{O(\eta)}} \varrho^{-(2t-\kappa)} |\overline{\Theta}|$ where we restrict the last element of the tuple to be $\theta_0'$.
Now for each $a \in \overline{A}_{z_0}$ define 
\[
\overline{\Theta}_a = \{ \theta \in \overline{\Theta}: \#\{ z: (a, \theta,\theta'_0) \in G_{z_0,z}\} \ge \d^{C_0 \eta} |Z({\bf q})|\}.
\]
For an appropriate choice of $C_0$ we then obtain 
\[
\sum_{a \in \overline{A}_{z_0}} |\overline{\Theta}_a| \gtrsim \d^{C_0 \eta}  |\overline{A}_{z_0}| |\overline{\Theta}|.
\]
Now we can define a set of tubes $\mc T_a$ by including a $\varrho^2 \times \varrho$ tube through $a$ and in direction $\theta$ if $\theta \in \overline{\Theta}_a$ and $|\theta - \theta_0'| > \d^{C_1 \eta/\kappa}$ for some large fixed $C_1$. Since $\overline{\Theta}$ is $([\varrho^2,\varrho],\kappa, \d^{-O(\eta)})$-AD-regular, the latter condition does not significantly affect the total number of included tubes. We conclude that on average $\mc T_a$ has cardinality $\d^{O(\eta)} |\overline{\Theta}|$ (when counted with multiplicity). Now we observe that $\mc T_a$ is a $(\varrho, \kappa-\gamma, \d^{-O(\eta)})$-set for every $a$ with $|\mc T_a| \gtrsim_{\d^{O(\eta)}} |\overline{\Theta}|$. Indeed, let $\tau \in [\varrho,1]$ and choose arbitrary $\tau \varrho \times \varrho$ tube $T_\tau$ through $a$. Then each $T \in \mc T_a[T_\tau]$ gives $\theta \in \overline{\Theta}$ such that the direction vector $\operatorname{Dir}(\theta-\theta'_0)$ lies in a $O(\tau)$-neighborhood of $\operatorname{Dir}(T_\tau)$.
Thus, for some $\sim \tau \times 1$ tube $T'$ we get $|\overline{\Theta} \cap T'| \ge |\mc T_a[T_\tau]|$. By the tube Katz-Tao property 
\[
|\overline{\Theta} \cap T'| \lesssim_{\d^{O(\eta)}} (\tau/\varrho)^{\kappa} |\overline{\Theta} \cap T'|_\tau \lesssim_{\d^{O(\eta)}}(\tau/\varrho)^{\kappa} \tau^{-\gamma} \sim_{\d^{O(\eta)}} \tau^{\kappa-\gamma} |\overline{\Theta}|,
\]
as desired. By the Furstenberg set estimate applied to a rescaled collection of tubes, we conclude that
\[
\left |\bigcup_{a \in \overline{A}_{z_0}} \mc T_a\right|_{\varrho} \gtrsim_{\d^{O(\eta)}} \varrho^{\zeta} \varrho^{-\min( 2t-\kappa, \frac{(\kappa-\gamma)+(2t-\kappa)}{2},1)} |\overline{\Theta}|.
\]
Here $|\cdot |_{\varrho}$ means that we are counting $\varrho^2\times \varrho$-tubes up to an essential equivalence. 
Moreover this estimate holds after an arbitrary $\d^{C\eta}$-refinement of the tube sets $\mc T_a$.
On the other hand, we can estimate this union as follows. Fix $\theta \in \overline{\Theta}$. For $T \in\mc T_a$ in direction $\theta$ the number of $z$ such that $a+(z-z_0) (\theta-\theta_0') \in N_{C\varrho^2}\overline{A}_{z_0}$ is at least $\d^{O(\eta)}|Z({\bf q})|$. Thus, $|T \cap \overline{A}_{z_0}| \gtrsim_{\d^{O(\eta)}} |Z({\bf q})|_{\varrho^2} \gtrsim_{\d^{O(\eta)}} \varrho^{\eps - s}$. Thus, the number of essentially distinct tubes in direction $\theta$ is upper bounded by $\d^{-C\eta}\varrho^{s-\eps}|\overline{A}_{z_0}|$. We conclude
\[
\left |\bigcup_{a \in \overline{A}_{z_0}} \mc T_a\right|_{\varrho}  \lesssim_{\d^{O(\eta)}} \varrho^{s-\eps}|\overline{A}_{z_0}| |\overline{\Theta}| \lesssim \varrho^{s-\eps -(2t-\kappa)} |\overline{\Theta}|.
\]
Thus,
\begin{align*}
    \varrho^{\zeta} \varrho^{-\min( 2t-\kappa, \frac{(\kappa-\gamma)+(2t-\kappa)}{2},1)} &\lesssim_{\d^{O(\eta)}} \varrho^{s-\eps -( 2t-\kappa)} \\
    -\zeta + \min(2t-\kappa, \frac{2t-\gamma}{2}, 1) &\le -s +\eps + 2t-\kappa + O(\eta/\mu)
\end{align*}

Take $\zeta,\eps$ small enough in $\eps_0$ and then $\eta$ small enough in $\mu$ (recall $\mu=\mu(\eps)$), then the minimum can only be achieved on the second or third term. If it is the second term, then we get $t-\gamma/2 \le -s + 2t-\kappa + \eps_0/2$ and
\[
\kappa \le t-s + \eps_0/2 + \gamma/2 \le t-s+\eps_0,
\]
as desired. If minimum is achieved on third term then $1 \le -s + 2t-\kappa + \eps_0/2$ and
\[
\kappa \le 2t-1-s +\eps_0/2.
\]
Recall that we assume in this section that $\kappa\le 1-s+\eps_0$ and so taking average of these two inequalities gives $\kappa \le t-s+\eps_0$, as desired.
This completes the proof.

\begin{remark}
    In this argument, we could have used the radial projections theorem to prove that $\mc T_a$ are $\kappa$-dimensional as long as $\gamma < \kappa-\eps_0$. 
    Since this range is already covered by Proposition \ref{prop: trilinearABC}, we don't need to do this. Note that the argument in Section \ref{sec:intermediate} is still necessary as it covers the case $s+\kappa>1$ while the bound in this section can be lossy there.
\end{remark}

\section{Fully coplanar case: proof of Proposition \ref{prop: very-coplanar}}\label{sec:very-coplanar}

In this section we consider the coplanar case $\gamma \ge \kappa-\eps_1$ and $s+\kappa < 1+\eps_0$. Here we can choose $\eps_1$ arbitrarily small compared to $\eps_0$ and our goal is to prove $\kappa \le \kappa_0(s,t)+2\eps_0$.
Note that the coplanar reduction argument from Section \ref{sec:global-grains} requires $s+\gamma\ge 1$ in order to guarantee that the sets $X_{y,z}$ parameterizing global grains are 1-dimensional. 
First, we show that after rescaling, we can obtain a global grains structure where the grains are filled in by at least $s+\kappa$-dimensional sets. This will allow us to run a modified version of the arguments in Sections \ref{sec:structured-slope}--\ref{sec:proof-of-coplanar-case}. Similarly to those sections we will consider the non-linear, linear and slow variation cases.
The proofs follow quite closely the arguments from the above sections, so we omit various technical details and focus on the required modifications.

\begin{prop}\label{prop:porous-global-grain}
    Fix $s,t,\kappa$ such that $s+\kappa\le 1+\eps_0$ and $\kappa \ge \max(\frac{2}{3}(t-s),0)+2\eps_0$. 
    Suppose that for every $\eps_1>0$ and $\eta_0$, $\d_0$ there are $\d\le \d_0$, $\eta\le \eta_0$ such that there is a $(\d, \eta, s,t, \kappa, \gamma)$-Furstenberg configuration with $\gamma \ge \kappa-\eps_1$ and whose set of tubes has $t$-Convex Wolff constant at most $\d^{-\eta}$.

    Then for arbitrarily small $\eta,\d$ one can construct a $(\d, \eta, s,t, \kappa, \kappa)$-Furstenberg configuration $(\TT, E,\{\TT_p\})$ with the following properties. 

    \begin{itemize}
    \item[(i)] We can write
    \[
    A_z = \{ (x, y+ f(z)x), ~ x \in X_{y,z},~ y\in Y_z \},
    \]
    for some $f: Z\to [-1,1]$. Here $Y_z$ is $(\d, b,\d^{-\eta})$-AD-regular and $X_{y,z}$ is $(\d, a, \d^{-\eta})$-AD-regular with 
    \[
    a+b = 2t-\kappa, \quad a \ge s+\kappa.
    \]
    Furthermore, for every $\varrho\times \tau$ tube $T$, we have $|A_z \cap T|_\varrho \lesssim \d^{-\eta} (\tau/\varrho)^{a}$.

    \item[(ii)] For every $p \in E_z$, there is a $(\d, \kappa, \d^{-\eta})$-AD-regular set $\Phi_p \subset [-1,1]$ such that for every $\varphi \in \Phi_p$ there is $T \in \TT$ containing $p$ such that $\theta(T) \in N_\d (\varphi, \xi_p + f(z) \varphi)$. We denote by $\TT_p$ the set of all such tubes. 
    
    The function $f:Z\to [-1,1]$ is dyadic Lipschitz, i.e. for every ${\bf z} = z^\varrho$ and $z_1,z_2\in {\bf z} \cap Z$ we have $|f(z)-f(z')| \lesssim \varrho$.

    \item[(iii)] The set $\TT$ is $t$-Convex Wolff, i.e. $C_{t-CW}(\TT) \le \d^{-\eta}$.
    \end{itemize}
\end{prop}

\begin{proof}
    Fix $\eta_0>0$, let $\eps_0\ll \eta_0$ and fix $\d_0$. Consider an $(\d, \eta, s, t, \kappa, \gamma)$-Furstenberg configuration such that $\gamma \ge \kappa-\eps_1$ with $\eta,\d$ chosen sufficiently small with respect to other parameters.

    By applying Lemma \ref{lem:finding-parallel-rectangles} to sets $A_z$ and passing to a refinement, we can find a pair of scales $\varrho, \tau \in [\d,1]$ such that the sets $A_z \cap B_\varrho$ rescale to nearly $(\varrho/\tau, 2t-\kappa, a, (\varrho/\tau)^{-\varrho})$-aligned sets. Here we fix some $\zeta\ll \eta_0$ and we have $\varrho/\tau \le \d^{\chi(\zeta)}$ and $a \in [0,2t-\kappa]$. Now we can apply the rescaling and thickening map $\TT \mapsto \TT_{\varrho}^{T_\tau}$ for a fixed $T_\tau \in \TT_\tau$ and define the corresponding Furstenberg configuration on scale $\tilde \d = \varrho/\tau$. By adjusting the new set $\tilde E$ below scale $\tilde \d$, we can ensure that the new sets $\tilde A$ are $(\tilde \d, 2t-\kappa, a, \tilde d^{-\zeta}\d^{-C\eta})$-aligned. 
    It follows that we can write $\tilde A_z = \{ (x, y+f(z)x, z), x\in \tilde X_{y,z},y\in \tilde Y_z\}$ where $\tilde X_{y,z}$ is $(\tilde \d, a, \tilde \d^{-O(\zeta)})$-AD-regular, $\tilde Y_z$ is $(\tilde \d, b, \tilde \d^{-O(\zeta)})$-AD-regular ($b= 2t-\kappa-a$) and $\tilde A_z$ is tube Katz-Tao. We can check that $a \ge s+\kappa$ holds by running the same argument as in Proposition \ref{prop:global-grains-1}. When we apply the ABC sum-product theorem, we use the fact that $\min(s+\gamma,1) \ge s+\kappa-\eps_0$ to conclude that $a \ge s+\kappa-\eps_0-\eps$. By adjusting the error terms we can turn this into $a\ge s+\kappa$, as desired. This gives (i) with $\eta_0$ in place of $\eta$ if we choose parameters sufficiently small.

    To establish (ii), we can run the proof of Proposition \ref{prop:global-grains-2} again: if the slope $\alpha_{\bf q}$ describing the sets $\Theta_p$, $p\in E({\bf q})$, differs significantly from $f(z)$, then we can estimate the $\varrho^2$-covering number of $A_z({\bf q})$ by considering orthogonal projections along vectors $(1,f(z))$ and $(1, \alpha_{\bf q})$. These projections have at least $s+\kappa$ dimensional fibers, so this leads to the inequality
    \[
    2t-\kappa \le 2(2t-\kappa - (s+\kappa)) + o(1),
    \]
    so $\kappa \le \frac{2}{3}(t-s)+o(1) < \frac{2}{3}(t-s)+2\eps_0$ for small enough parameters, as desired. The dyadic Lipschitz property of $f$ is established in the same way (see e.g. proof of Lemma \ref{lemma:f-Lipschitz}). 
\end{proof}

\subsection{$\Pi$-tuples.} Suppose we are given a Furstenberg configuration satisfying properties in Proposition \ref{prop:porous-global-grain}. Furthermore suppose that for $\varrho\in [\d^{1/4}, \d^\mu]$ we have defined sets $E({\bf q})$, $Z({\bf q})$ etc according to Lemma \ref{lem:pigeonhole-rho-ball}. Fix ${\bf q} = p_0^\varrho$ and consider the following set of configurations:
\[
\begin{aligned}
\Pi_{z,z'} = \{ (p_1, p_2, p_1', p_2', T_1, T_2) \in E_z({\bf q})^2 \times E_{z'}({\bf q})^2 \times \TT^2,\\ 
T_1 \in \TT_{p_1} \cap \TT_{p_1'}, \quad T_2 \in \TT_{p_2} \cap \TT_{p_2'} \\
p_1' \leftrightarrow p_2', \quad |\theta(T_1)-\theta(T_2)| \le \varrho \}    
\end{aligned}
\]
Recall here that $p_1' \leftrightarrow p_2'$ means that we can write $p'_i = (x_i', y' +f(z') x_i', z')$ for some $y' \in Y_{z'}$ and $x_1',x_2' \in X_{y,z}$. 

As usual, we have a counting lemma.

\begin{lemma}\label{lem:counting-pi-tuples}
    For any $\tilde Z\subset Z({\bf q})$ we have
    \[
    \sum_{z,z'\in \tilde Z} |\Pi_{z,z'}| \sim_{\d^{O(\eta)}} |E_{\tilde Z}({\bf q})| |\tilde Z| (\varrho/\d)^{a} (\varrho /\d^2)^{\kappa}
    \]
\end{lemma}

\begin{proof}
    We start by selecting `half' of the $\Pi$-tuple, i.e. $(p, p', T)$ so that $T\in \TT_p \cap \TT_{p'}$. We then fix the values of $\theta(T)^\varrho$ and $y'$ (y-coordinate of $p'$) and apply Cauchy--Schwarz. 
\end{proof}

\medskip 
\noindent{\em Additive relations.} Change coordinates so that $f(z_0) = 0$ for some fixed $z_0 \in Z({\bf q})$. 
Fix a tuple $(p_1, p_2, p_1', p_2', T_1, T_2) \in \Pi_{z,z'}$ and write
\[
p_i = (x_i, y_i+f(z)x_i, z), \quad p_i' = (x_i', y_i'+f(z')x_i', z')
\]
and 
\[
\theta(T_i) \in N_\d (\varphi_i, \xi_{p_i} + f(z)\varphi_i)
\]
for some $\varphi_i \in \Phi_{p_i}$. By definition, we have $|\varphi_1-\varphi_2|\lesssim \varrho$, $|\xi_{p_1} - \xi_{p_2}|\lesssim \varrho$. Thus, using
\[
p_i' = p_i + (z'-z) \theta(T_i) + O(\d),
\]
we conclude that
\[
\begin{aligned}
    x_i' &= x_i + (z'-z) \varphi_i +O(\d),\\
    y_i' + f(z') x_i' &= y_i + f(z) x_i + (z'-z) (\xi_{p_i} +f(z) \varphi_i) + O(\d).
\end{aligned}
\]
by subtracting $i=1$ from $i=2$ case, we get $x_2'-x_1' = x_2-x_1 + O(\varrho^2+\d)$ and 
using $y_1'=y_2'$, we obtain
\begin{equation}\label{eq:Pi-y-difference}
\begin{aligned}
y_2-y_1& = f(z') (x_2'-x_1') - f(z) (x_2-x_1) - (z'-z) (\xi_{p_2} -\xi_{p_1} + f(z)(\varphi_1-\varphi_2)) +O(\d) \\
&= (f(z')-f(z)) (x_2-x_1) - (z'-z) (\xi_{p_2}-\xi_{p_1}) + O(\varrho^3 +\d).
\end{aligned}
\end{equation}
Here we also used that $|f(z)| = |f(z)-f(z_0)|\lesssim \varrho$. 

Now for a fixed choice of $y_0 \in Y_{z_0}({\bf q})$, define a function $\xi: X_{y_0,z_0}\to [-1,1]$ by setting $\xi(x) = \xi_{(x, y_0,z_0)}$. By using the conclusion of Lemma \ref{lem:pigeonhole-rho-ball} on scale $\varrho^2$, we may pass to a refinement such that $\xi_p$ satisfies $|\xi_{p}-\xi_{p'}|\lesssim \varrho^2$ for $p,p'$ in the same dyadic $C\varrho^2$-cube. Furthermore, we may assume that the union of sets $X_{y,z_0}$ over all $|y-y_0| \lesssim \varrho^2$ is contained in a $C\varrho^2$-neighborhood of $X_{y_0,z_0}$. Indeed, the tube Katz--Tao axiom of the set $A_{z_0}$ implies
\[
\left|\bigcup_{y\in Y_{z_0}({\bf q}) \cap N_{C\varrho^2}y_0} X_{y, z_0}({\bf q}) \times \{y\} \right|_{\varrho^2} \lesssim_{\d^{O(\eta)}} \varrho^{-a} \lesssim |X_{y_0,z_0}({\bf q})|_{\varrho^2},
\]
where we use that the set on the left hand side is covered by a $\varrho^2\times \varrho$ tube. Thus, by choosing an average $y_0$ and pigeonholing, we get the claim.

Now suppose that the chosen $\Pi$-tuple satisfies $z=z_0$ and $y_1 = y_0+O(\varrho^2)$ (so that $x_1 \in N_{C\varrho^2} X_{y_0,z_0}$). Then we conclude from \eqref{eq:Pi-y-difference} and Lipschitz properties that $|y_2-y_1| \lesssim \varrho^2$ and so $x_2 \in N_{C\varrho^2} X_{y_0,z_0}$. Replace $x_i$ by an element $\tilde x_i\in X_{y_0,z_0}$, then we get that
\begin{equation}
    y_2 - y_1 = f(z') (\tilde x_2-\tilde x_1) - z' (\xi(\tilde x_2) - \xi(\tilde x_1)) + O(\varrho^3+\d).
\end{equation}

It follows from Lemma \ref{lem:counting-pi-tuples}, that there exists a graph $G \subset Y_{z_0}({\bf q})\times X_{y_0,z_0}({\bf q})^2 \times Z({\bf q})$ of density $\d^{O(\eta)}$ such that for every $(y,x_1,x_2,z)\in G$ we have 
\begin{equation}\label{eq:y-pi-additive-symmetry}
y+ f(z)(x_2-x_1) - z (\xi(x_2) - \xi(x_1)) \in N_{C\varrho^3} Y_{z_0}({\bf q}).    
\end{equation}
In the following sections we use this expression in cases when $f$ is linear and non-linear. 

\subsection{Non-linear case.}
When $f$ cannot be approximated by a linear function on $Z({\bf q})$, we will reduce to the following ABC sum-product type statement. 
\begin{lemma}\label{lem:ABC-dot-product}
    Fix $\beta\le \gamma \in (0,1]$ and let $\alpha \in (0, \frac{\beta+\gamma}{2})$. Let $\varepsilon_1>0$ be sufficiently small in $\alpha, \beta, \gamma$. Let $\varepsilon_2>0$ be arbitrary and let $\eta,\d$ be sufficiently small depending on previous parameters. 
    
    Let $B, C \subset [-1,1]^2$ and $A \subset [-1,1]$ be $\d$-separated sets such that (a), (b) and either (c) or (c') hold: 
    \begin{itemize}
        \item[(a)] $B$ is a $(\d, \beta, \d^{-\eta})$-set, $C$ is a $(\d, \gamma, \d^{-\eta})$-set,
        \item[(b)] for every $\d^{\varepsilon_1}\times 1$ tube $T$ we have 
        \[
        |B \cap T| \le \d^{\varepsilon_2} |B|,
        \]
        \item[(c)] there is a graph $G\subset A\times B\times B\times C$ such that $|G| \ge \d^{\min(\eta, \varepsilon_2)/10} |A||B|^2 |C|$ and
        \[
        \{ a+ (b-b')\cdot c, \quad (a,b,b',c) \in G \} \subset N_{C\d} A.
        \]
        \item[(c')] there is a graph $G\subset A\times B\times C\times C$ such that $|G| \ge \d^{\min(\eta, \varepsilon_2)/20} |A||B| |C|^2$ and
        \[
        \{ a+ b\cdot (c-c'), \quad (a,b,c,c') \in G \} \subset N_{C\d} A.
        \]
    \end{itemize}
    Then $|A| \ge \d^{-\alpha}$.

\end{lemma}

\begin{proof}
    Let us first assume that (a), (b) and (c) hold. 
    Fix some $\zeta>0$ such that $\varepsilon_1 \ll \zeta \ll_{\alpha, \beta, \gamma} 1$ and apply radial projections theorem (Theorem \ref{thm:radial-projections}) to the set $B$. We obtain a graph $H \subset B\times B$ such that $|H| \ge (1 - \d^{\min(\eta,\varepsilon_2/2)})|B|^2$ and so that for every $b_0 \in B$ with $|H|_{b_0}| \ge \nu |B|$, the set of radial directions $\frac{b-b_0}{|b-b_0|}$, $b\in H|_{b_0}$ forms a $(\d, \beta, \nu^{-1}\d^{-\zeta})$-set. By pruning, we may assume that $|b-b'| \ge \d^{C\eta/\beta}$ for every $(b,b')\in H$.

    It follows that for sufficiently small $\d$, we have $|G \cap (A\times H \times C)| \ge |G|/2$, thus we may restrict to quadruples $(a, b, b', c)$ with $(b, b') \in H$. 
    Now we can apply sharp Kaufman's exceptional set estimate \cite{ren2023furstenberg} (which is a consequence of the Furstenberg set estimate, Theorem \ref{thm:ren-wang}) to the set $C$ and the set of directions $\Lambda = \{ \frac{b-b'}{|b-b'|}, (b,b')\in H \}$ to conclude that for every $\chi>0$ and $\zeta\ll_{\chi} 1$ there exists $(b,b')\in B$ such that 
    \[
    |(b-b')\cdot C' |_\d \gtrsim \d^{\chi} \d^{-\min(\gamma, \frac{\gamma+\beta}{2}, 1)}
    \]
    for any $C' \subset C$, $|C'| \gtrsim \d^{\eta} |C|$. By selecting $C'=C_{a,b,b'}$ to be the set of $c$ such that $(a_0, b, b', c)\in G$ and averaging over the choice of $a_0$, it follows that 
    \[
    |A|_\d \gtrsim |(b-b') \cdot C|_\d \gtrsim \d^{\chi} \d^{-\min(\gamma, \frac{\gamma+\beta}{2}, 1)} > \d^{-\alpha}
    \]
    by the assumption on $\alpha$ and for $\chi < \frac{\beta+\gamma}{2}-\alpha$. This completes the proof.  

    Now let us consider the case when (c') holds instead of (c) and let $G \subset A\times B\times C \times C$ be the corresponding graph. By Cauchy-Schwarz, there are at least $\gtrsim \frac{|G|^2}{|A|_\d |C|^2}$ many tuples $(a_1, a_2, b_1, b_2, c, c')$ such that $(a_1, b_1, c, c'), (a_2, b_2, c,c') \in G$ and there exists $a' \in A$ such that $a_i + b_i(c-c') \in N_{C\d} a'$. By rearranging, $a_1 + (b_1-b_2) (c-c') \in N_{2C\d} a_2$. Recall $A$ is $\d$-separated, it follows, that there exists a graph $\tilde G \subset A \times B \times B \times C\times C$ of density $\gtrsim \d^{\min(\eta, \varepsilon_2)/10}$ such that $a+(b-b')(c-c') \in N_{2C\d}A$ for $(a,b,b',c,c')\in \tilde G$. Restrict the graph to a fiber above a typical element $c' \in C$ and replace $C$ by $C-c'$. This gives a new graph in $A\times B\times B\times C$ satisfying property (c), so the above proof gives the desired conclusion.
\end{proof} 

Applying this lemma with rescaled versions of the sets $A = N_{C\varrho^3}Y_{z_0}({\bf q}) \cap N_{C\varrho^2}y_0$, 
\[
B = \left\{\begin{pmatrix}
z \\ 
f(z)
\end{pmatrix},~~ z \in Z({\bf q}) \right\}, \quad C = \left\{\begin{pmatrix}
-\xi(x) \\ 
x
\end{pmatrix},~~ x \in X_{y_0,z_0}({\bf q}) \right\}.
\]
So we have $\beta = s$ since $Z({\bf q})$ is a rescaled $(\d/\varrho, s, \d^{-\eps})$-set and $\gamma = a \ge s+\kappa$ by the assumption on $X_{y,z}$. 

If the 2-ends assumption in Lemma \ref{lem:ABC-dot-product} applies, then using the (c') alternative we conclude 
\[
\varrho^{-b} \sim_{\d^{O(\eta)}} |Y_{z_0}({\bf q})\cap N_{C\varrho^2}y_0|_{\varrho^3} \gtrsim \varrho^\zeta \varrho^{ -\frac{s+a}{2} }
\]
where $\zeta>0$ is arbitrary. So 
\[
b \ge \frac{s+a}{2} -\zeta - O(\eta/\mu)
\]
and using $a+b=2t-\kappa$, $a \ge s+\kappa$, this implies $\kappa \le \frac{4}{5}(t-s) +o(1) \le \max(\frac{4}{5}(t-s),0)+\eps_0$, as desired. 

\subsection{Linear and slow variation case.} Now we may assume that $B$ concentrates in a $\d^{\varepsilon_1}\times 1$ tube. Here $\varepsilon_1$ depends only on the difference $\kappa-\kappa_0(s,t)$ and we can choose $\varepsilon_2\sim \eta$. 
In other words, for some $r, u\in \RR$ we have
\[
\#\{z \in Z({\bf q}):~ |f(z) - rz-u| \le \varrho^{\varepsilon_1+1}\} \gtrsim_{\d^{O(\eta)}} |Z({\bf q})|
\]
Let $\tilde Z \subset Z({\bf q})$ be the set of $z$ as above. Recall that we assume that $|f(z)|\lesssim \varrho$ for all $z\in Z({\bf q})$, so we have $|u| \lesssim \varrho$ and $|u| \lesssim 1$. 

By using $V$-tuples from Section \ref{subsec:W-tuples} we can construct a $\d^{O(\eta)}$-dense subgraph $G' \subset Y_{z_0}({\bf q}) \times \tilde Z\times \overline{\Phi}\times \overline{\Phi}$ such that for every $(y_1, z, \varphi_1, \varphi_2) \in G'$ there is $y_2 \in Y_{z_0}({\bf q})$ such that
\[
y_2 = y_1 -(z-z_0) (f(z)-f(z_0)) (\varphi_1-\varphi_2) + O(\varrho^3),
\]
see \eqref{eq:V-tuple-relation-Furstenberg}. Here $\overline{\Phi} \subset \varrho \ZZ \cap [-1,1]$ is the $(\varrho, \kappa, \d^{-O(\eta)})$-AD-regular set such that $\overline{\Phi}^\varrho = \Phi_{\bf q}^\varrho$. The linear approximation of $f$ implies that
\[
y_2 = y_1 - (z-z_0)^2 r (\varphi_1-\varphi_2) + O(\varrho^3+\varrho^{2+\varepsilon_1}).
\]
Suppose that $|r| > \varrho^{\varepsilon_1/10}$ holds. The set of $(z-z_0)^2$ over $z \in \tilde Z$ ranges over a rescaled $(\varrho, s, \varrho^{-\eps}\d^{-O(\eta)})$-set (after removing all $z$ which are too close to $z_0$) and the set of angles $\varphi_i$ ranges over a $(\varrho, \kappa, \d^{-O(\eta)})$-AD-regular set. 

Apply the ABC sum-product theorem (Theorem \ref{thm: ABC}) to the (rescaled) sets $\tilde A = Y_{z_0}({\bf q}) \cap N_{C\varrho^2 |r|}y_0$, $\tilde B = \{ (z-z_0)^2, ~z\in \tilde Z \}$ and $C = \{ \varphi_1-\varphi_2, \varphi_i\in \overline{\Phi} \}$ and the scale $\Delta = \frac{\varrho^{2+\varepsilon_1}}{C\varrho^2 |r|}$ . It follows that 
\[
| Y_{z_0}({\bf q}) \cap N_{C\varrho^2 |r|}y_0|_{C\varrho^2 |r|\Delta} \gtrsim \Delta^{\zeta} \Delta^{-\min(s+\kappa, 1)}.
\]
Using the AD-regularity of $Y_{z_0}$ with dimension $b \le 2t-\kappa-(s+\kappa)$, we conclude that 
\[
2t-\kappa \ge 2(s+\kappa) -o(1)
\]
giving the desired bound $\kappa \le \max(\frac{2}{3}(t-s),0)+\eps_0$. 

Thus, we may restrict to the case $|r| \le \varrho^{\varepsilon_1/10}$. In other words, the slope function $f$ varies slower than Lipschitz on a dense subset of $Z({\bf q})$. The proof of Lemma \ref{lem:slowly-rotating} can be applied essentially verbatim to deal with this situation. Note that the construction of $W$-tuples didn't need the fact that global grains have full dimension (this was only used in $L$-tuples). By applying the same ABC sum-product neighborhood expansion argument, we can prove the lower bound $|X_m|_{\tau\varrho} \gtrsim\varrho^{\zeta}\tau^{ - (s+\kappa)}$ (note that by definition $X_m$ is a $\varrho^2$-separated subset in a $\varrho$-interval). For $\sigma  = \max(\varrho^2, \varrho^{1+\varepsilon_1/10})$ this then gives us a convex $\sim_m\sigma \times \varrho \times 1$ box $U$ which satisfies
\[
|\TT_\sigma[U]| \gtrsim_{\d^{O(\eta)}} |X_m|_{\sigma} (\varrho/\sigma)^{\kappa} \gtrsim \varrho^\zeta (\varrho/\sigma)^{s+2\kappa}.
\]
On the other hand, the Convex Wolff axiom of $\TT$ implies that
\[
|\TT_\sigma[U]| \le \d^{-\eta} (\varrho/\sigma)^t
\]
so $t \ge s+2\kappa-o(1)$ giving the desired lower bound on $\kappa$. This covers all possible cases and completes the proof of Proposition \ref{prop: very-coplanar}.

\begin{remark}
    In this last step, it suffices to use a weaker Convex--Wolff axiom. Indeed, suppose that $C_{t;t'-CW}(\TT) \le \d^{-\eta}$ for some parameter $t'$ (see Definition \ref{def:convex-wolff-axiom}). Then to make the above argument work, we need to require $t' < s+2\kappa$. Recalling that $\kappa_0(s,t) = \max(0,t-s, 2t-2)$ and in this section $s+\kappa\le 1$, this leads to the following permissible $t'$:
    \[
    t'=     \begin{cases}
        \min(2t, s), \quad t \in (0,s],\\
        2t-s, \quad t \in [s,1].
    \end{cases}
    \]
    Note that if $t>1$ then $\kappa \ge t-s \ge 1-s$, so this contradicts the assumption $s+\kappa\le 1$.
    Thus, for Proposition \ref{prop: very-coplanar}, i.e. the case when $s+\kappa\le 1$ and $\gamma \ge \kappa-\eps_0$ it suffices to assume this weaker $(t;t')$-Convex Wolff axiom.
\end{remark}

\newpage 
\part{Sticky Kakeya}

\section{Setup for Sticky Kakeya}\label{sec:setup-kakeya}

We say that a shading $Y$ on $\TT$ is $\lambda$-dense, if $\sum_{T\in \TT} |Y(T)| \ge \lambda \sum_{T\in \TT}|T|$ (where $|\cdot|$ denotes the volume).

\begin{definition}
    Let $\operatorname{K}_d(\kappa)$ denote the following statement. 

    For any $\eps >0$ there exists $\eta = \eta(d, \eps)>0$ such that the following holds for $\d>0$ sufficiently small. Let $\TT$ be a $(d-1, \d^{-\eta})$-AD-regular set of $\d$-tubes in $\RR^d$ and $C_{CW}(\TT) \le \d^{-\eta}$. Suppose $Y$ is a $\d^\eta$-dense shading on $\TT$. Then 
    \[
    |U(\TT, Y)|_\d \ge \d^\eps \d^{-d+\kappa}.
    \]
\end{definition}


Let $\operatorname{Mat}(k, \ell)$ denote the set of $k\times \ell$ matrices with entries in $[-1,1]$.

\begin{definition}\label{def:kakeya-configuration}
    Let $d \ge 3$ and $\ell \in\{2, \ldots, d-1\}$. 
    A global $(\ell-1)$-grains  $(\d, \eta, \kappa)$-Kakeya configuration $(E, \TT, \{\TT_p\}_{p\in E})$ in $\RR^d$ is the following collection of data:
    \begin{enumerate}
        \item A $(d-1, \d^{-\eta})$-AD-regular set of $\d$-tubes $\TT$.
        \item A $\d$-separated set $Z \subset [-1,1]$ with $|Z| \ge \d^{\eta-1}$ and a $f: Z \to \operatorname{Mat}(d-\ell,\ell-1)$ satisfying ``dyadic Lipchitz continuity'': for any  $\varrho\in [\delta, 1]\cap \delta^{\eta\mathbb{Z}}$ and any  $z_1, z_2$ in the same dyadic $\varrho$-interval, $|f(z_1)-f(z_2)|\lesssim \varrho$, where $|\cdot|$ means $\ell^{\infty}$-norm over all entries of the matrix.  
        \item A $\d$-separated set $E \subset [-1,1]^d$ of the form 
        \[
        E = \bigsqcup_{z\in Z} E_z =\{ (x, y+ f(z) x, z),~ z \in Z,~ y\in Y_z, ~x \in X_{y,z} \},
        \]
        where $Y_z \subset [-1,1]^{d-\ell}$ is $(\d, d-\ell-\kappa, \d^{-\eta})$-AD-regular set and $X_{y,z} \subset [-1,1]^{\ell-1}$ is such that $|X_{y,z}|_\d \ge \d^{\eta-\ell+1}$. In addition, for every $\varrho \in [\delta, 1]\cap \delta^{\eta\mathbb{Z}}$ and a dyadic $\varrho$-interval $\mathbf{z}$ such that $Z\cap \mathbf{z} \neq\emptyset$, $\bigcup_{z\in \mathbf{z}\cap Z} Y_z$ is a $(\varrho, d-\ell-\kappa, \delta^{-\eta})$-AD-regular set. 
        
        \item For each $p \in E_z$ there is $\xi_p \in [-1,1]^{d-\ell}$ and $(\d, \kappa, \d^{-\eta})$-AD-regular set $\Phi_{p}\subset [-1,1]^{\ell-1}$ such that for every $\varphi \in \Phi_{p}$ there is a tube $T\in \TT$ through $p$ such that $\theta(T)=  (\varphi, \xi_p+f(z) \varphi)  +O(\delta)$. In addition, for any $\varrho \in [\delta, 1]\cap \delta^{\eta\mathbb{Z}}$, write $\mathbf{q}=p^{\varrho}$, then there is $\xi_{\mathbf{q}}\in [-1,1]^{d-\ell}$ and  a $(\varrho, \kappa, \delta^{-\eta})$-AD regular set $\Phi_{\mathbf{q}}\subset [-1, 1]^{\ell-1}$  such that $\Phi_{p'} \subset \Phi_{\mathbf{q}}^{\varrho}$ and $|\xi_{p'}-\xi_{\mathbf{q}}|\lesssim \varrho$ for any $p'\in E\cap \mathbf{q}$ (including $p$). 
    \end{enumerate}
\end{definition}

\begin{remark}\label{rem: configuration shading}
    Given a global $(\ell-1)$-grains $(\delta, \eta, \kappa)$-Kakeya configuration, we can also define a shading on the corresponding set of $\delta$-tubes $\TT$: for each $T\in \TT$, let $Y(T)$ be the union of dyadic $\d$-cubes $p^{(\d)}$ over $p\in E_z, z\in Z$ such that $\theta(T)= (\varphi, \xi_p+f(z)\varphi)+O(\delta)$. It follows from Definition \ref{def:kakeya-configuration} that this shading is $\gtrapprox \d^{C\eta}$-dense. 
\end{remark}

\begin{definition}[Slab Frostman, Tube Katz-Tao]
We say that a set of points $P \subset [-1,1]^d$ is $(\d, s, C)$-{\em slab Frostman} if for any $\varrho \in [\d, 1]$ and any $\varrho\times 1\times \ldots \times 1$ slab $S \subset \RR^d$ we have $|P \cap S| \le C \varrho^s |P|$. Similarly, let us say that $P$ is $(\d, s, C)$-{\em tube Katz-Tao} if for any $\varrho \in [\d, 1]$, $\tau \in [\varrho,1]$ and any $\varrho\times\cdots \times \varrho\times  \tau$ tube $T \subset \RR^d$ we have $|P \cap T|_\varrho \le C (\tau/\varrho)^{s}$.
\end{definition}

If, in addition, the set $\Phi_p$ is $(\d, \gamma, \d^{-\eta})$-slab Frostman for every $p \in E$, then we say that $(E, \TT, \{\TT_p\})$ is an $\ell$-linear $(\d,\eta,\kappa, \gamma)$-Kakeya configuration. 
If for every $z \in Z$ the set $Y_z$ is $(\d, 1-\gamma, \d^{-\eta})$-tube Katz Tao, then we say that $E$ is an $\ell$-planar $(\d,\eta,\kappa, \gamma)$-Kakeya configuration. If both of these conditions are satisfied then we call $E$ an $\ell$-planar $\ell$-linear $(\d,\eta,\kappa, \gamma)$-Kakeya configuration. 

These conditions ensure that the slope function $f$ gives a non-trivial description of the Kakeya configuration $E$. Indeed, if the Kakeya configuration is not $\ell$-linear, then the set of tubes though a point $\TT_p$ spans some $\ell'$-dimensional plane for $\ell'<\ell$. In this case, we can replace the $\ell$-dimensional global grains by $\ell'$-dimensional grains and study the resulting $\ell'$-linear Kakeya configuration. If the Kakeya configuration is not $\ell$-planar for some tiny $\gamma>0$, then this means that sets $Y_{z}$ concentrate a lot in tubes. Since $Y_z$ is AD-regular, the structure theorem for AD-regular sets (Proposition \ref{lem:finding-parallel-rectangles}) implies that after some rescaling $Y_z$ is a union of essentially fully filled out parallel tubes. So we can replace global $(\ell-1)$-grains by $\ell$-dimensional grains and study the resulting global $\ell$-grains configuration. 

In fact, we will need a more general definition of a $k$-planar $\ell$-linear Kakeya configuration for an arbitrary pair $k \ge \ell \in \{2, \ldots, d-1\}$. In this definition, $E$ has $(k-1)$-dimensional horizontal global grains that are parameterized by a slope function $\tilde f$ and the set of tubes passing through a fixed point $p\in E$ is contained in an $\ell$-dimensional plane parameterized by a slope function $f$. 

\begin{definition}\label{def:kakeya-configuration-2}
    Let $d\ge 3$ and $\ell\le k \in \{2, \ldots, d-1\}$. A $\ell$-linear $k$-planar  $(\d,\eta, \kappa,\gamma)$-Kakeya configuration $(E, \TT, \{\TT_p\}_{p\in E})$ in $\RR^d$ is
    a global $(\ell-1)$-grains $(\delta, \eta, \kappa)$-Kakeya configuration with  $Z, f, E$, $\Phi_p$ for each $p\in E$,  and the following additional information:

    \begin{enumerate}
        \item For $Z\subset [-1,1]$  and $f: Z \to \operatorname{Mat}(d-\ell,\ell-1)$ there is a slope function $\tilde f: Z \to \operatorname{Mat}(d-k,k-1)$ such that if 
        we further write 
        \[
         f =\begin{pmatrix}
            f_1 \\
            f_2 
        \end{pmatrix}, 
         \quad \tilde f = (\tilde f_1, \tilde f_2)
        \] 
        where $f_1 \in \operatorname{Mat}(k-\ell, \ell-1)$, $f_2, \tilde f_1 \in \operatorname{Mat}(d-k, \ell-1)$ and $\tilde f_2 \in \operatorname{Mat}(d-k, k-\ell)$, they obey the following consistency relation:
            \begin{equation}\label{eq:consistent-f}
      f_2(z) = \tilde f_1(z) + \tilde f_2(z) f_1(z).
    \end{equation}
        
        \item The $\d$-separated set $E \subset [-1,1]^d$  can be written as  
        \[
        E = \bigsqcup_{z\in Z} E_z =\{ (\tilde x, \tilde y+ \tilde f(z) \tilde x, z),~ z \in Z,~ \tilde y\in \tilde Y_z, ~\tilde x \in \tilde X_{y,z} \},
        \]
        where $\tilde Y_z \subset [-1,1]^{d-k}$ is $(\d, d-k-\kappa, \d^{-\eta})$-AD-regular set\footnote{We sometimes abuse the notation and reuse $Y_z$ for $\tilde Y_z$ to denote a different object than the one in Definition~\ref{def:kakeya-configuration}. It should be clear from context which one we mean and they rarely show up the same time. } and $\tilde X_{y,z} \subset [-1,1]^{k-1}$ is such that $|X_{y,z}|_\d \ge \d^{\eta-k+1}$.  Here we view $\tilde x,\tilde y$ as column vectors. 
        \item $\Phi_p \subset [-1, 1]^{\ell-1}$ is $(\d, \gamma, \d^{-\eta})$-slab Frostman for every $p \in E$. 
        \item $Y_z\subset [-1, 1]^{d-k}$ is $(\d, 1-\gamma, \d^{-\eta})$-tube Katz-Tao for every $z \in Z$. 
    \end{enumerate}
\end{definition}

Note that equation (\ref{eq:consistent-f}) reflects the fact that the set of directions of tubes in $\TT_p$ is contained in a translate of the $k$-dimensional plane at slope $\tilde f$. 

If $\ell=2$ then the set $\Phi_p$ is automatically $(\d, \kappa, \d^{-\eta})$-slab Frostman. Similarly, if $k = d-1$, then the set $Y_z$ is automatically $(\d,1- \kappa, \d^{-\eta})$-tube Katz-Tao.
In some situations we only need to keep track of only one of the two slope functions $\tilde f$ or $f$ in which case the simpler Definition \ref{def:kakeya-configuration} would suffice.

The proof of Theorem \ref{thm:kakeya-R4} splits into the following steps. Some of the steps can be generalized to the higher dimensional sticky Kakeya and so we state those for arbitrary dimension $d$.

First, we show that if $\kappa$ is the smallest value such that $\operatorname{K}_d(\kappa)$ is true then we can construct a $k$-planar $\ell$-linear $(\d, \eta, \kappa, \gamma)$-Kakeya configuration for some $\ell\le k \in \{2, \ldots, d-1\}$.

\begin{prop}\label{prop:kakeya-reduction}
    Let $\kappa$ be the infimum of the set of $\kappa'$ such that $\operatorname{K}_d(\kappa')$ holds and assume that $\kappa>0$. Let $g:(0,1)\to (0,1)$, $h:(0,1)^2\to (0,1)$ be continuous functions. 
    Then for any $\d_0, \eta_0$ there exist $\d\in (0,\d_0)$, $\eta \in (0,\eta_0), \gamma > 0$ and $\ell\le k \in \{2, \ldots, d-1\}$ so that there exists a $k$-planar $\ell$-linear $(\d, \eta, \kappa, \gamma)$-Kakeya configuration whose set $\TT$ of $\d$-tubes satisfies $C_{CW}(\TT) \le \d^{-\eta}$ and we have $\eta < g(\gamma)$ and $\d \le h(\eta, \gamma)$.
\end{prop}

A typical choice of the decreasing functions would be $h(\eta, \gamma) = \exp(-\eta^{-2})$ and $g(\gamma) = \gamma^2$. 

This proposition reduces the sticky Kakeya problem in $\RR^4$ to analyzing Kakeya configurations with $(k, \ell) = (3, 3)$, $(2, 2)$  and $(3,2)$. 
We start by dealing with the case of $(d-1)$-linear configurations in $\RR^d$. 

\begin{prop}\label{prop:codimension-1-kakeya}
    Let $\zeta>0, \gamma>0$ and $d\ge 3$. The following holds for $\eta \le \eta_{\ref{prop:codimension-1-kakeya}}(\zeta, \gamma, d)$ and $\d \le \d_0(\eta, \gamma, \zeta, d)$.
    Suppose that there exists a $(d-1)$-linear $(\d, \eta, \kappa, \gamma)$-Kakeya configuration in $\RR^d$ whose set of tubes $\TT$ satisfies $C_{CW}(\TT) \le \d^{-\eta}$. Then we have $\kappa \le \zeta$. 
\end{prop}

Next, we consider $2$-planar Kakeya configurations in $\RR^d$. 

\begin{prop}\label{prop:2-planar-kakeya}
    Let $\zeta>0$ and $d\ge 3$. The following holds for $\gamma >0$, $\eta \le \eta_{\ref{prop:2-planar-kakeya}}(\zeta, \gamma, d)$ and $\d \le \d_0(\eta, \gamma, \zeta, d)$. 
    Suppose that there exists a $2$-planar $(\d, \eta, \kappa, \gamma)$-Kakeya configuration in $\RR^d$ whose set of tubes $\TT$ satisfies $C_{CW}(\TT)\le \d^{-\eta}$. Then we have $\kappa\le \zeta$.
\end{prop}

Note that either of Propositions \ref{prop:codimension-1-kakeya} or \ref{prop:2-planar-kakeya} suffices to prove the sticky Kakeya conjecture in $\RR^3$. The proofs of these propositions generalize two different arguments outlines in Section \ref{sec:proof-overview}.
In $\RR^4$, it remains to consider $3$-planar $2$-linear Kakeya configurations.

\begin{prop}\label{prop:32-planar-kakeya}
    Let $\zeta>0$. The following holds for  $\gamma>0$, $\eta \le \eta_{\ref{prop:32-planar-kakeya}}(\zeta, \gamma)$ and $\d \le \d_0(\eta, \gamma, \zeta)$. 
    Suppose that there exists a 3-planar $2$-linear $(\d, \eta, \kappa, \gamma)$-Kakeya configuration in $\RR^4$ whose set of tubes $\TT$ satisfies $C_{CW}(\TT)\le \d^{-\eta}$. Then we have $\kappa \le \zeta$. 
\end{prop}

By combining these propositions we can deduce Theorem \ref{thm:kakeya-R4}. 

\begin{proof}[Proof of Theorem \ref{thm:kakeya-R4}]
    Note that Theorem \ref{thm:kakeya-R4} is equivalent to the statement $\operatorname{K}_4(0)$. 
    Let $\kappa$ be the infimum of all $\kappa'$ such that $\operatorname{K}_4(\kappa')$ holds. For the sake of contradiction suppose that $\kappa>0$. Let $\zeta = \kappa/2$ and let $g(\gamma) =  \min\{ \eta_{\ref{prop:codimension-1-kakeya}}(\zeta, \gamma, 4), \eta_{\ref{prop:2-planar-kakeya}}(\zeta, \gamma, 4), \eta_{\ref{prop:32-planar-kakeya}}(\zeta, \gamma)\}$ (without loss of generality the $\eta_0$-functions from these propositions are monotone decreasing). Similarly, let us define $h(\eta)$ as the minimum of functions $\d_0(\eta, \gamma, \zeta, d)$ appearing in the statements of propositions. 
    By Proposition \ref{prop:kakeya-reduction}, there is a pair $(k, \ell) \in \{ (2,2), (3, 2), (3,3) \}$ such that for any $\d_0, \eta_0$ there exist $\gamma>0, \d< \d_0, \eta<\eta_0$ and a $k$-planar $\ell$-linear $(\d, \eta, \kappa, \gamma)$-Kakeya configuration such that $\eta < g(\gamma)$ and the set of $\d$-tubes $\TT$ satisfies $C_{CW}(\TT) \le \d^{-\eta}$ and $\d$ is can be taken arbitrarily small compared to $\eta, \gamma$. Then for $(k, \ell) = (3,3)$ by Proposition \ref{prop:codimension-1-kakeya}, for $(k, \ell) = (2,2)$ by Proposition \ref{prop:2-planar-kakeya} and for $(k, \ell)=(3,2)$ Proposition \ref{prop:32-planar-kakeya} we conclude that $\kappa \le \zeta$.
\end{proof}

\section{Construction of Kakeya configurations: proof of Proposition \ref{prop:kakeya-reduction}}\label{sec:configuration-kakeya}

Let us restate Proposition \ref{prop:kakeya-reduction} for convenience.

\begin{prop}\label{prop: structure1}
        Let $\kappa$ be the infimum of the set of $\kappa'$ such that $\operatorname{K}_d(\kappa')$ holds. Let $g:(0,1)\to (0,1)$, $h:(0,1)^2\to (0,1)$ be monotone decreasing functions. 
    Then for any $\d_0, \eta_0$ there exist $\d\in (0,\d_0)$, $\eta \in (0,\eta_0), \gamma > 0$ and $\ell\le k \in \{2, \ldots, d-1\}$ so that there exists a $k$-planar $\ell$-linear $(\d, \eta, \kappa, \gamma)$-Kakeya configuration whose set $\TT$ of $\d$-tubes satisfies $C_{CW}(\TT) \le \d^{-\eta}$ and we have $\eta < g(\gamma)$ and $\d \le h(\eta, \gamma)$.
\end{prop}

We first prove a weaker version and then upgrade. 

\begin{prop}\label{prop: structure2}
    Let $\kappa$ be the infimum of the set of $\kappa'$ such that $\operatorname{K}_d(\kappa')$ holds. 
    Then for any $\d_0, \eta_0$ there exists $\d<\d_0$, $\eta \le \eta_0$, and $\ell\in \{2, \ldots, d-1\}$ so that there exists a global $(\ell-1)$-grains  $(\d, \eta, \kappa)$-Kakeya configuration whose set $\TT$ of $\d$-tubes satisfies $C_{CW}(\TT) \le \d^{-\eta}$. 

\end{prop}

\begin{proof}
    Since $\kappa$ is the infimum of all $\kappa'$ such that $\operatorname{K}_{d}(\kappa')$ holds. Then for any $\nu>0$, $\operatorname{K}_{d}(\kappa-\nu)$ does not hold. 
Therefore there exists $\eps>0$ such that for all sufficiently small $\eta>0$ and $\delta>0$, there exists  a $(d-1, \delta^{-\eta})$-AD-regular set of $\delta$-tubes with $C_{CW}(\TT)\leq \delta^{-\eta}$ and a $\delta^\eta$-dense shading $Y$ on $\TT$ and 
\begin{equation*}
|U(\TT, Y)|_{\delta}\leq \delta^{\epsilon}\delta^{-d+\kappa-\nu}\leq \delta^{-d+\kappa-\nu},
\end{equation*}
Since the properties of Kakeya configurations are monotone in $\eta$, we may replace the parameter $\eta$ by $\max(\eta, \nu)$, i.e. without loss of generality,
\begin{equation}\label{eq: Kextremal}
|U(\TT, Y)|_{\delta}\leq \d^{-\eta} \delta^{-d+\kappa}.
\end{equation}
Equivalently, 
\begin{equation}\label{eq: Kextremalmu}
    \mu(\TT, Y):=\frac{\sum_{T\in \TT}|Y(T)|_{\delta}}{|U(\TT, Y)|_{\delta}}\geq \delta^{-\kappa+\eta}\geq \delta^{-\kappa+\eta}.
\end{equation}
    
    For any $\eta_0>0$, let $\eta \ll \eps_1\ll \eta_d\ll \eta_{d-1} \ll \cdots \ll \eta_{1}\ll \eta_0$ be small constants 
    to be determined. 
    Let $\delta^{\eta_0} \ll \varrho_{d-1} \ll \cdots \ll \varrho_1 \ll \delta_0$ be small scales to be determined. 
We are going to show that for any $\delta_0$ and $\eta_0$, there exists $\tilde \delta \leq \delta_0$,  
$\tilde\eta < \eta_0,$  and $\ell\in \{2, \dots, d-1\}$ so that there exists a global $(\ell-1)$-grains $(\tilde \delta, \tilde \eta, \kappa)$-Kakeya configuration in $\mathbb{R}^d$.
We are going to take $\eta, \delta$ sufficiently small compared to $\eta_0, \delta_0$, respectively.

\medskip

\noindent \textit{Step 1. Pigeonhole $\ell$-linear points.}
\label{c}
After pigeonholing and 
passing to a $\gtrapprox 1$-refinement of $Y$, we may assume that  $\TT_Y(p)$ is $\approx \delta^{-\eta}$-uniform for  each $p\in U(\TT, Y)$:  $|\TT_Y(p)[T^\varrho]|$ is approximately the same for each  $T\in \TT_Y(p)$ and $\varrho \in [\delta, 1]\cap \delta^{\eta\mathbb{Z}}$ (here $T^\varrho$ means the dyadic $\varrho$-tube containing $T$).  In addition, the quantity $|\TT_Y(p)[T^\varrho]|$   is approximately the same independent of $p\in U(\TT, Y)$, and so  $|\TT_Y(p)|_\delta \sim \mu(\TT, Y)$.

For each $p\in U(\TT, Y)$, for $\ell=1, \ldots, d-1$, let $\delta_\ell$ be the smallest scale such that there exists an $\ell$-dimensional plane $\Sigma_\ell$ through $p$ such that 
\[
|\TT_Y(p)[N_{\delta_\ell} \Sigma_\ell]|_\delta \geq \delta_\ell^{\eta_\ell} |\TT_Y(p)|_\delta.
\]
Start from $\ell=1$, if $\delta_\ell\leq \varrho_\ell$, we say that $p$ is an $\ell$-linear point, otherwise we replace $\ell$ with $\ell+1$ and define $\delta_{\ell+1}$ accordingly. When $\ell$ increases to $d$, we stop and say that $p$ is a $d$-linear point (multilinear point) and write $\delta_d=\delta^{1/2}$.

Therefore an  $\ell$-linear point $p$ satisfies 
\begin{enumerate}
    \item  there exists $\delta_\ell\leq \varrho_\ell$ and  an $\ell$-plane $\Sigma_\ell$ through $p$ such that 
    \[
    |\TT_Y(p)[N_{\delta_\ell}\Sigma_\ell]|_{\delta}\geq \delta_\ell^{\eta_\ell} |\TT_Y(p)|_\delta. 
    \]
    \item for any $(\ell-1)$-plane $\Sigma_{\ell-1}$ through $p$ and any $r\leq \varrho_{\ell-1}$ (note that $\delta_\ell \leq \varrho_\ell \ll \varrho_{\ell-1}$), 
    \[
    |\TT_Y(p)[N_{r}\Sigma_{\ell-1}]|_\delta \leq r^{\eta_{\ell-1}} |\TT_Y(p)|_\delta.
    \]
\end{enumerate}

After pigeonholing, there exists $\ell\in [1, d]$ such that a $\gtrapprox 1$-fraction of points $x\in U(\TT, Y)$ is $\ell$-linear. Define  $Y_0(T)\subset Y(T)$ as the set of $\ell$-linear points such that $T\in \TT_Y(p)[N_{\delta_\ell}\Sigma_\ell]$, then $Y_0$ is a $\gtrapprox  \delta_\ell^{\eta_\ell}\delta^\eta$-dense shading on $\TT$. 

For now we have $\ell\in \{1, \dots, d\}$, we first rule out the case of $\ell=1$ and then rule out $\ell=d$ in Step 5.

Suppose that $\ell=1$, then $N_{\delta_1}\Sigma_1$ is just a $\delta_1$-tube. By item $(1)$ and since $\TT_Y(p)$ is $\approx 1$-uniform, $|\TT_Y(p)|_\delta\sim \mu(\TT, Y)$, we have 
\[
\mu(\TT, Y) \leq \delta_1^{-\eta_1}  \mu(\TT[2T_{\delta_1}], Y)
\]
for any $T_{\delta_1}\in \TT_{\delta_1}$. If $\delta_1\leq \delta^{1-\kappa/10d}$, we are done by the trivial bound $\mu(\TT[2T_{\delta_1}], Y)\lesssim (\delta/\delta_1)^{-d+1} \lesssim \delta^{-\kappa/10}$. 

Otherwise, $\delta/\delta_1\leq \delta^{\kappa/10d}$ is sufficiently small, apply anisotropic rescaling on $\TT[2T_{\delta_1}]$ to obtain a set  $\TT^{2T_{\delta_1}}$ of $\sim(\delta/\delta_1)$-tubes which is   $(d-1, \delta^{\eta})$-AD regular  and $C_{CW}(\TT^{2T_{\delta_1}})\leq \delta^{-\eta}$ (inherited from $C_{CW}(\TT)\leq \delta^{\eta}$). Therefore we can apply $\operatorname{K}_{d}(\kappa)$  on $\TT^{2T_{\delta_1}}$ to show 
\[
\mu(\TT, Y) \leq \delta_1^{-\eta_1}  \mu(\TT[2T_{\delta_1}], Y)\lesssim \delta_1^{-\eta_1}  (\delta/\delta_1)^{-\kappa-\eps}< \delta^{-\kappa}
\]
when $\delta_1^{-\eta_1+\kappa-\epsilon} \delta^{-\epsilon} \ll 1$ by taking $\delta_1 \leq \delta^{2\eps/\kappa}$.

\medskip

\noindent \textit{Step 2. Define point line configuration and obtain uniform refinement.}

Define $X\subset \Omega_d$ as the  set of $\delta$-separated  triples $(x,t, \theta)\in [-1,1]^{d-1}\times [-1, 1]\times [-1, 1]^{d-1}$ such that $(x,t)\in U(\TT, Y)^{(\d)}\cap (\delta\mathbb{Z})^{d}$ and $\theta=\theta(T)$ for some $T\in \TT_{Y_0}(p)$ with $p=(x,t)$. By Observation~\ref{obs:pigeonhole-slope}, we can assume all pairs $(p, T), T\in \TT_Y(p)$ can be represented in $X$. 

Apply Lemma~\ref{lem:uniform-subset} to $X$ with  $\Delta_j =\delta^{j\eta}$, $j=1, \dots, \eta^{-1}$ to find a subset $X'\subset X$, $|X'|_{\delta}\gtrapprox |X|_{\delta}$ and $X'$ is $\approx 1$-uniform on scales $\Delta_j$, $j=1\dots, \eta^{-1}$, which means for all $j \leq i+k$ and any $\omega_0=(x_0, t_0,  \theta_0)\in X'$, we have 
\[
|X' \cap R_{\Delta_i, \Delta_j, \Delta_k}(\omega_0)| \approx M_{X'}(\Delta_i, \Delta_j, \Delta_k):= \max_{R \in \mc R_{\Delta_i, \Delta_j, \Delta_k}} |X' \cap R|.
%
\]

\medskip

\noindent \textit{Step 3. Obtain a refined shading $Y'$ and show it is uniform.}
For each $T\in \TT$, define shading $Y'$ on $T$ as the union of  dyadic $\delta$-cubes containing points $p$ for which  $(p, \theta(T))\in X'$. 
Let $\TT'$ be the set of $T\in \TT$ such that there exists $p$ such that  $(p, T)\in X'$. The uniformity of $X'$ yields a lot of uniformity information on $(\TT', Y')$.

For any $a\in [\delta, 1]\cap \delta^{\eta\mathbb{Z}}$ and $T\in \TT', x\in Y'(T)$, 
\[
|\TT'[T^a]_{Y'}(x)|_{\delta}\approx  \M_{X'}(\delta, \delta,a).
\]

Observe that for any $p\in Y'(T)$ and $b\in (\delta,1)$, 
\begin{equation*}\label{eq: shadingsegment3}
|Y'(T)\cap B(p,b)|_{\delta} \approx  \M_{X'}(b, \delta, \delta).
\end{equation*}
Therefore $Y'(T)$ is $\approx \delta^{\eta}$-uniform. 

By the definition of $\TT'$ and $Y'$, for any $T\in \TT$, if $Y'(T)\neq \emptyset$, then $T\in \TT'$. For any $\varrho \in [\delta, 1]\cap \delta^{\eta\mathbb{Z}}$, define $\TT_\varrho'= \cD_\varrho (\TT')$.  Then for any $T_\varrho\in \TT_\varrho'$, 
\[
\sum_{T\in \TT'[T_\varrho]} |Y'(T)|_{\delta} \approx \M_{X'}(1, \varrho, \varrho).
\]
Since $Y'$ is a $\gtrapprox 1$-refinement of $Y_0$, which is a $\gtrapprox \delta^{\eta} \delta_\ell^{\eta_\ell}$-dense shading, and $\TT$ is $(d-1, \delta^{-\eta})$-AD-regular,   we have 
\begin{equation*}\label{eq: TT'AD2}
 \delta^\eta \delta_\ell^{\eta_\ell}\delta^{-1}(\delta/\varrho)^{-d+1} \lessapprox   \sum_{T\in \TT'[T_\varrho]} |Y'(T)|_{\delta} \lessapprox \delta^{-1-\eta} (\delta/\varrho)^{-d+1}. 
\end{equation*}
In particular, $\TT'$ is $(d-1, \delta_\ell^{-2\eta_\ell})$-AD-regular since we choose $\eta>0$ and $\varrho_\ell$ sufficiently small that $\delta^\eta \geq \varrho_\ell^{\eta_\ell} \geq  \delta_\ell^{\eta_\ell}$.

For any $T_\varrho\in \TT_\varrho'$, define  
\[
Y'(T_\varrho)= \cD_\varrho(U(\TT'[T_\varrho], Y')).
\] Then for any $b\in [\varrho, 1]\cap \delta^{\eta\mathbb{N}}$,  
\begin{equation}
     |Y'(T_\varrho)|_b \approx \frac{\sum_{T\in \TT'[T_\varrho]}|Y'(T)|_{\delta}}{\M_{X'}(b, \varrho, \varrho)} \approx \frac{ \M_{X'}(1,\varrho, \varrho)}{\M_{X'}(b, \varrho, \varrho)}. 
\end{equation}
This shows that $Y'(T_\varrho)$ is $\approx \delta^{\eta}$-uniform for any $\varrho \in [\delta, 1]$.

\medskip

\noindent \textit{Step 4. Apply $\operatorname{K}_{d}(\kappa)$ on $(\TT[T_\varrho], Y')$ and $(\TT_\varrho, Y')$ to show extremal multiplicities.}

Recall  $\eps_1\in (0, \eta_d)$ is a small parameter to be determined. By $\operatorname{K}_{d}(\kappa)$, the following inequalities hold for $\eta>0$ sufficiently small and $\varrho>0$ sufficiently small (in particular, $\varrho^{\eps_1/2}\leq \delta^\eta$), 
\begin{equation}
    |U(\TT_\varrho, Y)|_{\varrho} \geq  \varrho^{\eps_1/2} \varrho^{-d+\kappa}\geq \varrho^{\eps_1} \delta^{-\eta} |Y(T_\varrho)|_\varrho \varrho^{-d+1+\kappa}, 
\end{equation}
or equivalently 
\begin{equation}\label{eq: muthicktubes0K}
\mu(\TT_\varrho, Y) \lessapprox \varrho^{-\eps_1-\kappa}.
\end{equation}
Apply $\operatorname{K}_d(\kappa)$ to $\TT^{T_\varrho}$ when  $\delta/\varrho>0$ is sufficiently small, we get
\begin{equation}\label{eq: muthintubes0K}
    \mu(\TT[T_\varrho], Y)\lessapprox (\delta/\varrho)^{-\eps_1-\kappa}.
\end{equation}
Since $Y'$ is a $\gtrapprox \delta_\ell^{\eta_\ell}$-refinement of $Y$, 
\[
\mu(\TT, Y')\gtrapprox \delta_\ell^{\eta_\ell} \mu(\TT, Y)\gtrsim \delta^{-\kappa+\eta}\delta_\ell^{\eta_\ell}.
\]
On the other hand, since $\mu(\TT_\varrho, Y)\approx |\TT_{\varrho, Y}(x)|$ for each $x\in U(\TT, Y)$ and $|\TT_{\varrho, Y'}(x)|\leq |\TT_{\varrho, Y}(x)|$ pointwise, we have 
\begin{equation}
    \mu(\TT_\varrho, Y')\leq \mu(\TT_\varrho, Y)\lessapprox \varrho^{-\kappa-\eps_1}.
\end{equation}

Since $\mu(\TT, Y')\lesssim \mu(\TT[T_\varrho], Y') \mu(\TT_\varrho, Y')$, combining with \eqref{eq: muthicktubes0K} and  \eqref{eq: muthintubes0K}, we have 
\begin{equation}\label{eq: muthicktubesK}
     \delta^{\epsilon_1}\delta_\ell^{2\eta_\ell} \varrho^{-\kappa} \leq \mu(\TT_\varrho, Y') \leq \varrho^{-\eps_1-\eta-\kappa}
\end{equation}
and 
\begin{equation}\label{eq: muthintubesK}
        \delta^{\epsilon_1}\delta_\ell^{2\eta_\ell} (\delta/\varrho)^{-\kappa} \leq \mu(\TT[T_\varrho], Y')\leq (\delta/\varrho)^{-\eps_1-\eta-\kappa}.
\end{equation}
Recall that by the definition of $Y'$ and $\TT'$, if $T\in \TT$ has nonempty shading $Y'(T)$, then $T\in \TT'$. Choose $\eps_1$ sufficiently small so that $\delta^{\epsilon_1} \geq \delta_\ell^{\eta_\ell}$. 
Then  by \eqref{eq: muthicktubesK} and \eqref{eq: muthintubesK}, for any $\delta\leq r<\varrho  \leq 1$, and $T_\varrho\in \TT_\varrho'$, 
\begin{equation}\label{eq: mutubesK}
   \delta_\ell^{6\eta_\ell}  (r/\varrho)^{-\kappa} \leq \mu(\TT_r'[T_\varrho], Y')=\mu(\TT_r[T_\varrho], Y') \leq    \delta_\ell^{-6\eta_\ell}  (r/\varrho)^{-\kappa}. 
\end{equation}

\medskip

\noindent \textit{Step 5. Finding grains using extremal multiplicities.}

Let $E_0=U(\TT', Y')$. For each $p\in E_0$, let $\Theta_p$ be the set of directions $\theta(T)$ for $T\in \TT_{Y'}(p)$ and let $\Lambda_p\subset \Theta_p^\ell:= \Theta_p \times \cdots \Theta_p$  denote the set of $\ell$-tuples $(\theta_1, \dots, \theta_\ell)$ such that 
\begin{equation}\label{eq: l-linear}
|v_1\wedge \cdots \wedge v_\ell|\gtrsim \delta_{\ell}^{C_0\eta_\ell/\eta_{\ell-1}}
\end{equation}
with $v_j=(\theta_j, 1)$ and some fixed constant $C_0\ge 1$ (depending on $d$ only). By $(2)$ of  the definition of $\ell$-linear points, 
since $p$ is $\ell$-linear but not $(\ell-1)$-linear, 
if $(\theta_1, \dots, \theta_\ell)\in \Theta_p^{\ell}\setminus \Lambda_p$, then the corresponding $T_j\in \TT_Y(p)$ such that $\theta(T_j)=\theta_j$ lies in  $\delta_\ell^{2\eta_\ell/\eta_{\ell-1}}$-neighborhood of an $(\ell-1)$-plane $\Sigma$ and 
\[
|\TT_Y(p)[N_{\delta_\ell^{2\eta_\ell/\eta_{\ell-1}}}\Sigma]|_\delta \leq \delta_\ell^{2\eta_\ell} |\TT_Y(p)|_\delta. 
\]
Here we choose $\varrho_\ell \ll \varrho_{\ell-1}$ so that $\delta_\ell^{\eta_\ell/\eta_{\ell-1}}\leq \varrho_{\ell-1}$. Therefore
\begin{equation} \label{eq: l-lineardominates}
|\Theta_p^\ell\setminus \Lambda_p|_{\delta}\leq \delta_\ell^{\eta_\ell}|\Theta_p|_{\delta}^\ell.
\end{equation}
By \eqref{eq: mutubesK}, for any $\varrho \in [\delta_\ell, 1]$, 
\begin{equation}\label{eq: l-linear-rho}
    |\Lambda_p|_{\varrho} \gtrsim \frac{|\Theta_p|_\delta^\ell}{\mu(\TT[T_\varrho], Y')^\ell} \gtrsim  \delta_\ell^{12\eta_\ell} \varrho^{-\kappa \ell } \gtrsim \delta_\ell^{6(\ell+3)\eta_\ell} |\TT_{\varrho, Y'}(\bp)|_{\varrho}^\ell
\end{equation}
where $\bp=p^{\varrho}$ is the $\varrho$-dyadic cube containing $p$.

Let $\Delta_j =\delta^{\eta j}\in [\delta_\ell, 1]$ for $j=1, \dots, J$ with $\Delta_J \sim  \delta_\ell$ up to a factor of $\delta^\eta$.  We are going to show, after refinements, any $\Delta_j$-cube $\bp_j\in \cD_{\Delta_j}E_0$, $\bp_j\cap N_{\Delta_j^2}E_0$ is approximately a union of parallel $\ell$-dimensional slabs (the $\Delta_j^2$-neighborhood of $\ell$-planes intersecting $\bp_j$). In addition, the slabs at different scales are aligned, which implies the slope function $f$ we are going to define is dyadic Lipchitz.  

Start with $j=J$. For any $\Delta_j$-cube $\bp_j\in \cD_{\Delta_j}E_0$, let $\Lambda_{\bp_j}\subset \underset{p\in \bp_j\cap E_{0}}{\bigcup} \Lambda_p$ be the set of $\ell$-tuples  $\vec{\theta}=(\theta_1, \dots, \theta_\ell)$ such that 
\begin{equation} \label{eq: densityofcubes}
|E_{\bp_j, \vec{\theta}}|_{\delta_\ell^2}\gtrsim \delta_\ell^{20\ell\eta_\ell} |\bp_j\cap E_0|_{\delta_\ell^2}
\end{equation}
where  $E_{\bp_j, \vec{\theta}}\subset \bp_j\cap E_{0}$ is the set of points $p$ such that $\vec{\theta}\in \Lambda_p^{\Delta_j}$ (the union of dyadic $\Delta_j$-cubes in $[-1, 1]^{(d-1)\ell}$ intersecting $\Lambda_p$). Once we have pigeonholed $\ell$-linear points and defined $Y_0$, the scale $\delta_\ell^2$ is the smallest scale we are going to work with for the rest of the proof.

By double counting using \eqref{eq: mutubesK}, \eqref{eq: l-lineardominates} and \eqref{eq: l-linear-rho},
\begin{equation} \label{eq: numberktuplesK}
|\Lambda_{\bp_j}|_{\Delta_j} \gtrapprox \delta_\ell^{6(\ell+3)\eta_\ell}|\TT_{\Delta_j, Y'}(\bp_j)|_{\Delta_j}^\ell.
\end{equation}

Next step  we replace $j$ by $j-1$. 
Suppose we have defined $\Lambda_{\bp_{j+1}}$, then for any $\Delta_j$-cube $\bp_j\in \cD_{\Delta_j}E_0$, let $\Lambda_{\bp_j}\subset \underset{\bp_{j+1}\in \bp_j\cap \cD_{\Delta_{j+1}}E_0}{\bigcup} \Lambda_{\bp_{j+1}}$  be the set of $\ell$-tuples $\vec{\theta}$ such that  
\eqref{eq: densityofcubes} still holds with 
\begin{equation}\label{eq: defineEpjtheta}
E_{\bp_j, \vec{\theta}} =\bigcup_{\bp_{j+1}\in \bp_j\cap \cD_{\Delta_{j+1}}E_0}  \, \, \bigcup_{\vec{\theta}': \vec{\theta}'^{\Delta_j}=\vec{\theta}^{\Delta_j}} E_{\bp_{j+1}, \vec{\theta}'}.
\end{equation} 
Then \eqref{eq: numberktuplesK} remains true. Here the condition $\vec{\theta}'^{\Delta_j}=\vec{\theta}^{\Delta_j}$ is what makes our slope function $f$ defined later dyadic Lipchitz. 
Repeat the process until $j=1$.

Now start with $j=1$. For each $\bp_j\in \cD_{\Delta_j} E_0$, choose a $\ell$-tuple $\vec{\theta}_j =(\theta_{j1}, \cdots, \theta_{j\ell})\in \Lambda_{\bp_j}$. Let $E^0_{\bp_j, \vec{\theta_j}}=E_{\bp_j, \vec{\theta_j}}$ and for $i= 1, \dots, \ell$, let  $
\mathbb{L}_{bad, i}^j$ be the set of $\Delta_j^2 \times \cdots \times \Delta_j^2\times \Delta_j$-segments $L$ in direction $\theta_{ji}$ with 
\[
|Y'(L)\cap E_{\bp_j, \vec{\theta_j}}^{i-1}|_{\delta_\ell^2} \leq \delta_\ell^{30\ell \eta_\ell} |Y'(L)|_{\delta_\ell^2}, \quad \text{where } Y'(L):=\bigcup_{T_{\Delta_j^2}\in \TT_{\Delta_j^2}: \, L \subset T_{\Delta_j^2}} Y'(T_{\Delta_j^2}\cap L).
\]

Write \[ E_{\bp_j, \vec{\theta}_j}^i=E_{\bp_j, \vec{\theta}_j}^{i-1}\setminus  \bigcup_{L\in \mathbb{L}^j_{bad, i}} L, \]
it follows from \eqref{eq: densityofcubes} that 
\begin{equation} \label{eq: densityofcubes1}
|E_{\bp_j, \vec{\theta}_j}^\ell|_{\delta_\ell^2}\gtrsim |E_{\bp_j, \vec{\theta}_j}|_{\delta_\ell^2}.
\end{equation} Write $E_{\bp_j, \vec{\theta}_j}'=E_{\bp_j, \vec{\theta}_j}^\ell$.

We claim that $E_{\bp_j, \vec{\theta}_j}'^{\Delta_j^2}$ is  roughly a union of parallel $\ell$-slabs, where a $\ell$-slab $\Sigma_{\bp_j}$ means the $\Delta_j^2$-neighborhood of an $\ell$-plane intersecting $\bp_j$. More precisely, 
\begin{equation}\label{eq: grain1}
    E_{\bp_j, \vec{\theta}_j}'\subset \bigsqcup \Sigma_{\bp_j}, 
\end{equation}
where the union is over parallel $\ell$-slabs and for each such $\ell$-slab $\Sigma_{\bp_j}$, 
\begin{equation}\label{eq: grain2}
    |E_{\bp_j, \vec{\theta}_j}'\cap \Sigma_{\bp_j}|_{\Delta_j^2} \gtrsim \delta_\ell^{O(\eta_\ell/\eta_{\ell-1})} |\Sigma_{\bp_j}|_{\Delta_j^2}.
\end{equation}

To see this, for $i=1, \dots, \ell$, let $\Pi_{i}: \mathbb{R}^d \rightarrow \langle v_{j1}, \dots, v_{ji}\rangle^{\perp}$ be the orthogonal projection from $\mathbb{R}^d$ to $\mathbb{R}^{d-i}$, where $v_{ji}=(\theta_{ji}, 1)$. By the definition of $\mathbb{L}_{bad, i}^j$ and $E_{\bp_j, \vec{\theta_j}}^i$, for each $y\in \Pi_{i}(E_{\bp_j, \vec{\theta_j}}^i)$,  since $|v_{j1}\wedge \cdots \wedge v_{ji} | \geq \delta_\ell^{O(\eta_\ell/\eta_{\ell-1})}$ and $Y'(T_{\Delta_j^2})$ is $\approx \delta^{-\eta}$-uniform
\begin{equation}
    |\Pi_i^{-1}(y)\cap E_{\bp_j, \vec{\theta_j}}^i|_{\Delta_j^2} \geq \delta_\ell^{O( \eta_\ell/\eta_{\ell-1})} \Delta_j^{-i}. 
\end{equation}
We obtain \eqref{eq: grain2} by letting $i=\ell$ and $\Sigma_{\bp_j}$ be the $\ell$-slabs  of the form $N_{\Delta_j^2}(\Pi_\ell^{-1}(y))\cap \bp_j$), which are parallel to the $\ell$-plane spanned by $v_{j1}, \dots, v_{j\ell}$.

Now we proceed to $j+1$. For each $\bp_{j+1}\in \cD_{\Delta_{j+1}} E_{\bp_j, \vec{\theta}_j}'$ with 
\begin{equation}\label{eq: densityofcubesj}
|\bp_{j+1}\cap E_{\bp_j, \vec{\theta}_j}'|_{\delta_\ell^2}\gtrsim \delta_\ell^{20\ell\eta_\ell}|\bp_{j+1}\cap E_0|_{\delta_\ell^2}, 
\end{equation}
by the definition of $E_{\bp_j, \vec{\theta}_j}$, there exists $\vec{\theta}_{j+1}$ satisfying $\vec{\theta}_{j+1}^{\Delta_j}=\vec{\theta}_j^{\Delta_j}$ and 
\[E_{\bp_{j+1}, \vec{\theta}_{j+1}}\subset E_{\bp_j, \vec{\theta}_j}', \quad 
|E_{\bp_{j+1}, \vec{\theta}_{j+1}}|_{\delta_\ell^2} \gtrsim \delta_\ell^{20\ell \eta_\ell} |\bp_{j+1}\cap E_0|_{\delta_\ell^2}
\]
since \eqref{eq: densityofcubes} and \eqref{eq: densityofcubes1} imply that the set of $\bp_{j+1}$ satisfying \eqref{eq: densityofcubesj} covers a $\gtrsim 1$-fraction of $E_{\bp_j, \vec{\theta}_j}'$.
Define $\mathbb{L}^{j+1}_{bad, i}, i =1, \dots, \ell$ and $E_{\bp_{j+1}, \vec{\theta}_{j+1}}'$ similarly. 
Iterate the process until $j=J$. 

Define 
\begin{equation}\label{eq: E1}
    E_1= \bigcup_{\bp_1\in \cD_{\Delta_1}E_0}  \,\, \bigcup_{\bp_2\in \cD_{\Delta_2}E_{\bp_1, \vec{\theta}_1}'} \cdots \bigcup_{\bp_{J}\in \cD_{\Delta_{J}}E_{\bp_{J-1}, \vec{\theta}_{J-1}}'} E_{\bp_J, \vec{\theta}_J}'.
\end{equation}

Then 
\begin{equation} \label{eq: densityE1}
|E_1|_{\delta_\ell^2}\gtrsim \delta_\ell^{20\ell\eta_\ell}|E_0|_{\delta_\ell^2} \text{ and both  } E_1 \text{ and } E_0 \text{ are } \approx \delta_\ell^{20\ell\eta_\ell}\text{-uniform}. 
\end{equation}

Recall the definition of $Y_0$ in Step 1, since $Y'$ is a refinement of $Y_0$, for any $p\in E_{\bp_J, \vec{\theta_J}}'$, $p$ is $\ell$-linear and $\TT_{Y_0}(p)\subset N_{\delta_\ell}(\Sigma_\ell)$ for some $\ell$-dimensional plane $\Sigma_\ell$, we know that the tubes $T_i$ with $\theta(T_i)= \theta_{Ji}, \vec{\theta}_J=(\theta_{J1}, \dots, \theta_{J\ell})$  are contained in $N_{\delta_\ell}\Sigma_\ell$. In addition, the tubes $\{T_i\}$ span roughly $\Sigma_\ell$ up to error $\delta_\ell^{O(\eta_\ell/\eta_{\ell-1})}$. We are going to use this observation to show Item $(4)$ of Definition~\ref{def:kakeya-configuration}.

Define $Y''(T)=Y'(T)\cap E_1$, then $Y''$ is a $\gtrsim  \delta_\ell^{30\ell\eta_\ell}$-dense shading of $\TT'$.  

Now we rule out the case $\ell=d$. By \eqref{eq: grain2}, if $\ell=d$, then $\Sigma_{\bp_J}$ is the $\delta_d=\delta^{1/2}$-cube spanned by vectors $v_{j1}, \dots v_{jd}$. Therefore
\[
|U(\TT, Y)|\geq |U(\TT', Y')|\geq \delta_d^{O(\eta_{d-1})} |U(\TT'_{\delta_d}, Y')|\geq \delta^{O(\eta_{d-1})} \delta_{d}^{\eps} \delta^{\kappa/2}\gg \delta^{\kappa}
\]
which yields a contradiction.

When $\ell \leq d-1$, to ease the notation, let $\tilde \delta= \delta_\ell^{1/2}$. 


Let $T_{\tilde\delta}$ be a tube such that $Y''$ is a $\gtrsim \tilde\delta^{30\ell\eta_\ell}$-dense shading on $\TT[T_{\tilde\delta}]$.  After a rotation, we may assume $T_{\tilde\delta}$ is pointing in direction of the last coordinate $z$. Let 
$\phi: T_{\tilde\delta}\rightarrow T_1$ be the anisotropic rescaling of $T_{\tilde\delta}$ to the unit tube and let 
$\tilde\TT=\phi(\TT_{\tilde\delta^2}[T_{\tilde\delta}])$ and $\tilde Y'(\phi(T_{\tilde\delta^2}))$ be the union of $\tilde\delta$-cubes in $\cD_{\tilde\delta} \phi(Y''(T_{\tilde\delta^2}))$. Similarly, let $\tilde Y(\phi(T_{\tilde\delta^2}))$ be the union of $\tilde\delta$-cubes in  $\cD_{\tilde\delta}\phi(Y'(T_{\tilde\delta^2}))$.

When $\ell\leq d-1$, each $\ell$-slab  in \eqref{eq: grain1} with $j=J$ 
becomes a horizontal slice.

\medskip

\noindent \textit{Step 6. AD-regularity of horizontal slices.}

Let $Z'$ be the set of $z$ such that $\mathbb{R}^{d-1}\times \{z\} \cap U(\tilde \TT, \tilde Y')\neq \emptyset$. 
Then $E'=U(\tilde\TT, \tilde Y')$ is of the form 
\[E'=\cup_{z\in Z'} E_z', \quad E_z' = A_z'\times \{z\}\subset \mathbb{R}^{d-1}\times \{z\}.
\]
Let $E=U(\tilde\TT, \tilde Y)$ and let  $E_z= U(\tilde\TT, \tilde Y)\cap \mathbb{R}^{d-1}\times \{z\}$. Then $E_z'\subset E_z$. 
After a rigid motion, \eqref{eq: grain1} shows we can write $E_z'=\{ (x, y+f(z)x, z), y\in Y_z, x\in X_{y,z}\}$ with  $f: Z'\rightarrow \operatorname{Mat}(d-\ell, \ell-1)$,  $X_{y,z}\subset [-1, 1]^{\ell-1}$ and \eqref{eq: grain2} says  $|X_{y,z}|_{\tilde\delta} \geq \tilde\delta^{O(\eta_{\ell}/\eta_{\ell-1})-\ell+1}$. In addition, 
for any  $\varrho \in [\tilde{\delta}, 1]\cap \delta^{\eta\mathbb{Z}}$ and any $z_1, z_2\in Z'$ with  $z_1^\varrho=z_2^\varrho$, \eqref{eq: defineEpjtheta} shows \begin{equation}\label{eq: Lipschitzf}
    |f(z_1)-f(z_2)| \lesssim  \varrho. 
\end{equation}

To see Item $(4)$ of Definition~\ref{def:kakeya-configuration} when $\ell\leq d-1$,  for each $p\in E_1\subset U(\TT, Y')$, there is an $\ell$-plane $\Sigma_\ell$ such that $\TT_{Y'}(p)$ is contained in $N_{\delta_\ell}\Sigma_\ell$. By the definition of $E_1$ in \eqref{eq: E1}, $p\in E_{\bp_j, \vec{\theta}_j}'$ for some $\bp_j$ and $\Delta_j=\tilde\delta$. Write $\vec{\theta}_j = (\theta_{j1}, \dots, \theta_{j\ell})$ and let $v_i=(\theta_{ji}, 1)$, then $\{v_i\}_{i=1}^\ell$ spans an $\ell$-plane parallel to $\Sigma_\ell$ and $\Sigma_{\bp_j}$ in \eqref{eq: grain1}  up to angle difference $\lesssim \delta_\ell^{1-O(\eta_\ell/\eta_{\ell-1})}$. Since  we take  $\tilde\delta=\delta_{\ell}^{1/2}$ and   $T_{\tilde\delta}$ points in the direction of  the last coordinate $\phi: T_{\tilde\delta} \rightarrow T_1$ is the anisotropic rescaling that rescales all other coordinates by $\tilde\delta^{-1}$, the tubes in $\TT[T_{\tilde\delta}]_{Y'}(p)$, which are contained in $N_{\delta_\ell}\Sigma_\ell$, rescale to $\tilde\delta$-tubes contained in the $\tilde\delta$-neighborhood of an $\ell$-plane $\Sigma_\ell(p)$ that is parallel to the horizontal slice $\phi(\Sigma_{\bp_j})$ up to angle difference $\tilde\delta^{1-O(\eta_\ell/\eta_{\ell-1})}$.

Write $E_z=A_z\times \{z\}$ for some $A_z\subset [-1, 1]^{d-1}$. We first show $A_z$ is a $(d-1-\kappa, \tilde\delta^{-\eta_\ell})$-AD regular set, then we show that  for a typical $z\in Z'$, $A_{z}'$ is a significant fraction of $A_z$ and thus $A_{z}'$ contains a $(d-1-\kappa, \tilde\delta^{-O(\eta_\ell)})$-AD-regular set.

To see this, for any $z\in Z, \omega=(x,y)\in A_z$ and $r\in [\tilde\delta, 1]\cap \delta^{\eta\mathbb{N}}$, $(A_z\cap B(\omega, r))\times \{z\}$ is a subset of $E_z$, if we take its $\tilde\delta$-neighborhood and apply $\phi^{-1}$, then it is a union of $\tilde\delta^2$-tube segments of length $\tilde\delta$ contained in an $r\tilde\delta$-tube segment of length $\tilde\delta$. Therefore 
\begin{equation} \label{eq: AzregularK}
|A_z\cap B(\omega, r)|_{\tilde\delta} 
\approx \frac{ \M_{ X'}(\tilde\delta, r\tilde\delta, \tilde\delta) }{\M_{ X'}(\tilde\delta, \tilde\delta^2, \tilde\delta)}.
\end{equation}
(recall that $A_z$ is obtained from rescaling $(\TT_{\tilde\d^2}, Y'(T_{\tilde \d^2}))$ and $Y'$ is defined using the $\approx \d^\eta$-uniform configuration $X'$). More generally, for any $\varrho \leq r\in [\tilde\delta, 1]$, we have 
\begin{equation}\label{eq: AzregularKrho}
    |A_z\cap B(\omega, r)|_\varrho \approx \frac{\M_{X'}(\tilde\delta, r\tilde\delta, \tilde\delta)}{\M_{X'}(\tilde\delta, \varrho \tilde\delta, \tilde\delta)}. 
\end{equation}

We are not going to use the RHS of \eqref{eq: AzregularKrho} to calculate $|A_z\cap B(\omega,r)|_{\varrho}$, what is important for us is the fact that for any $\varrho\leq r\in [\tilde\delta, 1]\cap \delta^{\eta\mathbb{N}}$ \eqref{eq: AzregularK} holds  uniformly for   $z\in Z, \omega\in A_z$.

 Recall $\TT'$ is $(d-1, \tilde\delta^{-O(\tilde\eta)})$-AD-regular and $|Y'(T_\varrho)|_\varrho$ is approximately the same for each $T_\varrho \in \TT_\varrho'$, $\varrho \in [\tilde\delta, 1]\cap \delta^{\eta\mathbb{N}}$.  By \eqref{eq: mutubesK}  and $\mu(\tilde \TT, \tilde Y)=\mu(\TT_{\tilde\delta^2}'[T_{\tilde\delta}], Y')$ since $|A_z|_{\tilde\delta}$ is approximately the same for each $z\in Z$ (apply \eqref{eq: AzregularK} with $r=1$), 
 \begin{equation}\label{eq: AzupperK}
 |A_z|_{\tilde\delta} \approx \frac{|E|_{\tilde\delta}}{|Z|_{\tilde\delta}} \lesssim   \frac{|E|_{\tilde\delta}}{|\tilde Y(\tilde T)|_{\tilde\delta}}\approx  \frac{ \sum_{\tilde T\in \tilde\TT} |\tilde Y(\tilde T)|_{\tilde\delta} }{ \mu(\tilde \TT, \tilde Y)  |\tilde Y (\tilde T)|_{\tilde\delta}} \lessapprox  \tilde\delta^{-O(\eta_\ell)-d+1+\kappa}. 
 \end{equation}

Since $|\tilde Y (\tilde T)|_{\tilde\delta} \gtrsim \tilde\delta^{O(\eta_\ell)-1} \gtrsim  \tilde\delta^{O(\eta_\ell)} |Z|_{\tilde\delta}$, 
 \begin{equation}
      |A_z|_{\tilde\delta} \approx \frac{|E|_{\tilde\delta}}{|Z|_{\tilde\delta}} \gtrapprox \tilde\delta^{O(\eta_\ell)}  \frac{|E|_{\tilde\delta}}{|\tilde Y(\tilde T)|_{\tilde\delta}}\approx \tilde\delta^{O(\eta_\ell)} \frac{ \sum_{\tilde T\in \tilde\TT} |\tilde Y(\tilde T)|_{\tilde\delta} }{ \mu(\tilde \TT, \tilde Y)  |\tilde Y (\tilde T)|_{\tilde\delta}} \gtrapprox \tilde\delta^{O(\eta_\ell)-d+1+\kappa}. 
 \end{equation}

Similarly, for any $r\in [\tilde\delta, 1]$, 
\begin{equation}\label{eq: lowerbdAzrK}
|A_z|_{r} \gtrapprox r \, |E|_r \gtrsim r \, \frac{\sum_{\tilde T \in \tilde\TT_r}|\tilde Y(\tilde T)|_r}{\mu(\tilde\TT_r, \tilde Y)}  
\gtrapprox \tilde\delta^{O(\eta_\ell)} r^{-d+1+\kappa}. 
\end{equation}

Therefore by \eqref{eq: AzregularKrho} we have the following upper bound 
\begin{equation}\label{eq: upperAzxrK}
    |A_z\cap B(\omega, r)|_{\tilde\delta} \lessapprox \frac{|A_z|_{\tilde\delta}}{|A_z|_r} \lessapprox \tilde\delta^{-O(\eta_\ell)} (\tilde\delta/r)^{-d+1+\kappa}.
\end{equation}

For any $T_r\in \tilde \TT_r$, since $|Z|_{\tilde\delta/r} \lessapprox \tilde\delta^{-O(\eta_\ell)} |\tilde Y(\tilde T)|_{\tilde\delta/r}$ and $|A_z\cap B(\omega, r)|_{\tilde\delta}$ is approximately the same for all $z\in Z, \omega\in A_z$ and $r\in [\tilde\delta, 1]\cap \delta^{\eta N}$ (in particular, when averaging over $z\in Z$), 
\begin{equation}\label{eq: lowerAzxrK}
    |A_z\cap B(\omega,r)|_{\tilde\delta} \gtrsim (\tilde\delta/r) \, U(\tilde\TT[T_r], \tilde Y) \gtrsim  (\tilde\delta/r) \, \frac{\sum_{\tilde T\in \tilde \TT[T_r]} |\tilde Y (\tilde T)|_{\tilde\delta/r}}{\mu(\tilde \TT[T_r], \tilde Y) } \gtrapprox \tilde\delta^{O(\eta_\ell)} (\tilde\delta/r)^{-d+1+\kappa}. 
\end{equation}
Combining with \eqref{eq: lowerAzxrK} and \eqref{eq: upperAzxrK}, we have 
\begin{equation}\label{eq: AzxrK}
    \tilde\delta^{O(\eta_\ell)}  (\tilde\delta/r)^{-d+1+\kappa}\lessapprox     |A_z\cap B(\omega,r)|_{\tilde\delta} \lessapprox \tilde\delta^{-O(\eta_\ell)} (\tilde\delta/r)^{-d+1+\kappa}. 
\end{equation}
This shows that $A_z$ is $(d-1-\kappa, \tilde\delta^{-O(\eta_\ell)})$-AD regular. 

Since $A_z'\subset A_z$ and $|E'|_{\tilde\delta} \gtrsim \tilde\delta^{O(\eta_\ell)}|E|$, $|Z'|_{\tilde\delta} \gtrsim \tilde\delta^{O(\eta_\ell)}|Z|_{\tilde\delta}$, we can choose a $\approx 1$-uniform set $Z''\subset Z'$, $|Z''|_{\tilde\delta}\gtrapprox |Z'|_{\tilde\delta}$ and for each $z\in Z''$, $|A_z'|_{\tilde\delta} \gtrsim \tilde\delta^{O(\eta_\ell)}|A_z|_{\tilde\delta}$ and therefore we can choose a $(d-1-\kappa, \tilde\delta^{-O(\eta_\ell)})$-AD-regular set $A_z''\subset A_z'$. In particular, \eqref{eq: AzxrK} holds for $A_z''$ in place of $A_z$ with slightly worst constants in $O(\eta_\ell)$.

The next step is to pigeonhole so that after a refinement, $Y_z$ is $(d-\ell-\kappa, \tilde\delta^{-O(\eta_\ell/\eta_{\ell-1})})$-AD-regular. 

Recall that by \eqref{eq: grain1} and \eqref{eq: grain2},
for each $z\in Z'$,  $y\in Y_z$, $|X_{y,z}|_{\tilde\delta} \geq \tilde\delta^{O(\eta_\ell/\eta_{\ell-1})-\ell+1}$.  We also have trivially $|X_{y,z}|_{\tilde\delta} \lesssim \tilde\delta^{-\ell+1}$. 
After pigeonholing, we can find a $\tilde\delta^{O(\eta_\ell)}$-uniform set  $Y_z''\subset Y_z$ such that 
$|Y_z''|_{\tilde\delta} \gtrsim |Y_z|_{\tilde\delta}$ and for each $y\in Y_z''$, $X_{y,z}''\subset X_{y,z}$, $|X_{y,z}''|_{\tilde\delta} \gtrapprox \tilde\delta^{O(\eta_\ell/\eta_{\ell-1})-\ell+1}$ such that 
$A_z''':=\{(x, y+f(z)x): y\in Y_z'', x\in X_{y,z}''\}\subset A_z''$ and $|A_z'''|_{\tilde\delta}\gtrapprox |A_z''|_{\tilde\delta}$.


Apply \eqref{eq: AzxrK} with $r=1$ and $|A_z'''|_{\tilde\delta}\gtrapprox |A_z|_{\tilde\delta}$, and by the upper and lower bounds of $|X_{y,z}''|_{\tilde\delta}$, 
\begin{equation}\label{eq: Yz'}
 \tilde\delta^{O(\eta_\ell/\eta_{\ell-1})} \tilde\delta^{-d+\ell+\kappa} \lesssim   |Y_z''|_{\tilde\delta} \lessapprox \tilde\delta^{-O(\eta_\ell/\eta_{\ell-1})} \tilde\delta^{-d+\ell+\kappa}
\end{equation}

Apply \eqref{eq: AzxrK} with $r\in [\tilde\delta, 1]$, for $\omega$ such that $|A_z''' \cap B(\omega, r)|_{\tilde\delta}\gtrapprox |A_z\cap B(\omega, r)|_{\tilde\delta}$ (which happens for a $\gtrapprox 1$-fraction of $\omega\in A_z'''$),
\begin{equation}
  \left|\bigcup_{y\in Y_z'}X_{y,z}''\bigcap B(\omega, r)\right|_{\tilde\delta} \gtrapprox \tilde\delta^{O(\eta_\ell)} (\tilde\delta/r)^{-d+1+\kappa}. 
\end{equation}
 Since $|X_{y,z}''\cap B(x, r)|_{\tilde\delta}\lesssim (\tilde\delta/r)^{-\ell+1}$ for any $x\in X_{y,z}$, for  $\omega=(x,y)$ we have 
\begin{equation}
    |Y_z''\cap B(y, r)|_{\tilde\delta} \gtrapprox \tilde\delta^{O(\eta_\ell)} (\tilde\delta/r)^{-d+\ell+\kappa}. 
\end{equation}

On the other hand, since $A_z$ is $(d-1-\kappa, \tilde\delta^{-O(\eta_\ell)})$-AD regular and $|A_z'''|_{\tilde\delta}\gtrapprox |A_z|_{\tilde\delta}$,   \eqref{eq: lowerbdAzrK} holds for $A_z'''$ in place of $A_z$, and by $|X_{y,z}''|_r\lesssim r^{-\ell+1}$, 
\begin{equation}
    |Y_z''|_{r} \gtrapprox \tilde\delta^{O(\eta_\ell)} r^{-d+\ell +\kappa}. 
\end{equation}
Since  $Y_z''$ is $\tilde\delta^{O(\eta_\ell)}$-uniform, combine with  \eqref{eq: Yz'} we get
\begin{equation}
        |Y_z''\cap B(y, r)|_{\tilde\delta} \lessapprox \frac{|Y_z''|_{\tilde\delta}}{|Y_z''|_r}\lessapprox \tilde\delta^{-O(\eta_\ell/\eta_{\ell-1})} (\tilde\delta/r)^{-d+\ell+\kappa}. 
\end{equation}
Therefore $Y_z''$ is $(d-\ell-\kappa, \tilde\delta^{O(\eta_\ell/\eta_{\ell-1})})$-AD-regular.  

Let $\tilde\eta = C\eta_\ell/\eta_{\ell-1}$ where $C$ is a sufficiently large constant that dominates all implicit constants in previous $O(\eta_\ell/\eta_{\ell-1})$. 

It follows from Observation~\ref{obs:wolff-axiom-rescaling} that $C_{CW}(\tilde\TT)\leq \tilde\delta^{-\tilde\eta}$. This finishes the construction of a global $(\ell-1)$-grains $(\d, \eta, \kappa)$-Kakeya configuration.


\end{proof}

Now we prove the stronger version, Proposition~\ref{prop: structure1}, using Proposition~\ref{prop: structure2}.

\begin{proof} By Proposition~\ref{prop: structure2}, for any $\delta_0, \eta_0>0$ we can find $\eta\leq \eta_0$ and $\delta<\delta_0$ such that there exists a global $(\ell-1)$-grains $(\delta, \eta, \kappa)$-Kakeya configuration whose set $\TT$ of $\delta$-tubes satisfies $C_{CW}(\TT)\leq \delta^{-\eta}$. The global $(\ell-1)$-grains $(\delta, \eta, \kappa)$-Kakeya set is organized into a union of horizontal $(\ell-1)$-dimensional slices. These slices themselves might be organized into a union of higher dimensional slices, this is where the $k$-planar comes from.  

\medskip
\noindent \textit{Step 1 Apply structure of AD-regular set.}
Suppose $d-\ell > 1$, otherwise $d-\ell=1$ and we proceed to the next step.

Let $0<\eta_1<\eps_1<\eta_2< \eps_2<\cdots \eta_{d-\ell}<\eps_{d-\ell}$ be small parameters such that $\eta_j< g(\eps_j)$ and $\eps_j \ll \eta_{j+1}$ to be determined.  

For each $z\in Z$, recall $Y_z\subset[-1, 1]^{d-\ell}$ is a $(\delta, d-\ell-\kappa, \delta^{-\eta})$-AD regular and we would like to find product structure within it. To do so,  apply Lemma~\ref{lem:finding-parallel-rectangles} to 
$Y_z\subset [-1, 1]^{d-\ell}$ with $\zeta=\eta_1/2$ and obtain $\gamma_1 \in [0, 1]$ and two scales $\varrho < \tau \in [\delta, 1]$, $\tau/\varrho \geq \delta^{-\chi(\zeta)}$,  $Y_z'\subset Y_z$ such that $|Y_z'|_\delta \geq (\varrho/\tau)^{\zeta}|Y_z|_{\delta}$ and  for every $y\in Y_z'$, there exists an affine transformation $\phi_{y, \tau}$ such that $\phi_{y, \tau}(Y_z'\cap B(y, \tau))$ is $(\varrho/\tau, d-\ell-\kappa, 1-\gamma_1, (\varrho/\tau)^{-\zeta})$-aligned.  By pigeonholing, there exists a $(\delta, \delta^{-\eta})$-uniform set $Z_1\subset Z$, $|Z_1|_\delta \gtrapprox |Z|_\delta$ such that the corresponding $(\varrho, \tau)$ are comparable for every $z\in Z_1$.  Let $Y$ be defined as in Remark~\ref{rem: configuration shading} with respect to our global $(\ell-1)$-grains $(\delta, \eta, \kappa)$-Kakeya configuration.  For each $T\in \TT$, define $Y'(T)\subset Y(T)$ as the set of $p=(x,y,z)$ with $z\in Z_1$, $y\in Y_z'$ and $x\in X_{y,z}$. Then $Y'$ is a $\gtrapprox \delta^{2\eta} (\varrho/\tau)^{\zeta}$-dense shading on $\TT$.

By pigeonholing, there exists a tube $T_\tau\in \TT_\tau$ such that $Y'$ is $\gtrapprox \delta^{2\eta} (\varrho/\tau)^{\zeta}$-dense on $\TT[T_\tau]$.  Let $\TT_1= \TT_\varrho^{T_\tau}$, the anisotropic rescaling of $\TT_\varrho[T_\tau]$ to $\varrho/\tau$-tubes, and let $Y_1$ be the rescaling of $Y$.  We claim that  $E_1= U(\TT_1, Y_1)$ is a global $(\ell-1)$-grains $(\delta_1, \eta_1, \kappa)$-Kakeya configuration with $\delta_1=\varrho/\tau$ and $\eta_1= 2\zeta$, when $\eta$ is sufficiently small so that $\delta^{-\eta}\leq \delta^{-\xi(\zeta)\eta_1}$. In addition,  since $\phi_{y, \tau}(Y_z'\cap B(y, \tau))$ is $(\varrho/\tau, d-\ell-\kappa, 1-\gamma_1, (\varrho/\tau)^{-\zeta})$-aligned, $E_1$ can be written in two ways 
\begin{equation}\label{eq: E1first}
    E_1=\bigsqcup_{z\in Z_1} E_{1, z}= \{ (x, y +f_1(z)x, z), z\in Z_1, y\in Y_z^1, x\in X_{y,z}^1\},
\end{equation}
where $Y_z^1\subset [-1, 1]^{d-\ell}$ is $(\delta_1, d-\ell-\kappa, \delta_1^{-\eta_1})$-AD-regular set and $X_{y,z}^1\subset[-1, 1]^{\ell-1}$ is such that $|X_{y,z}^1|_{\delta_1}\geq \delta_1^{\eta_1-\ell+1};$
\begin{equation}\label{eq: E1second}
    E_1= \bigsqcup_{z\in Z_1} \tilde{E}_{1, z} = \{ (\tilde x, \tilde y +\tilde f_1(z)\tilde x, z), z\in Z_1, \tilde y\in \tilde Y_z^1, \tilde x\in \tilde X_{y,z}^1\},
\end{equation}
where $\tilde Y_z^1\subset [-1, 1]^{d-\ell-1}$ is $(\delta_1, d-\ell-1-\kappa+\gamma_1, \delta_1^{-\eta_1})$-AD-regular set and $\tilde X_{y,z}^1\subset [-1, 1]^{\ell}$  is such that $|\tilde X_{y,z}^1|_{\delta_1}\geq \delta_1^{\eta_1-\ell+\gamma_1}$. 
 In addition, $Y_z^1$ is $(\delta_1, 1-\gamma_1, \delta_1^{-\eta_1})$-tube Katz-Tao.  
 
If $\gamma_1 \geq \eps_1$, we stop and move to Step 2 (recall that $\eta_1< g(\eps_1)$).

 Otherwise, $\gamma_1< \eps_1$ and $\tilde Y_z^1$ is a $(\delta_1, d-\ell-1-\kappa, \delta_1^{-\eta_1-\eps_1})$-AD regular set and $|\tilde X_{y,z}^1|_{\delta_1}\geq \delta_1^{-\eta_1+\eps_1 -\ell}$. Thus, $E_1$ is a global $\ell$-grains $(\d_1, \eta_1, \kappa)$-Kakeya configuration.

 We iterate the previous argument with $E_1$ written in \eqref{eq: E1second}, and $\TT_1, \eps_1+\eta_1, \eta_2$ in place of $E, \TT, \eta, \eta_1$ and obtain the corresponding set $\TT_2$ of $\delta_2$-tubes with $\delta_2 \leq \delta_1^{\chi(\eta_2/2)}$ and global $(\ell+1)$-grains $(\delta_2, \eta_2, \kappa)$-Kakeya configuration provided $\eps_1$ and $\eta_1$ are sufficiently small compared to $\eta_2$. We stop until   $\gamma_j\geq \eps_j$ or until $\ell+j=d$.  If $\ell+j=d$, then we know that $|E_j|_{\delta_j} \geq \delta_j^{-d(1-\eta_j)}$ which is a  contradiction when $\eta_j\ll \kappa/d$. 

\medskip

\noindent \textit{Step 2. Obtain the $(\delta_j, \gamma_j, \delta_j^{-\eta_j})$-Slab Frostman condition. }

Suppose we stop at the $j$th step and $\gamma_{j} \geq \eps_{j}$ and let $k= \ell+j-1$. 
Let $Y_j$ be the corresponding shading for $E_j$, the global $(k-1)$-grains $(\delta_j, \eta_j, \kappa)$-Kakeya configuration in the $j$th step (see Remark~\ref{rem: configuration shading}). At this stage, we know for any $p\in U(\TT_{j}, Y_j)$, the set of $\delta_j$ tubes in $\TT_{j, Y_j}(p)$ is contained in a $\delta_j$-neighborhood of an $\ell$-plane $\Sigma_\ell(p)$, but it is not clear if a typical $\ell$-tuple of tubes in it spans $\Sigma_\ell(p)$. 

To obtain the Slab Frostman condition, we should pigeonhole as in the proof of Proposition~\ref{prop: structure2} to find $\ell'$-linear points and possibly replace $\TT_j$ by  $\TT_{j,\varrho}^{T_\tau}$ for some $\varrho<\tau\in [\delta_j, 1)$ and replace $Y_j$ by a $\geq \delta_j^{\tilde\eta_j}$-refinement $Y_{j}^{T_\tau}$ for some $\tilde\eta_j>0$.  We can still write  $E':=U(\TT_{j,\varrho}^{T_\tau}, Y_{j}^{T_\tau})$ as 
\begin{equation}\label{eq: E'}
 E' = \bigsqcup_{z\in Z'} E_{z}'=\{(x, y+f_j(z)x, z), z\in Z', y\in Y_z', x\in  X_{y,z}'\}
\end{equation}
where $Y_z'\subset[-1, 1]^{d-\ell-j+1}$ is a $(\varrho/\tau, d-\ell-j+1-\kappa, \delta_j^{-2\eta_j-\tilde\eta_j})$-AD-regular set and $X_{y,z}'\subset [-1, 1]^{\ell+j-2}$ is such that $|X_{y, z}'|_{\varrho/\tau}\gtrapprox \delta_j^{2\eta_j+\tilde\eta_j}(\varrho/\tau)^{-\ell-j+2}$. In addition, $Y_z'$ is $(\varrho/\tau, 1- \gamma_j, \delta_j^{-2\eta_j-\tilde\eta_j})$-tube Katz-Tao.  



Choose $\eta_j$ sufficiently small compared to $\eps_j$, we repeat  the proof of Proposition~\ref{prop: structure2} until the end of Step 1 with $0< \eta_j\ll \tilde\eta_\ell \ll \cdots \ll \tilde\eta_1 \ll \eps_1/10d$ and $\tilde\varrho_{\ell-1}\ll \cdots \tilde\varrho_1\ll \delta_0$,  small parameters satisfying $10 \tilde\eta_{\ell'}/\eps_j \leq g(\tilde \eta_{\ell'-1}), 1\leq \ell'\leq \ell$, determined the same way as the proof of Proposition~\ref{prop: structure2} with $(\eta_j, \delta_j, \eps_j/(10d), \delta_0)$ in place of $(\eta, \delta, \eta_0, \delta_0)$,  to find $\ell'\in \{2, \dots, \ell\} $ and a $\gtrapprox 1$-refinement $Y_j'$ of $Y_j$ such that each point in $E_j'=U(\TT_j, Y_j')$ is $\ell'$-linear.  Here we stop at $\ell$ instead of $d$ because through any $p\in U(\TT_j, Y_j)$, $\TT_{j, Y_j}(p)$ is contained in a $\delta_j$-neighborhood of an $\ell$-plane $\Sigma_\ell(p)$.

Note that $\ell'\leq \ell$ could be different than $\ell$. If $\ell'=\ell$, then obtain the $\ell$-linear $k$-planar $(\delta_j, \eta_j, \kappa, \tilde\eta_{\ell-1})$-Kakeya configuration with $E_j'$. The $(\delta_j, \tilde\eta_{\ell-1}, \delta_j^{-\eta_j})$-Slab Frostman condition comes from the fact that the points in   $E_j'$ are not $(\ell-1)$-linear.

If $\ell'<\ell$, then there is a scale $\tilde\delta_{\ell'}\leq \tilde\varrho_{\ell'}$ and  for each $p\in E_j'$, there is an $\ell'$-plane $\Sigma_{\ell'}(p)$ such that $|\TT_{j, Y_j'}(p)\cap N_{\tilde\delta_{\ell'}}\Sigma_{\ell'}(p)|\geq \tilde\delta_{\ell'}^{\tilde\eta_{\ell'}} |\TT_{j, Y_j'}(p)|$ and for any $(\ell'-1)$-plane $\Sigma$ and any $r<\tilde\varrho_{\ell'-1}$, 
\begin{equation}
|\TT_{j, Y_j'}(p)\cap N_{r}\Sigma|\leq r^{\tilde\eta_{\ell'-1}} |\TT_{j, Y_j'}(p)|.
\end{equation}

For each $T\in \TT_j$, define $Y_{\tilde\delta_{\ell'}}(T)$ as the set of $p$ such that $T\in \TT_{j, Y_j'}(p)\cap N_{\tilde\delta_{\ell'}}\Sigma_{\ell'}(p)$.  Then $Y_{\tilde\delta_{\ell'}}$ is a $\geq \tilde\delta_{\ell'}^{\tilde\eta_{\ell'}}$-refinement of $Y_j'$.

Continue the proof of Proposition~\ref{prop: structure2} to find a global $(\ell'-1)$-grains $(...)$-Kakeya configuration of the form $E'$ for a suitable choice of $\varrho$ and $\tau$ (note that $\ell'$-linear property of points might get destroyed by the rescaling map $\phi$). We then need to verify the Slab Frostman condition of this configuration. We pigeonhole as in the proof of Proposition~\ref{prop: structure2} to find $\ell''$-linear  points, if  $\ell'' = \ell'$, then we are done. Otherwise  $\ell''<\ell'$, iterate the process until $\ell''=2$ (as in the proof of Proposition~\ref{prop: structure2}, if $\ell''=1$ we reach a contradiction).

To ease the notation, let $\delta'$-denote the resulting scale and  $\TT'$ denote the resulting set of $\delta'$-tubes and let $Y'$ denote the corresponding shading and $E'=U(\TT', Y')$ is of the form \eqref{eq: E'} such that each point is $\ell'$-linear but not $(\ell'-1)$-linear (we continue to use $\ell'$ to denote $\ell''$). And this gives an $\ell'$-linear $(\delta', \eta', \kappa)$-Kakeya configuration whose $\Phi_p$ satisfies $(\delta', \gamma', \delta'^{-\eta'})$-Slab Frostman condition for some $\eta', \gamma'>0$ with $\eta' < g(\gamma')$, $\delta_j^{-\eta_j}\ll \delta'^{-\eta'}$ and $\gamma'< \gamma_j$. By choosing the initial scale $\d$ sufficiently small compared to the other parameters, we can additionally ensure that $\d' < h(\eta', \gamma')$ holds. 

\medskip

\noindent \textit{Step 3. Consistency relation.}







Therefore,  $E'$ can be written of the form 
\begin{equation} \label{eq: Ej'}
E'=\bigsqcup_{z\in Z'} E_{z}'=\{x', y'+ f^\dag (z)x', z): z\in Z', y'\in Y_z , x'\in X_{y',z}\}
\end{equation}
where $Y_z\subset [-1, 1]^{d-\ell'}$ and $X_{y',z}\subset [-1, 1]^{\ell'-1}$ is such that $|X_{y',z}|_{\delta'} \geq \delta'^{ 5\eta'/\eps_j -\ell'}$ and the slope function $f^\dag: Z' \rightarrow \operatorname{Mat}(d-\ell', \ell'-1)$. 

Compare $E'$ written in \eqref{eq: Ej'}  with $E'$ written in \eqref{eq: E'}:
 \[ E'=\bigsqcup_{z\in Z'} E_{z}= \{ x, y +f_j(z)x, z), z\in Z', y\in Y_z', x\in X_{y,z}'\},\]
 with $f_j: Z' \rightarrow \operatorname{Mat}(d-k, k-1)$, with $k=\ell+j-1$  (recall we view $x, y$ and $x', y'$ as column vectors). Recall that  $Y_{z}'$ is $(\delta', \gamma_j, \delta'^{-\eta'})$-tube Katz-Tao (by definition, $\gamma_j<1-\eps_j$), which then 
 implies 
\[
f^\dag_2 = f_{j1}+f_{j2}f^\dag_1 +O(\delta'^{1-5\eta'/\eps_j})
\]
where 
        \[
         f^\dag =\begin{pmatrix}
            f^\dag_1 \\
            f^\dag_2 
        \end{pmatrix}, 
         \quad  f_j = ( f_{j1}, f_{j2}),
        \]         
$f^\dag_1\in \operatorname{Mat}(k-\ell', \ell'-1), f^\dag_2, f_{j1}\in \operatorname{Mat}(d-k, \ell'-1)$ and $f_{j2}\in \operatorname{Mat}(d-k, k-\ell')$. 
Indeed, the the $(\ell'-1)$-grains and $(k-1)$-grains in the two representations above do not align up to angle at least $\delta'^{1-5\eta'/\eps_j}$, this would imply that $Y'_{z}$ essentially fills out a tube, contradicting tube Katz-Tao condition.

We can make the error term zero by thickening to scale $\delta'^{1-5\eta'/\eps_j}$ and adjusting the set $E'$ slightly. 
Then $E'$ is the $\ell'$-linear $k$-planar $(\delta', 10\eta'/\eps_j, \kappa, \gamma)$-Kakeya configuration for some  $\gamma=  \min \{ \eps_j, \tilde\eta_{\ell'-1}\}$ and $10\eta'/\eps_j < g(\gamma)$. 
\end{proof}


\section{Properties of slope functions}\label{sec:slope-func-properties}

In this section we prove some useful basic properties of $\ell$-linear $k$-planar Kakeya configurations. First, we show that if either of the slope functions $f$ or $\tilde f$ varies at a slower rate than Lipschitz, then the set of tubes $\TT$ violates Convex Wolff axioms. In the two proofs it is crucial that we use the Slab Frostman and the Tube Katz--Tao axioms (respectively) included in Definition \ref{def:kakeya-configuration-2}. The proof for $f$ will be similar to the proof of Lemma \ref{lem:slowly-rotating} of the corresponding Fustenberg statement, while the proof for $\tilde f$ is based in part on \cite[Lemma 6.3]{wang2026sticky}.

The last property we need is that both functions $f$ and $\tilde f$ may be assumed to be approximately Lipschitz. 
In the next section we define `$W$-tuples' -- they were already used for sticky Fustenberg in Section \ref{subsec:W-tuples}. Here the computations are largely the same, we include all details for readers convenience. 

\subsection{\texorpdfstring{$W$}{W}-tuples.}

\label{subsec:W-tuples-kakeya}
The contents of this section will be used in Sections
\ref{subsec:convex-wolff-kakeya} and
\ref{subsec:non-linear-case-Kakeya}.

Fix an $\ell$-linear $k$-planar
$(\d,\eta,\kappa,\gamma)$-Kakeya configuration in $\RR^d$ and a scale $\varrho\in[\d^{1/4},1)$. Fix a dyadic $\varrho$-cube
${\bf q}=p_0^\varrho$ intersecting with $E$ and write
$E({\bf q})=E\cap \bf q$, $Z({\bf q})$ as the set of $z$ such that $E_z\cap \bf q\neq \emptyset$,  $\Phi_{\bf q}$, and $\TT({\bf q})=\cup_{p\in E(\bf q)} \TT_p$ for the corresponding objects describing the structure of the Kakeya configuration inside ${\bf q}$.
We can write $E(\bf q)$ in the form
\[
E({\bf q})
=
\{(x,y+f(z)x,z):
z\in Z({\bf q}),\ y\in Y_z({\bf q}),\
x\in X_{y,z}({\bf q})\}.
\]
We have the following properties that follow from Definition \ref{def:kakeya-configuration} (3), (4):
\begin{equation}\label{eq:kakeya-rho-cube-counts}
\begin{split}
|E({\bf q})|
&\sim_{\d^{O(\eta)}}|Z({\bf q})|
(\varrho/\d)^{d-1-\kappa},\\
|\TT({\bf q})|
&\sim_{\d^{O(\eta)}}\varrho^{-\kappa}
(\varrho/\d)^{d-1},\\
\sum_{T\in\TT({\bf q})}|Z({\bf q},T)|
&\sim_{\d^{O(\eta)}}|Z({\bf q})|\,|\TT({\bf q})|.
\end{split}
\end{equation}


\medskip\noindent\emph{Definition of the $W$-tuples.}
For $z, z', z'' \in Z({\bf q})$ define the set of $W$-tuples to be the following 
\begin{align*}
W_{z, z', z''}({\bf q}) = \{ (p,p_1,p_2,p_1',p_2', T_1, T_2,T_1', T_2') \in E_z({\bf q} ) \times E_{z'}({\bf q}) \times E_{z'}({\bf q}) \times E_{z''}({\bf q}) \times E_{z''}({\bf q}) \times \TT({\bf q})^4:   \\
T_i \in \TT_p \cap \TT_{p_i}, \quad T_i' \in \TT_{p_i} \cap \TT_{p_i'}, \quad |\theta(T_1') -\theta(T_2')| \le \varrho 
\}    
\end{align*}
Thus a $W$-tuple is obtained by taking two copies of the same
two-edge path from level $z$ to level $z''$, with a common starting
point and $\varrho$-close tube directions.
See Figure \ref{fig:W-tuple-furstenberg} for an illustration.

\medskip\noindent\emph{Counting the $W$ tuples.}
The next lemma gives the expected lower bound.

\begin{lemma}\label{lem:W-tuple-many-Kakeya}
We have
\[
\sum_{z,z',z''\in Z({\bf q})}|W_{z,z',z''}({\bf q})|
\gtrsim_{\d^{O(\eta)}}
|E({\bf q})|\d^{-4\kappa}\varrho^\kappa
|Z({\bf q})|^2.
\]
\end{lemma}

\begin{proof}
For $z,z',z''\in Z({\bf q})$, let
\[
\begin{split}
P_{z,z',z''}
=\{&(p,p_1,p_1',T_1,T_1'):
p\in E_z({\bf q}),\ p_1\in E_{z'}({\bf q}),\
p_1'\in E_{z''}({\bf q}),~ T_1, T_1' \in \TT({\bf q})\\
&T_1\in\TT_p\cap\TT_{p_1},\
T_1'\in\TT_{p_1}\cap\TT_{p_1'}\}.
\end{split}
\]
So each $W$-tuple is composed of two 3-vertex paths in $P_{z,z',z''}$ starting at $p$.
The estimates
\eqref{eq:kakeya-rho-cube-counts} give
\begin{equation}\label{eq:P-path-count-kakeya}
\sum_{z,z',z''\in Z({\bf q})}|P_{z,z',z''}|
\sim_{\d^{O(\eta)}}
|E({\bf q})|\d^{-2\kappa}|Z({\bf q})|^2.
\end{equation}
Indeed, first select $p$ in $|E({\bf q})|$ ways, then $T \in \TT_p \cap \TT({\bf q})$, $p_1 \in Z(T, {\bf q})$, $T_1' \in \TT_{p_1}$ and $p_1' \in Z(T_1',{\bf q})$ can be selected in $\sim \d^{-\kappa}$, $\sim |Z({\bf q})|$ and $\sim \d^{-\kappa}$ ways, respectively.  

Now for $p \in E_z({\bf q})$ and a dyadic $c\varrho$-cube ${\bf \theta}$ define a subset $P_{z,z',z''}(p, {\bf \theta}) \subset P_{z,z',z''}$ consisting of tuples $(p,p_1,p_1',T_1,T_1')$ with $\theta(T_1') \in {\bf \theta}$. For every $p$, there are $\sim_{\d^{O(\eta)}} \varrho^{-\kappa}$-many choices for ${\bf \theta}$ that make this set non-empty. So using Cauchy--Schwarz it follows
\begin{align*}
\sum_{z,z',z''}\sum_{p, {\bf \theta}} |P_{z,z',z''}(p, {\bf \theta})|^2    &\gtrsim_{\d^{O(\eta)}}
\frac{
\bigl(|E({\bf q})|\d^{-2\kappa}|Z({\bf q})|^2\bigr)^2
}{
|E({\bf q})||Z({\bf q})|^2\varrho^{-\kappa}
}\\
&=
\d^{O(\eta)}
|E({\bf q})|\d^{-4\kappa}\varrho^\kappa
|Z({\bf q})|^2. 
\end{align*}

On the other hand, every element of $P_{z,z',z''}(p, {\bf \theta})^2$ gives rise to a $W$-tuple and every $W$-tuple is obtained this way in at most $\sim 1$-many ways. This gives the desired formula for the number of $W$-tuples.
\end{proof}

\medskip\noindent\emph{Additive relations.}
Fix a tuple in $W_{z,z',z''}({\bf q})$ and write
\[
\begin{aligned}
p&=(x,y+f(z)x,z),\\
p_i&=(x_i,y_i+f(z')x_i,z'),\\
p_i'&=(x_i',y_i'+f(z'')x_i',z'').
\end{aligned}
\]
Let $\xi=\xi_p$, $\xi_i=\xi_{p_i}$, and
$\xi_i'=\xi_{p_i'}$.  We denote the first $\ell-1$ direction
coordinates of $T_i$ and $T_i'$ by $\varphi_i$ and $\varphi_i'$,
respectively.  In particular,
$|\varphi_1'-\varphi_2'|\leq\varrho$.

Since $T_i$ passes through both $p$ and $p_i$, while $T_i'$ passes
through both $p_i$ and $p_i'$, Definition
\ref{def:kakeya-configuration} gives
\[
\begin{aligned}
\theta(T_i)
&=(\varphi_i,\xi+f(z)\varphi_i)+O(\d)
 =(\varphi_i,\xi_i+f(z')\varphi_i)+O(\d),\\
\theta(T_i')
&=(\varphi_i',\xi_i+f(z')\varphi_i')+O(\d)
 =(\varphi_i',\xi_i'+f(z'')\varphi_i')+O(\d).
\end{aligned}
\]
Consequently,
\begin{equation}\label{eq:xii}
\xi_i=\xi+(f(z)-f(z'))\varphi_i+O(\d),
\qquad
\xi_i'=\xi_i+(f(z')-f(z''))\varphi_i'+O(\d),
\end{equation}
and hence
\begin{equation}\label{eq:xi-diff-kakeya}
\xi_2'-\xi_1'
=(f(z)-f(z'))(\varphi_2-\varphi_1)
+O\bigl(\varrho|f(z'')-f(z')|+\d\bigr).
\end{equation}

We have
\[
p_i=p+(z'-z)(\theta(T_i),1)+O(\d),
\qquad
p_i'=p_i+(z''-z')(\theta(T_i'),1)+O(\d).
\]
Subtracting the relations with $i=1,2$ and using
$|z''-z'|\leq\varrho$ and
$|\theta(T_2')-\theta(T_1')|\leq\varrho$, we obtain
\[
p_2'-p_1'
=(z'-z)(\varphi_2-\varphi_1,f(z)(\varphi_2-\varphi_1),0)
+O(\varrho^2+\d).
\]
In particular,
\begin{equation}\label{eq:x-diff-kakeya}
x_2'-x_1'
=(z'-z)(\varphi_2-\varphi_1)+O(\varrho^2+\d).
\end{equation}
Combining this with \eqref{eq:xi-diff-kakeya} gives the vector identity
\begin{equation}\label{eq:vector-identitiy}
\begin{pmatrix}
x_2'-x_1'\\
\xi_2'-\xi_1'
\end{pmatrix}
=
\begin{pmatrix}
z'-z\\
f(z)-f(z')
\end{pmatrix}
(\varphi_2-\varphi_1)
+
\begin{pmatrix}
O(\varrho^2+\d)\\
O(\varrho|f(z'')-f(z')|+\d)
\end{pmatrix}.
\end{equation}

Comparing the second coordinate block of $p_2'-p_1'$ with the
representation of $p_i'$ at level $z''$ also yields
\begin{equation}\label{eq:y-differ-little}
\begin{split}
y_2'-y_1'
&=(z'-z)(f(z)-f(z''))(\varphi_2-\varphi_1)
+O(\varrho^2+\d)\\
&=O\bigl(
\varrho|f(z)-f(z')|
+\varrho|f(z')-f(z'')|
+\varrho^2+\d
\bigr).
\end{split}
\end{equation}
Finally, subtracting the direction formulas at the common intermediate
points gives
\begin{equation}\label{eq:diff-theta-prime-kakeya}
\begin{split}
\theta(T_1')-\theta(T_2')
&=(\varphi_1'-\varphi_2', f(z')(\varphi_1'-\varphi_2'))
 +(0,(f(z)-f(z'))(\varphi_1-\varphi_2))
+O(\d).
\end{split}
\end{equation}

\subsection{Small variation of slope functions.} \label{subsec:convex-wolff-kakeya}
In this section we show that if either of the slope functions $f$ and $\tilde f$ of a $k$-planar $\ell$-linear Kakeya configuration varies slower than Lipschitz, then tubes in $\TT$ concentrate in convex sets. 

\begin{lemma}\label{lem:f-slow}
    Let $(E, \TT, \{\TT_p\})$ be an $\ell$-linear $k$-planar $(\d, \eta, \kappa,\gamma)$-Kakeya configuration in $\RR^d$ with the slope function $f(z) \in \operatorname{Mat}(d-\ell,\ell-1)$ and $\eta<\gamma$.  Suppose that for some $\varrho \in [\d^{1/3}, 1]$ and $\sigma \in [\varrho^2, \varrho]$ and every $z_0 \in Z$, ${\bf z} = z_0^\varrho$, we have 
    \begin{equation}\label{eq: assumption18.1}
    \#\{z \in Z \cap {\bf z}:~ |f(z) - f(z_0)| \le \sigma \} \ge \d^{\varepsilon} |Z \cap {\bf z}|.
    \end{equation}
    Then we have $C_{CW}(\TT) \gtrsim\d^{O(\eta/\gamma+\varepsilon)} (\varrho/\sigma)^{\kappa}$.
\end{lemma}

\begin{proof}
Fix $p_0\in E_{z_0}$ and let $\mathbf{q}=p_0^\varrho$ and $Z(\mathbf{q})= \{ z\in Z\cap \mathbf{z}: |f(z)-f(z_0)|\leq \sigma\}$. 
    
Write \[
f(z)=\begin{pmatrix}
    f_1(z) \\
    \vdots\\
    f_{\ell-1}(z)
\end{pmatrix} \in  \operatorname{Mat}(d-\ell, \ell-1).
\]

    
    For $T\in \TT$ and $z \in Z$ let 
    \begin{align*}
    U(z, T) = N_\sigma T + [-\varrho, \varrho] \cdot (1, 0, \ldots,0, f_1(z),0 ) + [-\varrho,\varrho] \cdot (0,1,\ldots,0, f_2(z), 0) + \cdots \\ 
     + [-\varrho, \varrho] \cdot (0,\ldots, 0, 1, f_{\ell-1}(z), 0)    
    \end{align*}
    Note that $U(z, T)$ is a $\sim \underbrace{\sigma\times \ldots \times \sigma}_{d-\ell} \times \underbrace{\varrho\times \ldots \times \varrho}_{\ell-1}\times 1$ box. We will show that $\TT[C \cdot U(z, T)]$ is large for an appropriately large constant $C$ and a typical choice of $T \in \TT({\bf q})$ and $z \in Z({\bf q})$. Consider the set of $W$-tuples $W_{z,z',z''}({\bf q})$ as defined in Section \ref{subsec:W-tuples-kakeya}. 

    For a pair of tubes $T_1', T_2'$ and $z'' \in Z({\bf q})$ we write $(z'', T_1') \sim (z'', T'_2)$ if they can be extended to a tuple $(p,p_1',p_2', p_1'', p_2'', T_1, T_2, T'_1, T'_2) \in W_{z,z', z''}({\bf q})$ for some $z, z' \in Z({\bf q})$.

    \begin{claim}\label{claim:adjacent-tubes-same-box-kakeya}
        There exists a constant $C$ only depending on the dimension $d$ such that if $(z'', T_1') \sim (z'', T'_2)$ then $U(z'', T_1') \subset C \cdot U(z'', T_2')\subset C^2\cdot U(z'', T_1')$.
    \end{claim}

    \begin{proof}
    Consider a tuple $(p, p_1, p_2, p_1', p_2', T_1,T_2, T'_1, T'_2) \in W_{z, z', z''}({\bf q})$ and write $p = (x, y+ f(z) x, z)$, $p_i = (x_i, y_i + f(z') x_i, z' )$ and $p_i' = (x_i', y'_i+ f(z'') x'_i,z'')$. Let $\varphi_1,\varphi_2, \varphi_1', \varphi_2'$ denote the first $\ell-1$ coordinates of the slopes $\theta(T_1), \theta(T_2), \theta(T_1'), \theta(T_2')$. Then we have by (\ref{eq:diff-theta-prime-kakeya}) and $|f(z')-f(z'')|\leq \sigma$ for any $z', z''\in Z(\mathbf{q})$:
    \begin{align*}
    \theta(T_1') - \theta(T_2') = (\varphi_1'-\varphi_2')(1, f(z')) + (\varphi_1-\varphi_2) (0, f(z)-f(z')) + O(\d)     \\
    =  (\varphi_1'-\varphi_2')(1, f(z')) +O(\sigma) = O(\varrho) (1, f(z'')) + O(\sigma).
    \end{align*}
  This means that even though $\theta(T_1')-\theta(T_2')$ could be as  large as $O(\varrho)$, but the difference lies in 
  \begin{align*}
  [-C\varrho, C\varrho]\cdot (1, 0, \dots, 0, f_1(z'')) + [-C\varrho, C\varrho]\cdot (0, 1, \dots, 0, f_2(z''))+\\ \cdots + [-C\varrho, C\varrho] \cdot (0, \dots, 0, 1, f_{\ell-1}(z'')) +O(\sigma).    
  \end{align*}
    
    By (\ref{eq:x-diff-kakeya}), (\ref{eq:y-differ-little}):
    \begin{align*}
    x_1' - x_2' &= (z'-z) (\varphi_1-\varphi_2) + O(\varrho^2) = O(\varrho),\\
    y_1'-y_2' &=  O(\varrho \sigma+\d).
    \end{align*}
    Similarly, this says $p_1'-p_2'$ lies in $O(\varrho\sigma+\delta)$-neighborhood of 
    \begin{align*}
  [-C\varrho, C\varrho]\cdot (1, 0, \dots, 0, f_1(z''),0) + [-C\varrho, C\varrho]\cdot (0, 1, \dots, 0, f_2(z''),0)+\\ \cdots + [-C\varrho, C\varrho] \cdot (0, \dots, 0, 1, f_{\ell-1}(z''),0) .    
  \end{align*}
    
    From this we conclude that $p_2' \in U(z'', T_1')$ and $T_2' \subset C\cdot U(z'', T_1')$ which implies
    $U(z'', T'_1)$ and $U(z'', T'_2)$ define essentially the same box, concluding the proof.  Since the relation $(z'', T_1')\sim (z'', T_2')$ is symmetric, we also have $C\cdot U(z'', T_2')\subset C^2\cdot U(z'', T_1').$
    \end{proof}

    Let $\overline{\Phi}_{\bf q} \subset \varrho \ZZ^{\ell-1}$ be the set such that $\Phi_{\bf q}^\varrho = \overline{\Phi}_{\bf q}^\varrho$. Define a subset
    \[
    G(z'', T_1') \subset Z({\bf q}) \times Z({\bf q}) \times \overline{\Phi}_{\bf q}\times \overline{\Phi}_{\bf q}
    \]
    to be the set of tuples $(z, z', \overline{\varphi}_1, \overline{\varphi}_2)$ for which there exists a tuple $(p, p_1, p_2, p_1', p_2', T_1, T_2, T_1', T_2') \in W_{z, z',z''}({\bf q})$ such that $\theta(T_i) \in N_\varrho (\overline{\varphi}_i, \xi_p + f(z) \overline{\varphi}_i)$. By the lower bound on $|W_{z,z',z''}({\bf q})|$ we get
    \[
    \frac{1}{|Z({\bf q}) ||\TT({\bf q})|}\sum_{z'' \in Z({\bf q}), T'_1 \in \TT({\bf q})} |G(z'', T_1')| \gtrsim\d^{O(\eta)} |Z({\bf q})|^2 |\overline{\Phi}_{\bf q}|^2.
    \]
    Heuristically, the above inequality says the following: fix $z''\in Z(\mathbf{q}), T_1'\in \TT(\mathbf{q})$, this essentially fixes $p'_1$; then for a typical $z', \theta(T_1), z, \theta(T_2)$, the corresponding $p_2'$ exists and since $\Phi_{p_1'}^\varrho \approx \Phi_{p_2'}^\varrho$ (in the sense that they are both a significant fraction of $\Phi_{\mathbf{q}}^\varrho$), there exists $\theta(T_2')\in \Phi_{p_2'}$, $|\theta(T_2')-\theta(T_1')|\lesssim \varrho$.

    Let $\mc P$ be the set of all pairs $(z, T)$ with $T \in \TT({\bf q})$ and $z \in Z({\bf q}, T):=\{z\in Z(\mathbf{q}): \exists p\in E_z\cap \mathbf{q}, T\in \TT_p\}$. Let $\mc P' \subset \mc P$ be the subset of all $(z'', T'_1)$ with $|G(z'', T_1')| \ge \d^{C_1 \eta}|Z({\bf q})|^2 |\overline{\Phi}_{\bf q}|^2$ for some large constant $C_1$. It follows from the above that most of the $W$-tuples are captured by pairs $(z, T) \sim (z, T')$ with both $(z,T), (z, T') \in \mc P'$. Let $\mc G' \subset \mc P'\times \mc P'$ be the adjacency graph for the relation $(z, T) \sim (z, T')$.

    Fix some $(z, T) \in \mc P'$ and for an integer $m \ge 1$ consider the $m$-th neighborhood set $\mc N_m(z, T) \subset \mc P'$ in $\mc G$, i.e. this is the set of $(z, T') \in \mc P'$ which are connected to $(z, T)$ by a path with at most $m$ edges (note that the $z$-coordinate doesn't change along the edges of $\mc G'$). By Claim \ref{claim:adjacent-tubes-same-box-kakeya} all pairs $(z,T') \in \mc N_{m}(z, T)$ are all contained in the same convex set $Cm \cdot U(z, T)$. So it suffices to show that for a bounded $m$, the neighborhood $\mc N_m (z, T)$ is large.

    Now we observe that for $(z'', T'_1) \sim (z'', T_2')$, if we denote by $p_1' = (x_1', y_1'+f(z'') x_1', z'')$ and $p_2'=(x_2', y_2'+f(z'')x_2', z'')$ the points in $T_1' \cap E_{z''}$ and $T_2'\cap E_{z''}$, respectively, then by the computations we did in (\ref{eq:x-diff-kakeya}), we get that 
\[
x'_1-x'_2 = (z'-z) (\overline{\varphi}_1-\overline{\varphi}_2) + O(\varrho^2).
\]
 For $(z, T') \in \mc P'$ let $x(z, T')$ be the $x$-coordinate of the intersection point $T' \cap E_z({\bf q})$.
Now let $(z, T') \in \mc N_m(z, T)$ be arbitrary and let $(z_1,z_2, \overline{\varphi}_1, \overline{\varphi}_2) \in G(z, T')$. Then it follows from the above that there exists some $(z, T'') \in \mc N_{m+1}(z, T)$ such that
\[
x(z, T'') - x(z, T') = (z_2-z_1)(\overline{\varphi}_2-\overline{\varphi}_1) + O(\varrho^2).
\]
Let $X_{m}$ be a maximal $\varrho^2$-separated subset of points $x(z, T')$ over all $(z, T') \in \mc N_{m}(z, T)$. We obtain a graph 
\[
G_m \subset X_m \times Z({\bf q}) \times Z({\bf q}) \times \overline{\Phi}_{{\bf q}}\times \overline{\Phi}_{{\bf q}},
\]
so that for every $(x, z_1, z_2, \overline{\varphi}_1, \overline{\varphi}_2) \in G_m$ there exists $x' \in X_{m+1}$ such that
\[
x' - x = (z_2-z_1) (\overline{\varphi}_2-\overline{\varphi}_1)+ O(\varrho^2).
\]
Furthermore, we have the lower bound 
\[
|G_m| \ge \d^{O(\eta)}|X_m|  |Z({\bf q})|^2|\overline{\Phi}_{\bf q}|^2.
\]

By the assumption \eqref{eq: assumption18.1} we have $|Z({\bf q})|_{\varrho^2} \gtrsim_{\d^{O(\eta+\varepsilon)}} \varrho^{-1}$ and by the definition of an $\ell$-linear Kakeya configuration we have that $\overline{\Phi}_{\bf q}$ is $(\varrho, \gamma, \d^{-O(\eta)})$-slab Frostman set. Using this information and some pigeonholing, we can find $\ell-1$ pairs $(\overline{\varphi}_1^i, \overline{\varphi}_2^i)$, $i=1, \ldots, \ell-1$, so that  $|(\overline{\varphi}_1^1- \overline{\varphi}_2^1)\wedge \ldots \wedge (\overline{\varphi}_1^{\ell-1}- \overline{\varphi}_2^{\ell-1})| \gtrsim \d^{O(\eta/\gamma)}$ and so that $X_{\ell-1}$ has density $\d^{O(\eta)}$ in the set
\[
\left\{ \sum_{i=1}^{\ell-1} (z_i' - z_i) (\overline{\varphi}_1^i- \overline{\varphi}_2^i) \mid  z_i', z_i \in Z({\bf q}) \right\}.
\]
Take $m=\ell-1$, from this we we have the lower bound 
\[
|X_{\ell-1}|_{\varrho^2} \gtrsim \d^{O(\eta/\gamma+\varepsilon)} \varrho^{-(\ell-1)}. 
\]
Let $\{x_1,\ldots,x_n\} \subset X_m$ be a maximal $\sigma$-separated set (recall $\sigma\ge \varrho^2$). For each $x_i$ we have a corresponding pair $(z, T_i) \in \mc N_m(z, T)$ and a corresponding point $p_i \in E_z \cap T_i$. We know that $T_i \subset C m \cdot U(z,T)$. Furthermore, if we write $\theta(T_i) \in N_\d (\varphi_i, \xi_{p_i}+f(z) \varphi_i)$ for some $\xi_i \in \Xi_{p_i}$ and $\varphi_i \in \Phi_{p_i}$, then using AD-regularity of sets $\Phi_{p_i}$, we get $\sim_{\d^{O(\eta)}} (\sigma/\varrho)^{-\kappa}$-many $\sigma$-tubes $\tilde T_i \in \TT_{\sigma}$ which intersect the $z$-plane within $\sigma$-ball around $T_i \cap E_z$ and satisfy $\theta(\tilde T_i) = (\tilde\varphi_i, \xi_{p_i} + f(z) \tilde \varphi_i+O(\sigma) )$ with some $\tilde \varphi_i = \varphi_i+ O(\varrho)$. Indeed, we choose some $\varphi_i' \in \Phi_{p_i}$ within distance $\varrho$ from $\varphi_i$, then we choose a tube $T'_i \in \TT_{p_i}$ in direction $(\varphi_i', \xi_{p_i}+f(z)\varphi_i')$ and let $\tilde T_i = (T'_i)^{\sigma}$. For large enough constant $K$ we get
\[
|\TT_{\sigma}[ K\cdot U(z,T) ]| \gtrsim \d^{O(\eta/\gamma+\varepsilon)} n (\varrho/\sigma)^{\kappa} \gtrsim \d^{O(\eta/\gamma+\varepsilon)} (\varrho/\sigma)^{\ell-1+\kappa}.
\]
On the other hand, note that
\[
\frac{\operatorname{Vol}(U(z,T))}{\operatorname{Vol}(T_\sigma)} \sim (\varrho/\sigma)^{\ell-1}
\]
and so we conclude that
\[
C_{CW}(\TT) \gtrsim_{\d^{O(\eta)}} C_{CW}(\TT_\sigma) \gtrsim_{\d^{O(\eta/\gamma+\varepsilon)} } (\varrho/\sigma)^{\kappa},
\]
finishing the proof.
\end{proof}

\begin{lemma}\label{lem:f-tilde-slow}
    Let $(E, \TT, \{\TT_p\})$ be an $\ell$-linear $k$-planar $(\d, \eta, \kappa, \gamma)$-Kakeya configuration in $\RR^d$ with the slope function $\tilde f(z) \in \operatorname{Mat}(d-k,k-1)$. 
    Suppose that for some $\varrho \in [\d^{1/3}, 1]$ and $\sigma \in [\varrho^2, \varrho]$ and $z_0 \in Z$, ${\bf z} = z_0^\varrho$, we have 
    \[
    \#\{z \in Z \cap {\bf z}:~ |\tilde f(z) - \tilde f(z_0)| \le \sigma \} \ge \d^{\varepsilon} |Z \cap {\bf z}|.
    \]
    Then we have $C_{CW}(\TT) \gtrsim \d^{C (\varepsilon+\eta)/\gamma} (\varrho/\sigma)^{\kappa}$. 
\end{lemma}

\begin{proof}
We write 
\[
E = \{(x, y+\tilde f(z)x, z), ~ x \in X_{y,z}, y\in Y_z, z\in Z\}
\]
for sets $X_{y,z} \subset [-1,1]^{k-1}$, $Y_z \subset [-1,1]^{d-k}$ and $Z \subset [-1,1]$. In this proof we will not use the representation of $E$ associated with the $\ell$-linear slope function $f$.

For $z \in Z$, $y\in Y_z$ and $x \in X_{y,z}$ let $\xi(x,y,z) \in [-1,1]^{d-k}$ be a vector such that for $p = (x, y+\tilde f(z)x, z) \in E$ and any $T\in \TT_p$ we have $\theta(T) \in N_\d(\varphi, \xi(x,y,z) + \tilde f(z) \varphi)$. The function $\xi(x,y,z)$ is constructed from $\xi_p$ in Definition \ref{def:kakeya-configuration} (4) and the consistency relation in Definition \ref{def:kakeya-configuration-2} (1): every $T \in \TT_p$ satisfies
\[
\theta(T) \in N_\d (\varphi_1, \xi_p + f(z) \varphi_1), \quad \varphi_1 \in [-1,1]^{\ell-1},
\]
write $f = \begin{pmatrix}
    f_1\\
    f_2
\end{pmatrix}$ and $\xi_p = (\xi_{p,1}, \xi_{p,2}) \in [-1,1]^{k-\ell} \times [-1,1]^{d-k}$ then 
\[
\begin{aligned}
(\varphi_1, \xi_p + f(z) \varphi_1) &= (\varphi_1,\xi_{p,1}+f_1(z)\varphi_1, \xi_{p,2}+f_2(z)\varphi_1)    \\
&= (\varphi, \xi_{p,2}- \tilde f_2(z)\xi_{p,1} + \tilde f(z) \varphi)
\end{aligned}
\]
with $\varphi = (\varphi_1,\xi_{p,1}+f_1(z)\varphi_1)$. So we can define $\xi(x,y,z)=\tilde \xi_{p}=\xi_{p,2}- \tilde f_2(z)\xi_{p,1}$ for $p = (x, y+\tilde f(z)x, z)$. From the dyadic Lipschitz property of $\xi_p$ in Definition \ref{def:kakeya-configuration} (4) it follows that $\xi(x,y,z)$ is $\approx \d^{-\eta}$ Lipschitz in the variable $x$. 

We claim that there exists $z \in Z$ and $y \in Y_z$ and some $\xi \in [-1,1]^{d-k}$ such that 
\[
\#\{ x \in X_{y,z}:~ |\xi(x, y, z) - \xi| \le \sigma /\varrho \} \gtrsim\d^{C (\varepsilon+\eta)/\gamma} |X_{y,z}|.
\]
Let us first see why this suffices to prove the lemma. 
Indeed, write $\tilde f(z) = (\tilde f_1(z), \ldots, \tilde f_{k-1}(z))$ where $\tilde f_i(z) \in [-1,1]^{d-k}$ and define 
\begin{align*}
U = [-1,1]\cdot (1,0,\ldots,0, \tilde f_1(z),0) + [-1,1]\cdot (0,1,\ldots,0, \tilde f_2(z),0) + \ldots \\+ [-1,1]\cdot (0,0,\ldots,1, \tilde f_{k-1}(z),0) + [-1,1] \cdot (0, \xi, 1)  + B(\sigma/\varrho,0).
\end{align*}
Thus, $U$ is a $\underbrace{\sigma/\varrho \times \ldots \times \sigma/\varrho}_{d-k}\times \underbrace{1\times \ldots \times 1}_{k}$ box. 
We have the following lower bound on $\TT_{\sigma/\varrho}[C\cdot U]$:
\[
|\TT_{\sigma/\varrho}[C\cdot U]| \gtrsim\d^{C (\varepsilon+\eta)/\gamma} |X_{y,z}|_{\sigma/\varrho} (\varrho/\sigma)^{\kappa}.
\]
where we sum over a maximal $\sigma/\varrho$-separated set of $x \in X_{y,z}$ such that $|\xi(x, y, z) -\xi| \lesssim \sigma/\varrho$ and all tubes $T \in \TT_{\sigma/\varrho}$ that pass through the point $(x, y+\tilde f(z)x,z)\in E$. Every such tube has the property that $\theta(T) \in N_{C\sigma/\varrho}(\varphi, \xi+\tilde f(z) \varphi)$ which implies that $T$ is contained in a constant dilate of $U$. We thus conclude that
\[
C_{CW}(\TT) \gtrsim_{\d^{O(\eta)}} C_{CW}(\TT_{\sigma/\varrho}) \gtrsim\d^{C (\varepsilon+\eta)/\gamma} |X_{y,z}|_{\sigma/\varrho} (\varrho/\sigma)^{\kappa} \frac{\operatorname{Vol}(T_{\sigma/\varrho})}{\operatorname{Vol}(U)} \gtrsim \d^{C (\varepsilon+\eta)/\gamma} (\varrho/\sigma)^\kappa
\]
as desired. 

Thus, it remains to prove the assertion that $\xi(x,y,z)$  concentrate in a $\sigma/\varrho$ interval as we vary $x$. Fix a dyadic interval ${\bf z} = z_0^\varrho$.
Let $Z({\bf z}) \subset Z\cap {\bf z}$ be a subset of size $\gtrsim \d^\varepsilon |Z\cap {\bf z}|$ so that $|\tilde f(z) - \tilde f(z')| \le \sigma$ for all $z, z' \in Z({\bf z})$.
For a dyadic $\varrho$-cube ${\bf q} = p^\varrho$ with Z-coordinate ${\bf z}$ let $\tilde \xi_{\bf q}$ denote a vector such that $|\tilde \xi_{\bf q} - \tilde \xi_{p'}| \lesssim \varrho$ for all $p' \in E\cap {\bf q}$. It follows from Definition \ref{def:kakeya-configuration} (4) and the computation above that we can take $\tilde \xi_{\bf q} = \xi_{{\bf q},2}-\tilde f_2(z) \xi_{{\bf q},1}$ for some $z\in Z({\bf z})$ with $E_z \cap {\bf q}\neq\emptyset$ (Since $\sigma \le \varrho$, the choice of $z$ does not matter.)


For $z, z' \in Z({\bf z})$ define $S_{z,z'}$ to be the set of triples $(y, x_1, x_2)$ such that $y \in Y_{z}, x_1, x_2\in X_{y, z}$ and such that for points $p_i = (x_i, y+f(z)x_i, z) \in E_z$ there exists tubes $T_i \in \TT_{p_i}$ such that $E_{z'} \cap T_i \neq \emptyset$ for $i=1,2$. By double counting we have 
\[
\frac{1}{|Z({\bf z})|^2}\sum_{z,z' \in Z( {\bf z})} |S_{z,z'}| \gtrsim_{\d^{O(\eta)}} (1/\d )^{d-k-\kappa} (1/\d)^{2(k-1)},
\]
(note that the right hand side is the typical size of $|Y_z| |X_{y,z}|^2$). Now let $(y, x_1, x_2) \in S_{z,z'}$ and write $p_i = (x_i, y+ \tilde f(z) x_i, z)$ and $T_i \in \TT_{p_i}$ be tubes so that $z' \in Z(T_i)$. We can write $\theta(T_i) \in N_\d(\varphi_i, \tilde \xi_{p_i} + f(z) \varphi_i)$. It follows that $p_i'=p_i + (z'-z) \theta(T_i) \in N_{C\d} E_{z'}$. Let $t = z'-z$ and compute
\[
p'_i=p_i + t\theta(T_i) = (x_i + t \varphi_i, y + t\tilde \xi_{p_i} +\tilde  f(z) (x_i + t\varphi_i), z).
\]
On the other hand, by writing $p'_i = (x_i', y'_i + \tilde f(z') x_i', z')+O(\d)$ for some $y'_i \in Y_{z'}$ we see that $y + t\tilde \xi_{p_i} + (\tilde f(z)-\tilde f(z')) (x_i + t\varphi_i) \in N_{C\d} Y_{z'}$. Using the condition $|\tilde f(z)-\tilde f(z')| \le \sigma$ we conclude that $y + (z'-z)\tilde  \xi_{p_i} \in N_{C\sigma}Y_{z'}$. Now if we let ${\bf q}_i = p_i^\varrho$ then using $|\tilde \xi_{p_i} - \tilde \xi_{{\bf q}_i}|\lesssim \varrho$ and $\sigma \ge \varrho^2$, we get $y+ (z'-z) \tilde\xi_{{\bf q}_i} \in N_{C\sigma}Y_{z'}$. Let $y' \in Y_{z'}$ be such that $|y' - y - (z'-z) \tilde\xi_{{\bf q}_1}| \lesssim \sigma$. We then get $y' + (z'-z) (\tilde \xi_{{\bf q}_2} - \tilde \xi_{{\bf q}_1}) \in N_{C\sigma} Y_{z'}$. By averaging the lower bound on $S_{z,z'}$ over $z'$ and $x_1,x_2$ we get that there exists a choice of $z'$ such that
\[
\frac{1}{\#\{({\bf q}_1, {\bf q}_2)\}}\sum_{{\bf q}_1, {\bf q}_2}\#\{ (y', z) \in Y_{z'}\times Z({\bf z}):~  y'+ (z'-z)(\tilde \xi_{{\bf q}_2} - \tilde \xi_{{\bf q}_1}) \in N_{C\sigma} Y_{z'} \} \gtrsim_{\d^{O(\eta)}} |Y_{z'}| |Z({\bf z})|
\]
where the sum is taken over all pairs of dyadic $\varrho$-cubes ${\bf q}_1, {\bf q}_2$ covering $[-1,1]^{k-1}$. Now suppose that $y' \in Y_{z'}$ and ${\bf q}_1, {\bf q}_2$ is a pair of cubes for which there are at least $\d^{C\eta}|Z({\bf z})|$ many choices for $z$ satisfying the above. Recall that $|Z({\bf z})|_\d \gtrsim \d^{O(\varepsilon+\eta)} (\varrho/\d)$. Let $\tau = |\tilde \xi_{{\bf q}_1}-\tilde\xi_{{\bf q}_2}|$ and assume that $\tau \ge \sigma/\varrho$. Let $T$ be a $\sim \sigma \times \ldots \times \sigma \times \tau\varrho$ tube in $\RR^{k-1}$ passing through $y'$ in direction $\tilde\xi_{{\bf q}_1}-\tilde\xi_{{\bf q}_2}$. Then we conclude that 
\[
|Y_{z'} \cap T|_\sigma \gtrsim \d^{O(\varepsilon+\eta)} (\tau \varrho / \sigma ). 
\]
On the other hand, by the definition of a Kakeya configuration, the set $Y_{z'}$ is $(\d, \gamma, \d^{-O(\eta)})$-tube Katz-Tao. In particular, we have 
\[
|Y_{z'} \cap T| \lesssim_{\d^{O(\eta)}} (\varrho\tau/\sigma)^{1 - \gamma}.
\]
By combining the bounds we conclude that 
\[
(\varrho\tau/\sigma)^{\gamma} \lesssim \d^{-O(\varepsilon+\eta)},
\]
i.e. $\tau \lesssim  \d^{-C (\varepsilon+\eta)/\gamma} \sigma /\varrho$. We assumed that $\tau \ge \sigma/\varrho$, so in either case we conclude that $\tau$ satisfies the upper bound above. 
Since this holds for a typical pair of dyadic cubes ${\bf q}_1,{\bf q}_2$ we get that $\xi(x, y,z)$ concentrates in intervals of length $\d^{-O ((\varepsilon+\eta)/\gamma)}\sigma/\varrho$, as desired. 
\end{proof}

\subsection{Slope functions are Lipschitz}\label{subsec:lipschitz}

In this section we show that after pigeonholing both slope functions $f$ and $\tilde f$ can be assumed to be dyadic Lipschitz. It is part of the data of Definition \ref{def:kakeya-configuration} (4) that $f$ is Lipschtitz, however we need to verify separately that so is $\tilde f$. 

\begin{lemma}\label{lem:f-is-Lipschitz}
    Consider an $\ell$-linear $k$-planar  $(\d, \eta, \kappa, \gamma)$-Kakeya configuration and let $f, \tilde f$ be its slope functions. Let $\varrho \in [\d^{1/3}, 1]$. Then there exists $\tilde Z \subset Z$ with $|\tilde Z| \gtrsim \d^{O(\eta)} |Z|$ such that for every $z_0 \in \tilde Z$ and ${\bf z} = z_0^\varrho$ we have for every $z, z' \in \tilde Z \cap {\bf z}$ that
    \[
    |f(z) - f(z')| \le  \d^{-O(\eta/\gamma)} \varrho
    \]
    \[
    |\tilde f(z) - \tilde f(z')| \le  \d^{-O(\eta/\gamma)} \varrho.
    \] 
\end{lemma}

\begin{proof}

    By Item $(2)$ of the definition of $\ell$-linear Kakeya configuration (Definition~\ref{def:kakeya-configuration}), for any $z_1, z_2\in Z\cap \mathbf{z}$, 
    \[
    |f(z) - f(z')|  \lesssim \varrho. 
    \]

    Now pick $\theta \in \Theta_{\bf q}$ such that $\theta^\varrho$ intersects $\d^{O(\eta)}$-fraction of sets $\Theta_p$ for $p \in E({\bf q})$. By double counting, we then have 
    \begin{equation}\label{eq:average-shift}
    \frac{1}{|Z({\bf q})|^2}\sum_{z, z' \in Z({\bf q})}| (E_z({\bf q}) + (z'-z) \theta) \cap N_{C\varrho^2} E_{z'}({\bf q})| \gtrsim_{\d^{O(\eta)}} (\varrho/\d)^{d-1-\kappa}.    
    \end{equation}
    (The right hand side is the typical size of $E_z({\bf q})$). Now we write
    \[
    E_z({\bf q}) = \{ (x, y + \tilde f(z) x, z):~ y \in Y_z({\bf q}), ~ x \in X_{y,z}({\bf q}) \}
    \]
    Denote $\theta = (\varphi, \psi) \in [-1,1]^{k-1}\times [-1,1]^{d-k}$.
    Now fix a pair $z, z' \in Z({\bf q})$ which have at least average intersection in \eqref{eq:average-shift}. Let $p=(x, y+f(z)x,z) \in E_z({\bf q})$ be such that $p'=p+ (z'-z)\theta \in N_{C\varrho^2}E_{z'}({\bf q})$ holds, in other words, there is a $C\varrho^2$-tube in direction $(\theta, 1)$ that contains both $p'$ and $p$.  If we write $p' = (x', y'+f(z')x', z')$, then 
    \begin{align*}
    x' &= x + (z'-z) \varphi + O(\varrho^2),\\
    y'+ \tilde f(z')x' &= y + \tilde f(z) x  + (z'-z)\psi + O(\varrho^2),
    \end{align*}
    from which it follows that
    \[
    y' = y+ (\tilde f(z) - \tilde f(z')) x + (z'-z) (\psi - \tilde f(z') \varphi) + O(\varrho^2)
    \] 
    
    Let $v= y+ (\tilde f(z) - \tilde f(z')) x + (z'-z) (\psi - \tilde f(z') \varphi)$. Then by \eqref{eq:average-shift}, for $\d^{O(\eta)}$-fraction of $\tilde x \in X_{y, z}({\bf q})$ by using $\tilde p = (\tilde x, y+f(z)\tilde x, z)$ in place of $p$, we get that 
    \[
    y + (\tilde f(z) - \tilde f(z')) \tilde x + (z'-z) (\psi - \tilde f(z') \varphi) \in N_{C\varrho^2}Y_{z'}
    \]
    and, therefore,
    \[
    v + (\tilde f(z)- \tilde f(z')) (\tilde x - x) \in N_{C\varrho^2} Y_{z'}.
    \] 
    So using the fact that we can choose $\tilde x$ from a $\d^{O(\eta)}$-dense subset of a $(k-1)$-dimensional slice of ${\bf q}$, we get that $N_{C\varrho^2}Y_{z'}$ is $\d^{O(\eta)}$-dense in the image $v + (\tilde f(z)- \tilde f(z')) \cdot [-\varrho,\varrho]^{k-1} + B_{C\varrho^2}$. Let $U = (\tilde f(z)- \tilde f(z')) \cdot [-\varrho,\varrho]^{k-1} + B_{C\varrho^2}$. Denote $\tau = \|\tilde f(z)-\tilde f(z')\|$ and suppose that $\tau \ge \varrho$. Thus, the diameter of $U$ is $\sim \varrho \tau$. 
    Decompose $U$ into a finitely overlapping union of parallel $\varrho^2\times\ldots\times \varrho^2\times \varrho\tau$-tubes $T_i$. By the $(\d, \gamma, \d^{-O(\eta)})$-tube Katz-Tao axiom satisfied by $Y_{z'}({\bf q})$, we then conclude that 
    \[
    |Y_{z'}({\bf q}) \cap U|_{\varrho^2} \lesssim \sum_i |Y_{z'}({\bf q}) \cap T_i|_{\varrho^2} \lesssim \d^{-O(\eta)}\sum_i (\tau/\varrho)^{1-\gamma} \sim \d^{-O(\eta)} |U|_{\varrho^2} (\tau/\varrho)^{-\gamma}.
    \]
    Since $N_{\varrho^2} Y_z'$ is $\delta^{O(\eta)}$-dense in $U$,  this implies that $\tau \lesssim \d^{-O(\eta/\gamma)}\varrho$ holds, giving the desired Lipschitz property. 
\end{proof}

\section{\texorpdfstring{$(d-1)$}{d-1}-linear configurations: Proof of Proposition \ref{prop:codimension-1-kakeya}}
\label{sec:codimension-1}
In this section we consider $(d-1)$-planar, $(d-1)$-linear Kakeya configurations, which means $k = \ell =d-1$ and $f=\tilde f$. We show that any such configuration must be full-dimensional.


\subsection{\texorpdfstring{$V$}{V}-tuples and \texorpdfstring{$L$}{L}-tuples.}  \label{subsec:L-tuple-kakeya}

Fix a $k$-planar $k$-linear Kakeya configuration
$(\TT,E,\{\TT_p\})$, $p_0\in E$,  and a dyadic $\varrho$-cube
${\bf q}=p_0^\varrho$. 
Let $Z(\mathbf{q}), E(\mathbf{q})=E\cap \mathbf{q}$ and $\TT(\mathbf{q}) =\cup_{p\in E(\mathbf{q})} \TT_p$ be defined in the usual way. 

\medskip\noindent\emph{Basic $V$-tuples.}
For $z,z'\in Z({\bf q})$, let $V_{z,z'}({\bf q})$ consist of the
tuples
\[
(p,p_1,p_2,T_1,T_2)
\in E_z({\bf q})\times E_{z'}({\bf q})^2\times\TT({\bf q})^2
\]
such that
\[
T_i\in\TT_p\cap\TT_{p_i},
\qquad i=1,2.
\]
Their total number satisfies
\begin{equation}\label{eq:V-tuples-counting}
\sum_{z,z'\in Z({\bf q})}|V_{z,z'}({\bf q})|
\sim_{\d^{O(\eta)}}
|E({\bf q})|\d^{-2\kappa}|Z({\bf q})|.
\end{equation}
Indeed, the upper bound follows by selecting $p$ and then $T_1,T_2$ and $p_1,p_2$ in at most $|E({\bf q})|$, $\sim \d^{-2\kappa}$ and $\sim |Z({\bf q})|$ ways, respectively. For the lower bound, count the number of $(p, p_1, T_1)$ and $(p,p_2,T_2)$ and use Cauchy--Schwarz.


For a $V$-tuple, write
\[
p=(x,y+f(z)x,z),\qquad
p_i=(x_i,y_i+f(z')x_i,z'),
\]
and let $\varphi_i$ be the first $k-1$ coordinates of
$T_i$, so that $\theta(T_i) = (\varphi_i, \xi_p +f(z) \varphi_i) + O(\d)$.  Subtracting
$p_i=p+(z'-z)(\theta(T_i),1)+O(\d)$ for $i=1,2$ gives
\[
\begin{aligned}
x_1-x_2&=(z'-z)(\varphi_1-\varphi_2)+O(\d),\\
y_1-y_2+f(z')(x_1-x_2)
&=(z'-z)f(z)(\varphi_1-\varphi_2)+O(\d).
\end{aligned}
\]
Consequently,
\begin{equation}\label{eq:V-tuple-relation}
y_1-y_2
=(z'-z)(f(z)-f(z'))(\varphi_1-\varphi_2)+O(\d).
\end{equation}
We use this identity in the linear case. 

\medskip\noindent\emph{One-arm paths.}
For $z\in Z({\bf q})$ and $p,\widetilde p\in E_z({\bf q})$, write
$p\leftrightarrow\widetilde p$ if there is a $y\in Y_z$
such that
\[
p=(x,y+f(z)x,z),\qquad
\widetilde p=(\widetilde x,y+f(z)\widetilde x,z)
\]
for some $x,\widetilde x\in X_{y,z}$.  That is,
$p$ and $\widetilde p$ lie in the same horizontal $(k-1)$-grain.

For $z,z',z''\in Z({\bf q})$, let $J_{z,z',z''}$ be the set of
one-arm paths
\[
(p,p',\widetilde p',p'',T,T')
\in E_z({\bf q})\times E_{z'}({\bf q})^2
\times E_{z''}({\bf q})\times\TT({\bf q})^2
\]
satisfying
\begin{equation}\label{eq:J-incidences-kakeya}
T\in\TT_p\cap\TT_{p'},
\qquad
p'\leftrightarrow\widetilde p',
\qquad
T'\in\TT_{\widetilde p'}\cap\TT_{p''}.
\end{equation}

We have
\begin{equation}\label{eq:J-tuples-counting-kakeya}
\sum_{(z,z',z'')\in Z({\bf q})^3}|J_{z,z',z''}|
\sim_{\d^{O(\eta)}}
|E({\bf q})|\d^{-2\kappa}|Z({\bf q})|^2
(\varrho/\d)^{k-1}.
\end{equation}
To verify this, choose successively $p,~ T,~(z',p'),~ \widetilde p',~
 T',~(z'',p'')$ 
 and estimate the number of choices at each step using properties in Definition \ref{def:kakeya-configuration}. 

\medskip\noindent\emph{Definition of the $L$-tuples.}
Fix $\tau\in[\d,\varrho]$.  For $p\in E({\bf q})$ and
$\widetilde x\in(\tau\ZZ)^{k-1}$, let
$J_{z,z',z''}(p,\widetilde x)$ be the subset of tuples $(p,p',\widetilde p',p'',T,T')$ in
$J_{z,z',z''}$ whose initial point is $p$ and for which the
$x$-coordinate of $\widetilde p'$ lies in the $\tau$-neighborhood of
$\widetilde x$.  Define
\[
L_{z,z',z''}(p,\widetilde x)
\subset
J_{z,z',z''}(p,\widetilde x)^2
\]
to consist of ordered pairs $((p, p_1', \widetilde p_1', p_1'', T_1, T_1'), (p,  p_2', \widetilde p_2', p_2'', T_2, T_2'))$ whose terminal tubes $T_1', T_2'$ satisfy
$|\theta(T_1')-\theta(T_2')|\leq\tau$, and put
\[
L_{z,z',z''}
=\bigcup_{p,\widetilde x}
L_{z,z',z''}(p,\widetilde x).
\]
See Figure \ref{fig:L-tuple} for an illustration. 

\medskip\noindent\emph{Counting the $L$-tuples.} The following lemma gives an estimate on the total number of $L$-tuples for an appropriate choice of $\tau$.

\begin{lemma}\label{lem:number-of-Lambda-tuple-kakeya} Let $\tau\in [\delta, \rho]$ be the fixed parameter in the definition of $L$-tuples. 
Let $\widetilde{Z} \subset Z({\bf q})$ be a subset and suppose that, for every $z,z'\in \widetilde Z$, $p\in E({\bf q})$,
and $\varphi\in\Phi_p$, we have $|(f(z)-f(z'))\varphi|\leq c\tau/\varrho$
for a sufficiently small absolute constant $c>0$.  Then
\[
\sum_{(z,z',z'')\in \widetilde Z^3}|L_{z,z',z''}|
\sim_{\d^{O(\eta)}}
|E_{\widetilde Z}({\bf q})|\d^{-4\kappa}\tau^\kappa
(\varrho/\d)^{k-1}(\tau/\d)^{k-1}
|\widetilde Z|^2.
\]
\end{lemma}

\begin{proof}

Note that if $((p, p'_1, \tilde p'_1, p''_1, T_1, T'_1), (p, p'_2, \tilde p'_2, p''_2, T_2, T'_2)) \in L_{z,z',z''}$ then $(p, p_1', p_2', T_1, T_2)$ forms a $V$-tuple.
So using (\ref{eq:V-tuple-relation}) we get that $p_1'$ and $p_2'$ have their $y$-coordinates differ by $(z'-z) (f(z) - f(z')) (\varphi_1-\varphi_2) +O(\d) = O( \tau)$. Since $\tilde p'_i$ lies in the same global grain as $p_i'$ for $i=1,2$, $\tilde p_1'$ and $\tilde p_2'$ also have their $y$-coordinates differ by $O(\tau)$. So given a choice of $(p, p_1', p_2', T_1, T_2)$, by Cauchy--Schwarz there are typically $\d^{O(\eta)}(\varrho/\d)^{k-1} (\tau/\d)^{k-1}$ ways to select $\tilde p_1'$ and $\tilde p_2'$. Indeed, we first select $\tilde p_1'$ in $|X_{y,z}({\bf q})| \ge \d^{O(\eta)} (\varrho/\d)^{k-1}$ ways and then choose $\tilde p_2'$ on the global grain of $p_2'$ and within distance $\tau$ from $\tilde p_1'$ in $\d^{O(\eta)} (\tau/\d)^{k-1}$ ways (on average). 
    Since both $x$-coordinates of $\tilde p_1'$ and $\tilde p_2'$ are within $\tau$ from some fixed point $\tilde x$, we get $|\tilde p_1' - \tilde p_2'| = O(\tau)$.

    For a given $\tilde p \in E({\bf q})$ let $\TT({\bf q}, \tilde p) \subset \TT({\bf q})$ be the set of tubes $T$ such that $T \in \TT_{\tilde p'}\cap \TT({\bf q})$ for some $\tilde p' \in N_{C\tau}\tilde p$ (where $C$ is a fixed absolute constant). By using the property (4) of Definition \ref{def:kakeya-configuration} with $\tau$ in place of $\varrho$, we have
    \[
    \frac{1}{|E_{\tilde Z}({\bf q})|}\sum_{\tilde p \in E_{\tilde Z}({\bf q})} |\TT({\bf q}, \tilde p)|_\tau \sim_{\d^{O(\eta)}} \tau^{-\kappa}.
    \]
    It follows that when we average over all $(z, z', z'') \in \widetilde Z^3$, all $p, \xi, \tilde x$ and all pairs in  $J_{z,z', z''}(p, \xi, \tilde x)\times J_{z,z', z''}(p, \xi, \tilde x)$,  the intersection $ N_\tau \TT_{\tilde p_1'} \cap \TT_{\tilde p_2'}$ has $\tau$-covering number $\d^{O(\eta)}\tau^{-\kappa}$.
    So, on average, there are $\sim_{\d^{O(\eta)}} \tau^\kappa\d^{-2\kappa}$ ways to choose a pair of tubes $T_1' \in \TT_{\tilde p'_1} \cap \TT({\bf q})$ and $T_2' \in \TT_{\tilde p'_2} \cap \TT({\bf q})$ with $|\theta(T'_1) - \theta(T'_2)| \le \tau$. Once these tubes are selected, we can choose the points $p_1'', p_2''$ on the $z''$-slice in an essentially unique way. 
    Altogether, this gives us the claimed estimate on the sum of sizes of $L_{z,z',z''}$. To make the `on average' type of claims rigorous, we need to use the Cauchy--Schwarz inequality to lower bound the number of pairs of $J$-tuples similarly to how it was done in Lemma \ref{lem:W-tuple-many-Kakeya}. We omit the details.

\end{proof}

\medskip\noindent\emph{Additive relations.}
Fix an $L$-tuple and, for $i=1,2$, write
\[
\begin{aligned}
p_i'&=(x_i',y_i'+f(z')x_i',z'),\\
\widetilde p_i'
&=(\widetilde x_i',\widetilde y_i'
  +f(z')\widetilde x_i',z'),\\
p_i''&=(x_i'',y_i''+f(z'')x_i'',z''),\\
\theta(T_i)&=(\varphi_i,\xi_p+f(z)\varphi_i)+O(\d).
\end{aligned}
\]
The subtuple $(p, p_1', p_2', T_1, T_2)$ forms a $V$-tuple, so
\eqref{eq:V-tuple-relation} gives
\begin{equation}\label{eq:L-first-y-kakeya}
y_1'-y_2'
=(z'-z)(f(z)-f(z'))(\varphi_1-\varphi_2)+O(\d).
\end{equation}
Moreover, $\widetilde y_i'=y_i'$ because $p_i'$ and $\widetilde{p}_i'$ lie in the same global grain. 
We have
\[
p_i''=\widetilde p_i'
+(z''-z')(\theta(T_i'),1)+O(\d).
\]
Since $|z''-z'|\leq\varrho$ and
$|\theta(T_1')-\theta(T_2')|\leq\tau$, subtracting gives
\[
p_1''-p_2''
=\widetilde p_1'-\widetilde p_2'+O(\tau\varrho+\d).
\]
Consequently,
\[
\begin{aligned}
x_1''-x_2''
&=\widetilde x_1'-\widetilde x_2'+O(\tau\varrho+\d),\\
y_1''-y_2''
&=\widetilde y_1'-\widetilde y_2'
+(f(z')-f(z''))(\widetilde x_1'-\widetilde x_2')
+O(\tau\varrho+\d).
\end{aligned}
\]
Because the two $\widetilde x_i'$ lie in the same $\tau$-ball,
$|\widetilde x_1'-\widetilde x_2'|=O(\tau)$.  Combining the above
formulas with \eqref{eq:L-first-y-kakeya} gives
\begin{equation}\label{eq:y-axes-another-relation-kakeya}
\begin{split}
y_1''-y_2''
&=(z'-z)(f(z)-f(z'))(\varphi_1-\varphi_2)\\
&\quad
+O\bigl(\tau|f(z')-f(z'')|+\tau\varrho+\d\bigr).
\end{split}
\end{equation}
In particular, dyadic Lipschitz continuity of $f$ makes the error
$O(\tau\varrho+\d)$.  We will use this identity in the non-linear case below.

\subsection{Linear case.}\label{subsec:linear-case}
In the following lemma we view $\operatorname{Mat}(d-k, k-1)$ as a subset in $\RR^{(d-k)(k-1)}$ and write 
$|r| $ to denote the $\ell_\infty$-norm of entries of a matrix $r \in \operatorname{Mat}(d-k, k-1)$.

\begin{lemma}\label{lem:codimension-1-linear-case}
Let $\varepsilon_0, \mu, \gamma>0$ and $\varepsilon, \eta>0$ be sufficiently small depending on $\varepsilon_0, \mu, \gamma$.

    Consider a $k$-planar $k$-linear $(\d, \eta, \kappa, \gamma)$-Kakeya configuration in $\RR^d$ with the slope function $f: Z \to \operatorname{Mat}(d-k, k-1)$. 

    Let $\varrho \in [\d^{1/4},\d^{\mu}]$. 
    Suppose that there is $\tilde Z \subset Z$ with $|\tilde Z| \ge \d^{\varepsilon} |Z|$ such that for any $z_0 \in \tilde Z$ and ${\bf z} = z_0^{\varrho}$    
    there exist $r, u \in \operatorname{Mat}(d-k, k-1)$ such that
    \[
    |f(z) - r z - u| \le \varrho^{1+\varepsilon_0} |r|, \quad z \in \tilde Z \cap {\bf z}.
    \]
    Then we have $C_{CW}(\TT) \gtrsim \d^{O(\eta+\varepsilon)} \varrho^{-\kappa}$.
\end{lemma}

\begin{proof}
    First, we note that if $|r| \le \varrho$ holds for at least half of the $\varrho$-intervals ${\bf z}$ covering $\tilde Z$ then for any $z, z' \in \tilde Z \cap {\bf z}$ we have 
    \[
    |f(z) - f(z')| \lesssim \varrho^{2}.
    \]
    By Lemma \ref{lem:f-slow} we then conclude that 
    \[
    C_{CW}(\TT) \gtrsim \d^{O(\eta+\varepsilon)} \varrho^{- \kappa},
    \]
    giving the desired conclusion. 
    Thus, without loss of generality we may assume that $|r| \ge \varrho$ holds for at least half of the $\varrho$-intervals ${\bf z}$. By passing to a dense subset in $\tilde Z$ we may assume that this holds for all $\varrho$-intervals.

Let $p_0\in E_{z_0}\subset E$ and 
  ${\bf q} = p_0^\varrho$ and ${\bf z} = z_0^\varrho$. 
    Let $Z({\bf q}) = \tilde Z\cap {\bf z}$ and $E({\bf q})=\cup_{z\in \tilde Z} E_z \cap {\bf q}$. 
    Consider the $V$-tuples $V_{z, z'}({\bf q})$ defined using the sets $E({\bf q})$ and $Z({\bf q})$ in Section \ref{subsec:L-tuple-kakeya}. By (\ref{eq:V-tuples-counting}), there exists some $z' \in Z({\bf q})$ and $p_1 = (x_1, y_1 + f(z')x_1, z')$ such that defining $$V_{p_1}=\{(p, p_1, p_2, T_1, T_2)\in V_{z, z'}, z\in Z(\bf q)\},$$ we have $|V_{p_1}|\gtrsim \delta^{O(\eta)} \delta^{-2\kappa} |Z(\bf q)|$.  For the fixed $p_1$, by choosing $T_1, p, T_2$ consecutively, we also have $|V_{p_1}|\lesssim \delta^{-O(\eta)} \delta^{-2\kappa} |Z(\bf q)|$. 
    Let $\overline{\Phi}_{{\bf q}} \subset (\varrho \ZZ)^{k}$ be the set of points so that  $\overline{\Phi}_{{\bf q}}^\varrho = \Phi_{\bf q}^\varrho$. 
    For each $(p, p_1, p_2, T_1, T_2)\in V_{p_1}$, we associate it with the quadruple $(z, \overline{\varphi}_1, \overline{\varphi}_2, y_2)$, where $\overline{\varphi}_i \in \overline{\Phi}_{{\bf q}}$ are chosen so that $ \theta(T_i) \in N_\d (\varphi_i, \xi_p+f(z) \varphi_i)$ and we have $\varphi_i \in \overline{\varphi}_i^{\varrho}$ and $p_2 = (x_2, y_2 + f(z) x_2, z')$. By fixing $\overline{\varphi}_1$ we determine the slope $\varphi_1$ of $T_1$ up to precision $O(\varrho)$. So for a fixed choice of $p_1$ (with $\overline{\varphi}_1$ fixed)  there are at most $\lesssim_{\d^{O(\eta)}}(\varrho/\d)^\kappa$ ways to fix $T_1$ (since $\TT_{p_1}$ is almost AD-regular). Once $T_1$ is fixed, there are at most $|Z({\bf q})|$ choices for $z$ and $p \in T_1 \cap E_z({\bf q})$. After that given $\overline{\varphi}_2$ there are again at most $\lesssim_{\d^{O(\eta)}}(\varrho/\d)^\kappa$ ways to fix $T_2$.  Once we fixed $(p_1, z, T_1, T_2)$, the V-tuple $(p, p_1, p_2, T_1, T_2)$ is uniquely determined since $T_1, z$ uniquely determines $p$ and $p_1$ uniquely determines $z'$ and $T_2, z'$ uniquely determines $p_2$.  We conclude by double counting that the number of tuples $(z, \overline{\varphi}_1, \overline{\varphi}_2, y_2) \in Z({\bf q}) \times \overline{\Phi}_{{\bf q}}\times \overline{\Phi}_{{\bf q}} \times Y_{z'}({\bf q})$ associated with some V-tuples in $V_{p_1}$ is $\sim_{\d^{O(\eta)}}  |V_{p_1}| \, (\varrho/\delta)^{-2\kappa} \sim_{\d^{O(\eta)}} |Z({\bf q})| \varrho^{-2\kappa}$.

    Recall that (\ref{eq:V-tuple-relation}) implies that
    \[
    y_2 - y_1 = (z'-z) (f(z') - f(z)) (\varphi_1 - \varphi_2) + O(\d).
    \]
    We have $|f(z') - f(z)| = |(z'-z) r| + O(\varrho^{1+\varepsilon_0}|r|) = O(\varrho |r|)$
    and so for our choice of $\overline{\varphi}_i$ we get using the triangle inequality and $|z'-z|\leq \varrho$: 
    \[
     (z'-z) (f(z') - f(z)) (\varphi_1 - \varphi_2) = (z'-z) (f(z') - f(z)) (\overline{\varphi}_1 - \overline{\varphi}_2) + O(\varrho^3|r|+\d).
    \]
    We conclude using $|r| \ge \varrho$ and $\varrho \ge \d^{1/4}$ that for a fixed $p_1$, which determines $z'$ and $y_1$,  any tuple $(z, \overline{\varphi}_1, \overline{\varphi}_2)$ coming from quadruples $(z, \overline{\varphi}_1, \overline{\varphi}_2, y_2)$ associate to V-tuples in $V_{p_1}$, satisfies 
    \[
    y_1  + (z'-z) (f(z') - f(z)) (\overline{\varphi}_1 - \overline{\varphi}_2) \in N_{C\varrho^3|r| }Y_{z'}({\bf q}).
    \]
    We conclude that given $p_1$, 
    \begin{equation}\label{eq: ycondition}
    \#\{ (z, \overline{\varphi}_1, \overline{\varphi}_2) \in Z({\bf q}) \times \overline{\Phi}_{{\bf q}}\times \overline{\Phi}_{{\bf q}}: y_1 + (z'-z) ( f(z) - f(z') )(\overline{\varphi}_1 -\overline{\varphi}_2) \in N_{C\varrho^3 |r| }Y_{z'}({\bf q})  \} \gtrsim_{\d^{O(\eta)}}  |Z({\bf q}) | |\overline{\Phi}_{{\bf q}}|^2
    \end{equation}
    By the assumption on $f$ that $f(z)-f(z') = (z-z') r + O(\varrho^{1+\varepsilon_0} |r|)$. Recall our configuration is a $k$-planar $k$-linear $(\d, \eta, \kappa, \gamma)$-Kakeya configuration. This implies that the set $\overline{\Phi}_{\bf q} \subset (\varrho \ZZ)^{k} \cap [-1,1]^{k}$ is $(\varrho, \gamma, \d^{-O(\eta)})$-slab Frostman. Therefore, for all but $\d^{O(\eta)}$-fraction of pairs $(\overline{\varphi}_1, \overline{\varphi}_2)$ we have $|r (\overline{\varphi}_1-\overline{\varphi}_2)| \gtrsim \d^{O(\eta)} \d^{O(\eta / \gamma)} |r|$. To see this, suppose the $\ell^{\infty}$-norm of $r$ is achieved at the $ij$-th entry $|r|=|r_{ij}|$ and write $\lambda=(\lambda_1, \dots, \lambda_k)=\overline{\varphi}_1-\overline{\varphi}_2$, then the set of $\lambda$ such that $|\sum_{\ell} r_{i\ell}\lambda_\ell|\leq  \delta^{O(\eta/\gamma)} |r_ij|$ is contained in the $\delta^{O(\eta/\gamma)}$-neighborhood of a hyperplane. Since $\overline{\Phi}_{\bf q}$ is $(\varrho, \gamma, \delta^{-O(\eta)})$-slab Frostman, for each $\overline{\varphi}_1$, the number of such $\overline{\varphi}_2$ is at most a $\delta^{O(\eta)}$-fraction.

    Thus, by double counting, we can fix values of $ \overline{\varphi}_1, \overline{\varphi}_2$ for which there are $\d^{O(\eta)}|Z({\bf q})|$ choices for a tuple $(z, \overline{\varphi}_1, \overline{\varphi}_2)$ satsifying \eqref{eq: ycondition} and which satisfy $|r (\overline{\varphi}_1-\overline{\varphi}_2)| \gtrsim \d^{O(\eta/\gamma)} |r|$. Let $v  = r (\overline{\varphi}_1-\overline{\varphi}_2) \in \RR^{d-k}$. Note that then we have 
    \[
    (z'-z) ( f(z) - f(z') )(\overline{\varphi}_1 -\overline{\varphi}_2) = (z'-z)^2 v + O(\varrho^{2+\varepsilon_0}|r|)
    \]
    and so we conclude that for $\Delta  = C\varrho^3|r| + C \varrho^{2+\varepsilon_0}|r| \sim  \varrho^{2+\varepsilon_0}|r| $ we have
    \[
    \#\{ z \in Z({\bf q}):~ y_1 + (z'-z)^2 v \in N_{\Delta } Y_{z'}({\bf q}) \} \gtrsim_{\d^{O(\eta)}} |Z({\bf q})|
    \]
    Recall that  $Z({\bf q}) = \tilde Z \cap {\bf z}$. For a fixed $p_1$ (which determines $z'$ and $y_1$) and $\overline{\varphi}_1, \overline{\varphi}_2$,  Let $\overline{Z}$ be the set of $z\in Z({\bf q})$ such that $y_1 + (z'-z)^2v \in N_{\Delta } Y_{z'}({\bf q})$ and $|z-z'| \gtrsim \d^{\varepsilon+C'\eta}\varrho$. By choosing large enough $C'$, from the above we get the lower bound $|\overline{Z}| \gtrsim_{\d^{O(\eta)}} |Z({\bf q})|$.

    Let $T_{y_1, v} \subset \RR^{d-k}$ be a $C\Delta \times C\varrho^2|r|$ tube centered at $y_1$ and in direction $v$. 
    Since we have $|(z'-z)^2 v| \le C \varrho^2 |r|$, we conclude that we have the following lower bound on the covering number of $Y_{z'}({\bf q}) \cap T_{y_1, v}$:
    \begin{align*}
        |  Y_{z'}({\bf q}) \cap T_{y_1, v}|_\Delta \gtrsim \d^{O(\varepsilon+\eta)}  |\overline{Z}|_{\Delta/(\varrho|v|)} \gtrsim \d^{O(\varepsilon+\eta)} (\varrho^2|v|/\Delta)
    \end{align*}
    Here we use that if $|z_1-z_2|\gtrsim \Delta/(\varrho |v|)$ and $|2z'-(z_1+z_2)|\gtrsim \delta^{\varepsilon +O(\eta)} \varrho$, then $|(z'-z_1)^2 v - (z'-z_2)^2v|\gtrsim  \delta^{O(\varepsilon +\eta)}\Delta$.


    On the other hand, by Definition \ref{def:kakeya-configuration-2}, the set 
    $Y_{z'}({\bf q})$ is $(\d/\varrho, 1-\gamma, \d^{O(\eta)})$-tube Katz-Tao. 
    It follows that
    \[
    | Y_{z'}({\bf q}) \cap T_{y_1, v}|_\Delta  \lesssim_{\d^{O(\eta)}} (\varrho^2 |r| / \Delta)^{1-\gamma}.
    \]
    We conclude that 
    \begin{align*}
        (\varrho^2 |r| / \Delta)^{1-\gamma} &\gtrsim \d^{O(\varepsilon+\eta)} (\varrho^2|v|/\Delta) \\
        &\gtrsim \d^{O(\varepsilon+\eta/\gamma)} (\varrho^2|r|/\Delta)\\
        (\varrho^2 |r| / \Delta)^{\gamma} &\lesssim  \d^{O(\varepsilon+\eta/\gamma)}.
    \end{align*}
    Now recall that $\Delta  \sim \varrho^{2+\varepsilon_0}|r|$ and so we conclude that $\gamma \mu\varepsilon_0 \le O( \varepsilon+\eta/\gamma ) $. By choosing
    $\varepsilon, \eta$ sufficiently small in $\varepsilon_0, \mu, \gamma$ we arrive at a contradiction. 
\end{proof}

\subsection{Non-linear case.}\label{subsec:non-linear-case-Kakeya}

In this section we show that there are no $(d-1)$-linear Kakeya configurations whose slope function cannot be locally approximated by a linear function.

\begin{lemma}\label{lem:codimension-1-non-linear-case}
Let $\zeta, \kappa,\gamma, \varepsilon_0, \mu>0$ be parameters 
and $\varepsilon>0$ is  sufficiently small depending on $\zeta, \mu, \gamma$ and $\eta$ is sufficiently small depending on all previous parameters. Finally, $\delta>0$ is sufficiently small depending on all previous parameters. 

    Consider a $(d-1)$-planar $(d-1)$-linear $(\d, \eta, \kappa, \gamma)$-Kakeya configuration in $\RR^d$ with the slope function $f: Z \to \operatorname{Mat}(1, d-2)$. 

    Let $\varrho \in [\d^{1/4},\d^{\mu}]$. 
    Suppose that for any $z_0 \in Z$ and ${\bf z} = z_0^{\varrho}$ and for all $r, u \in \RR^{d-2}$ we have 
    \begin{equation}\label{eq:f-non-linear-codim-1}
    \#\{ z \in Z \cap {\bf z} ~\mid ~ |f(z) - r z - u| \le \varrho^{1+\varepsilon_0} \} \le \varrho^{\varepsilon} (\varrho/\d).    
    \end{equation}
    Then we have $\kappa \le \zeta$.
\end{lemma}

\begin{proof}[Proof of Lemma \ref{lem:codimension-1-non-linear-case}]

    Let $p_0 \in E_{{z_0}} \subset E$ be an arbitrary point and denote ${\bf q} = p_0^\varrho$ and ${\bf z} = z_0^\varrho$. 
    Let $E({\bf q}) = E \cap {\bf q}$ and $Z({\bf q}) = \tilde Z\cap {\bf z}$.
    Consider the $L$-tuples $L_{z, z', z''}$ 
    defined using the sets $E({\bf q})$ and $Z({\bf q})$ and $\tau \sim \d^{- \eta}\varrho^2$. Here $\tau$ was chosen so that we can apply Lemma~\ref{lem:number-of-Lambda-tuple-kakeya}. Let $\overline{\Phi}_{{\bf q}} \subset (\varrho \ZZ)^{d-1}$ be the set of points so that  $\overline{\Phi}_{{\bf q}}^\varrho = \Phi_{\bf q}^\varrho$. 
    By Lemma \ref{lem:number-of-Lambda-tuple-kakeya} and double counting we can find some $z'' \in Z({\bf q})$ such that
    \[
    \sum_{z,z' \in Z({\bf q})} |L_{z, z', z''}| \gtrsim_{\d^{O(\eta)}} |E_{z''}({\bf q})| \d^{-4\kappa} \tau^\kappa  (\varrho/\d)^{d-2} (\tau/\d)^{d-2}|Z({\bf q})|^2.
    \]
    Note that we have $E_{z''}(\bf q)$ instead of $E(\bf q)$ on the  RHS of the inequality. 
    Define a graph 
    \begin{equation}\label{eq: graphG}
        G \subset Y_{z''}({\bf q}) \times Z({\bf q})\times Z({\bf q}) \times \overline{\Phi}_{{\bf q}} \times \overline{\Phi}_{{\bf q}}
    \end{equation} in the following way. Consider a tuple 
    $((p, p'_1, \tilde p'_1, p''_1, T_1, T'_1), (p, p'_2, \tilde p'_2, p''_2, T_2, T'_2)) \in L_{z,z',z''}$ and write $p_i'' = (x_i'', y''_i + f(z_0) x_i'', z_0)$. Let $\varphi_i$ be the first coordinate of $\theta(T_i)$. Since $T_i\in \TT({\bf q})$, there exists $\overline{\varphi}_i \in \overline{\Phi}_{\bf q}$ such that $\varphi_i$     
    lies in the dyadic $\varrho$-cube $\overline{\varphi}_i^{\varrho}$ for some $\overline{\varphi}_i \in \overline{\Phi}_{\bf q}$.
    Then we include the tuple $(y_1'', z,z',\overline{\varphi}_1,\overline{\varphi}_2)$ in $G$. 

    By (\ref{eq:y-axes-another-relation-kakeya}), $|\varphi_i -\overline{\varphi}_i|\lesssim \varrho$, and  $f$ is dyadic Lipchitz (see Definition~\ref{def:kakeya-configuration}), we then get
    \[
    y_2'' - y_1'' = (z'-z) (f(z') - f(z)) (\overline{\varphi}_1- \overline{\varphi}_2) + O(\tau \varrho +\d)
    \]
    Note that $\varrho \tau \ge \varrho^3 \gg \d$. 
    In other words, for any $(y_1'', z,z',\overline{\varphi}_1,\overline{\varphi}_2) \in G$ it holds that
    \[
    y_1'' + (z'-z) (f(z') - f(z)) (\overline{\varphi}_1- \overline{\varphi}_2) \in N_{C \tau\varrho} Y_{z''}({\bf q}).
    \]
    By estimating the number of ways to recover an $L$-tuple from an edge in $G$, the bound in Lemma \ref{lem:number-of-Lambda-tuple-kakeya} implies that $|G| \gtrsim_{\d^{O(\eta)}} |Y_{z''}({\bf q})| |Z({\bf q})|^2 |\overline{\Phi}_{\bf q}|^2$, in other words $G$ is a $\d^{O(\eta)}$-dense graph. 

    Now cover $Y_{z''}({\bf q})$ by finitely overlapping intervals $I_j$ of length $C \tau$ so that sets $Y_j =  Y_{z_0}({\bf q}) \cap I_j$ have $\d$-covering number $\sim_{\d^{O(\eta)}} (\tau/\d)^{1-\kappa}$.
    Let $G_j$ be the graph $G$ restricted to $Y_j$. By the pigeonhole principle, we can find an index $j$ so that $|G_j| \gtrsim_{\d^{O(\eta)}} |Y_j| |Z({\bf q})|^2 |\overline{\Phi}_{\bf q}|^2$.

    Let $\mc A \subset [-1,1]$ be a rescaling of the set $Y_j$ mapping $I_j$ to the unit interval and let $\mc B = \{ \frac{z-z''}{\varrho}, ~z \in Z({\bf q}) \} \subset [-1,1]$.
    Define $F(z) = \frac{f(z)- f(z'')}{\tau/\varrho} \in \operatorname{Mat}(1, d-2)$. 
    Let $\mc G \subset \mc A\times \mc B\times \mc B \times \overline{\Phi}_{\bf q} \times \overline{\Phi}_{\bf q}$ be the rescaled version of the graph $G_j$. With these definitions, we have the following property:
    \begin{equation}\label{eq:non-expanding-kakeya}
        |\{ a + (b' - b) (F(b') - F(b)) (\overline{\varphi}_1-\overline{\varphi}_2), \quad (a, b,  b', \overline{\varphi}_1, \overline{\varphi}_2) \in \mc G  \}|_{\varrho} \lesssim_{\d^{O(\eta)}} |\mc A|_\varrho.
    \end{equation}

For each choice of $\overline{\varphi}_1, \overline{\varphi}_2 \in \overline{\Phi}_{\bf q}$ this gives us an instance of Theorem~\ref{thm:ABC-with-F} with the function $F_{\overline{\varphi}_1- \overline{\varphi}_2}(z) := F(z) (\overline{\varphi}_1- \overline{\varphi}_2)$ and $C=\{1\}$. For a large constant $K$ let $\mc P \subset  \overline{\Phi}_{\bf q}\times  \overline{\Phi}_{\bf q}$ be the set of pairs $(\overline{\varphi}_1, \overline{\varphi}_2)$ such that $|\mc G|_{\overline{\varphi}_1, \overline{\varphi}_2}| \ge \d^{K \eta} |\mc A| |\mc B|^2$. By choosing $K$ appropriately, we obtain that $|\mc P| \gtrsim \d^{O(\eta)} |\overline{\Phi}_{\bf q}|^2$. 
Recall that by the definition of a $(d-1)$-planar $(d-1)$-linear configuration, the set of directions $\overline{\Phi}_{\bf q}$ is $(\varrho, \gamma, \d^{-O(\eta)})$-slab Frostman. On the other hand, since $Z({\bf q}) \subset Z\cap {\bf z}$ is a subset of density $\d^{O(\eta)}$, the condition (\ref{eq:f-non-linear-codim-1}) implies that for all $r, u \in \RR^{d-2}$ we have
\[
\#\{ b\in \mc B:~ |F(b) - b r - u| \le \varrho^{\varepsilon_0} (\varrho^2/\tau) \} \le \varrho^{\varepsilon} (\varrho/\d).
\]
Recall that $\tau \sim \d^{- \eta} \varrho^{2 }$. 
So if we take $\eta$ small enough in $\varepsilon\mu$ and $\varepsilon_0$ then this implies 
\begin{equation}\label{eq:big-F-non-linear}
\#\{ b\in \mc B:~ |F(b) - b r - u| \le \varrho^{2\varepsilon_0}  \} \le \varrho^{\varepsilon/2} |\mc B|.    
\end{equation}

Let $\mc P' \subset \mc P$ be the subset of pairs $(\overline{\varphi}_1, \overline{\varphi}_2) \in \mc P$ which satisfy $|\overline{\varphi}_1-\overline{\varphi}_2| \gtrsim \d^{K'\eta}$. For an appropriate choice of $K'$ we obtain $|\mc P'| \ge |\mc P|/2$. 

Let $(\overline{\varphi}_1, \overline{\varphi}_2) \in \mc P'$. 
Define a graph $\mc H \subset \mc B\times \mc B$ by connecting $b$ and $b'$ by an edge if $|(\mc A\times \{(b, b')\}) \cap \mc G_{\overline{\varphi}_1, \overline{\varphi}_2}| \ge \d^{C_1\eta} |\mc A|$. Note that the definition of $\mc G$ is symmetric in $b$ and $b'$ and so we have $(b, b') \in \mc H$ if and only if $(b', b) \in \mc H$, so $\mc H$ is an unoriented graph. 
For an appropriate choice of $C_1$ it follows that $|\mc H| \ge \d^{C_1\eta} |\mc B|^2$. Now let $\mc B(\overline{\varphi}_1, \overline{\varphi}_2) \subset \mc B$ be a subset of vertices such that $\mc H|_{\mc B(\overline{\varphi}_1, \overline{\varphi}_2)}$ has minimal degree $\gtrsim \d^{O(\eta)} |\mc B|$ (one can find such subset by pruning vertices of small degree). The resulting subset $\mc B(\overline{\varphi}_1, \overline{\varphi}_2) \subset \mc B$ has the property that for every $b \in \mc B(\overline{\varphi}_1, \overline{\varphi}_2)$ we have $| (\mc A\times \{b\} \times \mc B(\overline{\varphi}_1, \overline{\varphi}_2)) \cap \mc G_{\overline{\varphi}_1, \overline{\varphi}_2}| \ge \d^{O(\eta)} |\mc A||\mc B|$. In particular, we must have $|\mc B(\overline{\varphi}_1, \overline{\varphi}_2)| \ge \d^{O(\eta)} |\mc B|$.

Now for each $(\overline{\varphi}_1, \overline{\varphi}_2) \in \mc P'$ let us define $r_{\overline{\varphi}_1, \overline{\varphi}_2}$ and $u_{\overline{\varphi}_1, \overline{\varphi}_2}$ to be vectors in $\RR^{d-2}$ which maximize the size of the set 
\[
\mc B_{\overline{\varphi}_1, \overline{\varphi}_2} := \{b \in \mc B(\overline{\varphi}_1, \overline{\varphi}_2):~ |F_{\overline{\varphi}_1- \overline{\varphi}_2}(b) - r_{\overline{\varphi}_1, \overline{\varphi}_2}b - u_{\overline{\varphi}_1, \overline{\varphi}_2}| \le \varrho^{3\varepsilon_0}\}.
\]

Let $\varepsilon' = \varepsilon/2d$ and suppose that there exists $(\overline{\varphi}_1, \overline{\varphi}_2) \in \mc P'$ such that $|\mc B_{\overline{\varphi}_1, \overline{\varphi}_2}| \le \varrho^{\varepsilon'} |\mc B|$.
In this case we can apply Theorem~\ref{thm:ABC-with-F} to the sets $\mc A$, $\mc B_{\overline{\varphi}_1, \overline{\varphi}_2}$ and the function $F_{\overline{\varphi}_1-\overline{\varphi}_2}$. We choose parameters in Theorem~\ref{thm:ABC-with-F} to be: $\varepsilon_1 = 2\varepsilon_0$, $\varepsilon_2 = \varepsilon'/2$ and $\eta = O(\eta/\kappa)$, $\overline{\eta} = O(\eta)$. By taking $\eta$ small enough in $\mu$, $\varepsilon$, $\varepsilon_0$, $\zeta$, $\kappa$ we get that the conditions of the lemma are satisfied. We conclude that for  $\zeta>0$ and appropriately small parameters we have $|\mc A|_\varrho \ge \varrho^{\zeta/2-1}$. On the other hand, $\mc A$ is obtained from rescaling the $(\d/\varrho, 1-\kappa, \d^{-O(\eta)})$-AD-regular set $Y_{z''}({\bf q})$. This implies the bound $\kappa \le \zeta/2 + O(\eta/\mu) \le \zeta$, as desired. 

It remains to deal with the case when $|\mc B_{\overline{\varphi}_1, \overline{\varphi}_2}| > \varrho^{\varepsilon'} |\mc B|$ for all $(\overline{\varphi}_1, \overline{\varphi}_2) \in \mc P'$. We are going to reach a contradiction with the non-concentration property (\ref{eq:big-F-non-linear}) and the slab Frostman property of $\overline{\Phi}_{\bf q}$. 
Indeed, by double counting and Cauchy--Schwarz we have the following estimate
\[
\sum_{(\overline{\varphi}^1_1, \overline{\varphi}^1_2), \ldots, (\overline{\varphi}^{d-2}_1, \overline{\varphi}^{d-2}_2)\in \mc P'} |\mc B_{\overline{\varphi}^1_1, \overline{\varphi}_2^1} \cap \mc B_{\overline{\varphi}^2_1, \overline{\varphi}_2^2} \cap \ldots \cap \mc B_{\overline{\varphi}^{d-2}_1, \overline{\varphi}_2^{d-2}}| \gtrsim \d^{(d-2) \varepsilon'} |\mc B| |\mc P'|^{d-2}
\]
Let $\mc S \subset \mc P'^{d-2}$ be the subset of tuples of pairs $(\overline{\varphi}^1_1, \overline{\varphi}^1_2), \ldots, (\overline{\varphi}^{d-2}_1, \overline{\varphi}^{d-2}_2)\in \mc P'$ so that this intersection has size at least $\d^{ (d-1)\varepsilon' }|\mc B|$. We conclude that then $|\mc S| \gtrsim \d^{(d-2)\varepsilon'}|\mc P'|^{d-2}$.

Now since $\overline{\Phi}_{\bf q}$ is $(\varrho, \gamma, \d^{O(\eta)})$-slab Frostman, the number of ways to choose  $(\overline{\varphi}_1^{j}, \overline{\varphi}^{j}_2) \in\mc P'$ for $j=1, \ldots, d-2$ so that 
\[
|(\overline{\varphi}_1^{1}- \overline{\varphi}_2^{1}) \wedge (\overline{\varphi}_1^{2}- \overline{\varphi}_2^{2})\wedge \ldots \wedge (\overline{\varphi}_1^{d-2}- \overline{\varphi}_2^{d-2})| \le \alpha
\]
is upper bounded by $\d^{O(d\eta )} \alpha^{\gamma} |\mc P'|^{d-2}$. So by taking $\alpha = \d^{ C d\varepsilon'/\gamma }$ we get that there exists a tuple of pairs $((\overline{\varphi}^1_1, \overline{\varphi}^1_2), \ldots, (\overline{\varphi}^{d-2}_1, \overline{\varphi}^{d-2}_2))\in \mc S$ with wedge product at least $\alpha$. Denote $\psi_j = \overline{\varphi}^{j}_1-\overline{\varphi}^{j}_2 \in [-2,2]^{d-2}$ and let $\tilde {\mc B} = \mc B_{\overline{\varphi}^1_1, \overline{\varphi}_2^1} \cap \mc B_{\overline{\varphi}^2_1, \overline{\varphi}_2^2} \cap \ldots \cap \mc B_{\overline{\varphi}^{d-2}_1, \overline{\varphi}_2^{d-2}}$. Furthermore, denote $r_j  = r_{\overline{\varphi}^{j}_1, \overline{\varphi}^{j}_2}$ and $u_j = u_{\overline{\varphi}^{j}_1, \overline{\varphi}^{j}_2}$. Then for every $b \in \tilde{\mc B}$ and $j=1, \ldots, d-2$ we have 
\[
|F(b) \psi_j - r_j b - u_j| \le \varrho^{\varepsilon_0'}.
\]
In other words, if we let $\Psi = (\psi_1, \ldots, \psi_j) \in \operatorname{Mat}(d-2, d-2)$, $R = (r_1, \ldots, r_{d-2}) \in \RR^{d-2}$ and $U = (u_1, \ldots, u_{d-2}) \in \RR^{d-2}$ then we get $\|F(b) \Psi - bR - U  \|_{\infty} \le \varrho^{\varepsilon_0'}$. 
By the choice of $\psi_j$, we have $|\det \Psi| \ge \alpha $  and $\|\Psi\|_{\infty}\leq 2$, and so 
\[
\| F(b) - b R \Psi^{-1} - U \Psi^{-1} \|_{\infty} \lesssim_d \alpha^{-1}\varrho^{\varepsilon_0'}.
\]
Recall $\varepsilon_0' = 3\varepsilon_0$ and $\alpha = \d^{ C d\varepsilon'/\gamma }$. So assuming that $Cd \varepsilon' < \gamma \mu \varepsilon_0$, we get $|F(b) - b r - u| \le \varrho^{2\varepsilon_0}$ for $r = R \Psi^{-1}$ and $u = U \Psi^{-1}$. Thus, we conclude that the number of $b \in \mc B$ such that $|F(b) - b r - u| \le \varrho^{\varepsilon_0}$ is lower bounded by $|\tilde{\mc B}| \ge \d^{(d-1)\varepsilon'}|\mc B| \gg \d^{\varepsilon/2}|\mc B|$ (recall that $\varepsilon' = \varepsilon/2d$). This however contradicts (\ref{eq:big-F-non-linear}). Thus, the case when all pairs $(\overline{\varphi}_1, \overline{\varphi}_2) \in \mc P'$ have the corresponding set $\mc B(\overline{\varphi}_1, \overline{\varphi}_2)$ concentrate in a tube is impossible. This completes the proof. 
\end{proof}

\subsection{Proof of Proposition \ref{prop:codimension-1-kakeya}} Suppose that we are given a $(d-1)$-planar $(d-1)$-linear $(\d, \eta, \kappa, \gamma)$-Kakeya configuration in $\RR^d$. Let $f: Z \to \operatorname{Mat}(1, d-2)$ denote its slope function. Let $\varrho = \d^{1/4}$. By passing to a subset we may assume that $|Z \cap {\bf z}| \gtrsim \d^{\eta} (\d/\varrho)$ for every ${\bf z} = z^\varrho$ and $z \in Z$. 
Let $\varepsilon_0>0$ be a small constant to be determined later. 
For each $\varrho$-interval ${\bf z}$ let $r_{{\bf z}}, u_{\bf z} \in \RR^{d-2}$ be vectors maximizing the size of the set
\[
Z({\bf z}) = \{z \in Z\cap {\bf z}:~ |f(z) - r_{\bf z} z - u_{\bf z}| \le \varrho^{1+\varepsilon_0}\}.
\]

We split into two cases. Let $\varepsilon>0$ be a small parameter to be determined. 
First, suppose that $\sum_{\bf z} |Z({\bf z})| \ge \d^{\varepsilon} |Z|$ and define $\tilde Z = \bigsqcup Z({\bf z})$. By Lemma \ref{lem:f-slow}, for most intervals ${\bf z}$ it has to be the case that $|r_{\bf z}| \gtrsim \d^{O(\varepsilon /\kappa)}$ as otherwise $\TT$ violates the convex Wolff axiom.
Thus, we can replace $\varrho^{1+\varepsilon_0}$ on the right hand side with $\varrho^{\varepsilon_0/2} |r_{\bf z}|$ and apply Lemma \ref{lem:codimension-1-linear-case}. Provided that $\eta, \varepsilon$ are sufficiently small in $\varepsilon_0, \mu=1/4, \gamma$ we then conclude that $C_{CW}(\TT) \ge \d^{O(\eta+\varepsilon)} \varrho^{-\kappa}$, which is a contradiction to the assumption $C_{CW}(\TT) \le \d^{-\eta}$. Thus we get $\kappa = O(\eta+\varepsilon) \ll \zeta$ for any $\zeta>0$ and $\eta, \varepsilon\ll \zeta$. Thus, we may assume that $\sum |Z({\bf z})| < \d^{\varepsilon} |Z|$ holds. By passing to a subset in $Z$ of density $1/2$ we may assume that $|Z({\bf z})| \lesssim \d^{\varepsilon} |Z|$ holds for every ${\bf z} = z^\varrho$ and $z \in Z$. Now assuming that $\varepsilon$ is sufficiently small in $\zeta$ and $\eta$ is sufficiently small in $\varepsilon_0, \varepsilon, \zeta$ we can apply Lemma \ref{lem:codimension-1-non-linear-case} to conclude that $\kappa\le \zeta$. This concludes the proof.

\section{2-planar configurations: Proof of Proposition \ref{prop:2-planar-kakeya}}\label{sec:2-planar}

In this section we prove Proposition \ref{prop:2-planar-kakeya} which asserts that there are no non-trivial 2-planar 2-linear Kakeya configurations. Lemma \ref{lem:codimension-1-linear-case} already rules out the case when the slope function is approximately linear, so it suffices to deal with the non-linear case. Namely, in this section we prove the following. 

\begin{lemma} \label{lem:non-linear-grains-2-planar}
    Let $\zeta, \kappa,\mu,\varepsilon, \varepsilon_0 >0$ be so that $\varepsilon_0$ is sufficiently small in $\zeta$. Then the following holds for $\eta < \eta(\zeta,\mu, \varepsilon_0, \varepsilon)$.

    Consider a $2$-linear $(\d, \eta, \kappa)$-Kakeya configuration in $\RR^d$ and let $f: Z \to \operatorname{Mat}(d-2,1)$ be the slope function. Suppose that for all $z_0\in Z$, ${\bf z} = z_0^\varrho$, $\varrho \in [\d^{1/4}, \d^\mu]$  and for any $r, u \in \RR^{d-2}$ we have
    \begin{equation}\label{eq: assumption}
    \#\{ z \in Z \cap {\bf z}:  ~ |f(z) - r z -u| < \varrho^{1+\varepsilon_0} \} \le \varrho^{\varepsilon} (\varrho/\d).
    \end{equation}
    Then we have $\kappa \le \zeta$. 
\end{lemma}

Note that in this lemma we do not need the Kakeya configuration to be $2$-planar. This case only appears in the linear case from the previous section. 

\begin{proof}[Proof of Proposition \ref{prop:2-planar-kakeya}]
    Suppose that we are given a $2$-planar $2$-linear $(\d, \eta, \kappa, \gamma)$-Kakeya configuration in $\RR^d$. Let $f: Z \to \operatorname{Mat}(1, d-2)$ denote its slope function. Let $\varrho = \d^{1/4}$. By passing to a subset we may assume that $|Z \cap {\bf z}| \gtrsim \d^{\eta} (\d/\varrho)$ for every ${\bf z} = z^\varrho$ and $z \in Z$. 
Let $\varepsilon_0>0$ be a small constant to be determined later. 
For each $\varrho$-interval ${\bf z}$ let $r_{{\bf z}}, u_{\bf z} \in \RR^{d-2}$ be vectors maximizing the size of the set
\[
Z({\bf z}) = \{z \in Z\cap {\bf z}:~ |f(z) - r_{\bf z} z - u_{\bf z}| \le \varrho^{1+\varepsilon_0}\}.
\]
We split into two cases. Let $\varepsilon>0$ be a small parameter to be determined. 
First, suppose that $\sum_{\bf z} |Z({\bf z})| \ge \d^{\varepsilon} |Z|$. Then we can define $\tilde Z = \bigsqcup Z({\bf z})$ and apply Lemma \ref{lem:codimension-1-linear-case}. Provided that $\eta, \varepsilon$ are sufficiently small in $\varepsilon_0, \mu=1/4, \gamma$ we then conclude that $C_{CW}(\TT) \ge \d^{O(\eta+\varepsilon)} \varrho^{-\kappa}$. On the other hand, $C_{CW}(\TT) \le \d^{-\eta}$ by the assumption. Thus, we get $\kappa = O(\eta+\varepsilon) \ll \zeta$ for any $\zeta>0$ and $\eta, \varepsilon\ll \zeta$. Thus, we may assume that $\sum |Z({\bf z})| < \d^{\varepsilon} |Z|$ holds. By passing to a subset in $Z$ of density $1/2$ we may assume that $|Z({\bf z})| \lesssim \d^{\varepsilon} |Z|$ holds for every ${\bf z} = z^\varrho$ and $z \in Z$. Now assuming that $\varepsilon$ is sufficiently small in $\zeta$ and $\eta$ is sufficiently small in $\varepsilon_0, \varepsilon, \zeta$ we can apply Lemma \ref{lem:non-linear-grains-2-planar} to conclude that $\kappa\le \zeta$. This concludes the proof. 
\end{proof}

\subsection{Proof of Lemma \ref{lem:non-linear-grains-2-planar}}

In the proof we will reduce to the following $ABC$ sum-product type theorem. 

\begin{theorem}\label{thm:2-dim-two-ends-ABC}
    Let $\kappa, \zeta >0$.
    Let $\varepsilon_1, \varepsilon_2>0$ be so that $\varepsilon_1$ is small enough in $\zeta$. 
    For $\eta < \eta_0(\kappa, \zeta,\varepsilon_1)$ and $\overline{\eta} < \overline{\eta}(\kappa, \zeta, \varepsilon_1, \varepsilon_2)$ and $\d$ small enough in all parameters we have the following.

    Let $B \subset [-1,1]^2$ be a $\d$-separated $(\d, 1, \d^{-\eta})$-set such that $|B \cap T| \le \d^{\varepsilon_2} |B|$ for any $\d^{\varepsilon_1}\times 1$ tube $T$. Let $C \subset [-1,1]$ be a $\d$-separated $(\d, \kappa, \d^{-\eta})$-set. Let $A \subset [-1,1]^{2}$ be a $\d$-separated set. Let $G \subset A\times B\times C$ be a subset of density at least $\d^{\overline{\eta}}$ and suppose that
    \[
    |\{ a + cb\mid (a, b, c) \in G \}|_\d \le \d^{-\overline{\eta}} |A|.
    \]
    Then we have $|A| \ge \d^{\zeta-1} |C|$.
\end{theorem}

\begin{proof}
    Apply radial projection, Theorem \ref{thm:radial-projections}, to $B$, there exists $\eta_0=\eta_0(1, \varepsilon_1, \zeta)>0$ such that if $\eta < \eta_0$ and $\d$ is small enough then there is a subgraph
    $H \subset B\times B$ such that 
    \begin{equation}\label{eq: H}
    |H| \ge (1-\d^{\min(\eta_0, \varepsilon_2/2)}) |B|^2
    \end{equation} and for every $(b, b') \in H$ we have $|T_{b, b'}(\varrho)| \le \d^{-\zeta} \varrho |B|$ for all $\varrho \in [\d, 1]$. Thus, the set of directions $\Lambda=\{ \frac{b-b'}{|b-b'|}, ~ (b, b') \in H \} $ contains a $(\d, 1, O(\d^{-\zeta}))$-set. 

    Let $\mc S_j \subset \d \ZZ$ be the set of numbers $s = \d j$ such that $\#\{ (b, c) \in B\times C:~ bc \in [s, s+\d]\} \sim 2^j$. By dyadic pigeonhole and Cauchy--Schwarz, we can find an index $j$ such that $2^j |\mc S_j| \gtrsim \d^{\overline{\eta}} |B||C|$ and 
    \[
    \#\{ (a, \tilde a, s, \tilde s) \in A^2 \times \mc S_j^2: |a + s - \tilde a-\tilde s| \le \d \} \gtrsim \d^{O(\overline{\eta})} |A| |\mc S_j|^2.
    \]
    Note that since $B$ is a $(\d, 1, \d^{-\eta})$-set, $|\mc S_j| \ge \d^{O(\eta+\overline{\eta})} \d^{-1}$. 
    
    Thus, by Lemma \ref{lem:asymmetric-bsg} (Balog-Szemeredi-Gowers), for any $\varepsilon' >0$ there is $K=K(\varepsilon')$ and a subset $A' \subset A$ and $F \subset B\times C$ such that $|A'| \gtrsim \d^{\varepsilon' + K\overline{\eta}}|A|$ and $|F| \gtrsim \d^{\varepsilon' + K\overline{\eta}} |B| |C|$ such that
    \begin{equation}\label{eq:A-prime}
    |A' + \{ b c, ~ (b, c) \in F  \} -  \{ b c, ~ (b, c) \in F  \}|_\d \le \d^{-O_{n,m}(\varepsilon'+K\overline{\eta})} |A|_\d.  
    \end{equation}

    By Cauchy--Schwarz, we can construct a subgraph $H' \subset B\times B$ such that $|H'| \ge \d^{O(\varepsilon' + K\overline{\eta})} |B|^2$ such that for every $(b, b') \in H'$ we have $|F|_{b} \cap F|_{b'}| \ge \d^{O(\varepsilon' + K\overline{\eta})} |C|$. Denote $C_{b,b'} = F|_{b} \cap F|_{b'} \subset C$. 

    So by (\ref{eq:A-prime}) for every $(b,b') \in H'$ we have
    \[
    |A' + (b-b') C_{b,b'} |_{\d} \le  \d^{-O(\varepsilon'+K\overline{\eta})} |A|_\d
    \]
    
 For each $(b, b') \in H'$ and $\lambda = \frac{b-b'}{|b-b'|}$ let $\pi_{\lambda}$ denote the orthogonal projection from $\mathbb{R}^2$ to $\lambda^{\perp}$. Then for some constant $C_0$ we have
\[
\#\{(a, \tilde a, c, \tilde c): |a+(b-b')c -(\tilde a +(b-b')\tilde c)|\leq \d , \, |a-\tilde a|\geq \delta^{C_0(\eta/\kappa +\varepsilon' +K \overline{\eta})}, \, |c-\tilde c|\geq \delta^{O(\eta/\kappa)} \} \gtrsim \delta^{O(\varepsilon'+K \overline{\eta})} |A| |C|^2.
\]
Since $C$ is a $(\d, \kappa, \d^{-\eta})$-set,   most pairs $(c, \tilde c)\in C\times C$ are $\d^{O(\eta/\kappa)}$-separated.
Therefore
    \[
    \#\{ (a, \tilde a) \in A'^2: ~ |a-\tilde a| \ge \d^{C_0(\eta/\kappa + \varepsilon'+K\overline{\eta})},~ |\pi_\lambda(a)-\pi_\lambda(\tilde a)|\le \d \} \ge \d^{O(\varepsilon'+K\overline{\eta})} |A'| |C|.
    \]

    On the other hand, by taking $\varepsilon'$ small enough in $\min(\eta_0, \varepsilon_2/2)$ and $\overline{\eta}$ small enough in $\varepsilon', K$, we get $|H \cap H'| \ge |H'| /2$. So the set $ \{ \frac{b-b'}{|b-b'|}, ~ (b, b') \in H\cap H' \}$ contains a $(\d, 1, \d^{-\zeta - O(\varepsilon'+K\overline{\eta})})$-set $\Lambda'$ and we get that by double counting and the Frostman property of $\Lambda'$:
    \begin{align*}
    \sum_{
    \lambda \in \Lambda'}    \#\{ (a, \tilde a) \in A'^2: ~ |a-\tilde a| \ge \d^{C_0(\eta/\kappa + \varepsilon'+K\overline{\eta})},~ |\pi_\lambda(a)-\pi_\lambda(\tilde a)|\le \d \} \\\lesssim |A'|^2 \d^{-\zeta- O(\varepsilon'+K\overline{\eta})}  (\d/\d^{C_0(\eta/\kappa + \varepsilon'+K\overline{\eta})})  |\Lambda'|
    \end{align*}
    Comparing the lower and the upper bounds leads to
    \[
    |A'| \gtrsim \d^{\zeta +O(\eta/\kappa+\varepsilon'+K\overline{\eta})}  \d^{-1} |C|.
    \]
    By replacing $\zeta$ with $\zeta/2$ and taking $\eta, \varepsilon', \overline{\eta}$ small enough in $\zeta$, we get the desired bound. 
\end{proof}

\begin{remark}
    Theorem \ref{thm:2-dim-two-ends-ABC} can be extended the case when $A, B \subset [-1,1]^d$ for some $d\ge3$. Indeed, by taking a generic $\RR^d\to \RR^2$ projection $\pi$, one can show that $\pi(B)$ still satisfies the two-ends tube non-concentration condition (with slightly worse constants $\varepsilon_1, \varepsilon_2$ and after mild pruning). 
\end{remark}

Now we ready to prove the main result of this section.

\begin{proof}[Proof of Lemma \ref{lem:non-linear-grains-2-planar}]
    Recall from Definition~\ref{def:kakeya-configuration}, 
    $|f(z)-f(z')|\le \d^{-\eta}\varrho$ for $z, z'$ in the same dyadic $\varrho$-interval.  Fix $z_0 \in Z$, ${\bf z} = z_0^\varrho$ and some $p_0\in E$. 
    

    Let $W_{z,z',z''}$ be the set of $W$-tuples as defined in Section \ref{subsec:W-tuples-kakeya}. 
    Let $\overline{\Phi}_{\bf q} \subset \varrho \ZZ \cap [-1,1]$ be the set of points such that $\overline{\Phi}_{\bf q}^{\varrho} = \Phi^{\varrho}_{\bf q}$.
    For $z \in Z({\bf q})$ and $y \in Y_{z}({\bf q})$ we define a subset $G_{y,z} \subset \varrho^2 \ZZ \times Z({\bf q}) \times Z({\bf q})\times \overline{\Phi}_{\bf q} \times \overline{\Phi}_{\bf q}$ consists of tuples $(x, z_1, z_2, \overline{\varphi}_1, \overline{\varphi}_2)$ such that there exists a tuple $(p, p_1, p_2, p_1', p_2', T_1, T_2, T_1', T_2') \in W_{z_1, z_2, z}$ such that for $i=1,2$, the first coordinate of $\theta(T_i)$ lies in the dyadic $\varrho$-interval $\overline{\varphi}_i^\varrho$, and for $p'_i = (x_i', y_i'+f(z) x_i', z)$ we have $x_1' \in x^{(\varrho^2)}$ and $y_i' \in y^{(\varrho^2)}$. 

    For fixed $y, z$ and a constant $C\ge 1$ consider the following set: 
    \[
    A = A_{y,z}^C= \left\{ \begin{pmatrix}
        x' \\
        \xi'
    \end{pmatrix}:~\exists  p' = (x', y' + f(z)x', z) \in E({\bf q}), ~|y'-y| \le C \varrho^2, ~ \xi' = \xi_{p'} \right\}.
    \]
    Recall for each $p'\in E$,  $\xi_{p'}$ is given in Definition~\ref{def:kakeya-configuration} $(4)$.

        Suppose that $(x, z_1, z_2, \overline{\varphi}_1, \overline{\varphi}_2) \in G_{y, z}$. If $C$ is a sufficiently large constant,  then by \eqref{eq:vector-identitiy} and \eqref{eq:y-differ-little}, for some $\xi \in \varrho^2 \ZZ^{d-2} \cap[-1,1]^{d-2}$ we  have 
        \begin{equation}\label{eq: A}
         \begin{pmatrix}
            x \\
            \xi
        \end{pmatrix} \in A + 
        \begin{pmatrix}
            O(\varrho^2)\\
            O(\varrho^2)
        \end{pmatrix},
        \quad  \begin{pmatrix}
            x \\
            \xi
        \end{pmatrix} +  \begin{pmatrix}
            z_2-z_1\\
            f(z_2)-f(z_1)
        \end{pmatrix}(\varphi_2-\varphi_1) \in A+
        \begin{pmatrix}
            O(\varrho^2)\\
            O(\d^{-\eta}\varrho^2)
        \end{pmatrix}
        \end{equation}

    Recall that the set of $\xi_{p'}$ for $p'$ in a dyadic $w=C_0\varrho^2$ box all lie in the same dyadic $w$-box. Once we fix $x$, $p=(x, y+f(z)x, z)$ is fixed (since $y, z$ were fixed previously), and so is $\xi=\xi_p$. Therefore 
    \[
    |A|_{\varrho^2 \times \varrho^2} \lesssim \varrho^{-1}.
    \]
    Define $\mc Z = \{ \frac{z-z_0}{\varrho}, ~z\in Z({\bf q}) \} \subset [-1,1]$ and $F: \mc Z\to [-1,1]^{d-2}$ denote the function
    \[
    F(z) = \frac{f(z)-f(z_0)}{\d^{-\eta}\varrho}. 
    \]
    Fix some $x_0 \in X_{y,z}({\bf q})$, let $\xi_0 = \xi_{(x_0, y+f(z)x_0, z)}$ and define a set $\mc A \subset \varrho \ZZ \times (\varrho \ZZ)^{d-2}$ to be the set of vectors $\begin{pmatrix}
        \tilde x \\
        \tilde\xi
    \end{pmatrix}$ such that there exists $\begin{pmatrix}
        x'\\
        \xi' 
    \end{pmatrix} \in A$ such that $x' \in (x_0 + \varrho \tilde x)^{(\varrho^2)}$ and $\xi' \in (\xi_0 + \varrho \tilde\xi)^{(\varrho^2)}$. It follows that $|\mc A| \lesssim \varrho^{-1}$.
    
    Consider the graph 
    \[
    \mc B = \left\{\begin{pmatrix}
        \tilde z\\
        F(\tilde z)
    \end{pmatrix}, \tilde z\in \mc Z\right\}. 
    \]
    Since $|Z({\bf q})| \ge \d^{O(\eta)} (\varrho/\d)$,  by \eqref{eq: assumption}, we have the following two-ends condition on $F$:
    \[
    \#\{ \tilde z \in \mc Z: |F(\tilde z) - r \tilde z - u| \le \d^\eta\varrho^{\varepsilon_0} \} \lesssim_{\d^{O(\eta)}} \d^{\varepsilon} |\mc Z|
    \]
    for all $r, u \in \RR^{d-2}$. From this it follows that for any $\d^{2\varepsilon_0} \times 1$ tube $T \subset \RR^{d-1}$ we have
    \begin{equation}\label{eq:F-two-ends}
        |\mc B \cap T| \le \d^{\min(\varepsilon_0, \varepsilon)/2} |\mc B|
    \end{equation}
    provided that $\eta$ is small enough in $\varepsilon, \varepsilon_0$.

        Suppose that $(x, z_1, z_2, \overline{\varphi}_1, \overline{\varphi}_2) \in G_{y, z}$. Let $b_1, b_2 \in \mc B$ be elements in $\mc B$ corresponding to $z_1, z_2$.
        Restating \eqref{eq: A} in the new notation,  there are $a = \begin{pmatrix}
            \tilde x\\
            \tilde \xi
        \end{pmatrix}, a'=\begin{pmatrix}
            \tilde x'\\
            \tilde \xi'
        \end{pmatrix} \in \mc A$ 
        so that $x \in (x_0 + \varrho \tilde x)^{(\varrho^2)}$ and we have
        \begin{equation}\label{eq:mcA}
            a' = a + (b_2 - b_1) (\overline{\varphi}_2-\overline{\varphi}_1) + O(\d^{-\eta}\varrho).
        \end{equation}


    Using Lemma \ref{lem:W-tuple-many-Kakeya}, we conclude that there is a graph $\mc G \subset \mc A \times \mc B \times \mc B \times \overline{\Phi}_{\bf q} \times \overline{\Phi}_{\bf q}$ such that $|\mc G| \gtrsim_{\d^{O(\eta)}} |\mc A| |\mc B|^2 |\overline{\Phi}_{\bf q}|^2$ such that for each $(a, b_1, b_2, \overline{\varphi}_1, \overline{\varphi}_2) \in \mc G$ we have (\ref{eq:mcA}) for some $a' \in \mc A$. In particular, this implies the covering number upper bound
    \begin{equation}\label{eq: AB}
    |\{ a + (b_2 - b_1) (\overline{\varphi}_2-\overline{\varphi}_1), ~ (a, b_1, b_2, \overline{\varphi}_1,\overline{\varphi}_2) \in \mc G \}|_\varrho \lesssim \d^{-O(\eta)} |\mc A|_\varrho.
    \end{equation}

    Using that $|\mc Z|\gtrsim_{\d^{O(\eta)}}\varrho^{-1}$ we get that $\mc B$ is a $(\varrho, 1, \d^{-O(\eta)})$-set which satisfies the two-ends property (\ref{eq:F-two-ends}). We also have $|\mc A| \lesssim \varrho^{-1}$ and $\overline{\Phi}_{\bf q}$ is a $(\varrho, \kappa, \d^{-O(\eta)})$-AD-regular set.
    Fix values $b_1, \varphi_1$ for which the restriction $\mc G|_{b_1, \varphi_1}$ has density $\d^{O(\eta)}$.

    Let $V$ be a $2$-plane in $\mathbb{R}^{d-1}$ and consider the orthogonal projection $P_V: \mathbb{R}^{d-1}\rightarrow V$. By Mattila's projection theorem \cite[Chapter 9]{mattila1999geometry}, for almost all $2$-planes in $V$, $P_V(\mathcal{B})$ is a $(\rho, 1, \delta^{-O(\eta)})$-set (after a mild refinement) and $|P_V(\mathcal{A})|\sim |\mathcal{A}|\lesssim \varrho^{-1}$. In addition, \eqref{eq: AB} still holds with $P_V(a)$ in place of $a$, $P_V(b_2)-P_V(b_1)$ in place of $b_2-b_1$.  Apply  Theorem \ref{thm:2-dim-two-ends-ABC} applied to $P_V(\mathcal{B}), P_V(\mathcal{A})$ and $\overline{\Phi}_{\bf q}$, 
we obtain the lower bound $|\mc A|_\varrho >\varrho^{\zeta/2} \varrho^{-1-\kappa}$ provided that $\eta$ is small enough in $\zeta, \mu,\varepsilon, \varepsilon_0$.
    Comparing to the upper bound on $|\mc A|_\varrho$, we conclude that $\kappa \le \zeta$, as desired.
\end{proof}

\section{3-planar 2-linear in \texorpdfstring{$\RR^4$}{R4} case: proof of Proposition \ref{prop:32-planar-kakeya}} \label{sec:32planar}

In this section we consider the $3$-planar $2$-linear Kakeya configurations, this is the last case of sticky Kakeya in $\RR^4$. The overall strategy is the same as in the previous sections, except for the fact that here we need to split into `quadratic vs non-quadratic' cases instead of `linear vs non-linear' as before.

\subsection{\texorpdfstring{$V$}{V}-tuples and \texorpdfstring{$L$}{L}-tuples.}
\label{subsec:L-tuple-kakeya-2}
This subsection is the $k$-planar, $\ell$-linear analogue of Section
\ref{subsec:L-tuple-kakeya}. We reproduce the definitions and calculations from that section here for completeness. 

Fix a $k$-planar $\ell$-linear Kakeya configuration
$(\TT,E,\{\TT_p\})$ with slope functions
\[
\widetilde f(z)\in\operatorname{Mat}(d-k,k-1),
\qquad
f(z)\in\operatorname{Mat}(d-\ell,\ell-1).
\]
By Section \ref{subsec:lipschitz}, we may assume that both slope
functions are $\d^{-\eta}$-Lipschitz.  Fix $p_0\in E_{z_0}\subset E$ and the
$\varrho$-cube ${\bf q}=p_0^\varrho$ and define sets
$E({\bf q}) = E\cap \bf q$, $Z({\bf q}) = Z\cap z_0^\varrho$, and $\TT({\bf q}) =\cup_{p\in E(\bf q)} \TT_p$ in the usual way. As in Definition \ref{def:kakeya-configuration-2}, write $f = \begin{pmatrix}
    f_1 \\
    f_2
\end{pmatrix}$ and $\widetilde f = (\widetilde f_1, \widetilde f_2)$ where
$f_1 \in \Mat(k-\ell, \ell-1)$, $f_2, \widetilde f_1 \in \Mat(d-k, \ell-1)$, $\widetilde f_2 \in \Mat(d-k, k-\ell)$. 
Recall that we have the relation 
\begin{equation}\label{eq:ff-relation}
f_2(z) = \tilde f_1(z) + \widetilde f_2(z) f_1(z)
\end{equation}
that identifies the $k$-grains and $\ell$-grains decompositions of $E$. 

\medskip\noindent\emph{Basic $V$-tuples.}
For $z,z'\in Z({\bf q})$, let $V_{z,z'}({\bf q})$ consist of the
tuples
\[
(p,p_1,p_2,T_1,T_2)
\in E_z({\bf q})\times E_{z'}({\bf q})^2\times\TT({\bf q})^2
\]
such that
\[
T_i\in\TT_p\cap\TT_{p_i},
\qquad i=1,2.
\]
Their total number satisfies
\begin{equation}\label{eq:V-tuples-counting-3}
\sum_{z,z'\in Z({\bf q})}|V_{z,z'}({\bf q})|
\sim_{\d^{O(\eta)}}
|E({\bf q})|\d^{-2\kappa}|Z({\bf q})|.
\end{equation}
Indeed, the upper bound follows by selecting $p$ and then $T_1,T_2$ and $p_1,p_2$ in at most $|E({\bf q})|$, $\sim \d^{-2\kappa}$ and $\sim |Z({\bf q})|$ ways, respectively. For the lower bound, count the number of one-arm paths $(p, p_1, T_1)$ and $(p,p_2,T_2)$ and use Cauchy--Schwarz. 

For a $V$-tuple, write
\[
p=(x,y+\widetilde f(z)x,z),\qquad
p_i=(x_i,y_i+\widetilde f(z')x_i,z'),
\]
where $x, x_i$ are $k-1$-dimensional vectors and $y,y_i$ are $d-k$-dimensional vectors from the $k$-grain decomposition of $E({\bf q})$. 
Let $\varphi_i$ be the vector of first $\ell-1$ coordinates of
$T_i$, so that $\theta(T_i) = (\varphi_i, \xi_p +f(z) \varphi_i) + O(\d)$. Note that
\begin{align*}
\theta(T_1)-\theta(T_2) =& (\varphi_1-\varphi_2, f(z)(\varphi_1-\varphi_2))+O(\d) \\
=& (\varphi_1-\varphi_2, f_1(z)(\varphi_1-\varphi_2), f_2(z)(\varphi_1-\varphi_2) ) + O(\d)
\end{align*}
Subtracting
$p_i=p+(z'-z)(\theta(T_i),1)+O(\d)$ for $i=1,2$ gives
\[
\begin{aligned}
x_1-x_2&=(z'-z)(\varphi_1-\varphi_2, f_1(z)(\varphi_1-\varphi_2))+O(\d),\\
y_1-y_2+\widetilde f(z')(x_1-x_2)
&=(z'-z)f_2(z)(\varphi_1-\varphi_2)+O(\d).
\end{aligned}
\]
Consequently,
\begin{align}
y_1-y_2 =& (z'-z) (f_2(z) -\widetilde f(z') (I_{\ell-1}, f_1(z))) (\varphi_1-\varphi_2) + O(\d) \nonumber \\
=& (z'-z)( f_2(z) - \widetilde f_1(z') - \widetilde f_2(z') f_1(z) )(\varphi_1-\varphi_2)+O(\d) \nonumber \\ 
=& (z'-z) ( f_2(z) - f_2(z') - \widetilde f_2(z') (f_1(z) - f_1(z')) )(\varphi_1-\varphi_2)+O(\d),
\label{eq:V-tuple-relation-2}
\end{align}
where to obtain the last line we use (\ref{eq:ff-relation}) to express $\widetilde f_1(z')$ in terms of other functions. This identity will be used in the quadratic case in Section \ref{subsec:quadratic-case} below. We now proceed to define $L$-tuples, which will be used in the non-quadratic case.

\medskip\noindent\emph{One-arm paths.}
For $z\in Z({\bf q})$ and $p,\widetilde p\in E_z({\bf q})$, write
$p\leftrightarrow\widetilde p$ if
\[
p=(x,y+\widetilde f(z)x,z),\qquad
\widetilde p=(\widetilde x,y+\widetilde f(z)\widetilde x,z)
\]
for some $y\in Y_z({\bf q})$ and
$x,\widetilde x\in X_{y,z}({\bf q})$.  That is, the two points lie in
the same horizontal $k$-grain.

For $z,z',z''\in Z({\bf q})$, let $J_{z,z',z''}$ consist of
\[
(p,p',\widetilde p',p'',T,T')
\in E_z({\bf q})\times E_{z'}({\bf q})^2
\times E_{z''}({\bf q})\times\TT({\bf q})^2
\]
such that
\begin{equation}\label{eq:J-incidences-kakeya-2}
T\in\TT_p\cap\TT_{p'},
\qquad
p'\leftrightarrow\widetilde p',
\qquad
T'\in\TT_{\widetilde p'}\cap\TT_{p''}.
\end{equation}
As in Section \ref{subsec:L-tuple-kakeya} we have
\begin{equation}\label{eq:J-tuples-counting-kakeya-2}
\sum_{(z,z',z'')\in S}|J_{z,z',z''}|
\sim_{\d^{O(\eta)}}
|E({\bf q})|\d^{-2\kappa}|Z({\bf q})|^2
(\varrho/\d)^{k-1}.
\end{equation}

\medskip\noindent\emph{Definition of the $L$-tuples.}
Fix $\tau\in[\d,\varrho]$.  For $p\in E({\bf q})$ and
$\widetilde x\in(\tau\ZZ)^{k-1}$, let
$J_{z,z',z''}(p,\widetilde x)$ be the subset  of tuples $(p,p',\widetilde p',p'',T,T')$ in
$J_{z,z',z''}$ whose initial point is $p$ and for which the
$x$-coordinate of $\widetilde p'$ lies in the $\tau$-neighborhood of
$\widetilde x$.  Define
\[
L_{z,z',z''}(p,\widetilde x)
\subset
J_{z,z',z''}(p,\widetilde x)^2
\]
to consist of ordered pairs whose terminal tubes satisfy
$|\theta(T_1')-\theta(T_2')|\leq\tau$, and put
\[
L_{z,z',z''}
=\bigcup_{p,\widetilde x}
L_{z,z',z''}(p,\widetilde x).
\]
See Figure \ref{fig:L-tuple} for an illustration.

\medskip\noindent\emph{Counting $L$-tuples.}
The following lemma gives an estimate on the total number of $L$-tuples for an
appropriate choice of $\tau$. Since $f$ and $\tilde f$ are dyadic Lipchitz, \eqref{eq:tau-defining-relation-2} holds for $\tau \leq \varrho^2$, later we will show it holds for even smaller $\tau$. 

\begin{lemma}\label{lem:number-of-Lambda-tuple-kakeya-2}
Let $\widetilde{Z} \subset Z({\bf q})$ be a subset and suppose that for every $z, z' \in \widetilde Z$ we have 
\begin{equation}\label{eq:tau-defining-relation-2}
|f_2(z) - f_2(z') - \tilde f_2(z') (f_1(z)- f_1(z'))| \le c \tau/\varrho.    
\end{equation}
Then
\[
\sum_{(z,z',z'')\in \widetilde Z^3}|L_{z,z',z''}|
\sim_{\d^{O(\eta)}}
|E_{\widetilde Z}({\bf q})|\d^{-4\kappa}\tau^\kappa
(\varrho/\d)^{k-1}(\tau/\d)^{k-1}
|\widetilde Z|^2.
\]
\end{lemma}

\begin{proof}
The proof is essentially identical to the proof of Lemma \ref{lem:number-of-Lambda-tuple-kakeya}. The only difference compared to that lemma is the restriction that we impose on $\tau$. This comes from the difference between the formulas for $y_1-y_2$ in a $V$-tuple in (\ref{eq:V-tuple-relation-2}) and (\ref{eq:V-tuple-relation}). Let us briefly outline the counting argument leading to the formula for the number of $L$-tuples. Start with a $V$-tuple $(p, p_1', p_2', T_1, T_2)$ with $p \in E_z({\bf q})$ and $p_1',p_2' \in E_{z'}({\bf q})$. By \eqref{eq:V-tuples-counting-3}, the number of such $V$-tuples is $\sim_{\delta^{O(\eta)}} |E_{\tilde{Z}}(\mathbf{q})| \delta^{-2\kappa}|\tilde{Z}|$.

Using $|z-z'| \le \varrho$, $|\varphi_1-\varphi_2| \le O(1)$, the formula (\ref{eq:V-tuple-relation-2}) and the assumption on $\tau$ in the statement of the lemma, we conclude that the $y$-coordinates of $p_1'$ and $p_2'$ differ by $O(\tau)$. We can then choose a pair $\widetilde p_1'$ and $\widetilde p_2'$ on the same respective $k$-grains as $p_1'$ and $p_2'$ so that their $x$-coordinates differ by at most $\tau$ in $\sim_{\delta^{O(\eta)}} (\varrho/\d)^{k-1} (\tau/\d)^{k-1}$ ways. It follows that $|\tilde p'_1 - \tilde p'_2| = O(\tau)$. Since $\tilde p'_1, \tilde p'_2$ are contained in the same $\sim \tau$ cube, there is a set of $\tau$-cubes of size $\sim_{\d^{O(\eta)}} \tau^{-\kappa}$ covering the direction sets of $\TT_{\tilde p'_1}, \TT_{\tilde p'_2}$. Using a Cauchy--Schwarz argument, we conclude that there are $\sim_{\d^{O(\eta)}} \tau^\kappa \d^{-2\kappa}$ ways to select a pair of tubes $T_1' \in \TT_{\widetilde p'_1}$ and $T_2' \in \TT_{\widetilde p'_2}$ such that $|\theta(T'_1) -\theta(T'_2)| \le \tau$. Once these tubes are fixed, we can pick the remaining points $p_1'' \in Z(T_1', {\bf q}), p_2''\in Z(T_2', {\bf q})$ in $\sim |\widetilde Z|$ many ways. 
\end{proof}

\medskip\noindent\emph{Additive relations.} 
Suppose that
\[
((p,p_1',\widetilde p_1',p_1'',T_1,T_1'),
 (p,p_2',\widetilde p_2',p_2'',T_2,T_2'))
\in L_{z,z',z''}.
\]
For $i=1,2$, write
\[
\begin{aligned}
p_i'&=(x_i',y_i'+\widetilde f(z')x_i',z'),\\
\widetilde p_i'
&=(\widetilde x_i',\widetilde y_i'
  +\widetilde f(z')\widetilde x_i',z'),\\
p_i''&=(x_i'',y_i''+\widetilde f(z'')x_i'',z''),\\
\theta(T_i)&=(\varphi_i,\xi_p+f(z)\varphi_i)+O(\d),
\qquad \varphi_i\in\Phi_p.
\end{aligned}
\]
The subtuple $(p, p_1', p_2', T_1, T_2)$ forms a $V$-tuple, so \eqref{eq:V-tuple-relation-2} gives 
\[
y_1'-y_2' = (z'-z) ( f_2(z) - f_2(z') - \widetilde f_2(z') (f_1(z) - f_1(z')) )(\varphi_1-\varphi_2)+O(\d)
\]
Since $\widetilde p'_i$ and $\widetilde p'_i$ lie in the same horizontal $k$-grain, we have $\widetilde y'_i = y'_i$. We have
\[
p_i''=\widetilde p_i'
+(z''-z')(\theta(T_i'),1)+O(\d).
\]
As in \eqref{eq:y-axes-another-relation-kakeya}, we use $|z''-z'|\le \varrho$ and $|\theta(T_1')-\theta(T_2')| \le \tau$ to conclude that 
\[
p_1'' - p_2'' = \widetilde p_1' - \widetilde p_2' + O(\tau \varrho +\d).
\]
Consequently,
\[
\begin{aligned}
x_1''-x_2''
&=\widetilde x_1'-\widetilde x_2'+O(\tau\varrho+\d),\\
y_1''-y_2''
&=\widetilde y_1'-\widetilde y_2'
+(\widetilde f(z')-\widetilde f(z''))(\widetilde x_1'-\widetilde x_2')
+O(\tau\varrho+\d).
\end{aligned}
\]
Because the two $\widetilde x_i'$ lie in the same $\tau$-ball,
$|\widetilde x_1'-\widetilde x_2'|=O(\tau)$. It follows from Lemma \ref{lem:f-is-Lipschitz} that $\widetilde f$ is $\d^{-\eta}$-Lipschitz. Thus, by combining these facts with (\ref{eq:V-tuple-relation-2}) we get
\begin{equation}\label{eq:y-axes-another-relation-kakeya-2}
 y_1''-y_2''
= (z'-z) ( f_2(z) - f_2(z') - \widetilde f_2(z') (f_1(z) - f_1(z')) )(\varphi_1-\varphi_2)+O(\d^{-\eta}\tau\varrho+\d)
\end{equation}



Just as \eqref{eq: graphG} in Section \ref{subsec:non-linear-case-Kakeya}, given a choice of $\widetilde Z$ and the corresponding sets of $L$-tuples $L_{z,z', z''}$,  we can find some $z''\in \widetilde Z$ and construct a  subgraph 
\begin{equation}\label{eq: G}
G \subset Y_{z''}({\bf q}) \times \widetilde Z^2 \times \overline{\Phi}_{\bf q}^2,
\end{equation}
where $\overline{\Phi}_{\bf q} \subset \varrho \ZZ^{\ell-1}$ is the set whose dyadic $\varrho$-neighborhood covers $\Phi_{\bf q}$.
Any tuple $(y, z, z', \overline{\varphi}_1, \overline{\varphi}_2) \in G$ has the property 
\begin{equation}\label{eq: Ytau}
    y + (z'-z)  (f_2(z') -f_2(z) - \tilde f_2(z') (f_1(z') - f_1(z)) ) (\overline{\varphi}_1- \overline{\varphi}_2) \in N_{O(\d^{-\eta}\tau\varrho +\delta)} Y_{z''}({\bf q}). 
\end{equation}
Using Lemma \ref{lem:number-of-Lambda-tuple-kakeya-2} we can ensure that  $G$ is $\delta^{O(\eta)}$-dense
\begin{equation}\label{eq: lowerbdG}
|G| \gtrsim_{\d^{O(\eta)}}| Y_{z_0}({\bf q})| |Z({\bf q})|^2 |\overline{\Phi}_{\bf q}|^2
\end{equation} for $\tau$ satisfying (\ref{eq:tau-defining-relation-2}).

\subsection{Kaufman and radial projections in \texorpdfstring{$\RR^3$}{R3}.}

In the next section we will require 3-dimensional versions of  Theorem \ref{thm:radial-projections} and Theorem \ref{thm:kaufman}. For the radial projections result, we only state the version for sets of points with dimension $2$. This is the only case that we need and it can be done by a relatively simple modification of the $\RR^2$ argument.

\begin{theorem}\label{thm:radial-projections-3d}
    Let $\varepsilon_1, \varepsilon_2, \zeta \in (0,1)$. Suppose that $\varepsilon_1 \le \zeta/100$ .Then there exists $ \delta_0 = \d_0(\varepsilon_1, \varepsilon_2, \zeta), \eta_0 = \eta_0( \varepsilon_1, \varepsilon_2, \zeta)>0$ such that the following holds for all $\delta < \delta_0$, $\eta < \eta_0$.
    Let $P \subset [0,1]^3$ be a $\delta$-separated $(\d, \d^{-\eta})$-uniform $(\delta, 2, \delta^{-\eta})$-set such that $|P \cap S| \le \d^{\varepsilon_2} |P|$ for any $\d^{\varepsilon_1}\times 1\times 1$-slab $S$.

    Then there exists a subset $G \subset P\times P$ such that $|G| \ge (1-\delta^{\eta/10}) |P|^2$ such that for every $(p,p') \in G$ we have $|T_{p,p'}(\varrho) \cap P| \le \delta^{-\zeta} \varrho^{2} |P|$ for all $\varrho\in [\delta, 1]$.

    In particular, for every $p\in P$ such that $|G|_p| \ge \nu |P|$,  the set of tubes $\{T_{p, p'}(\delta), ~(p, p')\in G\}$ is a $(\delta, 2, C\delta^{-\zeta} \nu^{-1})$-set.
\end{theorem}

The proof is given in Appendix \ref{appendix:radial}.

\begin{lemma}\label{lem:kaufman-Rd}
    Let $d\ge 2$ and let $\varepsilon>0$. The following holds for $\d\le \d_0$ and $\eta < \eta_0(d,\varepsilon)$. 
    
    Let $P \subset [-1,1]^d$ be a $\d$-separated $(\d, 1, \d^{-\eta})$-set and let $\Lambda \subset [-1,1]^{d-1}$ be a $\d$-separated $(\d, d-1, \d^{-\eta})$-set. Then there exists $\Lambda' \subset \Lambda$, $|\Lambda'| \ge (1-\d^{\eta})|\Lambda|$, and $P_\lambda \subset P$ for each $\lambda \in \Lambda'$, $|P_\lambda| \ge (1-\d^{\eta}) |P|$
    such that $\pi_{\lambda}(P_\lambda)$ is a $(\d, 1, \d^{-\varepsilon})$-set. Here $\pi_\lambda:\RR^{d}\to\RR$ is the linear projection $(x_1, \ldots, x_d) \to \lambda_1x_1+\ldots+\lambda_{d-1}x_{d-1}+x_d$.
\end{lemma}

\begin{proof}
    The proof is identical to the proof of Theorem \ref{thm:kaufman}, so we skip it.
\end{proof}

\subsection{Proof of Proposition \ref{prop:32-planar-kakeya}} Let $(\TT, E, \{\TT_p\})$ be a 3-planar 2-linear $(\d,\eta, \kappa, \gamma)$-Kakeya configuration such that $C_{CW}(\TT) \lesssim \d^{-\eta}$. Note that both sets $Y_z$ and $\Phi_p$ are subsets in $[-1,1]$. This means that the $\gamma$ parameter does not play a role in the definition of $E$ (we can redefine $\gamma = 1$). Let 
\[
f(z) = \begin{pmatrix}
    f_1(z)\\
    f_2(z)
\end{pmatrix} \in \operatorname{Mat}(2, 1), \quad \tilde f(z) = (\tilde f_1(z), \tilde f_2(z))  \in \operatorname{Mat}(1,2)
\] 

Let $\varrho \in [\d^{1/4}, \d^\mu]$ be an intermediate scale. We consider several cases depending on how $f$ and $\tilde f$ behave on $\varrho$-intervals. First, by Lemma \ref{lem:f-is-Lipschitz}, we can pass to a subset in $Z$ of density $\d^{O(\eta/\kappa)}$ so that $f$ and $\tilde f$ are $\d^{-\eta}$-Lischitz on scale $\varrho$.

By Lemma \ref{lem:non-linear-grains-2-planar}, we either get $\kappa \le \zeta$ or there is some $\varepsilon_0 = \varepsilon_0(\kappa)>0$ and a subset $Z' \subset Z$ of density at least $\d^{\eta'}$ such that for every $z_0 \in Z'$ and ${\bf z} = z_0^\varrho$ there exists a vector $r_{\bf z}=r=(r_1, r_2) \in \RR \times \RR$ so that 
\begin{equation}\label{eq: epsilon0}
|f(z) - f(z') - (z-z')r| \le \varrho^{1+\varepsilon_0}
\end{equation}
holds for all $z, z' \in Z' \cap {\bf z}$. Here we take $\eta'$ to be the $\varepsilon$ from the Lemma \ref{lem:non-linear-grains-2-planar} and assume that $\eta$ is sufficiently small in $\eta'$, $\zeta$, $\mu$. 
Furthermore, by Lemma \ref{lem:f-slow} and Lemma \ref{lem:f-is-Lipschitz} we may assume that \begin{equation} \label{eq: r}
|r| \in [\d^{O_\kappa(\eta)}, \d^{-O_\kappa(\eta)}].
\end{equation}

For $z_0 \in Z$ and $z \in {\bf z} = z_0^\varrho$ define $f_{z_0}(z) = f_2(z) - \tilde f_2(z_0) f_1(z)$ and $\tilde f_{z_0}(z') = \tilde f_2(z')-\tilde f_2(z_0)$.  Consider the expression in \eqref{eq:tau-defining-relation-2}, we now can write
\begin{align*}
&f_2(z) - f_2(z') - \tilde f_2(z') (f_1(z) - f_1(z')) \\
=& f_{z_0}(z) - f_{z_0}(z') - (\tilde f_2(z') - \tilde f_2(z_0)) (f_1(z) - f_1(z'))\\
=& f_{z_0}(z) - f_{z_0}(z') - \tilde f_{z_0}(z') r_1 (z-z') + O(\d^{-\eta}\varrho^{2+\varepsilon_0}).
\end{align*}
Here we used the Lipschitz property $|\tilde f_{z_0}(z')|=|\tilde f_2(z') - \tilde f_2(z_0)| = O(\d^{-\eta}\varrho)$ and the linear approximation of $f_1$.  
Note that after a change of coordinates we may assume that $f(z_0) = 0$ and $\tilde f(z_0) = 0$. In these coordinates we then have $f_{z_0}(z) = f_2(z)$ and $\tilde f_{z_0}(z') = \tilde f_2(z')$. So one should think of functions $f_{z_0}$ and $\tilde f_{z_0}$ as $f_2$ and $\tilde f_2$ written in the coordinates `axis-aligned' at the level $z_0$. 

We split the proof into two cases depending on whether $f_{z_0}(z)$ can be well approximated by a quadratic polynomial. By combining Lemma \ref{lem:quadratic-case} and Lemma \ref{lem:non-quadratic-case} stated in the following sections, we conclude that $\kappa \le \zeta$ for any $\zeta>0$, provided that $\eta$ is small enough. 

More precisely, we choose parameters as follows. Let $\varepsilon = \varepsilon(\zeta)$ be given by Lemma \ref{lem:quadratic-case}. Let $\varepsilon_1=\varepsilon
$ and let $\varepsilon_2 = \varepsilon_2(\zeta, \varepsilon_1)$ be given by Lemma \ref{lem:non-quadratic-case}. Let $\tilde \eta = \min(\varepsilon_2, \tilde \eta(\varepsilon))$ where $ \tilde \eta(\varepsilon)$ is given by Lemma \ref{lem:quadratic-case}. Then if $\eta$ is small enough in $\varepsilon$ and $\kappa>\zeta$, then by Lemma \ref{lem:quadratic-case} we conclude that $f_{z_0}(z)$ cannot be approximated by a quadratic polynomial $a_0+a_1z+a_2z^2$ with error $\varrho^\varepsilon (\varrho^2+\varrho|a_1| + \varrho^2|a_2|)$ for a $\d^{\tilde \eta}$-dense subset of $z \in Z$. This then means that the non-concentration condition in Lemma \ref{lem:non-quadratic-case} is satisfied with $\varepsilon_1=\varepsilon$ and $\varepsilon_2 =\tilde \eta$. So for $\eta$ small enough in $\varepsilon, \tilde \eta, \mu$ we conclude that $\kappa\le \zeta$. This completes the proof of Proposition \ref{prop:32-planar-kakeya} and so it remains to prove the two lemmas.

\subsection{Quadratic case.} \label{subsec:quadratic-case}

The first case we consider is when $f_{z_0}(z)$ can be well approximated by a quadratic polynomial on each $\varrho$-interval.

\begin{lemma}\label{lem:quadratic-case}
    Let $\tilde \eta,\varepsilon>0$ and $\tilde Z \subset Z$ be a subset of density at least $\d^{\tilde \eta}$. 
    Suppose that for every $z_0 \in \tilde Z$ and ${\bf z} = z_0^\varrho$ there are coefficients $a_0, a_1, a_2\in \RR$ such that
    \[
    f_{z_0}(z) = a_0+ a_1 (z-z_0) +a_2 (z-z_0)^2 + O(\varrho^{\varepsilon}(\varrho^2+\varrho |a_1| + \varrho^2 |a_2|))
    \] 
    holds for all $z \in \tilde Z \cap {\bf z}$.  

    Then $\kappa \le \zeta$, provided that $\varepsilon$ is small enough in $\zeta$ and $\eta, \tilde \eta$ is small enough in $\varepsilon$. 
\end{lemma}

\begin{proof}
    Without loss of generality assume that $\varepsilon\le \varepsilon_0/2$ for the $\varepsilon_0$ in \eqref{eq: epsilon0} and $\eta$ is small enough in $\varepsilon_0$ and $\mu$.  
    
    Denote $\Delta = \varrho^{\varepsilon}(\varrho^2+\varrho |a_1| + \varrho^2 |a_2|)$. 
    Using the approximation $ f_1(z)-f_1(z') = r_1(z-z') + O(\varrho^{1+\varepsilon_0}) $ and $|\tilde f_{z_0}(z')| \lesssim \d^{-\eta} \varrho$ we can write
    \begin{align*}
    f_2(z) - f_2(z') - \tilde f_2(z') (f_1(z) - f_1(z'))&
    = f_{z_0}(z)-f_{z_0}(z')-\tilde f_{z_0}(z') r_1(z-z') +O(\delta^{-\eta}\varrho^{2+\varepsilon_0})           \\
    &= a_1 (z-z') + a_2 ((z-z_0)^2-(z'-z_0)^2) -\tilde f_{z_0}(z') r_1(z-z') + O(\Delta)    \\
    &= (z-z') (a_1 + a_2(z+z'-2 z_0) - r_1\tilde f_{z_0}(z'))+ O(\Delta).
    \end{align*}

    Since $|r_1|\lesssim_{\d^{O(\eta)}} 1$ (we absorb  $\kappa$ into the constant),  we have
    \begin{equation}\label{eq: shiftr1}
    |r_1\tilde f_{z_0}(z'))|\lesssim \delta^{-O(\eta)} \varrho. 
    \end{equation}

    To apply identities on $L$-tuples we need to define the scale $\tau$ according to \eqref{eq:tau-defining-relation-2}. We make a preparatory refinement of $\tilde Z$ in order to make the definition of $\tau$ stable under passing to subsets of density $\d^{O(\eta)}$. 
    Let $\eta_j$, $j=1, \ldots$ be an increasing sequence of parameters starting with $\eta_0= \eta$ and $\eta_{j+1} = \eta_j^{1/2}$. Define $\tilde Z_j\subset \tilde Z$ to be the set which satisfies $|\tilde Z_j| \ge \d^{\eta_j} |\tilde Z|$ and so that 
    \[
    \tau_j:=\varrho \max_{z,z' \in \tilde Z_j, z'\in z^\varrho} |f_2(z) - f_2(z') - \tilde f_2(z') (f_1(z) - f_1(z'))|
    \]
    is smallest possible.  Our goal will be to show that $\tau_j = O(\varrho^{3+\varepsilon/2})$ holds for some $j$ depending on $\zeta$ only.    
    Let $m  \geq 20/\varepsilon$ be an integer  and assume that for the sake of contradiction we have $\tau_j > \varrho^{3+ \varepsilon/2}$ for all $j=1, \ldots,m$. Select smallest $j$ so that $\tau_{j+1} \ge \varrho^{2/m} \tau_j \geq \varrho^{\varepsilon/10} \tau_j$ holds. By the assumption, there exists such a $j\le m$.

    We can further prune $\tilde Z_j$ slightly so that the maximum $\tau$ is essentially the same on every $\varrho$-interval. Let $\eta' = \eta_j$. For each $p_0\in E_{z_0}\in E$, define ${\bf q}=  p_0^{\varrho}$ and $Z({\bf q})= \tilde Z_j\cap z_0^{\varrho}$, $E(\mathbf{q})= E_{\tilde Z_j} \cap {\bf q}$.  Let $G \subset Y_{z''}({\bf q})\times Z({\bf q})^2 \times \overline{\Phi}_{\bf q}^2$ be the graph in \eqref{eq: G} obtained from the $L$-tuples defined on the set $E({\bf q})$. We can choose 
    \begin{equation}\label{eq: tau}
    \tau \sim \varrho \max_{z,z' \in Z({\bf q})} \left((z-z') (a_1 + a_2(z+z'-2z_0) - r_1\tilde f_{z_0}(z'))+ O(\Delta)\right)
    \end{equation}
    and so for every tuple $(y, z, z', \overline{\varphi}_1, \overline{\varphi}_2)\in G$ we have \eqref{eq: Ytau}, which is
    \begin{equation}\label{eq: Yz''}
    y + (z-z')^2 (a_1 + a_2(z+z'-2z_0) - r_1\tilde f_{z_0}(z')) (\overline{\varphi}_1- \overline{\varphi}_2) + O(\varrho \Delta) \in N_{O(\d^{-\eta'}\tau \varrho +\d + \varrho \Delta)} Y_{z''}({\bf q}).
    \end{equation}


    By averaging and using the lower bound \eqref{eq: lowerbdG} on the size of $G$, we can fix variables $y$, $z'$, $\overline{\varphi}_1, \overline{\varphi}_2$ so that $|\overline{\varphi}_1- \overline{\varphi}_2| \gtrsim_{\d^{O(\eta')}} 1$ and  there are $\gtrsim_{\d^{O(\eta')}} |Z({\bf q})|$ choices of $z$ such that $(y, z, z', \overline{\varphi}_1, \overline{\varphi}_2) \in G$ and $|z-z'| \gtrsim_{ \d^{O(\eta')}}\varrho$. Denote this set of such $z$ by $\overline{Z} \subset Z({\bf q})$. Let 
    \[P(z) = (z-z')^2 (a_1 + a_2(z+z'-2z_0) - r_1\tilde f_{z_0}(z')) (\overline{\varphi}_1- \overline{\varphi}_2),\]
    note that this is a cubic polynomial in the variable $z$. By the definition of $\tau$, we have $|P(z)| \lesssim  \tau$ and equality is achieved up to $\d^{O(\eta')}$ if we vary both $z$ and $z'$.

    By construction $|\overline{Z}|_\d \gtrsim_{\d^{O(\eta')}}  |Z({\bf q})|_\d \gtrsim_{\d^{O(\eta')}} \varrho/\d$. An image of a  $\delta^{O(\eta')}$-dense set under a bounded degree polynomial is again  a $\delta^{O(\eta')}$-dense on the image (with density bound independent of its coefficients using Remez inequality). 
Definition~\ref{def:kakeya-configuration-2}  with $d=4, k=3$ implies $Y_{z''}({\bf q})$ is a rescaled $(\d/\varrho, 1-\kappa, \d^{-O(\eta')})$-AD-regular set.      
    It follows from  \eqref{eq: Yz''} that $y+ P(\overline{Z})\subset N_{O(\delta^{-\eta'}\tau \varrho+\delta+\varrho \Delta)} Y_{z''}(\bf q)$ and so
 the image $P(\operatorname{conv}(\overline{Z}))$ is contained in an interval $I$ of length $\lesssim \d^{-O(\eta'/\kappa)}  (\tau \varrho+\delta + \varrho \Delta)$.

Consider 
\[
P'(z)= (z-z')(2(a_1-r_1\tilde f_{z_0}(z')) + a_2(3z-z'-2z_0))(\overline{\varphi}_1-\overline{\varphi}_2).
\]
Then $|P'(z)|\sim_{\delta^{O(\eta')}} \varrho |2(a_1-r_1\tilde f_{z_0}(z')) + a_2(3z-z'-2z_0)|$ and for a $\gtrsim 1$-fraction of $z$ satisfying  $|z-z'|\gtrsim_{\delta^{O(\eta')}} \varrho$ and $|P'(z)|\gtrsim_{\delta^{O(\eta')}} |a_2|\varrho^2$. Therefore $|P(\operatorname{conv}(\overline{Z}))|\gtrsim_{\delta^{O(\eta')}} |a_2|\varrho^3$ and so 
\begin{equation}\label{eq: a_2}
    |a_2|\lesssim_{\delta^{O(\eta')}} \tau/\varrho^2 + \delta/\varrho^3+ \Delta/\varrho^{2}
\end{equation}
where we treat $1/\kappa$ as a constant and put it in $O(1)$.  

Plug in $P(z)$ we get
\[
P(z)= (z-z')^2 (a_1 - r_1\tilde f_{z_0}(z'))(\overline{\varphi}_1-\overline{\varphi}_2)+ O(\delta^{-O(\eta')})( \tau \varrho + \delta +\varrho \Delta)
\]

The polynomial $Q(z)= (z-z')^2(a_1-r_1\tilde f_{z_0}(z')(\overline{\varphi}_1-\overline{\varphi}_2)$ also takes $\overline{Z}$ to a $\delta^{O(\eta')}$-dense set of an interval of length $\sim |a_1-r_1\tilde f_{z_0}(z')|\varrho^2$, which implies 
\begin{equation}\label{eq: a1}
    |a_1-r_1\tilde f_{z_0}(z')|\lesssim_{\delta^{O(\eta')}} \tau/\varrho+\delta/\varrho^2+\Delta/\varrho. 
\end{equation}

 Now note that this relation holds for any $z'$ for which there exists a choice of well-separated angles $\overline{\varphi}_1, \overline{\varphi}_2$ and $y$ such that the fiber of $G$ above $(y, z', \overline{\varphi}_1, \overline{\varphi}_2)$ is sufficiently dense. By double counting we get that there are at least $\d^{O(\eta')}|Z({\bf q})|$-many such choices for $z'$. It follows that for any $(z, z')$ chosen as above we get 
    \begin{equation}\label{eq: taubound}
    |(z-z')^2 (a_1 + a_2(z+z'-2 z_0) - r_1 \tilde f_{z_0}(z'))| \lesssim_{\d^{O(\eta')}} \tau \varrho +\delta+\Delta\varrho.
    \end{equation}

    Let $\tilde Z'_j \subset \tilde Z_j$ be the subset of $z'$ such that \eqref{eq: a_2} and \eqref{eq: a1} hold for the corresponding $a_1, a_2$ with respect to their dyadic $\varrho$-interval.  Then     and for any $z', z'' \in \tilde Z_j'$ in the same dyadic $\varrho$-interval,
\begin{equation}\label{eq: tildefz0}
|\tilde f_{z_0}(z')-\tilde f_{z_0}(z'')|\lesssim_{\delta^{O(\eta')}} \tau /\varrho +\d/\varrho^2+\Delta /\varrho
\end{equation}
and \eqref{eq: taubound} holds with $z''$ in place of $z$. 
    
    
    For constant $j$ and small enough $\eta$, we have $\eta_{j+1} =\eta_j^{1/2} \gg \eta_{j}$. In particular we get $|\tilde Z_j'| \gg \d^{\eta_{j+1}} |\tilde Z|$. So we can take $\tilde Z'_j$ as a candidate for the set $Z_{j+1}$ and conclude that we have the following estimate:
    \[
    \varrho^{\varepsilon/10}   \tau_j \leq \tau_{j+1} \lesssim_{\d^{O(\eta')}} \tau \lesssim \tau_j,
    \]
    where the first inequality is by the choice of $j$. 
    

Plug in \eqref{eq: tau}, we get 
$\tau \lesssim_{\delta^{O(\eta')}} \varrho^{-\varepsilon/10} (\varrho \tau+\delta  + \varrho\Delta)$ and so 
\begin{equation}\label{eq: tau-delta-estimate}
    \tau \lesssim_{\delta^{O(\eta')}} \varrho^{-\varepsilon/10} (\delta  + \varrho\Delta)   
\end{equation}
Using $|\tilde f_{z_0}(z')|\lesssim \delta^{-\eta}\varrho$, $|r_1|\lesssim_{\d^{O(\eta)}} 1$ and \eqref{eq: a1}, we get 
\begin{equation}\label{eq: a1'}
|a_1|\lesssim_{\d^{O(\eta')}} \varrho +\tau/\varrho +\delta/\varrho^2+\Delta/\varrho.
\end{equation}
Combining \eqref{eq: a1'} with \eqref{eq: a_2}, we get  
\begin{align*}
    \Delta &= \varrho^{\varepsilon} (\varrho^2 +\varrho|a_1| +\varrho^2|a_2|) \\
    &\lesssim_{\d^{O(\eta')}} \varrho^{\varepsilon} (\varrho^2 + \tau+\d/\varrho + \Delta)
\end{align*}
which implies for $\eta' \ll \varepsilon$,
\[
\Delta \lesssim_{\d^{O(\eta')}} \varrho^{\varepsilon} (\varrho^2 + \tau+\d/\varrho).
\]
Thus, we can continue \eqref{eq: tau-delta-estimate}:
\[
\tau \lesssim_{\d^{O(\eta')}} \varrho^{-\varepsilon/10} (\delta + \varrho^{1+\varepsilon} (\varrho^2 + \tau +\delta/\varrho)) \lesssim \varrho^{-\varepsilon/10}\d + \varrho^{3+9\varepsilon/10} +  \varrho^{1+\varepsilon/2} \tau.
\]
When $\varrho \in [\delta^{1/4}, \delta^{\mu}]$ and $\eta'>0$ sufficiently small, we get 
\begin{equation}\label{eq: taubound-2}
    \tau \leq  \varrho^{3+\varepsilon/2}.
\end{equation}

Then \eqref{eq: tildefz0} becomes 
\[
|\tilde f_{z_0}(z')-\tilde f_{z_0}(z'')|\lesssim_{\delta^{O(\eta')}} \varrho^{2+\varepsilon/2},
\]
which  by the Definition of $\tilde f_{z_0}(z')$ implies 
\[
|\tilde f_2 (z')- \tilde f_2 (z'')|\lesssim_{\delta^{O(\eta')}} \varrho^{2+\varepsilon/2}. 
\]

Next we are going to use $|f_{z_0}(z')-f_{z_0}(z'')|$ to get a bound involving $f_2(z')-f_2(z'')$ and then combine it with  the consistency relation \eqref{eq:ff-relation} to bound $|\tilde f_1(z') -\tilde f_1(z'')|$. 

By \eqref{eq: a_2} and \eqref{eq: a1'}, we have 
\[
|f_{z_0}(z')-f_{z_0}(z'')|\lesssim_{\delta^{O(\eta')}} \varrho^2,
\]
which implies
\[
|f_2(z')-f_2(z'') - \tilde f_2(z_0) (f_1(z')-f_1(z''))|\lesssim_{\delta^{O(\eta')}} \varrho^2.
\]
By the consistency relation \eqref{eq:ff-relation}, 
\begin{align*}
&|\tilde f_1(z') -\tilde f_1(z'')| = |f_2(z')-f_2(z'') -\tilde f_2(z')f_1(z')+\tilde f_2(z'') f_1(z'')| + O(\delta)\\
=&  |f_2(z')-f_2(z'') - \tilde f_2(z_0) f_1(z')  + \tilde f_2(z_0) f_1(z'') 
+(\tilde f_2(z_0)-\tilde f_2(z')) f_1(z') 
-(\tilde f_2(z_0)-\tilde f_2(z''))f_1(z'')| + O(\delta)\\
\leq & |f_2(z')-f_2(z'') - \tilde f_2(z_0) f_1(z')  + \tilde f_2(z_0) f_1(z'') | + |(\tilde f_2(z_0)-\tilde f_2(z')) f_1(z')| + |(\tilde f_2(z_0)-\tilde f_2(z''))f_1(z'')| +O(\delta)\\
\lesssim&_{\delta^{O(\eta')}} \varrho^2 
\end{align*}

So we finally conclude for any $z', z''\in \tilde Z_{j+1}$, 
\[
|\tilde f(z') - \tilde f(z'')|\leq |\tilde f_1(z') -\tilde f_1(z'')| + |\tilde f_2(z') -\tilde f_2(z'')| \lesssim_{\delta^{O(\eta')}} \varrho^2. 
\]

 By applying Lemma \ref{lem:f-tilde-slow} and recalling the lower bound $|\tilde Z_{j+1} (\bf q)| \gtrsim_{\d^{O(\eta')}}(\varrho/\d) $ we therefore conclude that $C_{CW}(\TT) \gtrsim_{\d^{O(\eta')}} \varrho^{-\kappa}$. 

This completes the proof when $\eta>0$ is sufficiently small and so $\eta'$ is sufficiently small. 
\end{proof}

\subsection{Non-quadratic case.}\label{subsec:non-quadratic-case}

Now we consider the case when $f_{z_0}$ cannot be approximated by a quadratic polynomial. 

\begin{lemma}\label{lem:non-quadratic-case}
    Let $\varepsilon_1, \varepsilon_2 \in (0, \varepsilon_0/2)$ and suppose that for every $z_0 \in Z$, ${\bf z} = z_0^\varrho$ and any coefficients $a_0, a_1, a_2\in \RR$ and $\Delta = \varrho^{\varepsilon_1}(\varrho^2 + \varrho |a_1| + \varrho^2 |a_2|)$ we have
    \begin{equation}\label{eq: nonquadratic}
    \#\{ z \in Z\cap {\bf z}:~ |f_{z_0}(z) - a_0 -a_1(z-z_0)-a_2(z-z_0)^2| \le \Delta \} \le \varrho^{\varepsilon_2} |Z\cap {\bf z}|.
    \end{equation}
    Then, provided that $\varepsilon_2$ is small enough in $\zeta, \varepsilon_1$, and $\eta$ is small enough in $\varepsilon_1, \varepsilon_2, \mu$, we have $\kappa \le \zeta$.
\end{lemma}

\begin{proof}

\medskip
\noindent \textit{ Step 1. Refine $Z$ so that $\tau$ is stable under refinement.}   Similarly to the proof of Lemma \ref{lem:quadratic-case}, we make a preparatory refinement of $Z$ in order to make the definition of $\tau$ stable under passing to subsets of density $\d^{O(\eta)}$. Below we will define a certain continuous function $h:(0,1)\to (0,1)$ that depends only on the parameters $\zeta, \mu>0$ (think of something like $h(\eta) = \eta^2$).
Let $\eta_j$, $j=1, \ldots$ be an increasing sequence of parameters starting with $\eta_0= \eta$ and such that $h(\eta_{j+1}) \ge \eta_j$. By taking $\eta$ small enough we can ensure each $\eta_j$ is small enough for every fixed value of $j$. Define $Z_j\subset Z$ to be the set which satisfies $| Z_j| \ge \d^{\eta_j} |Z|$ and so that 
    \[
    \tau_j:=\varrho\max_{z,z' \in Z_j, z'\in z^\varrho} |f_2(z) - f_2(z') - \tilde f_2(z') (f_1(z) - f_1(z'))|
    \]
    is smallest possible.
    Fix a large constant $m\ge 1$. We will choose $\eta$ small enough so that $\eta_m$ is much smaller than $\varepsilon_2$.
    First, suppose that $\tau_j \le \varrho^{3+2\varepsilon_1}$ for some $j\le m$. Then for $z, z' \in Z_j$ we have 
    \[
    |f_2(z) - f_2(z) - \tilde f_2(z') (f_1(z) - f_1(z'))| \lesssim \varrho^{2+2\varepsilon_1}.
    \]
Since $f_{z_0}(z) = f_2(z)-\tilde f_{2}(z_0) f_1(z)$, we get 
\begin{equation}
    |f_{z_0}(z) + (\tilde f_2(z_0)-\tilde f_2(z'))f_1(z) - f_2(z') + \tilde f_2(z')f_1(z')|\leq \varrho^{2+2\varepsilon_1}
\end{equation}

Using the approximation $f_1(z)-f_1(z_0) = r_1(z-z_0) +O(\varrho^{1+\varepsilon_0})$ and $|\tilde f_2(z_0)-\tilde f_2(z')|\lesssim_{\delta^{O(\eta)}} \varrho$, we get 
\[
|f_{z_0}(z) + (\tilde f_2(z_0) -\tilde f_2(z'))(r_1 (z-z_0) +f_1(z_0)) -f_2(z') +\tilde f_2(z')f_1(z')|\lesssim   \varrho^{2+2\varepsilon_1}
\]
since $\varepsilon_1 \leq \varepsilon_0/2$. 

Taking $a_1= -(\tilde f_2 (z_0)-\tilde f_2(z') r_1,  \, a_0= -(\tilde f_2 (z_0)-\tilde f_2(z') f_1(z_0) + f_2(z') -\tilde f_2(z') f_1(z')$ and $a_2=0$, we get 
\[
\#\{z\in Z\cap \mathbf{z}: |f_{z_0}(z) - a_0 -a_1(z-z_0)|\lesssim_{\delta^{O(\eta)}} \varrho^{2+2\varepsilon_1}\} \gtrsim |Z_j|\gtrsim \delta^{O(\eta_j)} |Z\cap \mathbf{z}|, 
\]
which yields a contradiction with the assumption of the lemma if we choose $\eta$ sufficiently small depending on $\varepsilon_1, \varepsilon_2, \mu$. 
    
Thus, we may assume that $\tau_j > \varrho^{3+2\varepsilon_1}$ holds for all $j\le m$. 
By the pigeonhole principle we can find some $j \le m$ such that $\tau_{j+1} \ge \varrho^{1/m} \tau_j$ holds. Write $\eta' = \eta_{j}$.

    \medskip
\noindent \textit{Step 2. Sum-product information from $L$-tuples.}
Let $Z({\bf q}) = Z_j\cap z_0^\varrho$ and let $E({\bf q}) =E_{Z_j}\cap \mathbf{q}$. Let $G \subset Y_{z''}({\bf q})\times Z({\bf q})^2\times\overline{\Phi}_{\bf q}^2$ be a dense graph constructed using $L$-tuples $L_{z, z', z''}$ in the previous sections. Note that by construction, we have $\tau \sim_{\d^{O(\eta')}} \tau_j$ in the definition of $L$-tuples. 
    By construction, for every $(y, z, z', \overline{\varphi}_1,\overline{\varphi}_2) \in G$ we have 
    \[
    y + (z'-z) (f_2(z) - f_2(z') - \tilde f_2(z') (f_1(z)-f_1(z')))(\overline{\varphi}_1-\overline{\varphi}_2) \in N_{C\d^{-\eta'} \tau\varrho} Y_{z''}({\bf q})
    \]
    where we use $\delta \leq \delta^{-\eta'}\tau\varrho$.
    Cover $Y_{z''}({\bf q})$ with intervals $I_j$ of length $\sim \tau$. Since $Y_{z''}(\mathbf{q})$ is a rescaled $(\delta/\varrho, 1-\kappa, \delta^{-\eta})$-AD regular set, $Y_j = I_j \cap Y_{z''}({\bf q})$ has size $\sim_{\d^{O(\eta')}} (\tau / \d)^{1-\kappa} $ for every $j$. Let $G_j \subset Y_j \times  Z({\bf q})^2\times\overline{\Phi}_{\bf q}^2$ be the restriction of $G$ onto $Y_j$. Then there exists an index $j$ so that $|G_j| \gtrsim_{\d^{O(\eta')}} |Y_j|  |Z({\bf q})|^2|\overline{\Phi}_{\bf q}|^2$ holds. Fix a pair $(\overline{\varphi}_1, \overline{\varphi}_2)$ such that $|G_j|_{\overline{\varphi}_1, \overline{\varphi}_2}| \gtrsim_{\d^{O(\eta')}} |Y_j|  |Z({\bf q})|^2$ and $|\overline{\varphi}_1- \overline{\varphi}_2| \gtrsim_{\d^{O(\eta')}}1$.

    After a change of coordinates we may assume that $z_0 = 0$ and $f(z_0) = 0$ and $\tilde f(z_0) =0$. Recall that we have approximations $f_1(z) = r_1 z + O(\varrho^{1+\varepsilon_0})$ and $\tilde f_2(z) = O(\d^{-\eta}\varrho)$ for all $z \in Z({\bf q})$. 
    Define 
    \[
    w=C\delta^{-\eta'}(\varrho + \frac{\varrho^{3+\varepsilon_0}}{\tau} )
    \]
    for a sufficiently large constant $C$. By definition, it follows that for any $(y, z_1, z_2) \in G_j|_{\overline{\varphi}_1, \overline{\varphi}_2}$ we have $y\in Y_j$ and
    \begin{equation}\label{eq: Yz''w}
    y + (z'-z) (f_2(z) - f_2(z') -\tilde f_2(z')r_1(z-z')) (\overline{\varphi}_1-\overline{\varphi}_2) \in N_{\tau w}Y_{z''}({\bf q}) \cap C I_j.
    \end{equation}
    (Note that changing coordinates does not affect the validity of this expression). Let $\psi:[-1,1] \to I_j$ be the affine bijection and let $\mc A \subset w\ZZ \cap [-1,1]$ be the set of points such that $\psi(\mc A^{(w)})$ forms a covering of $Y_j$. Let $\mc B = \{\varrho^{-1}z\mid z \in Z({\bf q}) \} \subset [-1,1]$. Note that by the construction we have the covering lower bound $|\mc B|_w \gtrsim_{\d^{O(\eta')}} w^{-1}$. 
    Define rescaled functions $F, \tilde F: \mc B\to \RR$ by
    \[
    F(b) = \varrho \tau^{-1}f_2( \varrho b) (\overline{\varphi}_1-\overline{\varphi}_2), \quad \tilde F(b) = \varrho^2 \tau^{-1}\tilde f_2(\varrho b) r_1(\overline{\varphi}_1-\overline{\varphi}_2).
    \]
    From the definition of $\tau$, we have $|F(b) - F(b') - \tilde F(b') (b-b')| \lesssim_{\d^{O(\eta')}} 1$ for all $b, b' \in \mc B$ ($\tau$ was essentially defined so that this expression is about 1 for typical choice of $b, b'$). By swapping $b, b'$ and triangle inequality it follows that 
    \[
    |(\tilde F(b')-\tilde F(b))(b-b')|\lesssim_{\d^{O(\eta')}} 1.
    \]
    Since this holds for all $b, b' \in \mc B$ and $\mc B$ is a dense subset in $[-1,1]$, it follows that $\{\tilde F(b')\}$ is contained in an interval of length $\d^{O(\eta')}$. Recall that we chose coordinates so that $\tilde f_2(z_0) = 0$, which implies that the interval containing $\tilde F(b')$ must be around the origin. Using this we can now conclude that the same conclusion holds for $F(b)$: we have $|F(b)| \lesssim_{\d^{O(\eta')}}1$ for all $b \in \mc B$. 

    \medskip
\noindent \textit{Step 3. Large $\tau$ case.} Let $\zeta_1>0$ be a small constant that depends on $\zeta, \mu$ only, to be determined. First, we consider the case when $\tau > \varrho^{3-\zeta_1}$ holds. 
Using $|\tilde f_2(z')| \lesssim_{\d^{O(\eta')}} \varrho$ and $|r_1| \lesssim_{\d^{O(\eta')}} 1$, \eqref{eq: Yz''w} then implies 
\[
y+ (z'-z) (f_2(z) - f_2(z')) (\overline{\varphi}_1 -\overline{\varphi}_2) \in N_{2\tau w}Y_{z''}({\bf q}) \cap CI_{j}.
\]
Let $\mc G \subset \mc A \times \mc B\times \mc B$ be the graph rescaled from $G_j|_{\overline{\varphi}_1,\overline{\varphi}_2}$. Then we obtain 
\[
|\{a+ (b'-b) (F(b)-F(b')), (a, b, b') \in \mc G\}|_w \lesssim_{\d^{O(\eta')}} |\mc A|_w.
\]
We would like to apply Theorem \ref{thm:ABC-with-F} on scale $w$ to conclude that $|\mc A|_w \gtrapprox w^{-1}$. Note that by the assumption on $\tau$,
\[
w \lesssim \d^{-\eta'} \varrho^{\varepsilon_0},
\]
and $\varepsilon_0$ depends on $\kappa$ only, so assuming that $\kappa\ge \zeta$ we have $\varepsilon_0 \gtrsim_{\zeta} 1$.

The parameters in this application are $\alpha = 1-\zeta$, $\beta=1$, $\gamma=0$, let $\overline{\varepsilon}_1 \ll_{\zeta} 1$ be the constant $\varepsilon_1$ in Theorem \ref{thm:ABC-with-F} and let $\overline{\varepsilon}_2>0$ be a small constant to be determined. Suppose that there exist $u,v \in \RR$ such that
\[
|\{ b \in \mc B: ~ |F(b) - v b-u| \le w^{\overline{\varepsilon}_1}  \}| \ge w^{\overline{\varepsilon}_2} |\mc B|.
\]
Plugging $b = \varrho^{-1} z$ and $F(b)= \varrho \tau^{-1} f_2(z)$ gives
\[
|  f_2(z) - \tau \varrho^{-2} v z - \tau \varrho^{-1} u| \le  (\tau/\varrho) w^{\overline{\varepsilon}_1}
\]
Put $a_0 = \tau \varrho^{-1} u$, $a_1 = \tau \varrho^{-2} v$ and $a_2=0$ and let $\Delta = \varrho^{\varepsilon_1} (\varrho^2 +\varrho |a_1| + \varrho^2 |a_2|)$. If $(\tau/\varrho) w^{\overline{\varepsilon}_1} \le \Delta $ then by the assumption of the lemma we conclude that
\[
w^{\overline{\varepsilon}_2} |\mc B| \le \varrho^{\varepsilon_2} |Z \cap {\bf z}|.
\]
Suppose that $\overline{\varepsilon}_2$ satisfies $\varepsilon_0 \overline{\varepsilon}_2 \ll \varepsilon_2$ and $\eta'$ is small enough, then this gives a contradiction. Thus, we may assume 
\[
(\tau/\varrho) w^{\overline{\varepsilon}_1} > \Delta = \varrho^{\varepsilon_1} (\varrho^2 +\varrho |a_1| + \varrho^2 |a_2|) = \varrho^{\varepsilon_1} (\varrho^2 + (\tau/ \varrho) |v|).
\]
Suppose that $\varepsilon_1 < \overline{\varepsilon}_1 \varepsilon_0/2$ holds, this then implies 
\[
|v| \le \varrho^{\overline{\varepsilon}_1 \zeta_1 - \varepsilon_1} \le w^{\overline{\varepsilon}_1 /2} \le \varrho^{\varepsilon_1}.
\]
We conclude that
\[
\{b \in \mc B:~ |F(b) - u| \le w^{\overline{\varepsilon}_1/2}\} \ge w^{\overline{\varepsilon}_2} |\mc B|.
\]
\begin{equation}\label{eq:f2-concentrates-on-interval}
|\{z \in Z({\bf q}): ~ |f_2(z) - \tau \varrho^{-1} u | \le (\tau/\varrho)\varrho^{\varepsilon_1} \}| \gtrsim_{\d^{O(\eta')}} w^{\overline{\varepsilon}_2} |Z \cap {\bf z}|.    
\end{equation}
Recall that $\eta' = \eta_j$ and $\eta_{j+1}$ is defined so that $h(\eta_{j+1}) \ge \eta_j$ for some function $h$ that we need to specify. Recall the function $\overline{\eta}(\alpha, \beta, \gamma, \varepsilon_1, \varepsilon_2)$ from Theorem \ref{thm:ABC-with-F}. Let us suppose that $\overline{\varepsilon}_2$ is chosen so that 
\[
w^{\overline{\varepsilon}_2} \gg  \d^{ \eta_{j+1}/2 }
\]
and on the other hand
\[
\eta' \ll \overline{\eta}(1-\zeta, 1, 0, \overline{\varepsilon}_1, \overline{\varepsilon}_2).
\]
Note that the right hand size is a function of $\zeta, \varepsilon_0$ and $\overline{\varepsilon}_2$. Thus, we can define the function $h$ depending only on $\zeta,\mu$ to arrange this. 

By combining the conclusion \eqref{eq:f2-concentrates-on-interval} with the assumption $\tau \ge \varrho^{3-\zeta_1}$ and with $\tilde f_2(z') (f_1(z)-f_1(z')) \lesssim_{\d^{O(\eta')}} \varrho^{2}$ we conclude that on a subset of $Z({\bf q})$ of density at least $w^{\overline{\varepsilon}_2}$ we have
\[
|f_2(z) - f_2(z') - \tilde f_2(z') (f_1(z) - f_1(z'))| \lesssim (\tau/\varrho) \varrho^{\varepsilon_1}.
\]

By the definition of $\tau_{j+1}$ and $\tau_j\sim \tau$ we then conclude that
\[
\tau_{j+1} \lesssim \varrho^{\varepsilon_1}\tau_j.
\]
Note that $\varepsilon_1$ needs to satisfy $\varepsilon_1 \le \overline{\varepsilon}_1\varepsilon_0/2$ where $\overline{\varepsilon}_1$ depends only on $\zeta$ and $\varepsilon_0$ depends only on $\zeta$. This means that it suffices to assume $\varepsilon_1 \ll_{\zeta, \mu} 1$. Recall that we chose $\tau_j$ by pigeonholing $j$ over $j=1, \ldots, m$ for some parameter $m$ chosen in advance. We now select $m \sim \frac{1}{\varepsilon_1}$ so that $\tau_{j+1} \ge \varrho^{1/m} \tau_j$ contradicts the above conclusion. 

Note that the value of $m$ needs to satisfy the property that $\eta_m$ is sufficiently small depending on $\varepsilon_1, \varepsilon_2, \mu$. Clearly, this holds by taking $\eta\ll_{\zeta, \mu, \varepsilon_1, \varepsilon_2} 1$. 

With these choices of parameters the argument above leads to a contradiction. This means that the initial assumption that the function $F$ did not satisfy the non-concentration conditions with exponents $(\overline{\varepsilon}_1, \overline{\varepsilon}_2)$ is false. Since the graph $\mc G \subset \mc A \times \mc B\times \mc B$ is $\d^{O(\eta')}$-dense by assumption and we selected $\eta' \ll \overline{\eta}(1-\zeta, 1, 0, \overline{\varepsilon}_1, \overline{\varepsilon}_2)$, we can apply Theorem \ref{thm:ABC-with-F} to conclude that 
\[
|\mc A|_w \gtrsim w^{-(1-\zeta)}
\]
and so since by construction $\mc A$ is $(w,1-\kappa,\d^{-O(\eta')} )$-AD-regular set, the conclusion $\kappa \le \zeta$ follows (after replacing $\zeta$ by $\zeta/2$).

    \medskip
\noindent \textit{Step 4. Non-concentration for $F$.}
We can thus assume that $\tau \le \varrho^{3-\zeta_1}$ holds for small parameter $\zeta_1$ which is allowed to depend on $\zeta$ and $\mu$.
Recall that by our assumption \eqref{eq: nonquadratic}, the function $f_{z_0}(z) = f_2(z)$ cannot be approximated by a quadratic polynomial. 
    We claim that for any $a_0, a_1, a_2\in \mathbb{R}$ and $\overline{\Delta}= \varrho^{2\varepsilon_1+\zeta_1} (1+|a_1|+|a_2|)$ we have the following non-concentration condition on $F$:
    \begin{equation}\label{eq:quadratic-non-concentration}
    \#\{ b \in \mc B:~ |F(b) - a_0 - a_1 b- a_2 b^2| \le \overline{\Delta}  \}  < \varrho^{\varepsilon_2/2} |\mc B|.
    \end{equation}
    Indeed, assume the contrary and plug in $b=\varrho^{-1}z$ and $F(b) = \varrho \tau^{-1} (\overline{\varphi}_1-\overline{\varphi}_2)  f_2(z)$. Denote $\phi = \overline{\varphi}_1-\overline{\varphi}_2\sim_{\d^{O(\eta')}} 1$.
    Then for $\d^{C\eta'}\varrho^{\varepsilon_2/2} \gg \varrho^{\varepsilon_2}$ fraction of $z \in Z({\bf q})$ we obtain:
    \[
    | \varrho \tau^{-1}\phi f_2(z) - a_0 - a_1 \varrho^{-1} z  - a_2 \varrho^{-2} z^2| \le \overline{\Delta},
    \]
    \[
    | f_2(z) - (\tau/\varrho \phi)a_0 - (\tau/\varrho^2\phi)a_1 z - (\tau / \varrho^3\phi) a_2z^2 | \le  \overline{\Delta} (\tau /\varrho\phi),
    \]
    put $a'_i = (\tau/\varrho^{1+i}\phi)a_i$, then 
    \[
    | f_2(z) -a'_0 - a'_1 z- a'_2 z^2 | \le  \overline{\Delta} (\tau /\varrho\phi) = \varrho^{2\varepsilon_1+\zeta_1} ( \tau/\varrho\phi + \varrho |a_1'| + \varrho^2 |a_2'| ).
    \]
    Since $\tau \le \varrho^{3-\zeta_1}$ this implies that for $\Delta = \varrho^{\varepsilon_1} (\varrho^2 + \varrho |a_1'| + \varrho^2 |a_2'|)$ we have 
    \[
    | f_2(z) -a'_0 - a'_1 z- a'_2 z^2 | \le \varrho^{2\varepsilon_1} ( \varrho^2/\phi + \varrho |a_1'| + \varrho^2 |a_2'| ) \le \Delta,
    \]
    and so we obtain a contradiction with the non-concentration assumption on $f_2 = f_{z_0}$. 
    Here we add another factor of $\varrho^{\varepsilon_1}$ to subsume the factor $\phi^{-1}$ in the definition of $F$.

    \medskip
\noindent \textit{Step 5. Reduction to projection problems.}
     Let $\mc G \subset \mc A \times \mc B\times \mc B$ be the graph rescaled from $G_j|_{\overline{\varphi}_1,\overline{\varphi}_2}$. We then obtain $|\mc G| \gtrsim_{\d^{O(\eta')}} |\mc A| |\mc B|^2$ and \eqref{eq: Yz''w} translates to  the following non-expansion property:
    \begin{equation}\label{eq:F-non-expanding-32-planar}
    |\{  a + (b'-b) (F(b) - F(b') - \tilde F(b') (b-b')), ~~ (a, b, b') \in \mc G \}|_w     \lesssim_{\d^{O(\eta')}} |\mc A|_w.
    \end{equation}

    Define $S(b, b') =  (b'-b) (F(b) - F(b') - \tilde F(b') (b-b'))$. For $b_1, b_2, b_1', b_2' \in \mc B$ define the polarized expression:
    \begin{equation}\label{eq:Q-polarization-32-planar}
    Q(b_1, b_2, b_1', b_2') = S(b_1, b_1') - S(b_1, b_2') - S(b_2, b_1') + S(b_2, b_2'). 
    \end{equation}
    Furthermore, let $S_1(b, b') = (b'-b) (F(b) - F(b'))$ and $S_2(b, b') = (b'-b)^2 \tilde F(b')$ so that $S=S_1+S_2$ and let $Q=Q_1+Q_2$ be the corresponding decomposition of the polarization. Then by a direct computation we have:
    \[
    Q_1(b_1, b_2, b_1', b_2') = (b_1-b_2) (F(b_1') - F(b_2')) + (b_1' - b_2') (F(b_1) - F(b_2))
    \]
    and 
    \begin{align*}
        Q_2(b_1, b_2,  b_1', b_2') &= ((b_1-b_1')^2- (b_2-b'_1)^2)\tilde F(b_1') - ((b_1-b_2')^2 - (b_2-b_2')^2) \tilde F(b_2') \\
        &= (b_1 - b_2)(b_1+b_2-2b_1') \tilde F(b_1') - (b_1-b_2) (b_1+b_2-2b_2') \tilde F(b_2') \\
        &= (b_1^2-b_2^2) (\tilde F(b_1') - \tilde F(b_2')) - (b_1-b_2) (2b_1' \tilde F(b_1') - 2b_2' \tilde F(b_2')).
    \end{align*}

    Using this we can now express 
    \begin{equation}\label{eq:Q-scalar-product}
        Q(b_1, b_2, b_1',b_2') = (V(b_1) - V(b_2))\cdot (U(b_1') - U(b_2'))
    \end{equation}
    where $V(b), U(b') \in \RR^3$ are vectors given by
    \begin{equation}\label{eq:vectors-V-U}
    V(b) = \begin{pmatrix}
        b \\
        b^2 \\
        F(b)
    \end{pmatrix},\quad U(b') = \begin{pmatrix}
        F(b') - 2 b' \tilde F(b') \\
        \tilde F(b') \\
        b'
    \end{pmatrix}    
    \end{equation}

    \medskip
\noindent \textit{Step 6. Applying radial projections. }
    Recall that we assume that $\tau \le \varrho^{3-\zeta_1}$ holds. 
    In this range we can redefine the $\Delta$ in (\ref{eq:quadratic-non-concentration}) to be $\Delta = \varrho^{2\varepsilon_1+\zeta_1}$ (taking smaller $\Delta$ makes the condition weaker). So in other words, we may assume that
    \begin{equation}\label{eq:F-non-concentrated-32planar}
    \#\{b \in \mc B:|F(b) - a_0-a_1b-a_2b^2| < \varrho^{2\varepsilon_1+\zeta_1}\} \le \varrho^{\varepsilon_2/2} |\mc B|    
    \end{equation}
    for any choice of coefficients $a_0, a_1, a_2 \in \RR$.

    For $t\ge 1$ let us define $\mc S_t$ to be the set of $s \in w\ZZ$ such that we have $\#\{(b, b')\in \mc B\times \mc B:~ |\mc G|_{b, b'}| \ge \d^{C_0\eta'} |\mc A|, ~~ S(b, b') \in s^{(w)} \} \sim 2^t$. By dyadic pigeonhole and an appropriate choice of $C_0$, we can find $t$ so that $2^t |\mc S_t| \gtrsim_{\d^{O(\eta')}} |\mc B|^2$ holds. Using (\ref{eq:F-non-concentrated-32planar}) and a Cauchy--Schwarz argument similar to the one given in the proof of Claim \ref{claim:S_j-non-concentrated}, we can show that $|\mc S_t| \gtrsim w^{-c_0}$ holds for some constant $c_0>0$ depending on $\varepsilon_0,\varepsilon_1, \varepsilon_2, \mu$. 
    
    So we can apply the asymmetric Balog--Szemer\'edi--Gowers theorem (Lemma \ref{lem:asymmetric-bsg}) to the set $\mc A$ and $\mc S_j$. Thus, for any $\tilde\eta>0$ such that $\eta'$ is sufficiently small in $\tilde \eta$ and $c_0$, we can find subsets $\mc A' \subset \mc A$ and $\mc S' \subset \mc S_j$ such that $|\mc A'| \ge w^{O(\tilde \eta)} |\mc A|$, $|\mc S'| \ge w^{O(\tilde \eta)} |\mc S_j|$  and such that for any $n, m\ge 0$ we have 
    \begin{equation}\label{eq:asym-bsg-32planar}
    |\mc A' + n \cdot \mc S' - m \cdot \mc S'|_w \lesssim_{n,m} w^{-O_{n,m}(\tilde \eta)} |\mc A'|_w.    
    \end{equation}
    Let $\mc H \subset \mc B\times \mc B$ be the set of pairs $(b, b')$ such that $S(b, b') \in \mc S'^{(w)}$. By construction, we have the lower bound $|\mc H| \ge w^{O(\tilde \eta)} |\mc B|^2$. Let $\mc C \subset \mc B^4$  be the set of 4-cycles in $\mc H$, i.e. the set of tuples $(b_1, b_2, b_1', b_2')$ such that $(b_i, b_j') \in \mc H$ for all $i, j \in \{1,2\}$. By Cauchy--Schwarz, we have $|\mc C| \ge w^{O(\tilde \eta)} |\mc B|^4$. Using the definition (\ref{eq:Q-polarization-32-planar}) and (\ref{eq:asym-bsg-32planar}), we conclude that 
    \begin{equation}\label{eq:A-does-not-expand}
    |\mc A' + \{ Q(b_1, b_2, b_1', b_2')~\mid~ (b_1, b_2, b_1', b_2')\in \mc C \}|_w \lesssim w^{-O(\tilde \eta)} |\mc A'|_w.    
    \end{equation}
    Recall that by construction we have $|\mc A|_w \lesssim_{\d^{O(\eta')}} w^{\kappa-1}$.
    Thus, in order to finish the proof it remains to show that the set of $ Q(b_1, b_2, b_1', b_2')$ for $(b_1, b_2,b_1', b_2')$ ranging in a dense subset of $\mc B^4$ defines an essentially 1-dimensional subset in the interval $[-1,1]$. 

    Recall from (\ref{eq:Q-scalar-product}) that $Q(b_1, b_2, b_1', b_2')$ is equal to a scalar product of vectors $V(b_1)-V(b_2)$ and $U(b_1')-U(b_2')$ defined in (\ref{eq:vectors-V-U}). 
    We now set up to apply radial projection theorem in $\mathbb{R}^3$ (Theorem \ref{thm:radial-projections-3d}). 
    Let $\varepsilon_1', \varepsilon_2'$ be defined so that $w^{\varepsilon_1'} = \varrho^{100 \varepsilon_1}$ and $w^{\varepsilon_2'} = \varrho^{\varepsilon_2/100}$. Recall that by definition we have $w = C \d^{-\eta'} (\varrho + \varrho^{3+\varepsilon_0}/\tau)$ and we assume that $\tau \in [ \varrho^{3+2\varepsilon_1}, \varrho^{3-\zeta_1}]$. So it follows that $\varrho^{\varepsilon_0+\zeta_1} \lesssim_{\d^{O(\eta')}} w \lesssim_{\d^{O(\eta')}} \varrho^{\varepsilon_0-2\varepsilon_1} $.

    Consider the following set:
    \[
    \mc P = \{ V(b_1) - V(b_2)~\mid ~ b_1,b_2 \in \mc B, ~~ |b_1-b_2| \gtrsim \d^{C\tilde \eta} \} \subset \RR^3.
    \]
    (if there are overlaps, we treat it as multiset). Since $|\mc B|_w \gtrsim_{\d^{O(\eta')}} w^{-1}$ and $V(b) = \begin{pmatrix}
        b\\
        b^2\\
        F(b)
    \end{pmatrix}$, it follows from the first two coordinates that $\mc P$ is a $(w, 2, \d^{-C\eta'})$-set. 
    Let $S$ be a $w^{\varepsilon_1'}$-slab, i.e. we can write
    \[
    S = \{ v \in \RR^3: \langle n, v\rangle \in [t, t+w^{\varepsilon_1'}] \}
    \]
    for some unit vector $n = (n_1, n_2,n_3)$ and $t\in \RR$. Suppose that we have $\# \{(b_1, b_2) \in \mc B\times \mc B:~ \langle n, V(b_1) - V(b_2)\rangle \in [t, t+w^{\varepsilon_1'}]\} \ge w^{\varepsilon_2'} |\mc B|^2$. By averaging and fixing a value of $b_2$, we get that for some $t'$ we have $\#\{ b\in \mc B:~ \langle n, V(b) \rangle \in [t', t'+w^{\varepsilon_1'}]\} \ge w^{\varepsilon_2'} |\mc B|$. If $|n_3| \le w^{\varepsilon_1'/2}$ then $|(n_1, n_2)|\gtrsim 1$ and we can estimate
    \[
    \#\{ b\in \mc B:~ \langle n, V(b) \rangle \in [t', t'+w^{\varepsilon_1'}]\} \le \#\{b \in \mc B: ~ |n_1 b + n_2 b^2 - t' |\le \d^{-C\eta'} w^{\varepsilon_1'/2}\} \lesssim \d^{-C\eta'} w^{\varepsilon_1'/4} |\mc B|
    \]
    where we recall that $\mc B \subset [-1,1]$ is a set satisfying $|\mc B|_w \gtrsim_{\d^{O(\eta')}} w^{-1}$ and we may assume $\mc B$ to be uniform on scale $w$. So since $\varepsilon_1' > \varepsilon_2'/100$ (recall $\varepsilon_2$ is sufficiently small depending on $\varepsilon_1$ and $\varepsilon_i'$ is proportional to $\varepsilon_i$), we obtain a contradiction.

    So we may assume that $|n_3| \ge w^{\varepsilon_1'/2}$ holds. We then conclude that for $a_1 = -n_1/n_3$, $a_2 = -n_2/n_3$ and $a_0 = t' /n_3$, we get $\#\{b \in \mc B:~ |F(b) - a_0 - a_1 b - a_2 b^2| \le w^{\varepsilon_1'/2} \} \ge w^{\varepsilon_2'} |\mc B|$. Comparing with (\ref{eq:F-non-concentrated-32planar}) and recalling the choice of $\varepsilon_1',\varepsilon_2'$ we obtain a contradiction. We conclude that the set $\mc P$ satisfies the assumption of Theorem \ref{thm:radial-projections-3d}. We conclude that provided that $\d^{-O(\eta')} \ll w^{-\eta'_0}$ for some $\eta'_0 = \eta'_0(\varepsilon_1', \varepsilon_2', \alpha)$, for any $\tilde \eta' \le \eta_0'$ and $\varepsilon_1' \le  \alpha/8$ we can find a subset $\mc D \subset \mc P\times \mc P$ with the following properties. We have $|\mc D| \ge (1-w^{c \tilde \eta'}) |\mc P|^2$ and for every $(p,p') \in \mc D$ we have the thin tube condition $|T_{p,p'}(\varrho) \cap \mc P| \le w^{-\alpha} \varrho^2 |\mc P|$ for all $\varrho \in [w,1]$. Using that $\mc P$ is a $(w,2, \d^{-O(\eta')})$-set, we can slightly modify the graph $\mc D$ to ensure that $|p-p'| \ge w^{3 \tilde \eta'}$ holds for every $(p, p') \in\mc D$.

    \medskip
\noindent \textit{Step 7. Applying Kaufman projection. }
    Note that each pair $(p, p') \in \mc D$ has the form $p = V(b_1)-V(b_2)$, $p' = V(b_3) - V(b_4)$ and corresponds to the slope $\frac{V(b_1) - V(b_2) - V(b_3) + V(b_4)}{b_1-b_2-b_3+b_4}$.
    We now apply Lemma \ref{lem:kaufman-Rd} to the set of slopes determined by $\mc D$ and the graph $\mc U = \{U(b'), ~b' \in \mc B\} \subset \RR^3$. Note that the set of slopes forms a $(w, 2, w^{-\alpha-O(\tilde\eta')})$-set and the graph $\mc U$ forms a $(w, 1, \d^{-O(\eta')})$-set (since $\mc B$ is dense and the third coordinate of $U(b')$ equals to $b'$).
    We conclude that for any $\beta, \overline{\eta}>0$ such that $\eta',\tilde \eta', \alpha \ll \overline{\eta} \le \eta_0(3, \beta)$ there exists a subset $\mc J \subset \mc B\times \mc B\times \mc B\times \mc B$ and for each $(b_1, b_2, b_3, b_4) \in J$ there exists a subset $\mc B_{b_1, b_2, b_3, b_4} \subset \mc B$ such that 
    \begin{equation}\label{eq:JB}
    |\mc J|\ge (1-w^{\overline{\eta}}) |\mc B|^4, \quad |\mc B_{b_1, b_2, b_3, b_4}| \ge (1-w^{\overline{\eta}}) |\mc B|    
    \end{equation}
    so that 
    \begin{equation}\label{eq:Vbi-1set}
    \{  (V(b_1) - V(b_2) - V(b_3) +  V(b_4)) \cdot U(b')~\mid ~ b' \in \mc B_{b_1, b_2, b_3, b_4}  \}    
    \end{equation}
    is a $(w, 1, w^{-\beta})$-set. On the other hand, by (\ref{eq:A-does-not-expand}) and Plunnecke-Ruzsa inequality (e.g. \cite[Chapter 6, Corollary 6.29]{tao2006additive}), we have 
    \[
    |\{ Q(b_1,b_2, b_1', b_2') - Q(b_3, b_4, b_3', b_4')~\mid ~ (b_1, b_2, b_1' ,b_2'), (b_3, b_4, b_3', b_4') \in \mc C \}|_w \lesssim w^{O(\tilde \eta)} |\mc A'|_w \lesssim \d^{O(\eta')} w^{O(\tilde \eta)} w^{\kappa-1}.
    \]
    Let $\mc W$ be the set of tuples $(b_1, b_2, b_3, b_4, b_1', b_2') \in \mc B^6$ so that $(b_1, b_2, b_1', b_2'), (b_3, b_4, b_1', b_2') \in \mc C$, in other words $(b_i, b_j')\in \mathcal{H}$ for all $i\in \{1, 2, 3, 4\}, j\in \{1, 2\}$. By Cauchy--Schwarz, we have $|\mc W|  \gtrsim w^{O(\tilde \eta)} |\mc B|^6 $. 
    By restricting the above set of differences to the set of octuples satisfying $(b_1', b_2') = (b_3', b_4')$ and using the representation of $Q$ as scalar product of vectors $V$ and $U$ we then conclude that
    \[
    |\{ (V(b_1) - V(b_2) - V(b_3) + V(b_4))\cdot (U(b_1') - U(b_2'))~\mid ~ (b_1, b_2, b_3, b_4, b_1', b_2') \in \mc W  \}|_w \lesssim  \d^{O(\eta')} w^{O(\tilde \eta)} w^{\kappa-1}.
    \]
    Let us choose parameters so that $\tilde \eta \ll \overline{\eta} $. Then by (\ref{eq:JB}) we get that for at least a half of tuples $(b_1, b_2, b_3, b_4, b_1',b_2') \in \mc W$ we also have that $(b_1,  b_2, b_3, b_4) \in \mc J$ and $b_1', b_2' \in \mc B_{b_1, b_2, b_3, b_4}$. Let $\mc W' \subset \mc W$ be the set of such tuples. By pigeonhole principle, there exists $(b_1, b_2, b_3, b_4) \in \mc B^4$ and a $w$-interval $I$ such that 
    \[
    \#\{ (b_1', b_2') \in \mc W'|_{b_1, b_2, b_3, b_4}:~   (V(b_1) - V(b_2) - V(b_3) + V(b_4))\cdot (U(b_1') - U(b_2')) \in I \} \gtrsim \d^{O(\eta')} w^{O(\tilde \eta)} w^{1-\kappa} |\mc B|^2
    \]
    and so 
    \[
    \#\{(b_1', b_2') \in \mc B^2_{b_1, b_2, b_3, b_4}:~  (V(b_1) - V(b_2) - V(b_3) + V(b_4))\cdot (U(b_1') - U(b_2')) \in I \} \gtrsim \d^{O(\eta')} w^{O(\tilde \eta)}  w^{1-\kappa} |\mc B|^2.
    \]
    On the other hand, since (\ref{eq:Vbi-1set}) is $(w, 1, w^{-\beta})$-set, we conclude that 
    \[
    w^{-\kappa} \lesssim \d^{O(\eta')} w^{O(\tilde \eta)} w^{-\beta}.
    \]
    By taking $\beta < \zeta/2$ and other parameters accordingly, we conclude that $\kappa \le\zeta$. This concludes the proof of Lemma \ref{lem:quadratic-case}.  
\end{proof}

\bibliographystyle{amsplain0.bst}
\bibliography{main}

\appendix

\section{Radial projection theorems}\label{appendix:radial}

In this appendix we prove the radial projection theorems used in our proofs. The arguments here broadly follow the proof steps in the Orponen--Shmerkin--Wang paper \cite{orponen2024kaufman} though details are quite different. 

\subsection{Radial projections in 2D}
We repeat the statement of Theorem \ref{thm:radial-projections} here for the readers convenience.

\begin{theorem}\label{thm:radial-again}
    Let $t \in (0,2]$ and $\varepsilon_1, \varepsilon_2, \zeta \in (0,1)$. Suppose that $\varepsilon_1 \le \zeta/8$. Then there exist $\delta_0 = \d_0(t, \varepsilon_1, \varepsilon_2, \zeta)>0$ and $\eta_0 = \eta_0(t, \varepsilon_1, \zeta)>0$ such that the following holds for all $\delta < \delta_0$ and $\eta < \eta_0$.
    Let $P \subset [0,1]^2$ be a $\delta$-separated $(\d, \d^{-\eta})$-uniform  $(\delta, t, \delta^{-\eta})$-set such that $|P \cap T| \le \d^{\varepsilon_2} |P|$ for any $\d^{\varepsilon_1}\times 1$-tube $T$.

    Then there exists a subset $G \subset P\times P$ such that $|G| \ge (1-\delta^{\min(\eta, \varepsilon_2/2)}) |P|^2$ such that for every $(p,p') \in G$ we have $|T_{p,p'}(\varrho) \cap P| \le \delta^{-\zeta} \varrho^{\min(1, t)} |P|$ for all $\varrho\in [\delta, 1]$.

    In particular, for every $p\in P$ such that $|G|_p| \ge \nu |P|$,  the set of tubes $\{T_{p, p'}(\delta), ~(p, p')\in G\}$ is a $(\delta, \min(1, t), C\delta^{-\zeta} \nu^{-1})$-set.
\end{theorem}

The proof of the following lemma is a mild modification of arguments in
\cite[Section 2]{orponen2018dimension} and
\cite[Appendix B]{shmerkin2023non}. 

\begin{lemma}\label{lem:initial-thin-tubes}
    Let $t\in (0,2]$ and $\varepsilon_1,\varepsilon_2>0$. There is an
    absolute constant $c>0$ such that the following holds whenever
    $0<\chi<ct\varepsilon_1$ and
    $\d<\d_0(t,\varepsilon_1,\varepsilon_2)$.

    Let $P \subset [-1,1]^2$ be a $(\d, t, \d^{-\chi})$-set such that $|P\cap T| \le \d^{\varepsilon_2}|P|$ for every $\d^{\varepsilon_1} \times 1$-tube $T$. 
    Then there exists $G \subset P\times P$ such that $|G| \ge (1- \d^{\min(\varepsilon_1 t/2, \varepsilon_2/2) }) |P|^2$ and $|T_{p, p'}(\d^{4\varepsilon_1}) \cap P| \le \d^{\chi} |P|$ for all $(p, p') \in G$.
\end{lemma}

\begin{proof}
    Put
    \[
        \Delta_0:=\d^{4\varepsilon_1},\qquad
        \varrho_0:=\d^{2\varepsilon_1},\qquad
        \theta_0:=\d^{\varepsilon_1}.
    \]
    Thus $\Delta_0/\varrho_0=\varrho_0$ and
    $\varrho_0/\theta_0=\theta_0$.

    Let $\TT_0$ be the collection of all $\Delta_0$-tubes $T$ with
    $|P\cap T|\geq \d^\chi|P|$. Choose a maximal subcollection
    $\TT\subset\TT_0$ such that
    \[
        \operatorname{diam}(T_1\cap T_2)\leq \varrho_0,
        \qquad T_1,T_2\in\TT,\quad T_1\neq T_2.
    \]
    Since $T_1\cap T_2$ is contained in a ball of radius
    $O(\varrho_0)$, the Frostman hypothesis gives
    \begin{equation}\label{eq:initial-tube-overlap}
        |P\cap T_1\cap T_2|
        \lesssim \d^{-\chi}\varrho_0^t|P|,
        \qquad T_1\neq T_2.
    \end{equation}
    
    We first claim that
    \begin{equation}\label{eq:initial-tube-number}
        |\TT|\lesssim \d^{-\chi}.
    \end{equation}

    Indeed, set $I:=\sum_{T\in\TT}|P\cap T|$. Then
    $I\geq \d^\chi|P||\TT|$.
    By Cauchy--Schwarz and \eqref{eq:initial-tube-overlap},
    \[
        \frac{I^2}{|P|}
        \leq \sum_{T_1,T_2\in\TT}|P\cap T_1\cap T_2|
        \lesssim
        I+\d^{-\chi}\varrho_0^t|P||\TT|^2.
    \]
    If the first term on the right dominates, then $I\lesssim |P|$,
    which gives \eqref{eq:initial-tube-number}. If the second term
    dominates, then
    \[
        \d^\chi|P||\TT|
        \lesssim I
        \lesssim
        \d^{-\chi/2}\varrho_0^{t/2}|P||\TT|.
    \]
    This is impossible when $\chi\ll t\varepsilon_1$, since
    $\varrho_0^{t/2}=\d^{t\varepsilon_1}$. This proves the claim.

    By maximality, every $T\in\TT_0$ has
    $\operatorname{diam}(T\cap T')>\varrho_0$ for some $T'\in\TT$. 
    Note that $\operatorname{diam}(T\cap T')>\varrho_0$ implies that $\angle T, T' = O(\Delta_0/\varrho_0)=O(\varrho_0)$ and that  $T$ is
    contained in a $C\varrho_0$-tube coaxial with $T'$. Write
    $\TT_{\varrho_0}$ for the family of
    $C\varrho_0$-tubes corresponding to all tubes in $\TT_0$.

    Let $P_{\mathrm{tr}}\subset P$ be the subset of points which are covered by at least two members of $\TT_{\varrho_0}$ making angle at least $\theta_0$. 
    The intersection of two such tubes is contained in a
    ball of radius $O(\varrho_0/\theta_0)=O(\theta_0)$. Therefore,
    by \eqref{eq:initial-tube-number},
    \begin{equation}\label{eq:initial-transverse-set}
    \begin{split}
        |P_{\mathrm{tr}}|
        &\lesssim
        \sum_{\substack{T_1,T_2\in\TT_{\varrho_0}\\
                        \angle(T_1,T_2)\geq\theta_0}}
        |P\cap T_1\cap T_2|  \\
        &\lesssim
        \d^{-\chi}\theta_0^t|P|\,|\TT|^2
        \lesssim \d^{-3\chi+t\varepsilon_1}|P|
        \leq \d^{2t\varepsilon_1/3}|P|.
    \end{split}
    \end{equation}

    Fix $p\in P\setminus P_{\mathrm{tr}}$. All members of
    $\TT_{\varrho_0}$ through $p$ make pairwise angle
    $O(\theta_0)$, so they are contained in a single
    $C\theta_0$-tube. Denote this tube by $T(p)$. In particular, every $\d^\chi|P|$-rich
    $\Delta_0$-tube through $p$ is contained in $T(p)$. By the two-ends assumption on $P$, we have
    \[
        |P\cap T(p)|\lesssim \d^{\varepsilon_2}|P|.
    \]
    Finally, define
    \[
        G:=\{(p,p')\in P\times P:
             p\notin P_{\mathrm{tr}}\ \text{and}\ p'\notin T(p)\}.
    \]
    If $(p,p')\in G$ and $T_{p,p'}(\Delta_0)$ were $\d^\chi |P|$-rich, then it would
    be contained in $T(p)$, forcing $p'\in T(p)$, a contradiction.
    Hence
    \[
        |T_{p,p'}(\Delta_0)\cap P|<\d^\chi|P|,
        \qquad (p,p')\in G.
    \]
    Moreover,
    \[
        |G|
        \geq
        (|P|-|P_{\mathrm{tr}}|)
        (1-C\d^{\varepsilon_2})|P|
        \geq
        \bigl(1-\d^{\min(t\varepsilon_1/2,\varepsilon_2/2)}\bigr)
        |P|^2
    \]
    for sufficiently small $\d$.
\end{proof}

In the proof we will need a $\eps$-improvement Furstenberg estimate with a minimal spacing condition due to Orponen and Shmerkin \cite[Theorem 5.61]{orponen2023projections}. The following is a dual version, as stated in \cite[Theorem 2.11]{orponen2025large}.

\begin{theorem}\label{thm:two-ends}
    Fix $s  \in (0,1]$, $t \in(0,2]$ and $u \in (0, \min\{t, 2-t\}]$. There exists $\eta = \eta(s, t, u) >0$ such that the following holds for small enough $\d>0$. 

    Let $\TT$ be a family of essentially distinct  $\d$-tubes with $|\TT| = \d^{-t}$ and satisfying the following non-concentration condition on scale $w = \d |\TT|^{1/2}$: $|\TT[T_w]| \le \d^u |\TT|$ for any $w$-tube $T_w$. 

    Let $Y(T) \subset T$ be a $(\d, s, \d^{-\eta})$-set such that $|Y(T)|_\d \sim N$ for all $T\in \TT$. Then 
    \[
    |U(\TT, Y)|_\d \gtrsim \d^{-su/8}|Y(T)|_\d |\TT|^{1/2}.
    \]
\end{theorem}

\begin{proof}[Proof of Theorem \ref{thm:radial-again}]
    Write
    \[
        s:=\min\{t,1\},\qquad
        \gamma:=\min\{\zeta/100,s/2\},\qquad
        \sigma:=s-\gamma.
    \]
    In particular, $\sigma>0$ and
    \[
        a:=t-\sigma\geq\gamma.
    \]
    Fix $\chi=c_0t\varepsilon_1$, where $c_0>0$ is sufficiently
    small for Lemma \ref{lem:initial-thin-tubes}. We choose the
    parameters in the order
    \[
        t,\varepsilon_1,\zeta,\chi
        \quad\text{first, and then}\quad
        0<\eta\ll\eta'\ll\chi,\gamma,\zeta.
    \]
    More precisely, below we take $\eta'=C_1\eta$ with $C_1$ large
    in terms of $t$ and $\zeta$, and then take $\eta$ sufficiently
    small.

    Since $\eta\leq\chi$, the set $P$ is also a
    $(\d,t,\d^{-\chi})$-set. Lemma
    \ref{lem:initial-thin-tubes} gives a graph
    $G_0\subset P\times P$ such that
    \begin{equation}\label{eq:G0-size-2d}
        |G_0|
        \geq
        \bigl(1-\d^{\min(t\varepsilon_1/2,\varepsilon_2/2)}\bigr)
        |P|^2
    \end{equation}
    and, for every $(p,p')\in G_0$,
    \begin{equation}\label{eq:G0-thin-2d}
        |C_0T_{p,p'}(\d^{4\varepsilon_1})\cap P|
        \leq \d^\chi|P|.
    \end{equation}
    Here $C_0$ is any fixed sufficiently large constant (we can get this version with the dilated tube by a slight change of scales when we apply the lemma). 
    
    Put $r:=\d^{2\eta/t}$ and let $G_{\geq r}$ be the set of pairs
    $(p,p')$ with $|p-p'|\geq r$. For each $p\in P$, the Frostman
    bound gives
    \[
        |P\cap B(p,r)|
        \leq \d^{-\eta}r^t|P|
        =\d^\eta|P|.
    \]
    Hence
    \begin{equation}\label{eq:close-pairs-2d}
        |(P\times P)\setminus G_{\geq r}|
        \lesssim \d^\eta|P|^2.
    \end{equation}
    Set $G'_0:=G_0\cap G_{\geq r}$.

    Let $G_1\subset G'_0$ be the set of pairs $(p,p')$ for which
    there is a scale $\varrho\in[\d,1]$ such that
    \[
        |T_{p,p'}(\varrho)\cap P|
        \geq \d^{-\zeta}\varrho^s|P|.
    \]
    It remains to prove that
    \begin{equation}\label{eq:G1-goal-2d}
        |G_1|\lesssim \d^\eta|P|^2.
    \end{equation}

    \medskip
    \noindent\emph{Step 1: selecting the densest scale $\varrho$.}
    For $(p,p')\in G_1$ and a dyadic scale
    $q\in[\d,1]$, define
    \[
        W_{p,p'}(q)
        :=q^{-\sigma}|T_{p,p'}(q)\cap P|\,|P|^{-1}.
    \]
    Choose a dyadic scale $\varrho(p,p')$ at which $W_{p,p'}$ is
    maximal. If $q$ is a dyadic scale comparable to a scale witnessing
    $(p,p')\in G_1$, then
    \[
        W_{p,p'}(q)
        \gtrsim \d^{-\zeta}q^{s-\sigma}
        =\d^{-\zeta}q^\gamma
        \geq \d^{-(\zeta-\gamma)}.
    \]
    Since $W_{p,p'}(\varrho)\leq \varrho^{-\sigma}$, maximality gives
    \begin{equation}\label{eq:rho-upper-bound}
        \varrho(p,p')
        \lesssim \d^{(\zeta-\gamma)/\sigma}.
    \end{equation}

    Dyadically pigeonhole both $\varrho(p,p')$ and
    $|T_{p,p'}(\varrho(p,p'))\cap P|$. We obtain values $\varrho,m$ and a
    graph $G_2\subset G_1$ such that
    \begin{equation}\label{eq:G2-density-2d}
        |G_2|\gtrsim (\log(1/\d))^{-2}|G_1|,
        \qquad
        |T_{p,p'}(\varrho)\cap P|\sim m
    \end{equation}
    for $(p,p')\in G_2$. Moreover,
    \begin{equation}\label{eq:m-estimate}
        m
        \gtrsim
        \d^{-(\zeta-\gamma)}\varrho^\sigma|P|
        \geq
        \d^{-(\zeta-\gamma)}\varrho^s|P|.
    \end{equation}

    The point of maximizing $W_{p,p'}$ is the following
    non-concentration estimate at all larger widths. For
    $\tau\in[\varrho,1]$, maximality at the scale $\varrho/\tau$ gives
    \begin{equation}\label{eq:non-concentration-in-tau}
        |T_{p,p'}(\varrho/\tau)\cap P|
        \lesssim \tau^{-\sigma}m,
        \qquad (p,p')\in G_2.
    \end{equation}

    \medskip
    \noindent\emph{Step 2: annulus estimate.}
    Set $\eta':=C_1\eta$, with $C_1$ sufficiently large in terms
    of $t$ and $\zeta$. For a dyadic $\tau\in[\varrho,1]$, define a subgraph $G_2(\tau) \subset G_2$ as follows: for each $p \in P$ let
    $G_2(\tau)|_p$ consist of the pairs $(p,p')\in G_2$ for which
    \begin{equation}\label{eq:bad-annulus-p}
    \begin{split}
        &|T_{p,p'}(r^{-3}\varrho)
          \cap(B(p,2\tau)\setminus B(p,\tau))\cap P|\\
        &\hspace{45mm}\geq
          \d^{-\eta'}\tau^{t-\sigma}m.
    \end{split}
    \end{equation}
    We estimate the size of $G_2(\tau)$. Fix $p$. The annulus
    $B(p,2\tau)\setminus B(p,\tau)$ can be covered by
    $O(r^3\tau/\varrho)$ finitely overlapping
    $2r^{-3}\varrho\times2\tau$ rectangles $R_i$. By the
    $(\d,t,\d^{-\eta})$-set property, the number of indices satisfying
    \[
        |R_i\cap P|\geq
        \d^{-\eta'}\tau^{t-\sigma}m
    \]
    is at most
    \begin{equation}\label{eq:heavy-rectangles-2d}
        \frac{|P\cap B(p,2\tau)|}
             {\d^{-\eta'}\tau^{t-\sigma}m}
        \lesssim
        \d^{\eta'-\eta}\tau^\sigma\frac{|P|}{m}.
    \end{equation}

    For one such rectangle $R_i$, let $S_i$ be the set of points
    $p'\in G_2(\tau)|_p$ whose tube segment in
    \eqref{eq:bad-annulus-p} is contained in $R_i$. All points of
    $S_i$ lie in a $Cr^{-3}\varrho/\tau$-tube $T_i$. If $S_i$ is
    non-empty, choose $p'\in S_i$. Since $|p-p'|\geq r$, elementary
    geometry gives
    \[
        T_i\subset
        T_{p,p'}(Cr^{-4}\varrho/\tau).
    \]
    Applying \eqref{eq:non-concentration-in-tau} (and using the
    trivial bound $|S_i|\leq |P|$ if the displayed width exceeds
    $1$) yields
    \[
        |S_i|
        \lesssim r^{-4\sigma}\tau^{-\sigma}m.
    \]
    Combining this with \eqref{eq:heavy-rectangles-2d}, we obtain
    \begin{equation}\label{eq:G2tau}
        \bigl|G_2(\tau)|_{p}\bigr|
        \lesssim
        r^{-4\sigma}\d^{\eta'-\eta}|P|.
    \end{equation}
    Define $G_2'(\tau)\subset G_2$ in the same way, with the roles of the
    two endpoints $p,p'$ reversed. The same estimate holds for every
    $p'$-fibre. Consequently,
    \begin{equation}\label{eq:G2tau-total}
        |G_2(\tau)|+|G_2'(\tau)|
        \lesssim
        r^{-4\sigma}\d^{\eta'-\eta}|P|^2.
    \end{equation}
    Define
    \[
        \tau_0:=\d^{2\eta'/a},
        \qquad a=t-\sigma>0.
    \]
    By \eqref{eq:rho-upper-bound}, we have $\varrho\leq\tau_0$ once
    $\eta$ is sufficiently small. Define a new subgraph by removing the subgraphs $G_2(\tau), G_2'(\tau)$ over all dyadic scales:
    \begin{equation}\label{eq:G3-G2}
        G_3:=
        G_2\setminus
        \bigcup_{\substack{\tau\in[\varrho,\tau_0]\\ \tau\ {\rm dyadic}}}
        \bigl(G_2(\tau)\cup G_2'(\tau)\bigr).
    \end{equation}
    If $(p,p')\in G_3$, summing
    \eqref{eq:bad-annulus-p} over dyadic annuli gives
    \begin{equation}\label{eq:non-clustering-radius}
    \begin{split}
        &|P\cap T_{p,p'}(r^{-3}\varrho)
          \cap(B(p,\tau_0)\cup B(p',\tau_0))|\\
        &\hspace{30mm}\lesssim
        \d^{-\eta'}\tau_0^{a}m
        +\d^{-\eta}\varrho^t|P|
        \lesssim \d^{\eta'}m.
    \end{split}
    \end{equation}
    In the last step we used $\tau_0^a=\d^{2\eta'}$ and
    \eqref{eq:m-estimate}; the contribution of the innermost
    $\varrho$-balls is smaller after taking $\eta$ sufficiently small.
    Since $|T_{p,p'}(\varrho)\cap P|\sim m$, it follows that
    \begin{equation}\label{eq:mass-away-endpoints-2d}
        |P\cap T_{p,p'}(\varrho)
          \setminus(B(p,\tau_0)\cup B(p',\tau_0))|
        \gtrsim m.
    \end{equation}
    Write $\pi_p(q):=\frac{q-p}{|q-p|}$ for the radial projection centered at $p$.
    We claim that, for every $p\in P$,
    \begin{equation}\label{eq:radial-cover-G3}
        |\pi_p(G_3|_p)|_\varrho
        \lesssim \tau_0^{-1}\frac{|P|}{m}.
    \end{equation}
    To see this, choose one pair $(p,p')$ for every element of a
    maximal $\varrho$-separated subset of $\pi_p(G_3|_p)$. By
    \eqref{eq:mass-away-endpoints-2d}, the corresponding $\varrho$-tubes
    contain $\gtrsim m$ points of $P$ at distance at least $\tau_0$
    from $p$. At such distances, a point belongs to at most
    $O(\tau_0^{-1})$ of these angularly separated tubes. Summing their
    incidences with $P$ proves \eqref{eq:radial-cover-G3}.

    \medskip
    \noindent\emph{Step 3: prune the set of pairs and define a shaded collection of tubes.}
    Partition $G_3$ into classes according to the first endpoint $p$
    and the dyadic $\varrho$-interval containing $\pi_p(p')$. Also
    consider the transposed partition, according to $p'$ and the
    interval containing $\pi_{p'}(p)$. By
    \eqref{eq:radial-cover-G3}, the number of classes in either
    partition is at most
    \begin{equation}\label{eq:number-direction-classes}
        C\tau_0^{-1}\frac{|P|^2}{m}.
    \end{equation}
    Starting with $G_3$, repeatedly delete every class in either
    partition which currently contains fewer than
    \[
        2k,\qquad k:=\tau_0^2m,
    \]
    pairs. A fixed class is deleted at most once. Thus
    \eqref{eq:number-direction-classes} shows that the total number of
    deleted pairs is $O(\tau_0|P|^2)$. Denote the remaining graph by
    $G_4$. Every pair of $G_4$ now has at least $k$ neighbours in
    each of its two direction classes, with all of these neighbours
    still belonging to $G_4$. In particular,
    \begin{equation}\label{eq:G4-core-loss}
        |G_3\setminus G_4|\lesssim \tau_0|P|^2.
    \end{equation}

    Let $\TT_\varrho$ be a minimal family of essentially distinct
    $C\varrho$-tubes which covers $T_{p,p'}(\varrho)$ for
    $(p,p')\in G_4$. Let $\bar\d$ be a dyadic scale comparable to
    $C'r^{-1}\varrho$, with $C'$ sufficiently large, and thicken every
    member of $\TT_\varrho$ to a $\bar\d$-tube. Denote the resulting family by $\TT$. Note that $\TT$ is not necessarily essentially distinct, but this property can only fail up to multiplicity factor of $O(r^{-2})$.

    For $T\in\TT$, let $Y(T)$ consist of those $p\in P$ for which
    $T_{p,p'}(\varrho)$ is covered by $T$ for some $(p,p')\in G_4$.
    Fix one pair $(p,p')\in G_4$ whose $\varrho$-tube is covered by $T$. By the construction of $G_4$, we have at least $k=\tau_0^2m$ points $p''$ such that
    $(p'',p')\in G_4$ and
    $\pi_{p'}(p'')$ lies in the same $\varrho$-interval as
    $\pi_{p'}(p)$. Since both $|p-p'|$ and $|p''-p'|$ are at least
    $r$, all tubes $T_{p'',p'}(\varrho)$ lie in the
    $Cr^{-1}\varrho$-thickening of $T_{p,p'}(\varrho)$. Hence
    \begin{equation}\label{eq:shading-size-2d}
        |Y(T)|\gtrsim k=\tau_0^2m.
    \end{equation}

    The annular estimates used to define $G_3$ also give, for
    $p\in Y(T)$ and $\tau\in[\bar\d,\tau_0]$,
    \begin{equation}\label{eq:properties-of-shading}
        |Y(T)\cap B(p,\tau)|
        \lesssim
        r^{-O(1)}\d^{-\eta'}\tau^{a}m.
    \end{equation}
    Indeed, choose a pair in $G_4$ witnessing $p\in Y(T)$ and sum
    the corresponding annular estimates up to radius $\tau$; the
    thickening from $\varrho$ to $\bar\d$ only costs a power of $r^{-1}$.
    For $\tau\geq\tau_0$ we instead use the trivial bound
    $|Y(T)\cap B(p,\tau)|\leq |Y(T)|$; the factor $\tau_0^{-2}$
    below makes this sufficient. 
    Now 
    \eqref{eq:shading-size-2d}--\eqref{eq:properties-of-shading}
    imply that $Y(T)$ is a $(\bar\d,\bar s,r^{-O(1)}\tau_0^{-2}\d^{-\eta'})$-set with $\bar s:=\min\{a,1\}\geq\gamma$. 

    By \eqref{eq:non-concentration-in-tau}, every tube in $\TT$ covers at most $r^{-O(1)}m^2$ pairs in $G_4$.
    Consequently,
    \begin{equation}\label{eq:T-lower-bound}
        |\TT|\gtrsim r^{O(1)}\frac{|G_4|}{m^2}.
    \end{equation}
    We now suppose, for a contradiction, that
    \begin{equation}\label{eq:G4-large-assumption}
        |G_4|\geq \d^\eta|P|^2.
    \end{equation}
    Then \eqref{eq:T-lower-bound} gives
    \[
        |\TT|
        \gtrsim r^{O(1)}\d^\eta\frac{|P|^2}{m^2}.
    \]
    We also record the complementary upper bound for $|\TT|$. Fix
    $T\in\TT$ and a representative $(p,p')\in G_4$ for it. The construction of $G_4$ implies that there are $k$ pairs $(q,p') \in G_4$ whose tubes are
    covered by an $O(r^{-1})$ dilation of $T$. For every such $q$, the construction of $G_4$ again provides $k$ pairs $(q,q')$ covered by a $O(r^{-2})$-dilation of $T$. Thus an $Cr^{-2}$ dilation of $T$ contains the tubes
    associated with at least $k^2$ distinct pairs of $G_4$. Summing over $T$ gives
    \begin{equation}\label{eq:T-upper-energy-2d}
        \tau_0^4m^2|\TT|
        =k^2|\TT|
        \lesssim r^{-O(1)} |G_4|
        \leq r^{-O(1)} |P|^2.
    \end{equation}
    From \eqref{eq:T-upper-energy-2d} and
    \eqref{eq:m-estimate},
    \begin{equation}\label{eq:T-upper-bound-2d}
        |\TT|
        \lesssim r^{-O(1)}
        \tau_0^{-4}
        \d^{2(\zeta-\gamma)}\varrho^{-2\sigma}
        \leq r^{-O(1)}
        \tau_0^{-4}
        \d^{2(\zeta-\gamma)}\varrho^{-2}.
    \end{equation}

    \medskip
    \noindent\emph{Step 4: apply Furstenberg estimate.}
    As in Theorem \ref{thm:two-ends}, set
    \[
        w:=\bar\d|\TT|^{1/2}.
    \]
    Since $\bar\d\sim r^{-1}\varrho$, the preceding estimate gives
    \begin{equation}\label{eq:w-small-2d}
        w
        \lesssim
        r^{-O(1)}\tau_0^{-2}
        \d^{\zeta-\gamma}\varrho^{1-\sigma}
        \leq \d^{3\zeta/4},
    \end{equation}
    after taking $\eta$ sufficiently small. Since
    $4\varepsilon_1\leq\zeta/2$, this implies
    $w\leq\d^{5\varepsilon_1}$.

    We verify the single-scale spacing hypothesis in Theorem
    \ref{thm:two-ends}. Let $T_w$ be any $Cw$-tube. If
    $T\in\TT$ is contained in $T_w$, choose a pair
    $(p,p')\in G_4\subset G_0$ whose $\varrho$-tube is covered by $T$.
    Then 
    \[
        C r^{-2} \cdot T_w\subset C_0T_{p,p'}(\d^{4\varepsilon_1}),
    \]
    since $w\le \d^{5\varepsilon_1}$ and $\eta$ is small enough in $\varepsilon_1$. Therefore
    \eqref{eq:G0-thin-2d} gives
    $|P\cap (C r^{-2} \cdot T_w)|\lesssim\d^\chi|P|$. Repeating the pair-counting argument using the definition of $G_4$ argument from \eqref{eq:T-upper-energy-2d}, now
    restricted to the members of $\TT$ contained in $T_w$, yields
    \[
        \tau_0^4m^2|\TT[T_w]|
        \lesssim
        |P\cap (Cr^{-2}\cdot T_w)|^2
        \lesssim
        r^{-O(1)}\d^{2\chi}|P|^2.
    \]
    Comparing this with \eqref{eq:T-lower-bound} and
    \eqref{eq:G4-large-assumption}, and then taking $\eta$ small
    relative to $\chi$, gives
    \begin{equation}\label{eq:minimal-spacing-2d}
        |\TT[T_w]|\leq \d^\chi|\TT|.
    \end{equation} 
    We next check that the parameters of Theorem
    \ref{thm:two-ends} stay in a fixed compact range. Write
    $|\TT|=\bar\d^{-\bar t}$. Since $\bar\d\leq w\leq
    \d^{5\varepsilon_1}$, \eqref{eq:G0-thin-2d} gives $m\lesssim
    \d^\chi|P|$. The lower bound \eqref{eq:T-lower-bound} and
    \eqref{eq:G4-large-assumption} therefore imply
    \[
        |\TT|\gtrsim r^{O(1)}\d^{\eta-2\chi}
        \geq \bar\d^{-\chi/2}.
    \]
    On the other hand, \eqref{eq:T-upper-bound-2d}, the relation
    $\varrho\sim r\bar\d$, and $\sigma\leq1-\gamma$ give
    \[
        |\TT|\leq \bar\d^{-2+\gamma},
    \]
    after taking $\eta$ sufficiently small. Here and below the powers
    of $r^{-1}$ and $\tau_0^{-1}$ are absorbed by decreasing $\eta$;
    this is legitimate because \eqref{eq:rho-upper-bound} also gives
    $\bar\d\leq\d^{c(t,\zeta)}$. Consequently,
    \[
        \bar t\in[\chi/2,2-\gamma].
    \]
    Put
    \[
        \bar u:=\min\{\chi/4,\gamma/2\}>0.
    \]
    Then $\bar u\leq\min\{\bar t,2-\bar t\}$, and
    \eqref{eq:minimal-spacing-2d} implies
    \[
        |\TT[T_w]|
        \leq \bar\d^{\bar u}|\TT|,
    \]
    since $\d^\chi\leq\bar\d^{\bar u}$.

    We may now apply Theorem
    \ref{thm:two-ends} at scale $\bar\d$, with parameters
    $(\bar s,\bar t,\bar u)$. By choosing $\eta'$ sufficiently
    small, the Frostman constant
    $r^{-O(1)}\tau_0^{-2}\d^{-\eta'}$ of the shadings is sufficiently small to apply Theorem
    \ref{thm:two-ends}.
    Notice that $|Y(T)|_{\bar \d} \gtrsim \frac{\tau_0^2 m}{|P \cap B_{\bar \d}|}$. 
    Since $P$ is $(\d, \d^{-\eta})$-uniform, we obtain
    \begin{align*}
        |P|_{\bar \d}
        &\geq |U(\TT,Y)|_{\bar\d}\\
        &\gtrsim
        \bar\d^{-\bar s\bar u/8}
        \frac{\tau_0^2m}{|P \cap B_{\bar \d}|}|\TT|^{1/2}\\
        &\gtrsim
        \bar\d^{-\bar s\bar u/8}
        \tau_0^2r^{O(1)}\d^{2\eta}|P|_{\bar \d},
    \end{align*}
    where the last line uses \eqref{eq:T-lower-bound} and
    \eqref{eq:G4-large-assumption}. Since
    $\bar\d\leq\d^{c(t,\zeta)}$ and
    $\bar s\bar u$ is bounded below in terms of the fixed
    parameters, this is impossible once
    $\eta\ll\eta'\ll\chi,\gamma,\zeta$. Thus
    \begin{equation}\label{eq:G4-small-2d}
        |G_4|<\d^\eta|P|^2.
    \end{equation}

    We finish by tracing the discarded pairs. From
    \eqref{eq:G4-core-loss} and \eqref{eq:G4-small-2d},
    \[
        |G_3|\lesssim(\tau_0+\d^\eta)|P|^2
        \lesssim\d^\eta|P|^2.
    \]
    Summing \eqref{eq:G2tau-total} over the
    $O(\log(1/\d))$ dyadic scales in \eqref{eq:G3-G2} gives
    \[
        |G_2|
        \lesssim
        \bigl(\d^\eta+
        \log(1/\d)\,r^{-4\sigma}\d^{\eta'-\eta}\bigr)|P|^2
        \lesssim\d^\eta|P|^2,
    \]
    provided $C_1$ is sufficiently large. Finally,
    \eqref{eq:G2-density-2d} gives
    \[
        |G_1|
        \lesssim(\log(1/\d))^2\d^\eta|P|^2
        \lesssim\d^{\eta/2}|P|^2.
    \]
    The graph $G:=G'_0\setminus G_1$ therefore has the required
    cardinality, by \eqref{eq:G0-size-2d} and
    \eqref{eq:close-pairs-2d}, and every $(p,p')\in G$ satisfies
    \[
        |T_{p,p'}(\varrho)\cap P|
        \leq\d^{-\zeta}\varrho^{\min\{1,t\}}|P|,
        \qquad \varrho\in[\d,1].
    \]
    Renaming the internal parameter $\eta/2$ as $\eta$ gives the
    exponent in the statement.
\end{proof}

\subsection{Radial projections in 3D}
We now adapt the proof of Theorem \ref{thm:radial-projections} to prove Theorem \ref{thm:radial-projections-3d}. We restate it here for convenience.

\begin{theorem}\label{thm:radial-projections-3d-again}
    Let $\varepsilon_1, \varepsilon_2, \zeta \in (0,1)$. Suppose that $\varepsilon_1 \le \zeta/100$. Then there exist $\delta_0 = \d_0(\varepsilon_1, \varepsilon_2, \zeta)>0$ and $\eta_0 = \eta_0(\varepsilon_1, \varepsilon_2, \zeta)>0$ such that the following holds for all $\delta < \delta_0$ and $\eta < \eta_0$.
    Let $P \subset [0,1]^3$ be a $\delta$-separated  $(\d, \d^{-\eta})$-uniform $(\delta, 2, \delta^{-\eta})$-set such that $|P \cap S| \le \d^{\varepsilon_2} |P|$ for any $\d^{\varepsilon_1}\times 1\times 1$-slab $S$.

    Then there exists a subset $G \subset P\times P$ such that $|G| \ge (1-\delta^{\eta/10}) |P|^2$ such that for every $(p,p') \in G$ we have $|T_{p,p'}(\varrho) \cap P| \le \delta^{-\zeta} \varrho^{2} |P|$ for all $\varrho\in [\delta, 1]$.

    In particular, for every $p\in P$ such that $|G|_p| \ge \nu |P|$,  the set of tubes $\{T_{p, p'}(\delta), ~(p, p')\in G\}$ is a $(\delta, 2, C\delta^{-\zeta} \nu^{-1})$-set.
\end{theorem}

\begin{remark}
    We state the result only for 2-dimensional sets $P$ since it suffices for our applications and allows for a simpler (although somewhat ad hoc) proof. A general sharp radial projections theorem in all dimensions was proven by Ren \cite{ren2023discretized}. However, the statements of Ren's results do not suffice for our applications although it is plausible that his techniques can be modified to prove an appropriate generalization of Theorem \ref{thm:radial-projections-3d-again}.
    
    Note that unlike Theorem \ref{thm:radial-again}, in this statement the parameter $\eta_0$ is allowed to depend on $\varepsilon_2$. It might be possible to remove this dependence via some generalization of Lemma \ref{lem:initial-thin-tubes} but it seems to require significant modifications to the proof strategy below and so we do not pursue it here.
\end{remark}

\begin{proof}
    Put
    \[
        \gamma:=\zeta/100,\qquad
        \sigma:=2-\gamma,\qquad
        r:=\d^{2\eta}.
    \]
    All constants implicit in $\d^{O(\eta)}$ are allowed to depend
    on $\varepsilon_1,\varepsilon_2,\zeta$. We first remove the pairs
    at distance less than $r$. Indeed,
    \begin{equation}\label{eq:close-pairs-3d}
        |\{(p,p')\in P^2:|p-p'|<r\}|
        \lesssim \d^{-\eta}r^2|P|^2
        \leq \d^\eta|P|^2.
    \end{equation}

    Let $G_{\rm bad}$ be the set of pairs $(p,p')$ for which
    \[
        |T_{p,p'}(q)\cap P|
        \geq \d^{-\zeta}q^2|P|
    \]
    at some scale $q\in[\d,1]$. Suppose, towards a contradiction,
    that
    \begin{equation}\label{eq:Gbad-large-3d}
        |G_{\rm bad}\cap\{|p-p'|\geq r\}|
        \geq \d^{\eta/3}|P|^2.
    \end{equation}
    For such a pair and a dyadic $q$, set
    \[
        W_{p,p'}(q):=
        q^{-\sigma}|T_{p,p'}(q)\cap P|\,|P|^{-1},
    \]
    and choose a scale at which $W_{p,p'}$ is maximal. As in the
    two-dimensional proof, a scale witnessing badness gives
    \[
        W_{p,p'}(q)
        \gtrsim \d^{-\zeta}q^\gamma
        \geq \d^{-(\zeta-\gamma)}.
    \]
    Pigeonholing the maximizing scale and the tube mass gives a
    dyadic $\varrho$, a number $m$, and a graph $G_1$ such that
    \begin{equation}\label{eq:stopping-data-3d}
    \begin{split}
        |G_1|&\geq \d^{\eta/2}|P|^2,\\
        |T_{p,p'}(\varrho)\cap P|&\sim m
            \gtrsim \d^{-(\zeta-\gamma)}
                     \varrho^{2-\gamma}|P|,
            \qquad (p,p')\in G_1,\\
        |T_{p,p'}(q)\cap P|
            &\lesssim (q/\varrho)^{2-\gamma}m,
            \qquad q\in[\varrho,1].
    \end{split}
    \end{equation}
    Here logarithmic losses have been absorbed into the change from
    $\d^{\eta/3}$ to $\d^{\eta/2}$.

    A $q$-tube in $\mathbb R^3$ can be covered by $O(q^{-1})$
    balls of radius $q$. Consequently,
    \[
        |P\cap T|\lesssim \d^{-\eta}q|P|
    \]
    for every $q$-tube $T$. Comparing this estimate at $q=\varrho$
    with \eqref{eq:stopping-data-3d} gives
    \begin{equation}\label{eq:rho-upper-3d}
        \varrho^{1-\gamma}
        \lesssim \d^{\zeta-\gamma-\eta},
        \qquad\text{and hence}\qquad
        \varrho\leq \d^{c\zeta}.
    \end{equation}

    \medskip
    \noindent\emph{Step 1: finding many heavy slices.}
    Let $\eta_1=C\eta/\zeta$, where $C$ is a sufficiently large
    constant, and let $\Theta$ be a $c\varrho$-net in $S^2$. For
    $\theta\in\Theta$, partition a ball of constant radius around the origin into finitely overlapping $\varrho$-slabs normal to
    $\theta$. Retain the slabs $S$ satisfying
    \begin{equation}\label{eq:heavy-rho-slab}
        |P\cap S|\geq \varrho^{1+\eta_1}|P|,
    \end{equation}
    and denote the retained family by $\mathbb S_\theta$.
    Enlarge the members of $\mathbb S_\theta$ by a constant factor so that every
    $2\varrho$-slab satisfying \eqref{eq:heavy-rho-slab} is contained
    in a member of $\mathbb S_\theta$.

    The union of the discarded slabs in a fixed direction $\theta$ contains
    at most
    \[
        O(\varrho^{-1})\varrho^{1+\eta_1}|P|
        \lesssim \varrho^{\eta_1}|P|
    \]
    points. 
    
    For each $\theta\in \Theta$ we have the following estimate:
    \begin{equation}\label{eq:not-covered-Stheta}
    \#\{(p, p') \in P\times P: |\pi_\theta(p-p')| \le \varrho,~\{p, p'\} \not\subset S ~ \forall S\in \mathbb S_\theta\} \lesssim \varrho^{\eta_1} \varrho |P|^2.    
    \end{equation}
    Indeed, first choose $p \in P$ arbitrarily and consider the $2\varrho$-slab $S_p=\{x \in \RR^3: |\pi_\theta(p-x)| \le \varrho\}$. If $|S_p| \le \varrho^{1+\eta_1} |P|$ then there are at most $ \varrho^{1+\eta_1}|P|$ choices for $p'$. Otherwise, note that $S_p \cap [-1,1]^3$ is contained in a $2\varrho\times 10\times 10$ slab with normal direction $\theta$ and so $S_p \cap P$ is covered by some slab $S \in \mathbb S_\theta$. It follows that there are zero choices for $p'$. 
    
    For a fixed pair $(p,p')$ at distance at least $r$, there are
    $\gtrsim\varrho^{-1}$ elements $\theta\in\Theta$ for which
    $|\theta\cdot(p-p')|\leq\varrho$. Summing
    \eqref{eq:not-covered-Stheta} over $\theta$, and then reversing
    the order of summation, shows that all but
    $O(\varrho^{\eta_1}|P|^2)$ pairs of $G_1$ are contained in
    $\gtrsim\varrho^{-1}$ heavy slabs. Since
    $\varrho^{\eta_1}\ll\d^{\eta/2}$ by
    \eqref{eq:rho-upper-3d}, we obtain a graph
    $G_2\subset G_1$ such that
    \begin{equation}\label{eq:many-slices-3d}
        |G_2|\gtrsim |G_1|,
        \qquad
        \#\{S\in\mathbb S:p,p'\in S\}
        \gtrsim\varrho^{-1}
    \end{equation}
    for every $(p,p')\in G_2$, where
    $\mathbb S:=\bigcup_{\theta\in\Theta}\mathbb S_\theta$. For $S\in \mathbb S$ let $G_2[S] = \{(p, p')\in G_2:~ p, p'\in S\}$. Using \eqref{eq:many-slices-3d} and double counting, we get
    \begin{equation}\label{eq:G2-density-average}
    \sum_{S\in \mathbb S} |G_2[S]| \gtrsim \varrho^{-1} |G_2| \gtrsim \d^{\eta/2}\varrho^{-1} |P|^2.    
    \end{equation}    
    If $p,p'\in S$ and $|p-p'|\geq r$, elementary geometry gives
    \[
        T_{p,p'}(\varrho)\cap[-2,2]^3
        \subset N_\Delta(S),
        \qquad
        \Delta\sim r^{-1}\varrho.
    \]
    Write $\Delta=\d^\beta$. 
    By \eqref{eq:rho-upper-3d}, for sufficiently small $\eta$ we have $c\zeta\leq\beta\leq1$.

    \medskip
    \noindent\emph{Step 2: two ends inside most slices.}
    Put $\kappa_0:=1/50$ and
    \begin{equation}\label{eq:slice-parameters-3d}
        \varepsilon_1^*
        :=\frac{(1+\kappa_0)\varepsilon_1}{\beta}.
    \end{equation}
    Choose $\varepsilon_2^*>0$ sufficiently small in terms of
    $\varepsilon_1^*,\varepsilon_2,\zeta$, and later choose
    $\eta\ll\beta\varepsilon_2^*$. For an incidence
    $((p,p'),S)$ occurring in \eqref{eq:many-slices-3d}, call the
    incidence \emph{concentrated} if
    \begin{equation}\label{eq:concentrated-slice-3d}
        |P\cap N_\Delta(S)\cap
          T_{p,p'}(\Delta^{\varepsilon_1^*})|
        >
        \Delta^{\varepsilon_2^*}
        |P\cap N_\Delta(S)|.
    \end{equation}

    We claim that the pairs for which a positive proportion
    of the incidences in \eqref{eq:many-slices-3d} are concentrated
    form a graph of size at most $\d^{C_0\eta}|G_2|$, where $C_0$
    may be chosen arbitrarily large. 

    Suppose the contrary and let $G_2'\subset G_2$ be the graph of all such pairs. Set
    \[
        \alpha:=M_0\varepsilon_2^*/\varepsilon_1^*,
        \qquad M_0:=200.
    \]
    For every concentrated incidence choose
    $\tau\in[\Delta^{\varepsilon_1^*},1]$ maximizing
    \[
        q^{-\alpha}
        |P\cap N_\Delta(S)\cap T_{p,p'}(q)|.
    \]
    Dyadically pigeonholing $\tau$, we get a graph $G_3 \subset G_2'$ of
    density $\d^{O(\eta)}$ in $P^2$ such that every pair $(p,p') \in G_3$
    has $\gtrsim\varrho^{-1}/\log(1/\d)$ many slabs $S \in \mathbb S$ incident to it such that $q^{-\alpha} |P\cap N_\Delta(S)\cap T_{p,p'}(q)|$ is maximized at $\tau$. 
    Maximality and \eqref{eq:concentrated-slice-3d} imply
    \begin{equation}\label{eq:tau-mass-3d}
    \begin{split}
        &|P\cap N_\Delta(S)\cap
        (T_{p,p'}(\tau)\setminus T_{p,p'}(c_\alpha\tau))|\\
        &\hspace{20mm}\gtrsim
        (\tau/\Delta^{\varepsilon_1^*})^\alpha
        \Delta^{\varepsilon_2^*}
        |P\cap N_\Delta(S)|.
    \end{split}
    \end{equation}
    On the other hand, this expression is at most
    $|P\cap N_\Delta(S)|$. Therefore
    \begin{equation}\label{eq:tau-upper-3d}
        \tau
        \lesssim
        \Delta^{\varepsilon_1^*
        -\varepsilon_2^*/\alpha}
        =
        \Delta^{(1-1/M_0)\varepsilon_1^*}.
    \end{equation}

    The normals of the slabs through a fixed pair run along a
    $\varrho$-net of a great circle. Choose a
    $C\Delta/\tau$-separated subfamily. It has
    $\gtrsim\tau/\Delta$ members, and the portions appearing on the
    left of \eqref{eq:tau-mass-3d} have bounded overlap. Using
    \eqref{eq:heavy-rho-slab}, we obtain, for every
    $(p,p')\in G_3$,
    \begin{equation}\label{eq:rich-tau-tube-3d}
        |P\cap T_{p,p'}(\tau)|
        \gtrsim \nu\tau|P|,
        \qquad
        \nu:=
        \frac{c\,r\varrho^{\eta_1}
              \Delta^{\varepsilon_2^*}}
             {\log(1/\d)}.
    \end{equation}
    In particular,
    $\nu\geq\d^{O(\eta)}\Delta^{\varepsilon_2^*}$. Note that this means that the tubes $T_{p,p'}(\tau)$ have essentially maximum possible intersection with $P$ allowed by the Frostman property. 

    By applying a standard hairbrush argument, if $|G_3| \ge \lambda |P|^2$, then we can find a slab $S$ such that 
    \begin{equation}\label{eq:hairbrush-output-3d}
    \begin{split}
        \operatorname{width}(S_\tau)
        &\lesssim
        \d^{-C\eta}(\lambda\nu)^{-C}\tau,\\
        |P\cap S_\tau|
        &\gtrsim
        \d^{C\eta}(\lambda\nu)^C|P|.
    \end{split}
    \end{equation}
    Let us briefly outline the main steps needed to justify this. First, consider the maximal family $\TT$ of essentially distinct $\tau$-tubes which contain $\gtrsim \nu \tau$-fraction of $P$. By double counting, we show that $|\TT| = (\lambda \nu)^{O(1)} \tau^{-2}$.
    By Cauchy--Schwarz inequality, there are $(\lambda \nu)^{O(1)} |\TT|^2$ many pairs of tubes $T, T' \in \TT$ that have non-empty intersection. Let $\theta$ be the typical angle between a pair of intersecting tubes in $\TT$. Then $\TT$ splits into essentially disjoint families $\TT[T_\theta]$ for some collection of $\sim\theta$-tubes $T_\theta$. By using a two-ends argument on the subfamily $\TT[T_\theta]$, we then conclude that $|P\cap T_\theta| \gtrsim (\lambda \nu)^{O(1)} |P|$. On the other hand the Frostman property of $P$ implies that $|P\cap T_\theta|\lesssim \d^{-\eta} \theta |P|$. So it follows that $\theta \gtrsim_{\d^{O(1)}} (\lambda \nu)^{O(1)}$, i.e. the tubes in $\TT$ intersect roughly transversally. Now the classical Wolff hairbrush argument shows that a large fraction of tubes in $\TT$ must be contained in a $\approx \tau$ slab $S_\tau$, giving \eqref{eq:hairbrush-output-3d}. We refer to \cite[Lemma 3.5]{maldague2025heilbronn} for a rigorous argument proving a closely related statement. 

    Now we use \eqref{eq:hairbrush-output-3d} to reach a contradiction. Recall that by assumption, $\lambda\geq\d^{O(\eta)}$. Choose
    $\varepsilon_2^*$ sufficiently small, and then $\eta$
    sufficiently small, so that \eqref{eq:tau-upper-3d} and
    \eqref{eq:hairbrush-output-3d} give
    \[
        \operatorname{width}(S_\tau)
        \leq\d^{\varepsilon_1},
        \qquad
        |P\cap S_\tau|>\d^{\varepsilon_2}|P|.
    \]
    Indeed, by \eqref{eq:slice-parameters-3d},
    $\Delta^{\varepsilon_1^*}
    =\d^{(1+\kappa_0)\varepsilon_1}$, leaving room to absorb the
    factors $\Delta^{-C\varepsilon_2^*}\d^{-C\eta}$ in the slab
    width. This contradicts the hypothesis of the theorem and
    proves the claim.

    Let $G_4 \subset G_2$ be the set of pairs $(p,p')$ such that there are $\gtrsim c\varrho^{-1}$ many slabs $S \in \mathbb S$ such that $((p,p'), S)$ is not concentrated. That is, every $(p,p')\in G_4$ belongs to
    $\gtrsim\varrho^{-1}$ slabs $S$ satisfying
    \begin{equation}\label{eq:slice-two-ends-3d}
        |P\cap N_\Delta(S)\cap
          T_{p,p'}(\Delta^{\varepsilon_1^*})|
        \leq
        \Delta^{\varepsilon_2^*}
        |P\cap N_\Delta(S)|.
    \end{equation}
    By the above claim and \eqref{eq:stopping-data-3d}, \eqref{eq:many-slices-3d}, we get that $|G_4|\gtrsim|G_2|$ holds. Using \eqref{eq:G2-density-average} and the fact that any pair $(p,p') \in G_2\setminus G_2$ is contained in at most $\lesssim r^{-1}\varrho^{-1}$-many slabs $S\in \mathbb S$, we obtain
    \begin{equation}\label{eq:G4-dense}
        \begin{aligned}
            \sum_{S \in \mathbb S} |G_4[S]| &\ge \sum_{S \in \mathbb S} |G_2[S]| - |(G_2\setminus G_4)[S]| \\ &\gtrsim \varrho^{-1} \d^{\eta/2} |P|^2 - C r^{-1} \varrho^{-1} \d^{C_0\eta} |P|^2 \gtrsim \varrho^{-1} \d^{\eta/2} |P|^2.
        \end{aligned}
    \end{equation}
    In particular, 
    \begin{equation}\label{eq:PcapS-dense}
        \sum_{S\in \mathbb S}|P \cap S|^2 \gtrsim_{\d^{O(\eta)}} \varrho^{-1} |P|^2.
    \end{equation} 
    
    

    \medskip
    \noindent\emph{Step 3: selecting and regularizing one slice.}
    We now apply a Marstrand slicing argument to locate a slab $S$ so that a large subset of $P\cap N_\Delta S$ satisfies 1-dimensional Frostman condition.  For $S_\Delta=N_\Delta(S)$, write
    $\mu_S:=\mathbf 1_{P\cap S_\Delta}$ and define the truncated Riesz energy of $\mu_S$ by
    \[
        I_1^\Delta(\mu_S)
        :=
        \sum_{p,q\in P\cap S_\Delta}
        \frac{1}{\max\{|p-q|,\Delta\}}.
    \]
    A pair $p,q$ belongs to at most
    \[
        \frac{r^{-O(1)}}
             {\Delta\max\{|p-q|,\Delta\}}
    \]
    members of the family $\{S_\Delta:S\in\mathbb S\}$. Hence
    \begin{align}
        \sum_{S\in\mathbb S}I_1^\Delta(\mu_S)
        &\lesssim
        \frac{r^{-O(1)}}{\Delta}
        \sum_{p,q\in P}
        \frac{1}{\max\{|p-q|,\Delta\}^2}
        \notag\\
        &\lesssim
        \d^{-O(\eta)}\Delta^{-1}|P|^2.
        \label{eq:riesz-sum-3d}
    \end{align}
    The last inequality follows by splitting into dyadic distance
    ranges and using the two-dimensional Frostman estimate.

    Dyadically pigeonhole the set of slabs $\mathbb S$ so that $|G_4[S]|$ and $|P \cap S_\Delta|$ are roughly the same for all $S \in \mathbb S$ and \eqref{eq:G4-dense} still holds. 
    Using the lower bound $I_1^\Delta(\mu_S) \gtrsim |P\cap S_\Delta|^2$, \eqref{eq:riesz-sum-3d} implies $\sum_{S\in \mathbb S} |P\cap S_\Delta| \gtrsim_{\d^{O(\eta)}} \Delta^{-1}|P|^2$. Together with \eqref{eq:PcapS-dense} this implies that 
    \[
    |\mathbb S| |P\cap S_\Delta|^2 \sim_{\d^{O(\eta)}} \Delta^{-1} |P|^2
    \] 
    for all $S \in \mathbb S$. It follows from  \eqref{eq:G4-dense} that $|G_4[S]| \sim_{\d^{O(\eta)}} |P\cap S_\Delta|^2$ for every $S \in \mathbb S$. 
    After pruning, \eqref{eq:riesz-sum-3d} implies that for every $S \in \mathbb S$ we have 
    \[
    I_1^\Delta(\mu_S) \lesssim \d^{-O(\eta)} \Delta^{-1} |P|^2 |\mathbb S|^{-1} \lesssim \d^{-O(\eta)} |P\cap S_\Delta|^2.
    \]
    Let $K = \d^{-C_1\eta}$ for a sufficiently large constant $C_1$ and let $P_S \subset P\cap S_\Delta$ be a $(\Delta,1,K)$-set contained in $P\cap S_\Delta$ of the maximum possible size. We claim that for sufficiently large $C_1$, we have the estimate $|P_S| \ge (1-K^{-1/2})|P \cap S_\Delta|$. Indeed, for every dyadic scale $w \in [\Delta,1]$ consider the set $\mc Q_{w,j}$ of finitely overlapping $w$-cubes ${\bf q}$ satisfying $|P \cap S_\Delta \cap {\bf q}| \sim 2^j w |P \cap S_\Delta|$. By definition, we have 
    \[
    I_1^\Delta(\mu_S) \gtrsim |\mc Q_{w,j}| w^{-1} (2^jw |P \cap S_\Delta|)^2
    \]
    and so the upper bound on Riesz energy gives $|\mc Q_{w,j}| \lesssim \d^{-O(\eta)} 2^{-2j}w^{-1}$ and $|\bigcup_{{\bf q}\in \mc Q_{w,j}} {\bf q}\cap P\cap S_\Delta| \lesssim \d^{-O(\eta)} 2^{-j}|P \cap S_\Delta|$. Summing this over all dyadic $w$ and all $j \ge \log_2 K$, we get that the union of all rich cubes covers at most $\d^{-O(\eta)} K^{-1}$ fraction of points of $P\cap S_\Delta$. For $C_1$ large enough this fraction is less than $K^{-1/2}$ and removing it from the set leaves a $(\Delta,1,K)$-set $P_S$. 

    Observe that 
    \[
    \begin{aligned}
    |G_4[S]\cap (P_S\times P_S)| &\ge |G_4[S]| - 2|(P\cap S_\Delta)\setminus P_S| \cdot |P\cap S_\Delta|\\ &\gtrsim \d^{O(\eta)} |P\cap S_\Delta|^2 - K^{-1/2} |P\cap S_\Delta|^2 \\&\gtrsim \d^{O(\eta)} |P\cap S_\Delta|^2.    
    \end{aligned}
    \]

    \medskip
    \noindent\emph{Step 4: reduction to the planar
    argument.} 
    Let $G_S = G_4[S] \cap (P_S\times P_S)$.
    We conclude that $S \in \mathbb S$ and $S_\Delta = N_\Delta S$ have the following properties:

    \begin{itemize}
        \item[(i)] The set $P_{S}$ is $(\Delta, 1, K)$-set,
        \item[(ii)] $|G_S| \gtrsim \d^{O(\eta)} |P\cap S_\Delta|^2 \gtrsim \d^{O(\eta)} |P_S|^2$,
        \item[(iii)]  every $(p, p') \in G_S$ satisfies 
    \[
    |T_{p,p'}(\varrho) \cap P\cap S_\Delta| \ge \d^{-\zeta/2} \varrho^2 |P|
    \]
        \item[(iv)] every $(p,p') \in G_S$ satisfies the two-ends condition
    \[
    |T_{p, p'}(\Delta^{\varepsilon_1^*}) \cap P\cap S_\Delta| \le \Delta^{\varepsilon_2^*} |P \cap S_\Delta|.
    \] 
    \end{itemize}

    We want to use properties (i)--(iv) to reduce to the situation to the one in the proof of Theorem \ref{thm:radial-again}. To achieve this, we need to strengthen condition (iii). First, we estimate $M\sim |P \cap S_\Delta|$ from above. Let $(p,p') \in G_4$, there are $\sim \varrho^{-1}$ slabs $S\in \mathbb S$ incindent to $(p,p')$ so that \eqref{eq:slice-two-ends-3d} holds. In particular, $|P \cap S \setminus T_{p,p'}(\Delta^{\varepsilon_1^*})| \gtrsim M$ for every such $S$. For slabs through $p,p'$ with $C\varrho/\Delta^{\varepsilon_1^*}$-angularly separated normal directions, the corresponding sets $S \setminus T_{p,p'}(\Delta^{\varepsilon_1^*})$ are pairwise disjoint. Thus, we get
    \[
    M \frac{\Delta^{\varepsilon_1^*}}{\varrho} \lesssim |P|, \quad M \lesssim \Delta^{-\varepsilon_1^*} \varrho |P|. 
    \]
    Recall that $\varepsilon_1^* = \frac{(1+\kappa_0)}{\beta}\varepsilon_1$, $\kappa_0 = 1/50$, $\Delta=\d^\beta$ and $\varepsilon_1 \le \zeta/100$. It follows that $\Delta^{-\varepsilon_1^*} < \d^{-\zeta/10}$. So we conclude that $|P\cap S_\Delta| \lesssim \d^{-\zeta/10} \varrho|P|$ and (iii) implies 
    \begin{equation}\label{eq:TppP}
    |T_{p,p'}(\varrho) \cap P\cap S_\Delta| \ge \d^{-\zeta/3} \varrho |P \cap S_\Delta|.    
    \end{equation}
    
    Second, we claim that we can replace the set $P\cap S_\Delta$ in \eqref{eq:TppP} with $P_S\cap S$.  
    Let $\TT$ be a maximal subset of essentially distinct $\sim\varrho$-tubes covering the collection $T_{p,p'}(\varrho)$ over all $(p, p') \in G_S$. 
    We have 
    \begin{align*}
        |P\cap S_\Delta|^2 &\gtrsim_{\d^{O(\eta)}}  I_1^\Delta(\mu_S)= \sum_{p, p' \in P\cap S_\Delta} \frac{1}{\max(|p-p'|, \Delta)} \\
        &\gtrsim_{\d^{O(\eta)}} \sum_{p, p' \in P\cap S_\Delta} \#\{T \in \TT: p,p' \in \TT\} \\
        &= \sum_{T\in \TT} |T\cap P \cap S_\Delta|^2 \gtrsim |\TT| \d^{-\zeta/3} (\varrho |P \cap S_\Delta|)^2
    \end{align*}
    By rearranging we get the upper bound $|\TT| \lesssim_{\d^{O(\eta)}} \d^{\zeta/3}\varrho^{-2}$. Thus, on average, a tube $T \in \TT$ contains at least $\gtrsim_{\d^{O(\eta)}}\d^{\zeta/3} \varrho^2|G_S|$ many pairs $(p,p')\in G_S$. Since $G_S$ is supported on $P_S\cap S$, we in particular conclude that $|T\cap P_S\cap S| \gtrsim _{\d^{O(\eta)}}\d^{\zeta/6} \varrho|P_S|$ for every such tube $T$. So by passing to a dense subset in $G_S$ we can guarantee that 
    \[
    |T_{p,p'}(C\Delta) \cap P_S\cap S| \ge \d^{-\zeta/8} \Delta|P_S|,
    \]
    for every $(p,p') \in G_S$. 
    
    We can now repeat the proof of Theorem \ref{thm:radial-again} starting from Step 1 using the graph $G_S$ in place of $G_1$ (after projecting everything onto the plane parallel to $S$). More precisely, we have the following matching between parameters in the current setup and Theorem \ref{thm:radial-again}:
    \begin{table}[]
        \centering
        \begin{tabular}{c|c}
             $\Delta$ & $\d$ \\
             $\zeta/8\beta$ & $\zeta$ \\
             $\varepsilon_1^*$ & $\varepsilon_1$ \\
             $\varepsilon_2^*$ & $\varepsilon_2$ \\
             $ C_1 \eta/\beta $ & $ \eta$
        \end{tabular}
        \caption{Left column lists the parameters of the configuration $(P_S, G_S)$ that was constructed in the proof so far, the right column lists the corresponding parameters in the proof of Theorem \ref{thm:radial-again}}
        \label{tab:parameters}
    \end{table}

    The restriction $\varepsilon_1 \le \zeta/100$ in the statement of Theorem \ref{thm:radial-projections-3d-again} implies that $\varepsilon_1^* \le (\zeta /8\beta)/8$ holds. So for small enough $\eta$ and $\d$, the parameters above satisfy the setup of Theorem \ref{thm:radial-again} and so its proof applies and we eventually reach a contradiction with the assumption that $|G_S| \gtrsim_{\Delta^{O(\eta/\beta)}} |P_S|^2$. This implies that the initial assumption $|G_1| \ge \d^{\eta/2}|P|^2$ was false and completes the proof of the theorem.

\end{proof}

\end{document}